%% file: gdnre-jas.tex
\documentclass[]{interact}

\usepackage{epstopdf}%
\usepackage[caption=false]{subfig}%

\usepackage[numbers,sort&compress]{natbib}%
\bibpunct[, ]{[}{]}{,}{n}{,}{,}%

\usepackage[utf8]{inputenc}
\usepackage{tikz}
\usetikzlibrary{plotmarks}
\usepackage{pgfplots} 
\pgfplotsset{compat=newest} 
\pgfplotsset{plot coordinates/math parser=false}
\newlength\fwidth
\newlength\fheight
\usepackage[makeroom]{cancel}
\newcommand{\Expc}{\mathop{\textrm{E}}}

\theoremstyle{plain}%
\newtheorem{theorem}{Theorem}[section]
\newtheorem{lemma}[theorem]{Lemma}

\theoremstyle{definition}

\theoremstyle{remark}
\newtheorem{remark}{Remark}

\begin{document}

\articletype{}%

\title{General Distributions of Number Representation Elements}

\author{
\name{Félix Balado\thanks{CONTACT A.~F. Balado. Email: felix@ucd.ie}
  and Guénolé C.M. Silvestre}
\affil{School of Computer Science, University College Dublin, Dublin, Ireland}
}

\maketitle

\begin{abstract}
  We provide general expressions for the joint distributions of
  the~$k$ most significant $b$-ary digits and of the $k$~leading
  continued fraction coefficients of outcomes of an arbitrary
  continuous random variable. Our analysis highlights the connections
  between the two problems. In particular, we give the general
  convergence law of the distribution of the $j$-th significant digit,
  which is the counterpart of the general convergence law of the
  distribution of the $j$-th continued fraction coefficient
  (Gauss-Kuz'min law). We also particularise our general results for
  Benford and Pareto random variables. The former particularisation
  allows us to show the central role played by Benford variables in
  the asymptotics of the general expressions, among several other
  results, including the equivalent of Benford's law for continued
  fractions.  The particularisation for Pareto variables ---which
  include Benford variables as a special case--- is specially relevant
  in the context of pervasive scale-invariant phenomena, where Pareto
  variables occur much more frequently than Benford variables. This
  suggests that the Pareto expressions that we produce have wider
  applicability than their Benford counterparts in modelling most
  significant digits and leading continued fraction coefficients of
  real data. Our results may find practical application in all areas
  where Benford's law has been previously used.
\end{abstract}

\begin{keywords}
  Significant digits models; continued fraction coefficients models;
  Benford's law; Pareto distribution
\end{keywords}

\section{Introduction}
\label{sec:introduction}
Benford's law for the distribution of the most significant
digits~\cite{benford1938} has found practical use in numerous
statistical applications related to fraud detection, data quality
analysis and consistency, generation of synthetic data, floating point
operations error analysis, etc, in areas such as accountancy,
economics, signal processing, computer science, and
beyond~\cite{miller15:_benfords_law}. However, Benford's law is known
not to apply to many datasets, and therefore it would be desirable to
have completely general models of the distribution of significant
digits to avail of more powerful models in practical
applications. Furthermore, the positional numeral systems that lead to
most significant digits models have a little explored alternative:
real numbers can also be represented using continued fraction
expansions, and so leading continued fraction coefficients can play
a role analogous to that of significant digits.

Our main goal in this paper is to provide a general treatment of the
related problems of modelling probabilistically the most significant
digits and the leading continued fraction coefficients of outcomes of
an arbitrary continuous random variable, and to evince the
parallelisms between the two problems. Ever since the observations
made by Newcomb~\cite{newcomb81:_note} and Benford~\cite{benford1938},
studies of the distribution of most significant digits have largely
focused on Benford variables ---for an overview, see the introduction
by Berger and Hill~\cite{berger11:_basic_th_benford} and the book by
Miller~\cite{miller15:_benfords_law}. Some generalisations have been
pursued by Pietronero et al.~\cite{pietronero01:_explaining} and by
Barabesi and Pratelli~\cite{barabesi20:_generalized}, among others,
but a truly general approach to modelling significant digits has never
been presented. A lot less attention has been devoted to modelling
continued fraction (CF) coefficients. Except for the well-known
Gauss-Kuz'min asymptotic law~\cite{khinchin61:continued} and an
approximation due to Blachman~\cite{blachman84}, most of the work in
this area has been pioneered by Miller and
Takloo-Bighash~\cite{miller06:_invitation}. In any case, all existing
finite results for CF coefficients models are solely for the
particular case in which the fractional part of the data represented
by means of continued fractions is uniformly distributed. No general
approach has been investigated in this problem either.

This paper is organised as follows. In
Section~\ref{sec:gen-prob-msd} we give the general expression for
the joint distribution of the $k$ most significant $b$-ary digits of
outcomes drawn from an arbitrary positive-valued distribution. One
application of our analysis is a proof of the general asymptotic
distribution of the $j$-th most significant $b$-ary digit, which, as
we will discuss, is the near exact counterpart of the Gauss-Kuz'min
law for the general asymptotic distribution of the $j$-th continued
fraction coefficient. Our approach to modelling the $k$ most
significant $b$-ary digits through a single variable ---rather than
through $k$ separate variables as in previous works--- leads to
further contributions\footnote{Some preliminary results in this
  paper previously appeared in~\cite{balado21:_benford}.}  in the
particularisation of our general results in
Section~\ref{sec:particular-cases}.  Therein we produce a new
closed-form expression for the distribution of the $j$-th significant
$b$-ary digit of a Benford variable, and we give a short new proof of the
asymptotic sum-invariance property of these variables. We also show
that Benford's distribution is just a particular case of a more
general distribution based on Pareto variables, which must have wider
applicability in the pervasive realm of scale-invariant data ---a fact
first pointed out by Pietronero et al.~\cite{pietronero01:_explaining}
and then expanded upon by Barabesi and
Pratelli~\cite{barabesi20:_generalized}, who however did not give
results as complete as ours.

In Section~\ref{sec:gen-prob-cf-coeff} we give the general expression
for the joint distribution of the $k$ leading continued fraction
coefficients of outcomes drawn from an arbitrary distribution. This is
shown to be explicitly analogous to modelling the $k$ most significant
$b$-ary digits of the data when the continued fraction coefficients
correspond to the logarithm base $b$ of the data. Therefore, modelling
leading CF coefficients is a realistic practical alternative to
modelling most significant digits. Most of the results in
Section~\ref{sec:gen-prob-cf-coeff} are novel, and so are their
particularisations in Section~\ref{sec:particular-cases} ---in
particular the counterpart of Benford's law for continued fractions---
except for the special cases previously given by Miller and
Takloo-Bighash~\cite{miller06:_invitation}.

Additionally, we show in Section~\ref{sec:benf-vari-spec} the central
role played by the particular analysis for Benford variables in the
asymptotics of the general expressions ---both when modelling
significant digits and leading continued fraction coefficients--- and
we demonstrate this numerically in the case of Pareto variables in
Section~\ref{sec:benf-based-appr}.

Finally, we empirically verify all our theoretical results in
Section~\ref{sec:empirical-tests}, using both Monte Carlo experiments
and real datasets.

\textbf{Notation and preliminaries.}  Calligraphic letters are sets,
and $\vert\mathcal{V}\vert$ is the cardinality of
set~$\mathcal{V}$. Boldface Roman letters are row vectors. Random
variables (r.v.'s) are denoted by capital Roman letters, or by
functions of these. The cumulative distribution function (cdf) of
r.v.~$Z$ is $F_Z(z)=\Pr(Z\le z)$, where $z\in\mathbb{R}$. The
expectation of $Z$ is denoted by $\Expc(Z)$. If~$Z$ is continuous with
support~$\mathcal{Z}$, its probability density function (pdf)
is~$f_Z(z)$, where $z\in\mathcal{Z}$. A r.v.~$Z$ which is uniformly
distributed between $a$ and $b$ is denoted by~$Z\sim U(a,b)$. The
probability mass function (pmf) of a discrete r.v. $Z$ with support
$\mathcal{Z}$ is denoted by $\Pr(Z=z)$, where $z\in\mathcal{Z}$. The
unit-step function is defined as ~$u(z)=1$ if $z\ge 0$, and $u(z)=0$
otherwise. The fractional part of $z\in \mathbb{R}$ is
$\{z\}=z-\lfloor z\rfloor$.  Curly braces are also used to list the
elements of a discrete set, and the meaning of $\{z\}$ (i.e. either a
fractional part or a one-element set) is clear from the context. We
exclude zero from the set of natural numbers $\mathbb{N}$. We use
Knuth's notation for the rising factorial powers of $z\in\mathbb{R}$:
$z^{\overline{m}}=\Pi_{i=0}^{m-1}(z+i)=\Gamma(z+m)/\Gamma(z)$~\cite{graham99:_concrete}.

Throughout the manuscript, $X$ denotes a positive continuous~r.v. We
also define the associated r.v.
\begin{equation*}\label{eq:y}
  Y=\log_b X
\end{equation*}
for an arbitrary $b\in\mathbb{N}\backslash\{1\}$. The fractional part
of $Y$, i.e. $\{Y\}$, will play a particularly relevant role in our
analysis. The cdf of $\{Y\}$ is obtained from the cdf of $Y$ as
follows:
\begin{equation}\label{eq:cdffrcy}
  F_{\{Y\}}(y)=\sum_{i\in\mathbb{Z}}F_{Y}(y+i)-F_{Y}(i),
\end{equation}
for~$y\in[0,1)$. Because $\{Y\}$ is a fractional part, it always holds
that $F_{\{Y\}}(y)=0$ for $y\le 0$ and $F_{\{Y\}}(y)=1$ for $y\ge
1$. Also, ~$F_{\{Y\}}(y)$ is a continuous function of $y$ because $X$
is a continuous r.v.

\section{General Probability Distribution of the \lowercase{$k$} Most  Significant \lowercase{$b$}-ary Digits}\label{sec:gen-prob-msd}
In this section we will obtain the general expression for the joint
probability distribution of the~$k$ most significant digits of
a positive real number written in a positional base~$b$ numeral system,
where $b\in\mathbb{N}\backslash \{1\}$. Let us first define
\begin{equation*}
  \label{eq:support_set_bary_digits}
  \mathcal{A}=\{0,1,\ldots,b-1\}.
\end{equation*}
The $b$-ary representation of $x\in\mathbb{R}^+$ is formed by the
unique digits $a_i\in\mathcal{A}$ such that
$x=\sum_{i\in \mathbb{Z}} a_i\, b^i$ ---unicity requires ruling out
representations where $a_i=b-1$ for all $i<j$, where $j<0$ and
$a_j<b-1$ or $j=0$ and $a_j\in\mathcal{A}$. If we now let
$n=\lfloor \log_b x\rfloor$, the most significant $b$-ary digit of $x$
is~$a_n$. This is because the definition of~$n$ implies
$n\le \log_b x < n+1$, or, equivalently, $b^{n}\le x < b^{n+1}$. Using
$n$, the~$k$ most significant $b$-ary digits of~$x$ can be inferred as
follows: %
\begin{equation}
  a=\lfloor x\,b^{-n+k-1}\rfloor=\lfloor b^{\{\log_b x\}+k-1}\rfloor.\label{eq:a}
\end{equation}
By using $0\le \{\log_b x\}<1$ in~\eqref{eq:a} we can verify that $a$
belongs to the following set of
integers: %
\begin{equation}\label{eq:support}
  \mathcal{A}_{(k)}=\{b^{k-1},\ldots, b^{k}-1\}, %
\end{equation}
whose cardinality is $\vert\mathcal{A}_{(k)}\vert=b^{k}-b^{k-1}$. We
propose to call $a$ in~\eqref{eq:a} the $k$-th \textit{integer
  significand} of $x$. We must mention that for some authors the
significand of~$x$ is the integer $\lfloor x\,b^{-n}\rfloor=a_n$, or
even $\lfloor x\,b^{-n+k-1}\rfloor$ itself~\cite[page
7]{nigrini12:_benford}, but for some others the significand of $x$ is
the real
$x\,b^{-n}\in
[1,b)$~\cite{berger11:_basic_th_benford,miller15:_benfords_law} ---which is sometimes also called the normalised significand. In any
case, the advantages of consistently working with the $k$-th integer
significand will become clear throughout this paper.
To give an example of~\eqref{eq:a} and~\eqref{eq:support}, say that
$b=10$ and $x=0.00456678$. In this case $n=\lfloor\log_{10}
x\rfloor=-3$, so if we choose for instance $k=2$ then $a=\lfloor x\,
10^{4}\rfloor=\lfloor
10^{1.65961}\rfloor=45\in\mathcal{A}_{(2)}=\{10,11,12,\ldots,98,99\}$.

\begin{theorem}[General distribution of the $k$ most significant $b$-ary digits]\label{thm:msd}
  If~$A_{(k)}$ denotes the discrete r.v. that models the
  $k$ most significant $b$-ary digits (i.e. the
  $k$-th integer significand) of a positive continuous r.v.  $X$, then
  \begin{equation}\label{eq:pmfA}  
    \Pr(A_{(k)}=a)=F_{\{Y\}}\big(\log_b(a+1)-k+1\big)-F_{\{Y\}}\big(\log_ba-k+1\big),
     \rule[-1em]{0pt}{0pt}
   \end{equation}
where $a\in\mathcal{A}_{(k)}$ and $Y=\log_b X$. 
\end{theorem}

\begin{proof}
  Seeing~\eqref{eq:a}, the r.v. we are interested in is 
  defined as
  \begin{equation}
  A_{(k)}=\lfloor b^{\{\log_b    X\}+k-1}\rfloor.\label{eq:ak_definition}
\end{equation}
  From this definition, $A_{(k)}=a$ when
  $a\le b^{\{\log_b X\}+k-1}< a+1$, or, equivalently, when
\begin{equation}
  \log_b a-k+1\le \{\log_b X\} < \log_b (a+1)-k+1.\label{eq:ineq_mod}
\end{equation}
Using~\eqref{eq:ineq_mod} and the cdf of $\{Y\}=\{\log_b X\}$ we get~\eqref{eq:pmfA}.
\end{proof}

\begin{remark}
  It is straightforward to verify that the pmf~\eqref{eq:pmfA} adds up
  to one over its support, as
  \begin{equation}
    \label{eq:pmfA_verification}
    \sum_{a\in\mathcal{A}_{(k)}} \Pr(A_{(k)}=a)=F_{\{Y\}}\big(1\big)-F_{\{Y\}}\big(0\big)
  \end{equation}
  due to the cancellation of all consecutive terms in the telescoping
  sum on the left-hand side of~\eqref{eq:pmfA_verification} except for
  the two shown on the right-hand side. Because $F_{\{Y\}}(y)$ is the
  cdf of a r.v. with support~$[0,1)$, the right-hand side
  of~\eqref{eq:pmfA_verification} must equal one. %
\end{remark}

\subsection{Distribution of the $j$-th Most Significant  $b$-ary  Digit}
\label{sec:prob-distr-k}

Next, let us denote by $A_{[j]}$ the r.v. that models the $j$-th most
significant $b$-ary digit of~$X$. This variable can be obtained from
the variable $A_{(j)}$ that models the $j$-th integer significand as
follows:
\begin{equation}
  A_{[j]}=A_{(j)}\pmod{b}.\label{eq:aj_mod_definition}
\end{equation}
Obviously,
$A_{[1]}=A_{(1)}$. From~\eqref{eq:aj_mod_definition}, the pmf of
$A_{[j]}$ 
for $j\ge 2$ is
\begin{align}
  \label{eq:pmfAj}
  \Pr(A_{[j]}\!=\!a)&= \sum_{r\in\mathcal{A}_{(j-1)}} \Pr(A_{(j)}=r
                        b+a),%
\end{align}
where $a\in\mathcal{A}$.

\begin{remark}\label{rem:kth_significand}
  Observe that~\eqref{eq:pmfA} is also the joint pmf of
  $A_{[1]},\dots,A_{[k]}$. To see this we just have to write
  $a=\sum_{j=1}^{k}a_j b^{-j+k}$, which implies that
  $\Pr(A_{[1]}=a_1,\dots,A_{[k]}=a_k)=\Pr(A_{(k)}=a)$. Under this
  view,~\eqref{eq:pmfAj} is simply a marginalisation
  of~\eqref{eq:pmfA}.  However, the derivation of~\eqref{eq:pmfA} is
  simpler using the $k$-th integer significand variable $A_{(k)}$ than
  using $A_{[1]},\dots,A_{[k]}$. Further examples of the advantages of
  working with $k$-th integer significands are the following theorem
  and the results in Sections~\ref{sec:most-sign-digits}
  and~\ref{sec:most-sign-digits-1}.
\end{remark}

\begin{theorem}[General asymptotic distribution of the $j$-th most significant
  $b$-ary digit]\label{thm:asympt-aj}
  For any positive continuous random variable $X$, it holds that
  \begin{equation}
    \lim_{j\to\infty}\Pr(A_{[j]}=a)=  b^{-1}.\label{eq:asympt-aj}
  \end{equation}
\end{theorem}
\begin{proof}
  For any $\epsilon>0$ there exists $j_{\min}$ such that for all
  $j\ge j_{\min}$
  \begin{equation}
    \label{eq:epsilon}
    \log_b(rb+a)-\log_b(rb)<\epsilon
  \end{equation}
  for all $a\in\mathcal{A}$ and
  $r\in\mathcal{A}_{(j-1)}$. Specifically,
  $j_{\min}=\lceil 1+\log_b(\frac{b-1}{b^\epsilon-1})
  \rceil$. Therefore, inequality~\eqref{eq:epsilon} and the continuity
  of $F_{\{Y\}}(y)$ in~\eqref{eq:pmfA} imply that
  $\lim_{j\to\infty}\Pr(A_{(j)}=r b+a)-\Pr(A_{(j)}=r b)=0$ for all
  $a\in\mathcal{A}$. Thus, from~\eqref{eq:pmfAj}, we have
  that~\eqref{eq:asympt-aj} holds.
\end{proof}

\begin{remark}
  In general, the larger~$b$ is, the faster the convergence
  of~$A_{[j]}$ to a uniform discrete r.v. This is because $j_{\min}$
  is nonincreasing on $b$ for a given $\epsilon$.  Informally,
  Theorem~\ref{thm:asympt-aj} can be argued as follows: since
  $rb\ge b^{j-1}$ in~\eqref{eq:pmfAj}, for large~$j$ we have that
  $rb\gg a$, and therefore $rb+a\approx rb$. Consequently, when $j$ is
  large $\Pr(A_{(j)}=rb+a)\approx \Pr(A_{(j)}=rb)$ because of the
  continuity of the cdf of~$\{Y\}$. In such case~\eqref{eq:pmfAj} is
  approximately constant over~$a\in\mathcal{A}$, and so we have
  that~$A_{[j]}$ is asymptotically uniformly distributed.
\end{remark}

\section{General Probability Distribution of the \lowercase{$k$}
  Leading Continued Fraction Coefficients}
\label{sec:gen-prob-cf-coeff}

Continued fraction expansions are an alternative to positional
base~$b$ numeral systems for representing real numbers.  In this
section we will obtain the general expression for the joint
probability distribution of the $k$ leading coefficients in the simple
CF of a real number. %
Let $y_0=y\in \mathbb{R}$ and
define the recursion $y_{j}=\{y_{j-1}\}^{-1}$ based on the
decomposition $y_{j-1}=\lfloor y_{j-1}\rfloor +\{y_{j-1}\}$. By letting
$a_j=\lfloor y_j\rfloor$ we can express $y$ as the following
simple~CF:
\begin{equation}\label{eq:continued_y0}
  y=a_0+\cfrac{1}{a_1+\cfrac{1}{a_2+\cdots}} %
\end{equation}
A CF is termed \textit{simple} (or regular) when all the numerators of
the nested fractions equal~$1$. %
For typographical convenience, we will write the CF representation of $y$
in~\eqref{eq:continued_y0} using the standard notation
\begin{equation*}\label{eq:continued_y0_2}
  y=[a_0;a_1,a_2,\ldots].
\end{equation*}
From the construction of the simple CF we have that $a_0\in\mathbb{Z}$,
whereas $a_j\in\mathbb{N}$ for $j\ge 1$.  The recursion stops if
$\{y_j\}=0$ for some $j$, which only occurs when
$y\in\mathbb{Q}$; otherwise the CF is infinite (for an in-depth introduction to
continued fractions see~\cite{khinchin61:continued}).
Our goal in this section is to model probabilistically the~$a_j$
coefficients for $j\ge 1$. To this end we will assume that~$y$ is
drawn from a continuous r.v.~$Y$. Because~$\mathbb{Q}$ is of measure
zero then $\Pr(Y\in\mathbb{Q})=0$, and so we may assume that the CF
of~$y$ drawn from $Y$ is almost surely infinite ---and thus that the
$a_j$ coefficients are unique. In practical terms, this means that,
letting $Y_0=Y$, we can define the continuous r.v.'s
\begin{equation*}
  Y_{j}=\{Y_{j-1}\}^{-1}
\end{equation*}
with support $(1,\infty)$ for all $j\ge 1$. Therefore, the simple CF
coefficients we are interested in are modelled by the discrete r.v.'s
\begin{equation}
  A_j=\lfloor Y_j\rfloor\label{eq:aj_rvs}
\end{equation}
with support $\mathbb{N}$, for all $j\ge 1$.

\textbf{Additional notation.} In order to streamline the presentation
in this section we define the $k$-vector
\begin{equation*}
  \mathbf{A}_{k}=[A_1,\ldots,A_k]\label{eq:ak_cf}
\end{equation*}
comprising the first $k$ r.v.'s defined by~\eqref{eq:aj_rvs}, i.e. the
r.v.'s modelling the $k$ leading CF coefficients of $Y$. A realisation
of $\mathbf{A}_k$ is
$\mathbf{a}_k=[a_1,\ldots,a_k]\in\mathbb{N}^k$. Also, $\mathbf{e}_k$
denotes a unit $k$-vector with a one at the $k$-th position and zeros
everywhere else, i.e. $\mathbf{e}_k=[0,\ldots,0,1]$. A~vector symbol
placed within square brackets denotes a finite CF; for example,
$[\mathbf{a}_k]$ denotes $[a_1;a_2,\ldots,a_k]$. Observe that we can
write
$[\mathbf{a}_k]^{-1}=[0;\mathbf{a}_k]=[0;a_1,a_2,\ldots,a_k]$. Finally,
the subvector of consecutive entries of $\mathbf{a}_k$ between its
$m$-th entry and its last entry is denoted by
$\mathbf{a}_k^{m:k}=[a_m,\ldots,a_k]$. When $m>1$ the amount
$[\mathbf{a}_k^{m:k}]$ is called a remainder of
$[\mathbf{a}_k]$~\cite{khinchin61:continued}.

As a preliminary step, we will next  prove a lemma which we will
then use as a stepping stone in the derivation of the joint pmf
of~$\mathbf{A}_k$ in Theorem~\ref{thm:jpmfcf}.
\begin{lemma}\label{thm:lemma}
  The following two sets of inequalities hold for the $(j-1)$-th order
  convergent $[\mathbf{a}_j]=[a_1;a_2,\dots,a_j]$ of the infinite
  simple continued fraction $[a_1;a_2,a_3,\dots]$ with
  $a_i\in\mathbb{N}$:
  \begin{equation}
    \label{eq:ineq1}
    (-1)^{j-1}[\mathbf{a}_j]< (-1)^{j-1}[\mathbf{a}_j+\mathbf{e}_j]
  \end{equation}
  and
  \begin{equation}
    \label{eq:ineq2}
    (-1)^{j-1}[\mathbf{a}_{j-1}+\mathbf{e}_{j-1}]\le (-1)^{j-1}[\mathbf{a}_j]<(-1)^{j-1}[\mathbf{a}_{j+1}+\mathbf{e}_{j+1}],
    \rule[-1em]{0pt}{0pt}
  \end{equation}
  where the lower bound in~\eqref{eq:ineq2} requires $j>1$.
\end{lemma}
\begin{proof}
  Consider the function
  $\varphi(\mathbf{a}_j)=[\mathbf{a}_j]$. Taking $a_j$ momentarily to be
   a continuous variable, we can obtain the partial derivative
  of $\varphi(\mathbf{a}_j)$ with respect to $a_j$. For $j\ge 2$,
  applying the chain rule $j-1$ times yields
\begin{equation}\label{eq:deriv}
\frac{\partial\varphi(\mathbf{a}_j)}{\partial a_j}=(-1)^{j-1}\prod_{r=2}^j[\mathbf{a}_j^{r:j}]^{-2},
\end{equation}
whereas when $j=1$ we have that $d\varphi(\mathbf{a}_1)/d a_1=1$. As the product
indexed by $r$ is positive, the sign of~\eqref{eq:deriv} only depends
on~$(-1)^{j-1}$.
Consequently, if $j$ is odd then $\varphi(\mathbf{a}_j)$
is strictly increasing on $a_j$, and if~$j$ is even then
$\varphi(\mathbf{a}_j)$ is strictly decreasing on $a_j$. Thus,
when $j$ is odd $[\mathbf{a}_j]<[\mathbf{a}_j+\mathbf{e}_j]$ and when
$j$ is even $[\mathbf{a}_{j}]>[\mathbf{a}_j+\mathbf{e}_j]$, which
proves inequality~\eqref{eq:ineq1}.

Let us now prove the two inequalities in~\eqref{eq:ineq2} assuming
first that~$j$ is odd. Considering again~\eqref{eq:deriv}, the upper
bound can be obtained by seeing that
$[\mathbf{a}_j]<[\mathbf{a}_{j-1},a_j+\epsilon]$ for $\epsilon>0$, and
then choosing $\epsilon=1/(a_{j+1}+1)$. The lower bound, which
requires $j>1$, is due to
$[\mathbf{a}_j]\ge
[\mathbf{a}_{j-1},1]=[\mathbf{a}_{j-1}+\mathbf{e}_{j-1}]$. To conclude
the proof, when $j$ is even the two inequalities we have just
discussed are reversed due to the change of sign in~\eqref{eq:deriv}.
\end{proof}

\begin{remark}
  Lemma~\ref{thm:lemma} is closely related to the fact that the $j$-th
  order convergent of a continued fraction is smaller or larger than
  the continued fraction it approximates depending on the parity of
  $j$~\cite[Theorem 4]{khinchin61:continued}.
\end{remark}

\begin{theorem}[General distribution of the $k$ leading continued fraction coefficients]\label{thm:jpmfcf}
  For any continuous r.v. $Y$ represented by the simple continued
  fraction $Y=[A_0;A_1,A_2,\dots]$ it holds that
  \begin{equation}\label{eq:cf_joint_pmf}
    \Pr(\mathbf{A}_{k}=\mathbf{a}_k)=(-1)^{k}\big(F_{\{Y\}}([0;\mathbf{a}_{k}+\mathbf{e}_k])-F_{\{Y\}}([0;\mathbf{a}_k])\big),
  \end{equation}
  where $\mathbf{a}_k\in\mathbb{N}^k$.
\end{theorem}

\begin{proof}
  For all $j\ge 2$, if $a_m\le Y_m< a_m+1$ for all $m=1,\dots,j-1$ then
  we have that $Y_1=[\mathbf{a}_{j-1} , Y_j]$.  In these conditions it
  holds that
  $\{a_j\le Y_j < a_j+1\}=\{(-1)^{j-1}[\mathbf{a}_{j}]\le (-1)^{j-1}
  Y_1 <(-1)^{j-1}[\mathbf{a}_{j}+\mathbf{e}_j]\}$, where the reason
  for the alternating signs is~\eqref{eq:ineq1}. Therefore
  \begin{align}\label{eq:cf_joint_pmf2}
    \Pr(\mathbf{A}_{k}=\mathbf{a}_k)&=\Pr(A_1=a_1,\dots,A_k=a_k)\nonumber\\
                                    &=\Pr(\cap_{j=1}^k \{a_j\le Y_j<
                                      a_j+1\})\nonumber\\
                                    &=\Pr(\cap_{j=1}^k
                                        \{(-1)^{j-1}[\mathbf{a}_{j}]\le
                                        (-1)^{j-1} Y_1<(-1)^{j-1}[\mathbf{a}_{j}+\mathbf{e}_j]\}).
  \end{align}
  From~\eqref{eq:ineq2} we have that the
  lower bounds on~$Y_1$ in~\eqref{eq:cf_joint_pmf2} are related as
  \begin{equation}\label{eq:shrinklow}
    [\mathbf{a}_1]<[\mathbf{a}_2+\mathbf{e}_2]\le
    [\mathbf{a}_3]<[\mathbf{a}_4+\mathbf{e}_4]\le [\mathbf{a}_5]\cdots
  \end{equation}
  whereas the upper bounds on $Y_1$ in the same expression are related~as
  \begin{equation}\label{eq:shrinkhigh}
    [\mathbf{a}_1+\mathbf{e}_1]\ge [\mathbf{a}_2]>
    [\mathbf{a}_3+\mathbf{e}_3]\ge[\mathbf{a}_4]> [\mathbf{a}_5+\mathbf{e}_5]\cdots
  \end{equation}
  Hence, except for possible equality constraints (which are anyway
  immaterial in probability computations with continuous random
  variables), the intersection of the $k$ events
  in~\eqref{eq:cf_joint_pmf2} equals the $k$-th event, and thus
  \begin{align}\label{eq:cf_joint_pmf3}
    \Pr(\mathbf{A}_{k}=\mathbf{a}_k)&=\Pr(
                                        (-1)^{k-1}[\mathbf{a}_{k}]<
                                        (-1)^{k-1} Y_1<(-1)^{k-1}[\mathbf{a}_{k}+\mathbf{e}_k])\nonumber\\
    &=(-1)^{k-1} \big(F_{Y_1}([\mathbf{a}_{k}+\mathbf{e}_k])-F_{Y_1}([\mathbf{a}_k])\big).
  \end{align}
Finally, using
\begin{equation*}\label{eq:FY1}
  F_{Y_1}(y)=\Pr(Y_1\le y)=\Pr(\{Y\}\ge y^{-1})=1-F_{\{Y\}}(y^{-1})
\end{equation*}
in~\eqref{eq:cf_joint_pmf3} we get~\eqref{eq:cf_joint_pmf}.
\end{proof}

\begin{remark}\label{rem:th2}
  Observe that if we choose $Y=\log_b X$, then
  both~\eqref{eq:cf_joint_pmf} and~\eqref{eq:pmfA} depend solely on
  the same variable $\{Y\}$, which is the reason why we have used the
  notation $Y$ rather than~$X$ in this section.  With this choice of
  $Y$, the general expression~\eqref{eq:cf_joint_pmf} models the~$k$
  leading CF coefficients of $\log_b X$ (with the
  exception of $A_0$), and becomes analogous to the general
  expression~\eqref{eq:pmfA} that models the~$k$ most significant
  $b$-ary digits of~$X$. The reason why we have left $A_0$ out of the
  joint distribution~\eqref{eq:cf_joint_pmf} is because, unlike the
  rest of variables (i.e. $A_j$ for $j\ge 1$), it cannot be put as a
  sole function of~$Y_1$. Moreover, it is not possible to model $A_0$ in
  one important practical scenario ---see
  Section~\ref{sec:lead-cf-coeff}. %

We can also verify that the joint pmf~\eqref{eq:cf_joint_pmf} adds up
to one over its support, namely~$\mathbb{N}^k$. Let us first add the
joint pmf of $\mathbf{A}_k$ over $a_k\in\mathbb{N}$ assuming $k>1$. As
this infinite sum is a telescoping series, in the computation of the
partial sum $S_k^{(n)}=\sum_{a_k=1}^n\Pr(\mathbf{A}_{k}=\mathbf{a}_k)$
all consecutive terms but two are cancelled, and so
\begin{align*}
  S_k^{(n)}%
           &=(-1)^k\big(F_{\{Y\}}([0;\mathbf{a}_{k-1},n+1])-F_{\{Y\}}([0;\mathbf{a}_{k-1},1])\big).
\end{align*}
Now, as $\lim_{n\to\infty} [0;\mathbf{a}_{k-1},n+1]=[0;\mathbf{a}_{k-1}]$
and
  $[0;\mathbf{a}_{k-1},1]=[0;\mathbf{a}_{k-1}+\mathbf{e}_{k-1}]$, 
we then have that
\begin{align*}
  \lim_{n\to \infty}
  S_k^{(n)}&=(-1)^{k-1}\big(F_{\{Y\}}([0;\mathbf{a}_{k-1}+\mathbf{e}_{k-1}])-F_{\{Y\}}([0;\mathbf{a}_{k-1}])\big)\\
  &=\Pr(\mathbf{A}_{k-1}=\mathbf{a}_{k-1}).
\end{align*}
The continuity of the cdf~$F_{\{Y\}}(y)$ allows
writing $\lim_{n\to\infty}F_{\{Y\}}(g(n))=F_{\{Y\}}(\lim_{n\to\infty} g(n))$,
which justifies the limit
above. %
In view of this result, it only remains to verify that
the pmf of $\mathbf{A}_1=A_1$ adds up to one. The partial sum up to $n$ is
\begin{equation*}
  S_1^{(n)}=F_{\{Y\}}(1)-  F_{\{Y\}}(1/(n+1)),
\end{equation*}
and therefore $\lim_{n\to\infty} S_1^{(n)}=1$ for the same reason that
makes~\eqref{eq:pmfA_verification} equal to one. Incidentally, observe
that it would have been rather more difficult to verify the fact
that~\eqref{eq:cf_joint_pmf} adds up to one by summing out the random
variables in~$\mathbf{A}_k$ in an order different than the decreasing
order $A_k,A_{k-1},\ldots,A_1$ that we have used above.
\end{remark}

\subsection{Distribution of the $j$-th CF Coefficient}
\label{sec:distribution-j-th}
Just like in Section~\ref{sec:prob-distr-k}, we can marginalise the
joint pmf of $\mathbf{A}_j$ to obtain the distribution of the $j$-th
CF coefficient $A_j$ of $Y$. Although we already know that
$A_1=\mathbf{A}_1$, the main obstacle to explicitly getting the
distribution of $A_j$ for $j>1$ is that in this case marginalisation
involves $j-1$ infinite series, rather than a single finite sum as
in~\eqref{eq:pmfAj}. In general, it is difficult to carry out the
required summations in closed form. Moreover, the order of evaluation
of these series may influence the feasibility of the computation,
which is connected to the comment in the very last sentence of the
previous paragraph.

However, under the sole assumption that~$\{Y\}$ is a continuous
r.v. with support $[0,1)$, the Gauss-Kuz'min theorem furnishes the
general asymptotic distribution of~$A_j$~\cite[Theorem
34]{khinchin61:continued}:
\begin{equation}
  \label{eq:gauss-kuzmin}
  \lim_{j\to\infty} \Pr(A_j=a)= \log_2\Big(1+\frac{1}{a(a+2)}\Big).
\end{equation}

\begin{remark}
  Observe that Theorem~\ref{thm:asympt-aj}, which gives the general
  asymptotic behaviour of $A_{[j]}$ (the $j$-th most significant
  $b$-ary digit), is the near exact counterpart of the Gauss-Kuz'min
  theorem~\eqref{eq:gauss-kuzmin}, which gives the general asymptotic
  behaviour of $A_j$ (the $j$-th continued fraction coefficient).  The
  only essential difference is the requirement that the support of
  $\{Y\}$ be precisely $[0,1)$
  in~\eqref{eq:gauss-kuzmin}~\cite[Theorem
  33]{khinchin61:continued}, whereas this condition is not required
  in~\eqref{eq:asympt-aj} ---i.e. the support of $\{Y\}$ may be a
  subset of $[0,1)$ in
  Theorem~\ref{thm:asympt-aj}.  %
\end{remark}

\section{Particular Cases}
\label{sec:particular-cases}

In this section we will particularise the general expressions in
Sections~\ref{sec:gen-prob-msd} and~\ref{sec:gen-prob-cf-coeff} for two
especially relevant distributions of $X$. As it is clear
from~\eqref{eq:pmfA} and~\eqref{eq:cf_joint_pmf}, we just need the cdf
$F_{\{Y\}}(y)$ of the r.v. $\{Y\}=\{\log_b X\}$ in order to achieve
our goal.

\subsection{Benford Variables}
\label{sec:benford}
We consider in this section a r.v.~$X$ for which $\{Y\}\sim
U(0,1)$. We call such a r.v. a \textit{Benford variable}, although we
must note that some authors call it a \textit{strong} Benford variable
instead. At any rate, this is the archetypal case in which a model of
the $k$ most significant $b$-ary digits has been widely used and
discussed ---i.e. Benford's law~\cite{benford1938}. The cdf
of~$\{Y\}$ for a Benford variable $X$ is simply
\begin{equation}
  F_{\{Y\}}(y)=y\label{eq:cdfYbenford}
\end{equation}
for $y\in[0,1)$.

\subsubsection{Most Significant $b$-ary Digits of $X$}
\label{sec:most-sign-digits}
For a Benford variable, applying~\eqref{eq:cdfYbenford} to~\eqref{eq:pmfA} yields
\begin{align}
  \Pr(A_{(k)}=a)%
  &=\log_b\left(1+\frac{1}{a}\right),\label{eq:benfordk}
\end{align}
which is the well-known Benford distribution for the $k$ most
significant $b$-ary digits. This distribution has almost always been
expressed in previous works as the joint pmf of
$A_{[1]},\dots,A_{[k]}$ rather than as the pmf of $k$-th integer
significand $A_{(k)}$ (see for
example~\cite{berger11:_basic_th_benford}). As evinced in
Theorems~\ref{thm:msd} and~\ref{thm:asympt-aj}, and as it will become
clear in the remainder of this section, working with the $k$-th
integer significand is not just an aesthetical notation choice
---although it does make for simpler expressions.

Let us obtain next the pmf of $A_{[j]}$ when $j\ge 2$ (i.e. the
distribution of the $j$-th most significant $b$-ary digit).
From~\eqref{eq:pmfAj} and~\eqref{eq:benfordk} we have that
\begin{align}
  \Pr(A_{[j]}\!=\!a)&= \sum_{r\in\mathcal{A}_{(j-1)}}\log_b\!\left(1+\frac{1}{rb+a}\right)\label{eq:margbendford}\\
                    &=\log_b\!\Bigg(\prod_{r\in\mathcal{A}_{(j-1)}}\frac{(a+1)b^{-1}+r}{ab^{-1}+r}\Bigg)\label{eq:rising}\\
                    &=\log_b\left(\frac{\Gamma\left((a+1)b^{-1}+b^{j-1}\right)\Gamma\left(a b^{-1}+b^{j-2}\right)}{\Gamma\left((a+1)b^{-1}+b^{j-2}\right)\Gamma\left(a b^{-1}+b^{j-1}\right)}\right). \label{eq:benfordk_j}
\end{align} 
The last equality is due to the fact that the argument of the
logarithm in~\eqref{eq:rising} can be expressed as a fraction whose
numerator and denominator are the rising factorial powers
$((a+1)b^{-1}+b^{j-2})^{\overline{\vert\mathcal{A}_{(j-1)}\vert}}$ and
$(a b^{-1}+b^{j-2})^{\overline{\vert\mathcal{A}_{(j-1)}\vert}}$, respectively.

We can also explicitly restate the general result in Theorem~\ref{thm:asympt-aj} for a Benford
variable by relying on~\eqref{eq:benfordk_j}. Invoking the continuity of
the logarithm and
using %
$\lim_{z\to\infty}z^{w-v}\Gamma(v+z)/\Gamma(w+z)=1$~\cite{abramowitz72:_handbook} %
in~\eqref{eq:benfordk_j} twice ---with $z=b^{j-1}$ and $z=b^{j-2}$,
respectively--- yields %
  \begin{equation*}\label{eq:limpaj}
    \lim_{j\to\infty}\Pr(A_{[j]}=a)=\log_b\lim_{j\to\infty}
    \frac{b^{(j-1)b^{-1}}}{b^{(j-2)b^{-1}}} =b^{-1},
  \end{equation*}
  a fact that was originally pointed out by Benford~\cite{benford1938}
  through the marginalisation of~\eqref{eq:benfordk}.

\begin{remark}
  The closed-form analytic expression~\eqref{eq:benfordk_j} for the
  pmf of $A_{[j]}$ deserves some comments, as it appears that it was
  never given in studies of Benford's distribution previous
   to~\cite{balado21:_benford}: only the equivalent
  of~\eqref{eq:margbendford} was previously published.
  This is another sensible reason for working with the pmf of the
  $j$-th integer significand variable $A_{(j)}$ instead of the joint
  pmf of $A_{[1]},\dots,A_{[j]}$.  The former approach makes the
  obtention of closed-form distributions for $A_{[j]}$ more feasible:
  if we use $A_{(j)}$ we just have to evaluate one single sum
  [i.e. \eqref{eq:pmfAj}], whereas if we use $A_{[1]},\dots,A_{[j]}$
  we have to evaluate $j-1$ separate sums ---which obscures the
  result.  This appears to be the reason why previous works never
  produced~\eqref{eq:benfordk_j}.
\end{remark}

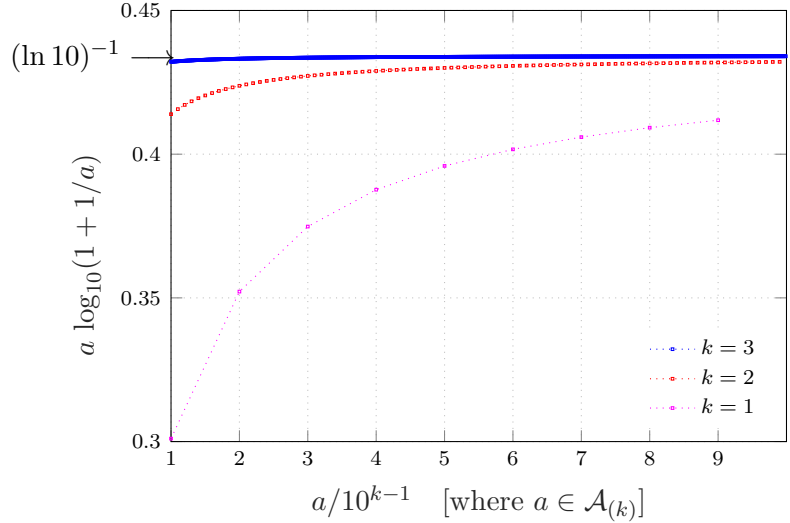
\begin{figure}[t!]
    \setlength\fwidth{.6\textwidth}
    \setlength\fheight{.4\textwidth} 
    \pgfplotsset{every tick label/.append style={font=\scriptsize}}
    \centering
     \input{fig2/si.tex}
    \caption{Illustration of the asymptotic sum-invariance property of a Benford variable
      for~$b=10$.}
    \label{fig:sum_invariance}
\end{figure}

\textbf{Asymptotic sum-invariance property.} A further advantage of
working with the $k$-th integer significand $A_{(k)}$ is that it
allows for an uncomplicated statement and proof of the asymptotic
sum-invariance property of a Benford
variable~\cite{berger11:_basic_th_benford,nigrini12:_benford}. In the
literature, this property has simply been called the ``sum-invariance
property''. Here we prefer to stress the fact that its validity is
only asymptotic when one considers a finite number $k$ of most
significant digits of $X$ ---i.e. the $k$-th integer significand
$A_{(k)}$--- which in fact originally motivated the empirical
definition of the sum-invariance property by
Nigrini~\cite{nigrini92:_phd}.

\begin{theorem}[Asymptotic sum-invariance property]\label{thm:sum-inv}
  If $X$ is a Benford variable, then it holds that
  \begin{equation}
    \label{eq:sum-invariance-prop}
    \lim_{\stackrel{k\to\infty,}{a\in\mathcal{A}_{(k)}}} a\Pr(A_{(k)}=a)=(\ln b)^{-1}.
  \end{equation}
\end{theorem}
\begin{proof}
  We just need to see that
  $\lim_{k\to\infty, a\in\mathcal{A}_{(k)}}
  a\log_b\left(1+\frac{1}{a}\right)=\lim_{v\to\infty}v\log_b\left(1+\frac{1}{v}\right)$
  due to~\eqref{eq:support}. The proof is completed by using either
  L'Hôpital's theorem, or the continuity of the logarithmic function
  and the definition of Euler's number as
  $\lim_{v\to\infty}(1+1/v)^v$.
\end{proof}

\begin{remark}
  Informally, Theorem~\ref{thm:sum-inv} tells us that in a large set
  of outcomes from a Benford variable the sum of
  $a\in\mathcal{A}_{(k)}$ over all those outcomes whose $k$-th integer
  significand is equal to~$a$ is roughly invariant over $a$ when $k$
  is large (i.e. the sum-invariance property as defined by Nigrini).
  Convergence speed to the limit~\eqref{eq:sum-invariance-prop} is
  exponential on~$k$ ---the faster the larger~$b$
  is. %
  Figure~\ref{fig:sum_invariance} shows that the sum-invariance
  property holds approximately when $k=3$ already, for
  $b=10$. Theorem~\ref{thm:sum-inv} also implies that
  \begin{equation*}
    \lim_{k\to\infty}
    \frac{\Expc(A_{(k)})}{\vert\mathcal{A}_{(k)}\vert}=(\ln b)^{-1}.\label{eq:lim_normalised_EAk}
  \end{equation*}
  The corresponding approximation
  $\Expc(A_{(k)})\approx (b^k-b^{k-1}) (\ln b)^{-1}$ improves with
  $k$, but it never achieves strict equality for finite $k$ ---in
  fact, $\Expc(A_{(k)})<(b^k-b^{k-1}) (\ln b)^{-1}$ for all $k$.
\end{remark}

\subsubsection{Leading CF Coefficients of\hspace{.13cm}$\log_b X$}
\label{sec:lead-cf-coeff}
For a Benford variable the application of~\eqref{eq:cdfYbenford}
to~\eqref{eq:cf_joint_pmf} yields
\begin{equation}\label{eq:cf_joint_pmf_benford}
  \Pr(\mathbf{A}_{k}=\mathbf{a}_k)=(-1)^{k}\big([0;\mathbf{a}_{k}+\mathbf{e}_k]-[0;\mathbf{a}_k]\big).
\end{equation}
According to our discussion in Remark~\ref{rem:th2}, this distribution
of the $k$ leading CF coefficients of $\log_b X$ is the counterpart of
Benford's distribution of the $k$ most significant $b$-ary digits of
$X$. Therefore~\eqref{eq:cf_joint_pmf_benford} can be seen as
Benford's law for continued fractions. In particular, any real dataset
that complies with~\eqref{eq:benfordk} will also comply
with~\eqref{eq:cf_joint_pmf_benford}.

By transforming the subtraction of fractions into a single
fraction,~\eqref{eq:cf_joint_pmf_benford} can also be written as the
inverse of the product of $[\mathbf{a}_k]$,
$[\mathbf{a}_k+\mathbf{e}_k]$ and all of their remainders, i.e.
\begin{equation}\label{eq:cf_identity}
  \Pr(\mathbf{A}_{k}=\mathbf{a}_k)  =\prod_{j=1}^k
      \big[0;\mathbf{a}_k^{j:k}\big]\big[0;\mathbf{a}_k^{j:k}+\mathbf{e}_{k-j+1}\big],
\end{equation}
which, apart from showing at a glance
that~\eqref{eq:cf_joint_pmf_benford} cannot be negative, may be more
suitable for log-likelihood
computations. %
The equivalent of expression~\eqref{eq:cf_identity} was previously
given by Miller and Takloo-Bighash~\cite[Lemma
10.1.8]{miller06:_invitation} in their exploration of the distribution
of CF coefficients---called digits by these authors. However, unlike
our result above, the version of~\eqref{eq:cf_identity} given by
Miller and Takloo-Bighash is not explicit, as it is presented in terms
of CF convergents. Therefore,~\eqref{eq:cf_identity}
or~\eqref{eq:cf_joint_pmf_benford} are clearly more useful when it
comes to practical applications ---furthermore, we show in
Section~\ref{sec:empirical-tests} the empirical accuracy of
our expressions using both synthetic and real data, something that was
not attempted by Miller and Takloo-Bighash. Lastly,
Blachman also gave the following
explicit approximation for uniform $\{Y\}$~\cite[equation (9)]{blachman84}:
\begin{equation}
  \label{eq:blachman}
  \Pr(\mathbf{A}_{k}=\mathbf{a}_k)\approx \left|\log_2\left(\frac{1+[0; \mathbf{a}_k]}{1+[0; \mathbf{a}_k+\mathbf{e}_k]}\right)\right|.
\end{equation}
The reader should be cautioned that this expression was given
in~\cite{blachman84} with an equal sign, although the author
unequivocally produced it as an approximation. Using
$\ln(1+z)\approx z$, which is accurate for $|z|\ll 1$, we can see
that~\eqref{eq:blachman} is roughly off by a factor of $(\ln 2)^{-1}$
with respect to the exact expression~\eqref{eq:cf_joint_pmf_benford}.

Let us now look at the marginals, that is to say, the distributions of
individual $A_j$ coefficients. When $k=1$ expression~\eqref{eq:cf_joint_pmf_benford} gives the distribution
of~$A_1=\mathbf{A}_1$ straightaway:
\begin{equation}\label{eq:a1_benford}
  \Pr(A_1=a)=a^{-1}-(a+1)^{-1}. %
\end{equation}
This pmf, also previously given by Miller and
Takloo-Bighash~\cite[page 232]{miller06:_invitation}, can be rewritten
as $\Pr(A_1=a)=a^{-1}(a+1)^{-1}$, which is the form
that~\eqref{eq:cf_identity} takes in this particular
case. Incidentally, observe that $\Expc(A_1)=\infty$ because of the
divergence of the harmonic series. It is also instructive to
particularise Blachman's approximation~\eqref{eq:blachman} for
$A_1$~\cite[equation (10)]{blachman84}: this renders the asymptotic
Gauss-Kuz'min law~\eqref{eq:gauss-kuzmin} instead of the exact
pmf~\eqref{eq:a1_benford}.

Recalling our discussion at the start of
Section~\ref{sec:distribution-j-th}, the Benford case is probably
unusual in the fact that we can also obtain the distribution of $A_2$
in closed form by marginalising~\eqref{eq:cf_joint_pmf_benford} for
$k=2$. Summing
$\Pr(\mathbf{A}_2=\mathbf{a}_2)=1/(a_1+(a_2+1)^{-1})-1/(a_1+a_2^{-1})$
over $a_1\in\mathbb{N}$, and using the digamma function defined as
$\psi(1+z)=-\gamma+\sum_{n=1}^\infty
z/(n(n+z))$~\cite{abramowitz72:_handbook} ---which is applicable
because the range of validity $z\notin \mathbb{Z}^-$ of this
definition always holds here---
one finds that
\begin{equation}\label{eq:a2_benford}
  \Pr(A_2=a)=\psi\left(1+a^{-1}\right)-\psi\left(1+(1+a)^{-1}\right).
\end{equation} %

It does not seem possible to obtain a closed-form exact expression for
the distribution of a single CF coefficient~$A_j$ when~$j>2$ in the
Benford case.  However it is possible to explicitly produce the
Gauss-Kuz'min law~\eqref{eq:gauss-kuzmin} by pursuing an approximation of $\Pr(A_j=a_j)$ for
all $j\ge 2$. To see this, consider first the sum
\begin{equation}\label{eq:sum}
  \sum_{x=1}^\infty  \left(\frac{1}{x+b}-\frac{1}{x+c}\right)=\psi(1+c)-\psi(1+b)
\end{equation}
for some $b,c>0$, which is just a generalisation
of~\eqref{eq:a2_benford}, and its integral approximation
\begin{equation}\label{eq:int_approx}
  \int_1^\infty \left(\frac{1}{x+b}-\frac{1}{x+c}\right) dx=
  \ln(1+c)-\ln(1+b),
\end{equation}
which attests to the intimate connection between the digamma function
and the natural logarithm~\cite[see Exercise 8.2.20 and equation
(8.51)]{arfken05:_mathematical}. Now, in the marginalisation that
leads to $\Pr(A_j=a_j)$ the summation on $a_1$ is of the
form~\eqref{eq:sum}, and so we may approximate it by the
integral~\eqref{eq:int_approx}:
\begin{align}
  \Pr(A_{j}=a_j)&=(-1)^j\sum_{a_{j-1}=1}^\infty\dots\sum_{a_{1}=1}^\infty\big([0;\mathbf{a}_{j}+\mathbf{e}_j]-[0;\mathbf{a}_j]\big)  \label{eq:aj_marginalisation}\\
                &\approx (-1)^{j}\sum_{a_{j-1}=1}^\infty\dots\sum_{a_{2}=1}^\infty\int_{1}^\infty\big([0;\mathbf{a}_{j}+\mathbf{e}_j]-[0;\mathbf{a}_j]\big)\, da_1 \nonumber\\
  &= (-1)^{j}\sum_{a_{j-1}=1}^\infty\dots\sum_{a_2=1}^\infty(-1)\big(\ln(1+[0;\mathbf{a}^{2:j}_{j}+\mathbf{e}_{j-1}])-\ln(1+[0;\mathbf{a}^{2:j}_j])\big).\label{eq:ln}
\end{align}
As $[0;\mathbf{a}^{i:j}_j]\ll 1$ nearly always for
$\mathbf{a}^{i:j}_j=[a_i,\cdots,a_j]\in \mathbb{N}^{j-i+1}$, then we may use
$\ln(1+z)\approx z$ in~\eqref{eq:ln} to obtain
\begin{align}\label{eq:ln1p}
    \Pr(A_{j}=a_j)&\approx
    (-1)^{j+1}\sum_{a_{j-1}=1}^\infty\dots\sum_{a_2=1}^\infty\big([0;\mathbf{a}^{2:j}_{j}+\mathbf{e}_{j-1}])-[0;\mathbf{a}^{2:j}_j]\big). 
\end{align}
Remarkably,~\eqref{eq:ln1p} has the exact same form
as~\eqref{eq:aj_marginalisation} ---but with one less infinite
summation. Therefore we can keep sequentially applying the same
approximation procedure described above to the summations on $a_2$,
$a_3$,\ldots,$a_{j-1}$. In the final summation on $a_{j-1}$ we do not
need the approximation in~\eqref{eq:ln1p} anymore, and thus in the
last step we have that
\begin{align}
  \Pr(A_{j}=a_j) &\approx (-1)^{2j-2}\int_{1}^\infty\big([0;\mathbf{a}^{j-1:j}_{j}+\mathbf{e}_{2}])-[0;\mathbf{a}^{j-1:j}_j]\big) \,da_{j-1}\nonumber\\
                &=(-1)^{2j-2}\int_{1}^\infty \left(\cfrac{1}{a_{j-1}+\cfrac{1}{a_j+1}}-\cfrac{1}{a_{j-1}+\cfrac{1}{a_j}}\right)da_{j-1}\nonumber\\
                &=(-1)^{2j-1}\left(\ln\left(1+\frac{1}{a_{j}+1}\right)-\ln\left(1+\frac{1}{a_{j}}\right)\right)\label{eq:gateway2a1}\\
  &=\ln\frac{(a_j+1)^2}{a_j(a_j+2)}.\label{eq:unnormalised_gk}
\end{align}
Due to the successive approximations~\eqref{eq:unnormalised_gk} is not
necessarily a pmf, and so we need to normalise it. The normalisation
factor is
\begin{equation*}
  \sum_{a_j=1}^\infty
  \ln\frac{(a_j+1)^2}{a_j(a_j+2)}=\ln\prod_{a_j=1}^\infty\frac{(a_j+1)^2}{a_j(a_j+2)}=\ln\frac{2
\cdot \cancel{2}}{1\cdot \cancel{3}} \cdot\frac{\cancel{3}\cdot \cancel{3}}{\cancel{2}\cdot \cancel{4}}\cdot\frac{\cancel{4}\cdot
\cancel{4}}{\cancel{3}\cdot \cancel{5}}\cdots=\ln 2.
\end{equation*}
Applying this factor to~\eqref{eq:unnormalised_gk} we finally obtain
\begin{equation}\label{eq:approx_aj_gk}
  \Pr(A_{j}=a_j)\approx \log_2\frac{(a_j+1)^2}{a_j(a_j+2)}=\log_2\left(1+\frac{1}{a_j(a_j+2)}\right),
\end{equation}
for $j\ge 2$.%

\begin{remark}
  Like in the verification of the general joint
  pmf~\eqref{eq:cf_joint_pmf} in Remark~\ref{rem:th2}, the right order
  of evaluation of the marginalisation sums is again key for us to be
  able to produce approximation~\eqref{eq:approx_aj_gk}.  Also, had we
  used $\ln (1+z)\approx z$ one last time in~\eqref{eq:gateway2a1}
  then we would have arrived at~\eqref{eq:a1_benford} instead of at
  the Gauss-Kuz'min law as the final approximation.  This shows that the pmf
  of the first CF coefficient and the asymptotic law are close
  already, which was also mentioned by Miller and
  Takloo-Bighash~\cite[Exercise 10.1.1]{miller06:_invitation}. Since
  the convergence of the distribution of $A_j$ to the asymptotic
  distribution is exponentially fast on $j$, it is unsurprising that
  the pmf~\eqref{eq:a2_benford} of the second CF coefficient turns out
  to be even closer to~\eqref{eq:gauss-kuzmin}, as suggested
  by~\eqref{eq:approx_aj_gk} ---see empirical validation in
  Section~\ref{sec:empirical-tests}.  %

  Although beyond the goals of this paper, it should be possible to
  refine the approximation procedure that we have given to get
  \eqref{eq:unnormalised_gk} in order to to explicitly obtain the
  exponential rate of convergence to the Gauss-Kuz'min law, by
  exploiting the expansion of the digamma function in terms of the
  natural logarithm and an error term series~\cite[equation
  (8.51)]{arfken05:_mathematical}.

  Finally, see  that although~$A_0$ is not included in the joint
  pmf~\eqref{eq:cf_joint_pmf_benford}, this variable cannot be modelled
  anyway when the only information that we have about $X$ is its
  ``Benfordness''.

\end{remark}

\subsubsection{Benford Variables and the Asymptotics of the General  Analysis}
\label{sec:benf-vari-spec}
To conclude Section~\ref{sec:benford} we examine the role played by
the particular analysis for a Benford variable
[i.e. \eqref{eq:benfordk} and~\eqref{eq:cf_joint_pmf_benford}] in the
general analysis [i.e. Theorems~\ref{thm:msd} and~\ref{thm:jpmfcf}]
when $k$ is large. Let us start by looking at the asymptotics
of~\eqref{eq:pmfA}. For any $\epsilon>0$ there exists~$k_{\min}$ such
that $\log_b(1+a^{-1})<\epsilon$ for all $k\ge k_{\min}$ and
$a\in\mathcal{A}_{(k)}$. Explicitly, this minimum index is
$k_{\min}=\lceil -\log_b(b^{\epsilon}-1)+1\rceil$. This inequality and
the continuity of $F_{\{Y\}}(y)$ allow us to
approximate~\eqref{eq:pmfA} for large~$k$ using the pdf of $\{Y\}$ as
\begin{equation}
  \label{eq:ak_pmf_asympt}
  \Pr(A_{(k)}=a)\approx f_{\{Y\}}(\log_b a -k+1)\,\log_b\left(1+\frac{1}{a}\right).
\end{equation}
We now turn our attention to the asymptotics
of~\eqref{eq:cf_joint_pmf}. Similarly as above, for any $\epsilon>0$
\eqref{eq:shrinklow} and~\eqref{eq:shrinkhigh} guarantee that there
exists~$k_{\min}$ such that
$(-1)^k([0;\mathbf{a}_k+\mathbf{e}_k]-[0;\mathbf{a}_k])<\epsilon$ for
all $k\ge k_\text{min}$. Invoking again the continuity of
$F_{\{Y\}}(y)$ we can approximate~\eqref{eq:cf_joint_pmf} for
large~$k$ using again the pdf of $\{Y\}$ as
\begin{equation}
  \label{eq:cf_joint_pmf_asympt}
  \Pr(\mathbf{A}_{k}=\mathbf{a}_k)\approx f_{\{Y\}}([0;\mathbf{a}_k])\,(-1)^{k}\big([0;\mathbf{a}_{k}+\mathbf{e}_k]-[0;\mathbf{a}_k]\big).
\end{equation}
The key point that we wish to make here is that the Benford
expressions~\eqref{eq:benfordk} and~\eqref{eq:cf_joint_pmf_benford}
appear as factors in the general asymptotic
approximations~\eqref{eq:ak_pmf_asympt}
and~\eqref{eq:cf_joint_pmf_asympt}, respectively, which illustrates
the special place that Benford variables take in the modelling of
significant digits and leading continued fraction coefficients. Of
course, for Benford~$X$ the pdf of~$\{Y\}$ is $f_{\{Y\}}(y)=1$ for
$y\in[0,1)$, and so in this case
approximations~\eqref{eq:ak_pmf_asympt}
and~\eqref{eq:cf_joint_pmf_asympt} coincide with their exact
counterparts.

\subsection{Pareto Variables}
\label{sec:pareto}

In this section we let~$X$ be a Pareto r.v. with minimum
value $x_\text{m}$ and shape parameter~$s$, whose pdf is
\begin{equation*}\label{eq:pareto}
  f_X(x)= s\, x_\text{m}^s\,x^{-(s+1)},\quad 0<x_\text{m}\le x,\; s>0.
\end{equation*}
The main motivation for considering the Pareto distribution is its
pervasiveness in natural phenomena, which is reflected in the fact
that Pareto variables are able to model a wealth of scale-invariant
datasets. According to Nair et al.~\cite{nair22:_fundam_heavy_tails}
heavy-tailed distributions are just as prominent as the Gaussian
distribution, if not more. This is a consequence of the Central Limit
Theorem (CLT) \textit{not} yielding Gaussian distributions ---but
heavy-tailed ones--- in common scenarios where the variance of the
random variables being added is infinite (or does not
exist). Furthermore, heavy-tailed distributions appear when the CLT is
applied to the logarithm of variables emerging from multiplicative
processes. In this context, the relevance of the Pareto distribution
owes to the fact that the tails of many heavy-tailed distributions
follow the Pareto law. Additionally, the Pareto distribution is the
only one that fulfils exactly the relaxed scale-invariance criterion
\begin{equation}\label{sec:scale_invariance_2}
  f_X(x)= \alpha^{s+1} f_X(\alpha\, x) %
\end{equation}
for any scaling factor $\alpha>0$, where $s>0$.

Let us firstly obtain the cdf of $\{Y\}$ in this case. The cdf of a
Pareto r.v. $X$ is $F_X(x)=1-x_\text{m}^s x^{-s}$ for
$x\ge x_\text{m}$, and thus the cdf of $Y=\log_b X$ is
$F_Y(y)=F_X(b^y)=1-x_\text{m}^sb^{-s y}$ for $y\ge \log_b
x_\text{m}$. Letting
\begin{equation*}
  \rho=\{\log_b x_\text{m}\}  \label{eq:frac_logbxm}
\end{equation*}
and using~\eqref{eq:cdffrcy}, we have that the cdf of $\{Y\}$ for a
Pareto r.v. $X$ is
\begin{equation}
  \label{eq:cdf_frac_logbx}
  F_{\{Y\}}(y)= b^{s(\rho-1)}\, \frac{1-b^{-sy}}{1-b^{-s\phantom{y}}}+u\big(y-\rho\big)\left(1-b^{-s(y-\rho)}\right)
\end{equation}
for $y\in[0,1)$, where $u(\cdot)$ is the unit-step function.

\begin{remark}\label{rem:pareto_benford}
  By application of l'Hôpital's rule, it can be verified
  that~\eqref{eq:cdf_frac_logbx} tends to~\eqref{eq:cdfYbenford} as
  $s\!\downarrow\! 0$, and so a Pareto variable becomes asymptotically
  Benford as its shape parameter $s$ vanishes ---for any value of
  $\rho$. Because~\eqref{eq:pmfA} and~\eqref{eq:cf_joint_pmf} only
  depend on $\{Y\}$, the distributions that we will produce in this
  section generalise their counterparts in the previous section
  [i.e. \eqref{eq:benfordk},~\eqref{eq:benfordk_j}
  and~\eqref{eq:cf_joint_pmf_benford} are particular cases
  of~\eqref{eq:paretok_general}, \eqref{eq::aj_pareto_general}
  and~\eqref{eq:cf_joint_pmf_pareto}, respectively, when
  $s\!\downarrow\! 0$]. The fact that Benford variables can appear as
  a particular case of Pareto variables is a likely reason for the
  sporadic emergence of Benford's distribution~\eqref{eq:benfordk} in
  scale-invariant scenarios. Finally, observe that, asymptotically
  as~$s\!\downarrow\! 0$, the relaxed scaled invariance
  property~\eqref{sec:scale_invariance_2} becomes \textit{strict},
  i.e. $f_X(x)=\alpha\, f_X(\alpha\, x)$. Strict scale invariance is
  in turn a property that drives the appearance of Benford's
  distribution~\cite{balado21:_benford}.

  An interesting line of research beyond the scope of this paper would
  entail pursuing analytical insights about the probability
  distribution of the $s$ parameter itself in scale-invariant
  scenarios. If this distribution could be found, perhaps under
  constraints yet to be specified, it would determine the frequency of
  emergence of Benford variables in those scenarios. In any case, it
  can be empirically verified that scale-invariant datasets are far
  more often Paretian rather than just Benfordian (see some examples
  in Figure~\ref{fig:msd_real_datasets}). Thus, the expressions that
  we will give in this section may have wider practical applicability
  than the ones in Section~\ref{sec:benford} in the context of
  scale-invariant datasets ---with the caveat that two parameters ($s$
  and $\rho$ or $x_\text{m}$) must be estimated when using the Pareto
  distribution
  results. %
\end{remark}

\subsubsection{Most Significant $b$-ary Digits of $X$}
\label{sec:most-sign-digits-1}
Combining~\eqref{eq:pmfA} and~\eqref{eq:cdf_frac_logbx}, and
letting
\begin{equation*}
  \xi=\rho+k-1\label{eq:xi}
\end{equation*}
yields the Paretian generalisation of~\eqref{eq:benfordk}:
\begin{align}
  \Pr(A_{(k)}=a)&= \frac{b^{s(\xi-1)}}{1-b^{-s}}\,\big(a^{-s}-(a+1)^{-s}\big)\nonumber\\
                &+u\big(a+1-b^{\xi}\big)\big(1-b^{s\,\xi}(a+1)^{-s}\big)%
                -u\big(a-b^{\xi}\big)\big(1-b^{s\,\xi}\,a^{-s}\big).
  \label{eq:paretok_general}
\end{align}
Let us obtain the distribution of $A_{[j]}$ for $j\ge 2$ next. For
this single purpose we make two definitions:
$\eta_v=\lceil b^{\xi-1}-v b^{-1}\rceil$ and
\begin{equation*}
  \label{eq:taus}
  \tau_s(v)=\left\{
    \begin{array}{l}
      -\psi(v),\quad s=1\\
      \zeta(s,v),\quad s\neq 1
    \end{array}\right.
\end{equation*}
where $\psi(\cdot)$ is again the digamma function and
$\zeta(s,v)=\sum_{n=0}^\infty (n+v)^{-s}$ is Hurwitz's zeta 
function~\cite{apostol76:_intro_analytic_nt}. Now, combining~\eqref{eq:pmfAj}
and~\eqref{eq:paretok_general} and using the two previous definitions
it is tedious but straightforward to show that the Paretian
generalisation of~\eqref{eq:benfordk_j} is
\begin{align} %
  \label{eq::aj_pareto_general}
  \Pr(A_{[j]}=a)
                          &= \frac{b^{s(\xi-2)}}{1-b^{-s}}\left(\tau_s(ab^{-1}+b^{j-2})-\tau_s((a+1)b^{-1}+b^{j-2})\right)\nonumber\\
                       &-\frac{b^{s(\xi-1)}}{1-b^{-s}}\left(\tau_s(ab^{-1}+b^{j-1})-\tau_s((a+1)b^{-1}+b^{j-1})\right)\nonumber\\
                       &+b^{s(\xi-1)}\Big(\tau_s(ab^{-1}+\eta_a)-\tau_s((a+1)b^{-1}+
                          \eta_{a+1})\Big)\nonumber\\
  &+\eta_a-\eta_{a+1}.
\end{align}

\begin{remark}
  Like in Section~\ref{sec:most-sign-digits}, we have been able to
  obtain a closed-form expression for the pmf of $A_{[j]}$ thanks to
  the use of the $j$-th integer significand. Of particular interest is
  the distribution of $A_{(k)}$~\eqref{eq:paretok_general}, which had
  only been published before our own work~\cite{balado21:_benford} for
  the special case in which the fractional part of the minimum of the
  Pareto distribution is zero, i.e. $\rho=0$ and thus $\xi=k-1$. In
  this case~\eqref{eq:paretok_general} becomes
  \begin{equation}
    \Pr(A_{(k)}=a) =\frac{a^{-s}-(a+1)^{-s}}{b^{-s(k-1)}-b^{-sk}}.\label{eq:paretok_dtp}
  \end{equation} %
  The case $k=1$ of~\eqref{eq:paretok_dtp} was first given by
  Pietronero et al.~\cite{pietronero01:_explaining} in the course of
  their investigation on the generalisation of Benford's distribution
  to scale-invariant phenomena. Barabesi and
  Pratelli~\cite{barabesi20:_generalized} then extended Pietronero et
  al.'s result and obtained~\eqref{eq:paretok_dtp} itself. As we will
  empirically verify in Section~\ref{sec:empirical-tests}, the fact
  that~\eqref{eq:paretok_general} can handle the general case
  $\rho> 0$ is not a minor detail, but a major factor in terms of that
  expression being able to model real data that cannot be modelled
  by~\eqref{eq:paretok_dtp} alone.

  Interestingly,~\eqref{eq:paretok_dtp} was first identified as a new
  distribution only a few years ago by Kozubowski et
  al.~\cite{kozubowski15:_pareto}, who called it the \textit{discrete
    truncated Pareto} (DTP) distribution. Kozubowski et al. also
  noticed that the DTP distribution generalises Benford's
  distribution, but they landed on this fact solely because of the
  mathematical form of~\eqref{eq:paretok_dtp}. In fact, their
  practical motivation was far removed from the distribution of most
  significant digits: it was a biological problem involving the
  distribution of diet breadth in Lepidoptera.  Another striking fact
  is that Kozubowski et al. arrived at the DTP distribution through
  the quantisation of a \textit{truncated} Pareto variable, instead of
  through the discretisation of the fractional part of the logarithm
  of a \textit{standard} Pareto variable ---i.e. the procedure that we
  have followed to get to~\eqref{eq:paretok_dtp}, which is the
  ultimate reason why the DTP distribution is connected with Benford's
  distribution. A Pareto variable must surely be the only choice for
  which two such remarkably different procedures yield the very same
  outcome.  The reason for this serendipitous coincidence is that the
  complementary cdf of the variable to be quantised or discretised,
  respectively, turns out to be a negative exponential function in
  both cases. To end this remark, Kozubowski et al. rightly point out
  that that the shape parameter~$s$ in~\eqref{eq:paretok_dtp} can be
  taken to be negative in terms of its validity as a pmf. However
  observe that $s$ must be strictly positive
  for~\eqref{eq:paretok_dtp} to have physical meaning in terms of
  modelling a distribution of most significant digits.
\end{remark}

\subsubsection{Leading CF Coefficients of\hspace{.13cm}$\log_b X$}
\label{sec:lead-cf-coeff-1}

Applying~\eqref{eq:cdf_frac_logbx} to~\eqref{eq:cf_joint_pmf} yields
the Paretian generalisation of~\eqref{eq:cf_joint_pmf_benford}:
\begin{align}\label{eq:cf_joint_pmf_pareto}
  \Pr(\mathbf{A}_{k}=\mathbf{a}_k)&=(-1)^{k}\Big(b^{s(\rho-1)}\,\frac{b^{-s[0;\mathbf{a}_k]}-b^{-s[0;\mathbf{a}_k+\mathbf{e}_k]}}{1-b^{-s}}\nonumber\\
  &+\,u\big([0;\mathbf{a}_k+\mathbf{e}_k]-\rho\big)\,\big(1- b^{-s([0;\mathbf{a}_k+\mathbf{e}_k]-\rho)}\big)\nonumber\\
  &-\,u\big([0;\mathbf{a}_k]-\rho\big)\,\big(1- b^{-s([0;\mathbf{a}_k]-\rho)}\big)\Big).
\end{align}
The special case of~\eqref{eq:cf_joint_pmf_pareto} for $\rho=0$ yields
the counterpart of the DTP distribution~\eqref{eq:paretok_dtp} in the
CF setting:
\begin{equation}
  \label{eq:dpt_CF}
  \Pr(\mathbf{A}_k=\mathbf{a}_k)=(-1)^k\left(\frac{b^{-s[0;\mathbf{a}_k]}-b^{-s[0;\mathbf{a}_k+\mathbf{e}_k]}}{1-b^{-s}}\right).
\end{equation}
Expression~\eqref{eq:cf_joint_pmf_pareto} is clearly not amenable to
analytic marginalisation beyond $A_1=\mathbf{A}_1$. An interesting
particular case of $A_1$ is given by specialising~\eqref{eq:dpt_CF} for $k=1$:
\begin{equation}\label{eq:a1_pareto_zero_rho}
  \Pr(A_1=a)=\frac{b^{-\frac{s}{a+1}}-b^{-\frac{s}{a}}}{1-b^{-s}}.
\end{equation}
Recalling
Remark~\ref{rem:pareto_benford},~\eqref{eq:a1_pareto_zero_rho} tends
to~\eqref{eq:a1_benford} as $s\downarrow 0$.

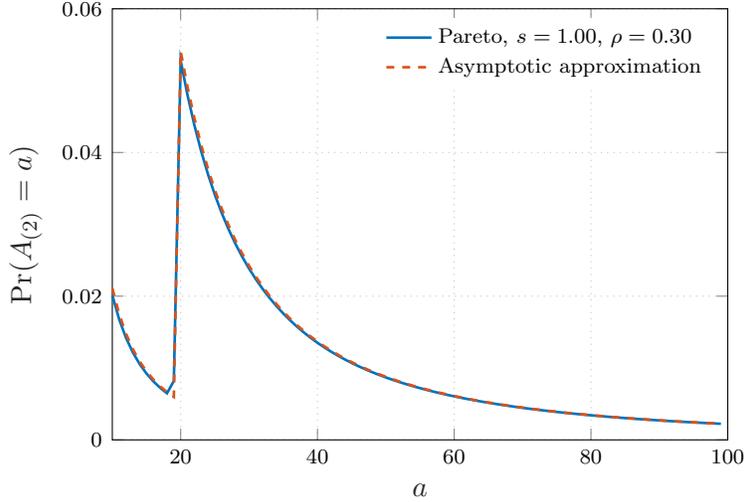
\begin{figure}[t!]
  \setlength\fwidth{.6\textwidth}
  \setlength\fheight{.4\textwidth} 
  \pgfplotsset{every tick label/.append style={font=\scriptsize}}
  \centering
   \input{fig2/general_pareto_vs_benford_approx_s1_00_rho0_30_mincha_291122_202751.tex}%
   \caption{Theoretical distribution of the two most significant
     decimal digits of Pareto~$X$~\eqref{eq:paretok_general} versus
     theoretical Benford-based asymptotic
     approximation~\eqref{eq:ak_pmf_asympt}. The lines join
     probability mass points for clarity.}
  \label{fig:bb_asympt_msd_pareto}
\end{figure}

\begin{figure}[t!]
  \setlength\fwidth{.6\textwidth}
  \setlength\fheight{.4\textwidth} 
  \pgfplotsset{every tick label/.append style={font=\scriptsize}}
  \centering
   \input{fig2/cf_general_pareto_vs_benford_approx_s1_00_rho0_30_mincha_271122_110309.tex}
   \caption{Theoretical joint pmf of the first two CF coefficients of
     $\log_{10}X$ for Pareto~$X$ with
     $s=1$ and $\rho=0.3$ [solid lines, \eqref{eq:cf_joint_pmf_pareto}]
     versus theoretical Benford-based asymptotic approximation [dashed
     lines, \eqref{eq:cf_joint_pmf_asympt}]. The lines join
     probability mass points corresponding to equal $a_2$ for
     clarity.}
  \label{fig:bb_asympt_lcf_pareto}
\end{figure}
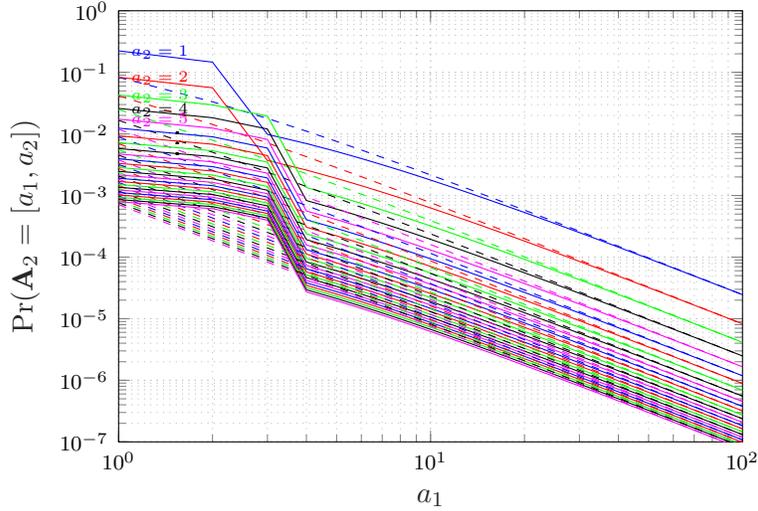

\subsubsection{Comparison with Benford-based Asymptotic Approximations}
\label{sec:benf-based-appr}
Now that we have particularised~\eqref{eq:pmfA}
and~\eqref{eq:cf_joint_pmf} for a non-Benfordian variable, we are in a
position to examine how well the Paretian
expressions~\eqref{eq:paretok_general}
and~\eqref{eq:cf_joint_pmf_pareto} are approximated by the
Benford-based asymptotic expressions discussed in
Section~\ref{sec:benf-vari-spec}. In order to evaluate
approximations~\eqref{eq:ak_pmf_asympt}
and~\eqref{eq:cf_joint_pmf_asympt} we need the pdf of $\{Y\}$,
$f_{\{Y\}}(y)$. This is obtained by
differentiating~\eqref{eq:cdf_frac_logbx}, which yields
$f_{\{Y\}}(y)= s (\ln
b)\,b^{-s(y-\rho)}\left(b^{-s}/(1-b^{-s\phantom{y}})+u\big(y-\rho\big)\right)$
for $y\in[0,1)$. We now compare in
Figures~\ref{fig:bb_asympt_msd_pareto}
and~\ref{fig:bb_asympt_lcf_pareto} the exact Paretian expressions and
their Benford-based approximations for some arbitrary values of $s$
and $\rho$ and for $k=2$, which we have chosen because the
visualisation of a joint distribution is not simple for $k>2$ in the
CF case. Even though the Benford-based approximations were obtained
assuming $k$ to be large, we can see that they are close to the exact
expressions already. This is markedly true for the most significant
digits model in Figure~\ref{fig:bb_asympt_msd_pareto}, where the
approximation is accurate for all $a\in\mathcal{A}_{(2)}$ already
---in fact, it can be verified that the asymptotic approximation is
also acceptable for $k=1$ in this case. Regarding the CF coefficients
model, the approximation becomes accurate when
$a_1+a_2\gtrapprox 10^2$.

\section{Empirical Tests}
\label{sec:empirical-tests}

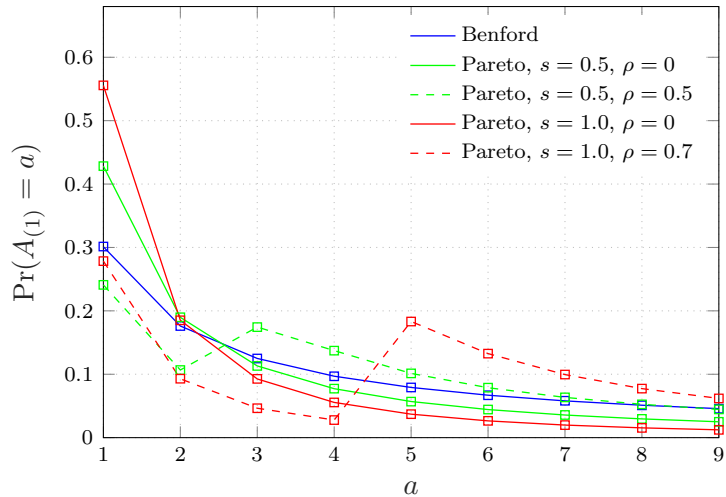
\begin{figure}[t!]
  \setlength\fwidth{.6\textwidth}
  \setlength\fheight{.4\textwidth} 
  \pgfplotsset{every tick label/.append style={font=\scriptsize}}
  \centering
   \input{fig2/b-p.tex}
  \caption{Distributions of the most significant decimal digit
    of~$X$. The theoretical pmf's (solid and dashed lines) are
    \eqref{eq:benfordk} and~\eqref{eq:paretok_general}, and the
    empirical frequencies (\raisebox{.8pt}{{\tiny $\square$}}) correspond to $p=10^7$
    pseudorandom outcomes in each case.}
  \label{fig:pr_ak_1}
\end{figure}

\begin{figure}[t!]
  \setlength\fwidth{.6\textwidth}
  \setlength\fheight{.4\textwidth} 
  \pgfplotsset{every tick label/.append style={font=\scriptsize}}
  \centering
  \input{fig2/b-p-2.tex}
  
  \caption{Distributions of the two most significant decimal digits
    of~$X$. The theoretical pmf's (solid and dashed lines) are
    \eqref{eq:benfordk} and~\eqref{eq:paretok_general}, and the
    empirical frequencies (\raisebox{.8pt}{{\tiny $\square$}}) correspond to $p=10^7$
    pseudorandom outcomes in each case.}
  \label{fig:pr_ak_2}
\end{figure}
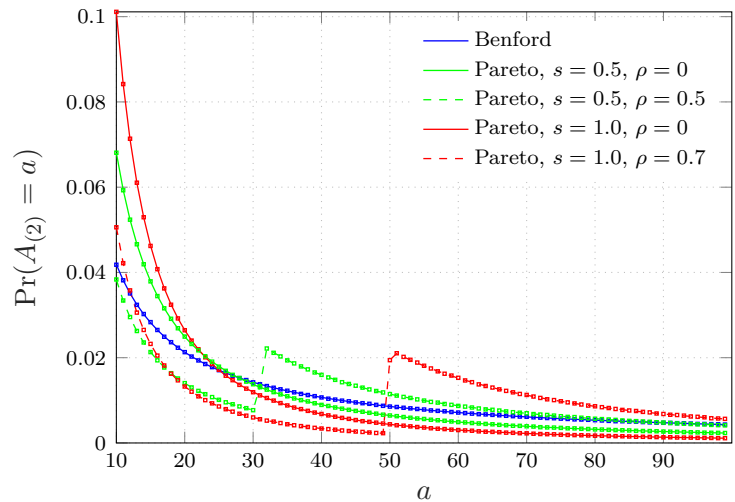

In this section the theoretical expressions given in
Sections~\ref{sec:benford} and~\ref{sec:pareto} are verified. In all
plots, solid or dashed lines represent theoretical probabilities
(joining discrete probability mass points) whereas square symbols
represent empirical frequencies obtained using a
dataset~$\{x_1,x_2,\ldots,x_p\}$. For simplicity, we use the maximum
likelihood (ML) estimators $\hat{x}_\text{m}=\min_i x_i$ and
$\hat{s}=\big(\frac{1}{p}\sum_i \ln(x_i/\hat{x}_\text{m})\big)^{-1}$
to drive the Paretian expressions with real datasets, but be aware
that better estimation approaches are
possible (see for instance~\cite{nair22:_fundam_heavy_tails}).

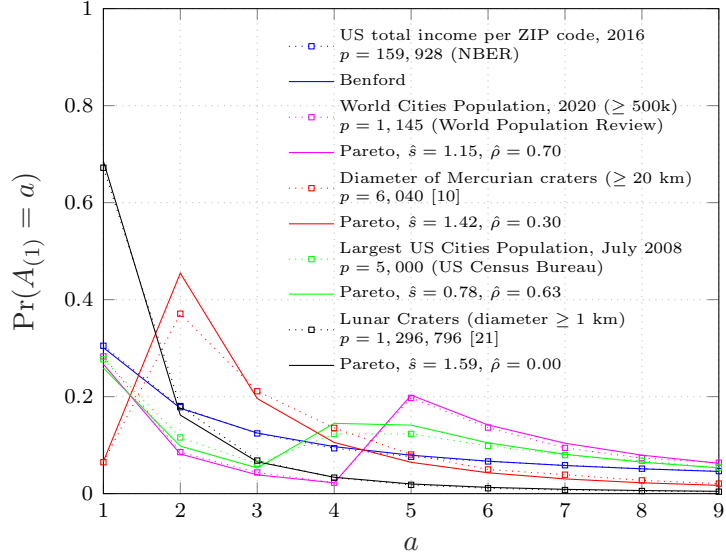
\begin{figure}
  \setlength\fwidth{.6\textwidth}
  \setlength\fheight{.45\textwidth} 
  \pgfplotsset{every tick label/.append style={font=\scriptsize}}
  \centering
  \input{fig2/msd_real_datasets.tex}
  \caption{Distributions of the most significant decimal digit of~$X$
    for real Benfordian and Paretian datasets. The theoretical pmf's
    (solid lines) are \eqref{eq:benfordk}
    and~\eqref{eq:paretok_general}, and the empirical frequencies
    (\raisebox{.8pt}{{\tiny $\square$}}) are joined by dotted lines for clarity.}
    \label{fig:msd_real_datasets}\phantom{\cite{fassett11:_mercury}}
\end{figure}

We start with distributions of the most significant digits of
$X$. Figures~\ref{fig:pr_ak_1} and~\ref{fig:pr_ak_2} present the
distributions of the most significant decimal digit $A_{(1)}$ and of
the two most significant decimal digits $A_{(2)}$, respectively. The
Benford results, which are a particular case
of~\eqref{eq:paretok_general}, are well known. The Paretian cases with
$\rho=0$ are covered by~\eqref{eq:paretok_dtp}, as previously shown by
Barabesi and Pratelli~\cite{barabesi20:_generalized}. However, it is
essential to use the general expression~\eqref{eq:paretok_general}
when $\rho>0$. In this case the pmf's of Paretian significant digits
do not behave anymore in a monotonically decreasing way (i.e. like
Benford's pmf) but rather feature a peak midway along the support
of~$A_{(k)}$.  Therefore modelling real Paretian datasets requires
being able to take a general value of $\rho$ into account, as there is
no special reason why $\rho$ should be zero in practice ---observe
some examples in
Figure~\ref{fig:msd_real_datasets}. %

Next, Figure~\ref{fig:Aj_digits_pareto_theoretical_montecarlo} shows
distributions of the $j$-th most significant decimal digit $A_{[j]}$.
Again, peaks can be seen in the distributions when $\rho>0$, but in
general these are less pronounced than in the distribution of
$A_{(k)}$ due to~\eqref{eq:asympt-aj}. An illustration of the
asymptotic behaviour proved in Theorem~\ref{thm:asympt-aj} is the fact
that~$A_{[4]}$ is nearly uniformly distributed for all three
distributions of~$X$ considered in
Figure~\ref{fig:Aj_digits_pareto_theoretical_montecarlo}.

\begin{figure}[t!]
  \setlength\fwidth{.6\textwidth}
  \setlength\fheight{.4\textwidth} 
  \pgfplotsset{every tick label/.append style={font=\scriptsize}}
  \centering
  \input{fig2/paretoj_general.tex}
  \caption{Distributions of the $j$-th most significant decimal digit
    of~$X$. The theoretical pmf's (solid and dashed lines)
    are~\eqref{eq:benfordk_j} and~\eqref{eq::aj_pareto_general}, and
    the empirical frequencies (\raisebox{.8pt}{{\tiny $\square$}}) correspond to $p=5\times 10^7$
    pseudorandom outcomes in each case.}
  \label{fig:Aj_digits_pareto_theoretical_montecarlo}
\end{figure}
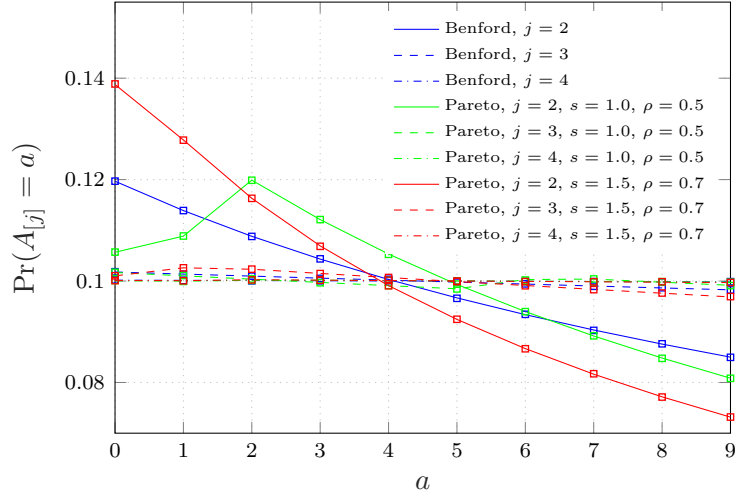

\begin{figure}[t!]
  \setlength\fwidth{.6\textwidth}
  \setlength\fheight{.4\textwidth} 
  \pgfplotsset{every tick label/.append style={font=\scriptsize}}
  \centering
  \input{fig2/cf_general_benford_pseudo_n1e+08_mincha_110522_092109.tex}
  \caption{Joint distribution of the first two CF coefficients of
    $\log_{10} X$ for Benford~$X$. The theoretical joint pmf (dashed
    lines) is~\eqref{eq:cf_joint_pmf_benford} and the empirical
    frequencies (\raisebox{.8pt}{{\tiny $\square$}}) correspond to $p=10^8$ pseudorandom
    outcomes.}
  \label{fig:cf_joint_A2_benford_theoretical_montecarlo}
\end{figure}
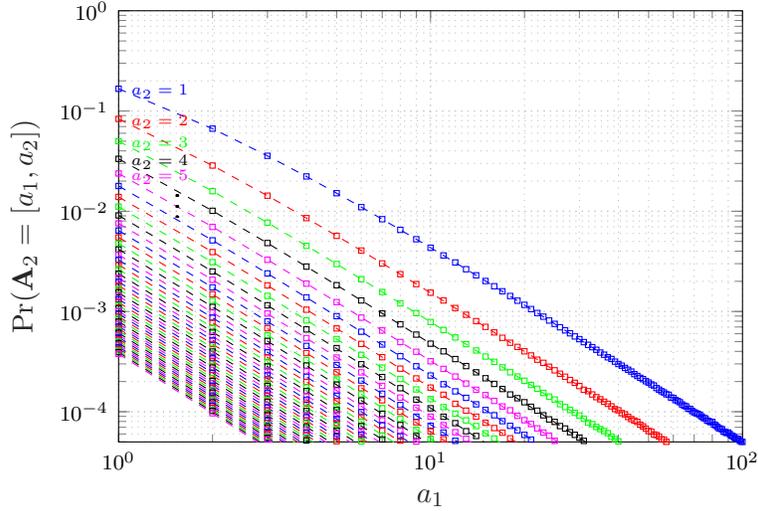

\begin{figure}[t!]
  \setlength\fwidth{.6\textwidth}
  \setlength\fheight{.4\textwidth} 
  \pgfplotsset{every tick label/.append style={font=\scriptsize}}
  \centering
  \input{fig2/cf_general_pareto_pseudo_n1e+08_s1_50_rho0_48_mincha_110522_102323.tex}
  \caption{Joint distribution of the first two CF coefficients of
    $\log_{10} X$ for Pareto~$X$, $s=1.5$, $\rho=0.48$. The
    theoretical joint pmf (dashed lines)
    is~\eqref{eq:cf_joint_pmf_pareto} and the empirical frequencies
    (\raisebox{.8pt}{{\tiny $\square$}}) correspond to $p=10^8$ pseudorandom outcomes.}
  \label{fig:cf_joint_A2_pareto_theoretical_montecarlo}
\end{figure}
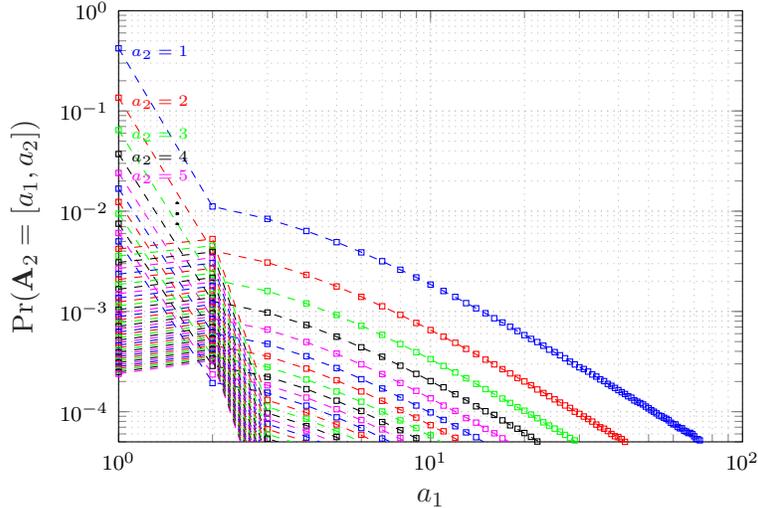

\begin{figure}[t!]
  \setlength\fwidth{.6\textwidth}
  \setlength\fheight{.4\textwidth} 
  \pgfplotsset{every tick label/.append style={font=\scriptsize}}
  \centering
 \input{fig2/cf_gk_a1_a2_n1e+08_mincha_110522_101623.tex}
 \caption{Distributions of the $j$-th CF coefficient of $\log_{10}
   X$. The theoretical pmf's (solid lines)
   are~\eqref{eq:gauss-kuzmin},~\eqref{eq:a1_benford},~\eqref{eq:a2_benford},
   and~\eqref{eq:cf_joint_pmf_pareto} [with $k=1$]. The empirical
   frequencies (\raisebox{.8pt}{{\tiny $\square$}}) correspond to $p=10^8$ pseudorandom
   outcomes in each case.}
\label{fig:cf_gk_a1_a2}
\end{figure}
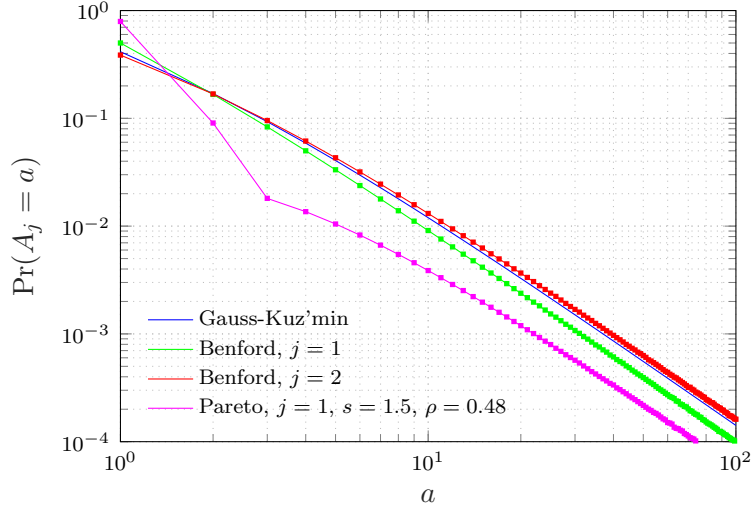

\begin{figure}[t!]
  \setlength\fwidth{0.6\textwidth}
  \setlength\fheight{0.4\textwidth}
 \pgfplotsset{every tick label/.append style={font=\scriptsize}}
   \centering
  \input{fig2/cf_general_benford_US_total_income_per_ZIP_code,_2016_n159924_mincha_110522_092226.tex}
  \caption{Joint distribution of the first two CF coefficients of
    $\log_{10} X$ for a real Benfordian dataset (US total income per
    ZIP code, National Bureau of Economic Research, 2016,
    $p=159,928$). The theoretical joint pmf (dashed lines)
    is~\eqref{eq:cf_joint_pmf_benford}, whereas the symbols
    (\raisebox{.8pt}{{\tiny $\square$}}) represent empirical
    frequencies. }
\label{fig:A2_cf_benford_income_dataset}
\end{figure}
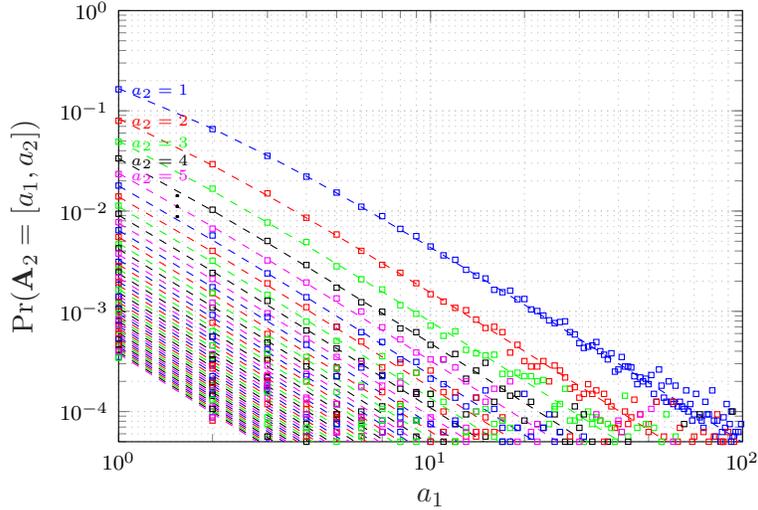

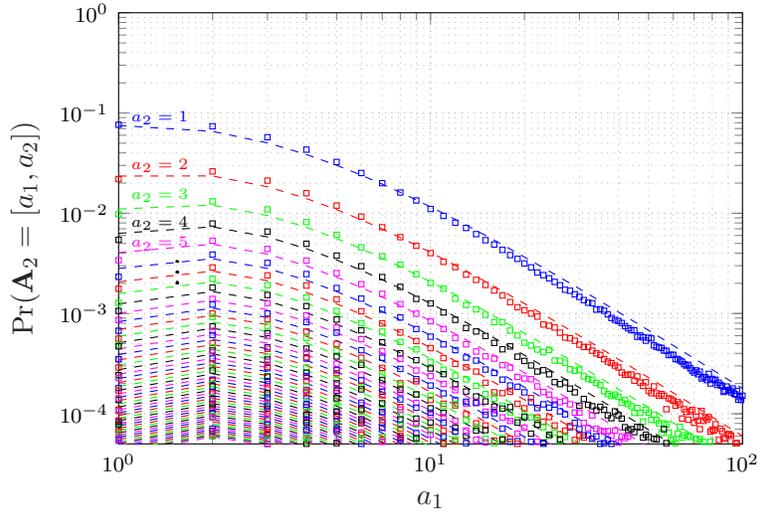
\begin{figure}[t!]
  \setlength\fwidth{0.6\textwidth}
  \setlength\fheight{0.4\textwidth}
  \pgfplotsset{every tick label/.append style={font=\scriptsize}}
  \centering
  \input{fig2/cf_general_pareto_Lunar_Craters_diameter_gt_1km_n1296785_s1_59_rho0_00_mincha_110522_094958.tex}
  \caption{Joint distribution of the first two CF coefficients of
    $\log_{10} X$ for a real Paretian dataset (diameter of Lunar
    craters $\ge 1$ km, $p=1,296,796$~\cite{robbins18}). The
    theoretical joint pmf (dashed lines)
    is~\eqref{eq:cf_joint_pmf_pareto}, driven by $\hat{s}=1.59$ and
    $\hat{\rho}=0.00$, whereas the symbols (\raisebox{.8pt}{{\tiny $\square$}}) represent  empirical frequencies.}
\label{fig:A2_cf_pareto_lunar_dataset}
\end{figure}

We move on next to distributions of the leading CF coefficients of
$\log_{10}
X$. Figures~\ref{fig:cf_joint_A2_benford_theoretical_montecarlo}
and~\ref{fig:cf_joint_A2_pareto_theoretical_montecarlo} verify the
validity of the joint distributions~\eqref{eq:cf_joint_pmf_benford}
and~\eqref{eq:cf_joint_pmf_pareto} of the two leading CF coefficients
$\mathbf{A}_2$ when $X$ is Benford and Pareto, respectively. In
Figure~\ref{fig:cf_gk_a1_a2} we show the distributions of the
marginals $A_1$ and $A_2$ for Benford $X$ [i.e. \eqref{eq:a1_benford}
and~\eqref{eq:a2_benford}], and we compare them to the Gauss-Kuz'min
law~\eqref{eq:gauss-kuzmin}. As we can see, the distribution of $A_j$
converges really fast to Gauss-Kuz'min: the distribution of $A_1$ is
close to it already, as remarked by Miller and
Takloo-Bighash~\cite{miller06:_invitation}, but that of $A_2$ is even
closer as expected
from~\eqref{eq:approx_aj_gk}. Figure~\ref{fig:cf_gk_a1_a2} also
depicts the distribution of~$A_1$ for Paretian~$X$,
using~\eqref{eq:cf_joint_pmf_pareto}. As we know, for Pareto $X$ the
distribution of $A_j$ must also converge exponentially fast to the
Gauss-Kuz'min law, although this is not graphically illustrated in
Figure~\ref{fig:cf_gk_a1_a2} due to the lack of a theoretical
expression for~$A_2$ in this case. Finally, we show in
Figures~\ref{fig:A2_cf_benford_income_dataset}
and~\ref{fig:A2_cf_pareto_lunar_dataset}
how~\eqref{eq:cf_joint_pmf_benford} and~\eqref{eq:cf_joint_pmf_pareto}
correctly model real Benfordian and Paretian datasets, respectively.

\section{Conclusions}
\label{sec:conclusions}
We have provided a general theoretical analysis of the distributions
of the most significant digits and the leading continued fraction
coefficients of the outcomes of an arbitrary random variable, which
highlights the connections between the two subjects. Empirical
verification for two relevant particularisations of our results
(Benford and Pareto variables, respectively) also supports the
accuracy of our results in practice. Our analysis reveals novel facts
---specially, but not only, concerning modelling continued fraction
coefficients--- and provides simpler proofs and new closed-form
expressions for already known ones.  In particular, we have shown that
the use of what we propose to call $k$-th integer significands
considerably simplifies modelling significant digits, allowing for
uncomplicated finite and asymptotic analyses. We have also shown the
parallelism between the general asymptotics of the probabilistic
models for the $j$-th significant $b$-ary digit and of the $j$-th
continued fraction coefficient ---i.e. between~\eqref{eq:asympt-aj}
and the Gauss-Kuz'min law~\eqref{eq:gauss-kuzmin}---, and the role
played by the Benford variables in the asymptotics of the general
analyses. Our results may find application in all areas where
Benford's law has been previously used.

\providecommand{\bysame}{\leavevmode\hbox to3em{\hrulefill}\thinspace}
\providecommand{\MR}{\relax\ifhmode\unskip\space\fi MR }
\providecommand{\MRhref}[2]{%
  \href{http://www.ams.org/mathscinet-getitem?mr=#1}{#2}
}
\providecommand{\href}[2]{#2}

\end{document}

%% file: fig2/si.tex
%
%
\definecolor{mycolor1}{rgb}{1.00000,0.00000,1.00000}%
\begin{tikzpicture}

\begin{axis}[%
width=0.951\fwidth,
height=\fheight,
at={(0\fwidth,0\fheight)},
scale only axis,
clip=false,
xmin=1,
xmax=10,
xlabel style={font=\color{white!15!black}},
xlabel={$a/10^{k-1}\quad [\text{where }a\in\mathcal{A}_{(k)}]$},
ymin=0.3,
ymax=0.45,
ylabel style={font=\color{white!15!black}},
ylabel={$a\, \log_{10} (1+1/a)\qquad$},
axis background/.style={fill=white},
xmajorgrids,
ymajorgrids,
legend style={at={(0.97,0.03)}, anchor=south east, legend cell align=left, align=left, draw=white!15!black},
ytick={.30,.35,.40,.45,.50},grid style={dotted},legend style={font=\scriptsize,draw=none},mark size=0.5pt,xtick={1,2,3,4,5,6,7,8,9}
]
\addplot [color=blue, dotted, mark=square, mark options={solid, blue}]
  table[row sep=crcr]{%
1	0.432137378264258\\
1.01	0.432158595906774\\
1.02	0.432179400211972\\
1.03	0.432199803141639\\
1.04	0.432219816200401\\
1.05	0.432239450457376\\
1.06	0.432258716566581\\
1.07	0.432277624786181\\
1.08	0.432296184996787\\
1.09	0.432314406718551\\
1.1	0.432332299127561\\
1.11	0.432349871071179\\
1.12	0.432367131082667\\
1.13	0.432384087394975\\
1.14	0.432400747953857\\
1.15	0.432417120430279\\
1.16	0.432433212232205\\
1.17	0.432449030515756\\
1.18	0.432464582195833\\
1.19	0.432479873956193\\
1.2	0.432494912259029\\
1.21	0.432509703354074\\
1.22	0.432524253287258\\
1.23	0.432538567908964\\
1.24	0.432552652881846\\
1.25	0.432566513688311\\
1.26	0.432580155637638\\
1.27	0.43259358387276\\
1.28	0.432606803376716\\
1.29	0.432619818978828\\
1.3	0.432632635360572\\
1.31	0.432645257061209\\
1.32	0.432657688483142\\
1.33	0.432669933897007\\
1.34	0.432681997446592\\
1.35	0.432693883153533\\
1.36	0.432705594921737\\
1.37	0.432717136541681\\
1.38	0.432728511694484\\
1.39	0.432739723955867\\
1.4	0.432750776799861\\
1.41	0.432761673602392\\
1.42	0.432772417644756\\
1.43	0.432783012116868\\
1.44	0.432793460120433\\
1.45	0.43280376467202\\
1.46	0.432813928705896\\
1.47	0.43282395507685\\
1.48	0.432833846562867\\
1.49	0.43284360586768\\
1.5	0.432853235623225\\
1.51	0.432862738392072\\
1.52	0.432872116669589\\
1.53	0.432881372886231\\
1.54	0.43289050940958\\
1.55	0.432899528546361\\
1.56	0.432908432544452\\
1.57	0.432917223594655\\
1.58	0.432925903832562\\
1.59	0.432934475340254\\
1.6	0.432942940147995\\
1.61	0.432951300235772\\
1.62	0.432959557534954\\
1.63	0.43296771392963\\
1.64	0.432975771258178\\
1.65	0.432983731314557\\
1.66	0.432991595849683\\
1.67	0.432999366572684\\
1.68	0.433007045152193\\
1.69	0.433014633217468\\
1.7	0.433022132359583\\
1.71	0.433029544132555\\
1.72	0.433036870054398\\
1.73	0.433044111608149\\
1.74	0.433051270242885\\
1.75	0.433058347374697\\
1.76	0.433065344387593\\
1.77	0.433072262634466\\
1.78	0.433079103437848\\
1.79	0.433085868090913\\
1.8	0.433092557858118\\
1.81	0.433099173976149\\
1.82	0.43310571765455\\
1.83	0.43311219007658\\
1.84	0.433118592399823\\
1.85	0.43312492575696\\
1.86	0.433131191256374\\
1.87	0.433137389982827\\
1.88	0.43314352299808\\
1.89	0.433149591341535\\
1.9	0.433155596030731\\
1.91	0.433161538062011\\
1.92	0.433167418411031\\
1.93	0.433173238033293\\
1.94	0.433178997864639\\
1.95	0.433184698821816\\
1.96	0.433190341802903\\
1.97	0.433195927687814\\
1.98	0.43320145733875\\
1.99	0.433206931600627\\
2	0.433212351301526\\
2.01	0.433217717253122\\
2.02	0.433223030251008\\
2.03	0.433228291075234\\
2.04	0.433233500490529\\
2.05	0.433238659246816\\
2.06	0.433243768079456\\
2.07	0.433248827709668\\
2.08	0.433253838844828\\
2.09	0.433258802178833\\
2.1	0.433263718392412\\
2.11	0.433268588153395\\
2.12	0.433273412117098\\
2.13	0.433278190926519\\
2.14	0.433282925212707\\
2.15	0.433287615595006\\
2.16	0.433292262681292\\
2.17	0.433296867068353\\
2.18	0.433301429341953\\
2.19	0.433305950077244\\
2.2	0.43331042983899\\
2.21	0.433314869181671\\
2.22	0.433319268649889\\
2.23	0.433323628778477\\
2.24	0.433327950092737\\
2.25	0.43333223310865\\
2.26	0.433336478333117\\
2.27	0.433340686264156\\
2.28	0.433344857390989\\
2.29	0.433348992194409\\
2.3	0.433353091146826\\
2.31	0.433357154712486\\
2.32	0.433361183347687\\
2.33	0.43336517750087\\
2.34	0.433369137612861\\
2.35	0.433373064117013\\
2.36	0.433376957439358\\
2.37	0.433380817998712\\
2.38	0.433384646206922\\
2.39	0.433388442468928\\
2.4	0.433392207182965\\
2.41	0.433395940740656\\
2.42	0.43339964352719\\
2.43	0.433403315921398\\
2.44	0.433406958295937\\
2.45	0.433410571017437\\
2.46	0.433414154446513\\
2.47	0.433417708937974\\
2.48	0.43342123484099\\
2.49	0.433424732499021\\
2.5	0.433428202250133\\
2.51	0.433431644426995\\
2.52	0.433435059357002\\
2.53	0.433438447362393\\
2.54	0.433441808760345\\
2.55	0.433445143863068\\
2.56	0.433448452977914\\
2.57	0.433451736407453\\
2.58	0.433454994449591\\
2.59	0.433458227397629\\
2.6	0.433461435540378\\
2.61	0.433464619162237\\
2.62	0.433467778543247\\
2.63	0.433470913959241\\
2.64	0.433474025681884\\
2.65	0.433477113978671\\
2.66	0.433480179113181\\
2.67	0.433483221345033\\
2.68	0.433486240929941\\
2.69	0.433489238119829\\
2.7	0.43349221316298\\
2.71	0.433495166303892\\
2.72	0.433498097783592\\
2.73	0.433501007839514\\
2.74	0.433503896705657\\
2.75	0.433506764612643\\
2.76	0.433509611787716\\
2.77	0.433512438454885\\
2.78	0.433515244834913\\
2.79	0.433518031145436\\
2.8	0.433520797600986\\
2.81	0.43352354441303\\
2.82	0.433526271790011\\
2.83	0.433528979937529\\
2.84	0.433531669058211\\
2.85	0.43353433935185\\
2.86	0.433536991015511\\
2.87	0.433539624243453\\
2.88	0.43354223922731\\
2.89	0.433544836155985\\
2.9	0.433547415215835\\
2.91	0.4335499765907\\
2.92	0.433552520461808\\
2.93	0.433555047008023\\
2.94	0.433557556405662\\
2.95	0.433560048828791\\
2.96	0.433562524449024\\
2.97	0.433564983435748\\
2.98	0.433567425955982\\
2.99	0.433569852174604\\
3	0.433572262254286\\
3.01	0.433574656355478\\
3.02	0.43357703463661\\
3.03	0.433579397253956\\
3.04	0.433581744361769\\
3.05	0.433584076112213\\
3.06	0.433586392655583\\
3.07	0.433588694140138\\
3.08	0.433590980712241\\
3.09	0.433593252516363\\
3.1	0.43359550969509\\
3.11	0.43359775238924\\
3.12	0.433599980737775\\
3.13	0.433602194877888\\
3.14	0.433604394945087\\
3.15	0.433606581073041\\
3.16	0.433608753393868\\
3.17	0.433610912037933\\
3.18	0.433613057134022\\
3.19	0.433615188809213\\
3.2	0.433617307189159\\
3.21	0.433619412397768\\
3.22	0.433621504557572\\
3.23	0.433623583789493\\
3.24	0.433625650212967\\
3.25	0.433627703945989\\
3.26	0.433629745105142\\
3.27	0.4336317738055\\
3.28	0.433633790160815\\
3.29	0.433635794283449\\
3.3	0.433637786284303\\
3.31	0.433639766273123\\
3.32	0.433641734358153\\
3.33	0.433643690646459\\
3.34	0.43364563524377\\
3.35	0.43364756825459\\
3.36	0.433649489782176\\
3.37	0.433651399928527\\
3.38	0.433653298794456\\
3.39	0.433655186479624\\
3.4	0.433657063082481\\
3.41	0.433658928700323\\
3.42	0.43366078342934\\
3.43	0.433662627364564\\
3.44	0.43366446059994\\
3.45	0.433666283228367\\
3.46	0.433668095341587\\
3.47	0.433669897030405\\
3.48	0.433671688384432\\
3.49	0.433673469492412\\
3.5	0.433675240441961\\
3.51	0.433677001319744\\
3.52	0.433678752211437\\
3.53	0.433680493201726\\
3.54	0.433682224374406\\
3.55	0.433683945812284\\
3.56	0.433685657597225\\
3.57	0.433687359810176\\
3.58	0.433689052531236\\
3.59	0.433690735839547\\
3.6	0.433692409813435\\
3.61	0.433694074530305\\
3.62	0.433695730066737\\
3.63	0.433697376498467\\
3.64	0.433699013900411\\
3.65	0.433700642346641\\
3.66	0.433702261910377\\
3.67	0.433703872664183\\
3.68	0.433705474679728\\
3.69	0.433707068027862\\
3.7	0.433708652778821\\
3.71	0.433710229001945\\
3.72	0.433711796765916\\
3.73	0.433713356138625\\
3.74	0.433714907187271\\
3.75	0.433716449978313\\
3.76	0.433717984577553\\
3.77	0.433719511050035\\
3.78	0.433721029460154\\
3.79	0.43372253987162\\
3.8	0.433724042347458\\
3.81	0.433725536950085\\
3.82	0.433727023741165\\
3.83	0.433728502781772\\
3.84	0.433729974132443\\
3.85	0.433731437852879\\
3.86	0.433732894002389\\
3.87	0.433734342639504\\
3.88	0.43373578382219\\
3.89	0.433737217607867\\
3.9	0.433738644053363\\
3.91	0.433740063214851\\
3.92	0.433741475148039\\
3.93	0.433742879907932\\
3.94	0.433744277549138\\
3.95	0.433745668125571\\
3.96	0.433747051690705\\
3.97	0.433748428297407\\
3.98	0.433749797998025\\
3.99	0.433751160844441\\
4	0.433752516887957\\
4.01	0.433753866179384\\
4.02	0.433755208769039\\
4.03	0.433756544706678\\
4.04	0.43375787404169\\
4.05	0.433759196822846\\
4.06	0.433760513098534\\
4.07	0.433761822916598\\
4.08	0.433763126324434\\
4.09	0.433764423369031\\
4.1	0.433765714096822\\
4.11	0.433766998553894\\
4.12	0.433768276785787\\
4.13	0.433769548837649\\
4.14	0.433770814754219\\
4.15	0.433772074579779\\
4.16	0.433773328358164\\
4.17	0.43377457613279\\
4.18	0.433775817946739\\
4.19	0.433777053842548\\
4.2	0.433778283862493\\
4.21	0.433779508048339\\
4.22	0.433780726441513\\
4.23	0.433781939082974\\
4.24	0.433783146013441\\
4.25	0.433784347273152\\
4.26	0.433785542901893\\
4.27	0.433786732939268\\
4.28	0.433787917424333\\
4.29	0.433789096395924\\
4.3	0.433790269892391\\
4.31	0.433791437951784\\
4.32	0.433792600611859\\
4.33	0.433793757909913\\
4.34	0.433794909882942\\
4.35	0.43379605656767\\
4.36	0.433797198000402\\
4.37	0.433798334217172\\
4.38	0.43379946525358\\
4.39	0.433800591145021\\
4.4	0.433801711926473\\
4.41	0.433802827632717\\
4.42	0.433803938298127\\
4.43	0.433805043956761\\
4.44	0.433806144642441\\
4.45	0.433807240388592\\
4.46	0.433808331228374\\
4.47	0.433809417194771\\
4.48	0.433810498320257\\
4.49	0.433811574637199\\
4.5	0.433812646177566\\
4.51	0.433813712973128\\
4.52	0.433814775055303\\
4.53	0.433815832455237\\
4.54	0.433816885203836\\
4.55	0.433817933331772\\
4.56	0.433818976869361\\
4.57	0.433820015846661\\
4.58	0.433821050293549\\
4.59	0.433822080239574\\
4.6	0.433823105714092\\
4.61	0.433824126746055\\
4.62	0.433825143364384\\
4.63	0.43382615559752\\
4.64	0.433827163473906\\
4.65	0.433828167021513\\
4.66	0.433829166268176\\
4.67	0.433830161241547\\
4.68	0.433831151968918\\
4.69	0.433832138477428\\
4.7	0.433833120794005\\
4.71	0.433834098945244\\
4.72	0.433835072957651\\
4.73	0.433836042857367\\
4.74	0.433837008670453\\
4.75	0.433837970422644\\
4.76	0.433838928139485\\
4.77	0.43383988184635\\
4.78	0.43384083156839\\
4.79	0.433841777330492\\
4.8	0.433842719157405\\
4.81	0.433843657073557\\
4.82	0.433844591103356\\
4.83	0.433845521270866\\
4.84	0.43384644759996\\
4.85	0.433847370114415\\
4.86	0.433848288837688\\
4.87	0.433849203793146\\
4.88	0.433850115003892\\
4.89	0.433851022492894\\
4.9	0.433851926282863\\
4.91	0.433852826396381\\
4.92	0.433853722855874\\
4.93	0.433854615683551\\
4.94	0.433855504901382\\
4.95	0.433856390531229\\
4.96	0.433857272594777\\
4.97	0.433858151113532\\
4.98	0.43385902610885\\
4.99	0.433859897601793\\
5	0.433860765613457\\
5.01	0.433861630164594\\
5.02	0.433862491275847\\
5.03	0.433863348967722\\
5.04	0.433864203260587\\
5.05	0.433865054174557\\
5.06	0.433865901729649\\
5.07	0.433866745945725\\
5.08	0.433867586842461\\
5.09	0.433868424439403\\
5.1	0.433869258755945\\
5.11	0.433870089811298\\
5.12	0.433870917624584\\
5.13	0.433871742214715\\
5.14	0.433872563600451\\
5.15	0.433873381800497\\
5.16	0.433874196833275\\
5.17	0.433875008717182\\
5.18	0.433875817470485\\
5.19	0.433876623111143\\
5.2	0.433877425657162\\
5.21	0.433878225126345\\
5.22	0.433879021536311\\
5.23	0.433879814904608\\
5.24	0.433880605248658\\
5.25	0.433881392585636\\
5.26	0.433882176932784\\
5.27	0.433882958307009\\
5.28	0.43388373672522\\
5.29	0.433884512204143\\
5.3	0.433885284760436\\
5.31	0.433886054410516\\
5.32	0.433886821170801\\
5.33	0.433887585057566\\
5.34	0.433888346086838\\
5.35	0.433889104274747\\
5.36	0.433889859637082\\
5.37	0.433890612189613\\
5.38	0.433891361948034\\
5.39	0.433892108927846\\
5.4	0.433892853144483\\
5.41	0.433893594613266\\
5.42	0.433894333349353\\
5.43	0.433895069367821\\
5.44	0.433895802683644\\
5.45	0.433896533311649\\
5.46	0.433897261266672\\
5.47	0.433897986563285\\
5.48	0.433898709216086\\
5.49	0.433899429239396\\
5.5	0.433900146647638\\
5.51	0.433900861455009\\
5.52	0.43390157367562\\
5.53	0.433902283323524\\
5.54	0.433902990412562\\
5.55	0.433903694956579\\
5.56	0.433904396969312\\
5.57	0.433905096464362\\
5.58	0.433905793455274\\
5.59	0.433906487955386\\
5.6	0.433907179978153\\
5.61	0.433907869536705\\
5.62	0.433908556644261\\
5.63	0.433909241313775\\
5.64	0.433909923558271\\
5.65	0.43391060339058\\
5.66	0.433911280823497\\
5.67	0.433911955869656\\
5.68	0.43391262854172\\
5.69	0.433913298852133\\
5.7	0.433913966813292\\
5.71	0.433914632437604\\
5.72	0.433915295737201\\
5.73	0.433915956724367\\
5.74	0.433916615411124\\
5.75	0.433917271809415\\
5.76	0.433917925931177\\
5.77	0.433918577788225\\
5.78	0.433919227392309\\
5.79	0.433919874755148\\
5.8	0.433920519888189\\
5.81	0.433921162803074\\
5.82	0.433921803511099\\
5.83	0.433922442023668\\
5.84	0.433923078352069\\
5.85	0.433923712507467\\
5.86	0.433924344500939\\
5.87	0.433924974343568\\
5.88	0.433925602046348\\
5.89	0.433926227620037\\
5.9	0.433926851075614\\
5.91	0.43392747242366\\
5.92	0.433928091674972\\
5.93	0.433928708840041\\
5.94	0.433929323929486\\
5.95	0.433929936953683\\
5.96	0.433930547923071\\
5.97	0.43393115684794\\
5.98	0.433931763738549\\
5.99	0.433932368605018\\
6	0.433932971457543\\
6.01	0.433933572306103\\
6.02	0.433934171160725\\
6.03	0.43393476803128\\
6.04	0.433935362927571\\
6.05	0.433935955859466\\
6.06	0.433936546836673\\
6.07	0.433937135868748\\
6.08	0.433937722965386\\
6.09	0.433938308136061\\
6.1	0.433938891390186\\
6.11	0.433939472737273\\
6.12	0.433940052186551\\
6.13	0.433940629747345\\
6.14	0.433941205428921\\
6.15	0.433941779240383\\
6.16	0.433942351190798\\
6.17	0.433942921289259\\
6.18	0.433943489544702\\
6.19	0.433944055966138\\
6.2	0.433944620562339\\
6.21	0.433945183342117\\
6.22	0.433945744314295\\
6.23	0.4339463034875\\
6.24	0.433946860870367\\
6.25	0.43394741647155\\
6.26	0.433947970299502\\
6.27	0.433948522362776\\
6.28	0.433949072669732\\
6.29	0.433949621228709\\
6.3	0.433950168048141\\
6.31	0.43395071313618\\
6.32	0.433951256501098\\
6.33	0.433951798151027\\
6.34	0.433952338094046\\
6.35	0.433952876338271\\
6.36	0.433953412891635\\
6.37	0.433953947762148\\
6.38	0.433954480957727\\
6.39	0.433955012486206\\
6.4	0.433955542355336\\
6.41	0.433956070572971\\
6.42	0.433956597146754\\
6.43	0.433957122084414\\
6.44	0.433957645393444\\
6.45	0.433958167081513\\
6.46	0.433958687156134\\
6.47	0.433959205624761\\
6.48	0.433959722494751\\
6.49	0.433960237773647\\
6.5	0.433960751468629\\
6.51	0.4339612635871\\
6.52	0.433961774136257\\
6.53	0.433962283123226\\
6.54	0.433962790555319\\
6.55	0.433963296439592\\
6.56	0.433963800783033\\
6.57	0.433964303592803\\
6.58	0.433964804875788\\
6.59	0.433965304638957\\
6.6	0.433965802889219\\
6.61	0.433966299633464\\
6.62	0.433966794878472\\
6.63	0.433967288631005\\
6.64	0.433967780897845\\
6.65	0.433968271685681\\
6.66	0.433968761001075\\
6.67	0.433969248850797\\
6.68	0.433969735241325\\
6.69	0.433970220179246\\
6.7	0.433970703671\\
6.71	0.433971185723081\\
6.72	0.433971666341908\\
6.73	0.433972145533821\\
6.74	0.433972623305269\\
6.75	0.433973099662411\\
6.76	0.433973574611688\\
6.77	0.43397404815922\\
6.78	0.433974520311176\\
6.79	0.433974991073817\\
6.8	0.433975460453209\\
6.81	0.433975928455453\\
6.82	0.433976395086587\\
6.83	0.433976860352617\\
6.84	0.433977324259606\\
6.85	0.433977786813368\\
6.86	0.433978248019953\\
6.87	0.433978707885134\\
6.88	0.433979166414749\\
6.89	0.433979623614716\\
6.9	0.433980079490708\\
6.91	0.433980534048557\\
6.92	0.433980987293854\\
6.93	0.433981439232354\\
6.94	0.4339818898697\\
6.95	0.433982339211487\\
6.96	0.433982787263317\\
6.97	0.433983234030721\\
6.98	0.433983679519215\\
6.99	0.433984123734206\\
7	0.43398456668125\\
7.01	0.433985008365755\\
7.02	0.433985448793091\\
7.03	0.433985887968627\\
7.04	0.433986325897704\\
7.05	0.43398676258553\\
7.06	0.433987198037549\\
7.07	0.433987632258797\\
7.08	0.433988065254669\\
7.09	0.433988497030198\\
7.1	0.433988927590621\\
7.11	0.433989356941056\\
7.12	0.433989785086517\\
7.13	0.433990212032207\\
7.14	0.433990637783042\\
7.15	0.433991062344078\\
7.16	0.43399148572031\\
7.17	0.433991907916677\\
7.18	0.43399232893806\\
7.19	0.433992748789425\\
7.2	0.433993167475613\\
7.21	0.433993585001469\\
7.22	0.43399400137184\\
7.23	0.433994416591424\\
7.24	0.433994830665105\\
7.25	0.433995243597527\\
7.26	0.43399565539342\\
7.27	0.433996066057514\\
7.28	0.433996475594453\\
7.29	0.433996884008837\\
7.3	0.43399729130534\\
7.31	0.433997697488471\\
7.32	0.43399810256283\\
7.33	0.43399850653293\\
7.34	0.433998909403312\\
7.35	0.433999311178433\\
7.36	0.433999711862704\\
7.37	0.434000111460729\\
7.38	0.434000509976704\\
7.39	0.434000907415216\\
7.4	0.434001303780493\\
7.41	0.434001699076856\\
7.42	0.434002093308805\\
7.43	0.434002486480405\\
7.44	0.434002878596156\\
7.45	0.434003269660081\\
7.46	0.434003659676528\\
7.47	0.434004048649754\\
7.48	0.434004436583789\\
7.49	0.434004823482953\\
7.5	0.434005209351282\\
7.51	0.434005594192896\\
7.52	0.434005978011859\\
7.53	0.434006360812326\\
7.54	0.434006742598325\\
7.55	0.434007123373827\\
7.56	0.434007503142888\\
7.57	0.434007881909429\\
7.58	0.434008259677517\\
7.59	0.434008636451067\\
7.6	0.434009012233904\\
7.61	0.434009387030042\\
7.62	0.434009760843311\\
7.63	0.434010133677628\\
7.64	0.434010505536762\\
7.65	0.434010876424539\\
7.66	0.434011246344821\\
7.67	0.434011615301335\\
7.68	0.434011983297799\\
7.69	0.434012350338119\\
7.7	0.434012716425848\\
7.71	0.434013081564743\\
7.72	0.434013445758523\\
7.73	0.434013809010828\\
7.74	0.43401417132527\\
7.75	0.434014532705544\\
7.76	0.43401489315524\\
7.77	0.43401525267789\\
7.78	0.434015611277144\\
7.79	0.43401596895651\\
7.8	0.43401632571954\\
7.81	0.434016681569707\\
7.82	0.434017036510589\\
7.83	0.434017390545597\\
7.84	0.434017743678231\\
7.85	0.434018095911951\\
7.86	0.434018447250133\\
7.87	0.434018797696241\\
7.88	0.434019147253622\\
7.89	0.434019495925681\\
7.9	0.434019843715751\\
7.91	0.434020190627184\\
7.92	0.434020536663372\\
7.93	0.434020881827531\\
7.94	0.434021226122952\\
7.95	0.434021569553008\\
7.96	0.434021912120864\\
7.97	0.434022253829795\\
7.98	0.434022594683041\\
7.99	0.434022934683818\\
8	0.43402327383525\\
8.01	0.434023612140601\\
8.02	0.434023949602979\\
8.03	0.434024286225569\\
8.04	0.434024622011479\\
8.05	0.434024956963823\\
8.06	0.434025291085721\\
8.07	0.434025624380219\\
8.08	0.434025956850366\\
8.09	0.43402628849937\\
8.1	0.434026619330121\\
8.11	0.43402694934566\\
8.12	0.43402727854901\\
8.13	0.434027606943177\\
8.14	0.434027934531133\\
8.15	0.434028261315893\\
8.16	0.434028587300341\\
8.17	0.434028912487458\\
8.18	0.434029236880072\\
8.19	0.434029560481237\\
8.2	0.434029883293748\\
8.21	0.434030205320516\\
8.22	0.434030526564398\\
8.23	0.434030847028251\\
8.24	0.43403116671484\\
8.25	0.434031485627134\\
8.26	0.434031803767804\\
8.27	0.434032121139764\\
8.28	0.434032437745745\\
8.29	0.434032753588511\\
8.3	0.434033068670756\\
8.31	0.434033382995394\\
8.32	0.434033696564975\\
8.33	0.434034009382295\\
8.34	0.434034321450049\\
8.35	0.434034632770959\\
8.36	0.434034943347661\\
8.37	0.434035253182839\\
8.38	0.434035562279108\\
8.39	0.434035870639221\\
8.4	0.434036178265632\\
8.41	0.434036485161108\\
8.42	0.434036791328122\\
8.43	0.434037096769467\\
8.44	0.434037401487485\\
8.45	0.434037705484842\\
8.46	0.434038008764159\\
8.47	0.43403831132789\\
8.48	0.434038613178526\\
8.49	0.434038914318714\\
8.5	0.434039214750886\\
8.51	0.434039514477539\\
8.52	0.43403981350111\\
8.53	0.434040111824176\\
8.54	0.434040409449074\\
8.55	0.434040706378429\\
8.56	0.434041002614451\\
8.57	0.434041298159701\\
8.58	0.434041593016632\\
8.59	0.434041887187536\\
8.6	0.434042180674823\\
8.61	0.434042473480963\\
8.62	0.434042765608272\\
8.63	0.434043057059018\\
8.64	0.434043347835697\\
8.65	0.434043637940563\\
8.66	0.434043927375899\\
8.67	0.434044216144132\\
8.68	0.434044504247513\\
8.69	0.434044791688316\\
8.7	0.434045078468857\\
8.71	0.434045364591433\\
8.72	0.434045650058172\\
8.73	0.434045934871477\\
8.74	0.434046219033543\\
8.75	0.434046502546537\\
8.76	0.434046785412791\\
8.77	0.434047067634399\\
8.78	0.434047349213648\\
8.79	0.434047630152706\\
8.8	0.434047910453801\\
8.81	0.434048190118983\\
8.82	0.434048469150476\\
8.83	0.434048747550441\\
8.84	0.434049025321071\\
8.85	0.434049302464432\\
8.86	0.43404957898265\\
8.87	0.434049854877831\\
8.88	0.434050130152073\\
8.89	0.434050404807478\\
8.9	0.43405067884616\\
8.91	0.434050952270253\\
8.92	0.434051225081662\\
8.93	0.434051497282497\\
8.94	0.434051768874914\\
8.95	0.434052039860855\\
8.96	0.434052310242319\\
8.97	0.434052580021448\\
8.98	0.434052849200101\\
8.99	0.434053117780376\\
9	0.434053385764271\\
9.01	0.434053653153764\\
9.02	0.434053919950826\\
9.03	0.434054186157323\\
9.04	0.434054451775378\\
9.05	0.434054716806835\\
9.06	0.434054981253619\\
9.07	0.434055245117815\\
9.08	0.434055508401145\\
9.09	0.434055771105662\\
9.1	0.434056033233228\\
9.11	0.434056294785768\\
9.12	0.434056555765094\\
9.13	0.434056816173161\\
9.14	0.434057076011802\\
9.15	0.434057335282895\\
9.16	0.434057593988358\\
9.17	0.434057852129967\\
9.18	0.434058109709613\\
9.19	0.43405836672903\\
9.2	0.434058623190136\\
9.21	0.434058879094765\\
9.22	0.434059134444655\\
9.23	0.434059389241615\\
9.24	0.434059643487512\\
9.25	0.434059897184091\\
9.26	0.43406015033313\\
9.27	0.43406040293635\\
9.28	0.434060654995571\\
9.29	0.434060906512528\\
9.3	0.434061157488944\\
9.31	0.434061407926695\\
9.32	0.434061657827348\\
9.33	0.434061907192682\\
9.34	0.434062156024408\\
9.35	0.434062404324236\\
9.36	0.434062652093958\\
9.37	0.43406289933515\\
9.38	0.434063146049615\\
9.39	0.434063392238906\\
9.4	0.434063637904766\\
9.41	0.434063883048832\\
9.42	0.434064127672891\\
9.43	0.434064371778408\\
9.44	0.434064615367135\\
9.45	0.434064858440643\\
9.46	0.43406510100066\\
9.47	0.434065343048778\\
9.48	0.434065584586613\\
9.49	0.434065825615781\\
9.5	0.434066066137877\\
9.51	0.43406630615445\\
9.52	0.434066545667163\\
9.53	0.434066784677584\\
9.54	0.434067023187255\\
9.55	0.434067261197848\\
9.56	0.434067498710772\\
9.57	0.4340677357277\\
9.58	0.434067972250175\\
9.59	0.434068208279762\\
9.6	0.434068443817843\\
9.61	0.434068678866138\\
9.62	0.434068913426123\\
9.63	0.434069147499282\\
9.64	0.434069381087165\\
9.65	0.43406961419127\\
9.66	0.43406984681301\\
9.67	0.434070078954058\\
9.68	0.434070310615754\\
9.69	0.434070541799633\\
9.7	0.434070772507209\\
9.71	0.434071002739848\\
9.72	0.434071232499202\\
9.73	0.434071461786523\\
9.74	0.434071690603376\\
9.75	0.434071918951141\\
9.76	0.434072146831348\\
9.77	0.434072374245271\\
9.78	0.434072601194551\\
9.79	0.4340728276805\\
9.8	0.434073053704538\\
9.81	0.434073279268155\\
9.82	0.434073504372591\\
9.83	0.434073729019369\\
9.84	0.434073953209876\\
9.85	0.434074176945516\\
9.86	0.434074400227579\\
9.87	0.434074623057481\\
9.88	0.434074845436637\\
9.89	0.434075067366419\\
9.9	0.434075288848114\\
9.91	0.434075509883174\\
9.92	0.434075730472888\\
9.93	0.434075950618623\\
9.94	0.434076170321697\\
9.95	0.434076389583422\\
9.96	0.434076608405166\\
9.97	0.434076826788211\\
9.98	0.434077044734001\\
9.99	0.434077262243712\\
};
\addlegendentry{$k=3$}

\addplot [color=red, dotted, mark=square, mark options={solid, red}]
  table[row sep=crcr]{%
1	0.413926851582251\\
1.1	0.415674169783397\\
1.2	0.417145275110543\\
1.3	0.418400883828216\\
1.4	0.419485127284205\\
1.5	0.420430854003653\\
1.6	0.421263019557586\\
1.7	0.422000923325546\\
1.8	0.422659725291412\\
1.9	0.423251499511892\\
2	0.423785981398762\\
2.1	0.424271107854027\\
2.2	0.424713414298506\\
2.3	0.425118330962302\\
2.4	0.425490407050359\\
2.5	0.425833482469509\\
2.6	0.426150820892404\\
2.7	0.426445213947261\\
2.8	0.426719063588633\\
2.9	0.426974447800485\\
3	0.427213173438309\\
3.1	0.427436819054632\\
3.2	0.427646769855408\\
3.3	0.427844246424132\\
3.4	0.428030328472696\\
3.5	0.428205974595405\\
3.6	0.428372038789477\\
3.7	0.42852928434316\\
3.8	0.428678395568185\\
3.9	0.428819987757063\\
4	0.428954615670923\\
4.1	0.429082780804764\\
4.2	0.429204937630813\\
4.3	0.42932149898384\\
4.4	0.429432840722875\\
4.5	0.429539305780366\\
4.6	0.429641207690598\\
4.7	0.429738833673877\\
4.8	0.429832447340468\\
4.9	0.429922291067752\\
5	0.430008588095878\\
5.1	0.430091544380001\\
5.2	0.430171350231472\\
5.3	0.430248181775511\\
5.4	0.43032220224887\\
5.5	0.430393563157609\\
5.6	0.430462405312293\\
5.7	0.430528859755417\\
5.8	0.430593048593999\\
5.9	0.430655085748466\\
6	0.430715077627403\\
6.1	0.430773123736698\\
6.2	0.430829317230325\\
6.3	0.430883745409243\\
6.4	0.430936490173978\\
6.5	0.43098762843585\\
6.6	0.431037232491212\\
6.7	0.431085370362463\\
6.8	0.43113210610929\\
6.9	0.431177500113106\\
7	0.43122160933729\\
7.1	0.431264487565715\\
7.2	0.431306185621494\\
7.3	0.431346751567983\\
7.4	0.431386230893568\\
7.5	0.431424666681851\\
7.6	0.431462099768477\\
7.7	0.431498568885883\\
7.8	0.431534110796958\\
7.9	0.431568760418673\\
8	0.431602550936492\\
8.1	0.431635513910424\\
8.2	0.431667679373292\\
8.3	0.431699075922041\\
8.4	0.431729730802529\\
8.5	0.431759669988372\\
8.6	0.43178891825437\\
8.7	0.431817499244862\\
8.8	0.43184543553749\\
8.9	0.431872748702677\\
9	0.431899459359184\\
9.1	0.43192558722601\\
9.2	0.431951151170948\\
9.3	0.431976169256009\\
9.4	0.432000658780019\\
9.5	0.432024636318463\\
9.6	0.432048117760941\\
9.7	0.432071118346249\\
9.8	0.432093652695392\\
9.9	0.432115734842561\\
};
\addlegendentry{$k=2$}

\addplot [color=mycolor1, dotted, mark=square, mark options={solid, mycolor1}]
  table[row sep=crcr]{%
1	0.301029995663981\\
2	0.352182518111362\\
3	0.3748162098249\\
4	0.387640052032226\\
5	0.395906230238124\\
6	0.401680737783679\\
7	0.405943628843807\\
8	0.40922017957905\\
9	0.411817415046076\\
};
\addlegendentry{$k=1$}

\node[right, align=left]
at (axis cs:-1.5,0.434) {$(\ln 10)^{-1} \longrightarrow$};
\end{axis}
\end{tikzpicture}%

%% file: fig2/general_pareto_vs_benford_approx_s1_00_rho0_30_mincha_291122_202751.tex
%
%
\definecolor{mycolor1}{rgb}{0.00000,0.44700,0.74100}%
\definecolor{mycolor2}{rgb}{0.85000,0.32500,0.09800}%
\begin{tikzpicture}

\begin{axis}[%
width=0.951\fwidth,
height=\fheight,
at={(0\fwidth,0\fheight)},
scale only axis,
xmin=10,
xmax=100,
xlabel style={font=\color{white!15!black}},
xlabel={$a$},
ymin=0,
ymax=0.06,
ylabel style={font=\color{white!15!black}},
ylabel={$\Pr(A_{(2)}=a)$},
axis background/.style={fill=white},
legend style={legend cell align=left, align=left, draw=white!15!black},
grid style={dotted},grid = major,legend style={font=\scriptsize,draw=none},mark size=1pt,yticklabel style={/pgf/number format/fixed,/pgf/number format/precision=5},scaled ticks=false
]
\addplot [color=mycolor1, line width=1.0pt]
  table[row sep=crcr]{%
10	0.0201541647976654\\
11	0.0167951373313879\\
12	0.0142112700496359\\
13	0.0121810886139736\\
14	0.0105569434654438\\
15	0.00923732553226333\\
16	0.00815058135199706\\
17	0.00724496120177516\\
18	0.00648233370685146\\
19	0.00820294285172688\\
20	0.0527847173272191\\
21	0.0479861066611081\\
22	0.0438134017340553\\
23	0.0401622849228841\\
24	0.0369493021290533\\
25	0.0341070481191261\\
26	0.031580600110302\\
27	0.0293248429595662\\
28	0.0273024399968374\\
29	0.0254822773303815\\
30	0.023838259438099\\
31	0.0223483682232177\\
32	0.0209939216642348\\
33	0.0197589850957505\\
34	0.0186299002331361\\
35	0.0175949057757397\\
36	0.0166438297878618\\
37	0.0157678387463954\\
38	0.0149592316311958\\
39	0.0142112700496357\\
40	0.0135180373642878\\
41	0.0128743212993218\\
42	0.0122755156574927\\
43	0.0117175376730614\\
44	0.0111967582209253\\
45	0.0107099426461023\\
46	0.0102542004058428\\
47	0.00982694205559931\\
48	0.00942584237986066\\
49	0.00904880868466607\\
50	0.0086939534421302\\
51	0.00835957061743298\\
52	0.00804411512243541\\
53	0.00774618493271551\\
54	0.00746450548061683\\
55	0.0071979159991663\\
56	0.00694535754305514\\
57	0.00670586245536364\\
58	0.00647854508399526\\
59	0.00626259358119552\\
60	0.006057262644107\\
61	0.00586186707494241\\
62	0.00567577605669023\\
63	0.00549840805491852\\
64	0.00532922626861354\\
65	0.00516773456350383\\
66	0.00501347383026507\\
67	0.00486601871761028\\
68	0.00472497469680988\\
69	0.00458997541975814\\
70	0.00446068033751146\\
71	0.00433677255035847\\
72	0.00421795686404713\\
73	0.00410395802988384\\
74	0.00399451914908677\\
75	0.00388940022411088\\
76	0.00378837684166644\\
77	0.00369123897393153\\
78	0.00359778988598369\\
79	0.00350784513883418\\
80	0.00342123167861619\\
81	0.00333778700352783\\
82	0.00325735840103314\\
83	0.00317980224862768\\
84	0.00310498337218934\\
85	0.00303277445655703\\
86	0.00296305550353282\\
87	0.00289571333299787\\
88	0.0028306411232677\\
89	0.00276773798719498\\
90	0.00270690858088296\\
91	0.00264806274216811\\
92	0.00259111515631516\\
93	0.00253598504660635\\
94	0.00248259588773025\\
95	0.00243087514006934\\
96	0.00238075400316062\\
97	0.00233216718676965\\
98	0.00228505269814794\\
99	0.00223935164418509\\
};
\addlegendentry{Pareto, $s=1.00$, $\rho=0.30$}

\addplot [color=mycolor2, dashed, line width=1.0pt]
  table[row sep=crcr]{%
10	0.0211298677773864\\
11	0.0175364163112079\\
12	0.014787610945285\\
13	0.0126380208612893\\
14	0.0109253076570956\\
15	0.00953861326589765\\
16	0.00840014050483144\\
17	0.00745398883125785\\
18	0.006659153646646\\
19	0.00598500452058972\\
20	0.0540828755046737\\
21	0.0491109175069127\\
22	0.0447944074199641\\
23	0.0410229926518984\\
24	0.0377086052259867\\
25	0.0347802715204207\\
26	0.0321802795635642\\
27	0.029861312715368\\
28	0.0277842809531541\\
29	0.0259166622518993\\
30	0.0242312214048964\\
31	0.022705011217347\\
32	0.0213185871204372\\
33	0.0200553846339102\\
34	0.0189012221983565\\
35	0.0178439013305631\\
36	0.0168728829217129\\
37	0.0159790235463409\\
38	0.0151543593957347\\
39	0.0143919282531828\\
40	0.0136856220443595\\
41	0.0130300641053437\\
42	0.0124205065436595\\
43	0.0118527440188713\\
44	0.0113230410079142\\
45	0.0108280701975774\\
46	0.0103648601003323\\
47	0.0099307503484761\\
48	0.00952335340674502\\
49	0.0091405216714269\\
50	0.008780319106987\\
51	0.00844099671884637\\
52	0.00812097128059964\\
53	0.00781880683134879\\
54	0.00753319853843436\\
55	0.00726295858617237\\
56	0.00700700380502395\\
57	0.00676434480012229\\
58	0.00653407637501116\\
59	0.00631536907719371\\
60	0.00610746171777628\\
61	0.00590965473902153\\
62	0.00572130432171907\\
63	0.00554181713954242\\
64	0.0053706456804617\\
65	0.00520728406622217\\
66	0.00505126431020054\\
67	0.00490215296187767\\
68	0.00475954809294368\\
69	0.0046230765858544\\
70	0.00449239169064118\\
71	0.00436717082006798\\
72	0.00424711355692527\\
73	0.00413193985045268\\
74	0.00402138838164888\\
75	0.00391521507963468\\
76	0.0038131917733243\\
77	0.00371510496448398\\
78	0.00362075470984653\\
79	0.00352995360134307\\
80	0.00344252583473171\\
81	0.0033583063579734\\
82	0.00327714009164482\\
83	0.00319888121450862\\
84	0.00312339250808994\\
85	0.00305054475475464\\
86	0.00298021618435594\\
87	0.00291229196502132\\
88	0.00284666373410127\\
89	0.0027832291657002\\
90	0.00272189157156449\\
91	0.00266255953241954\\
92	0.00260514655712944\\
93	0.00254957076730514\\
94	0.00249575460521289\\
95	0.00244362456303547\\
96	0.00239311093172157\\
97	0.00234414756781951\\
98	0.00229667167683904\\
99	0.00225062361181506\\
};
\addlegendentry{Asymptotic approximation}

\end{axis}
\end{tikzpicture}%

%% file: fig2/cf_general_pareto_vs_benford_approx_s1_00_rho0_30_mincha_271122_110309.tex
%
%
\definecolor{mycolor1}{rgb}{1.00000,0.00000,1.00000}%
\begin{tikzpicture}

\begin{axis}[%
width=0.964\fwidth,
height=\fheight,
at={(0\fwidth,0\fheight)},
scale only axis,
xmode=log,
xmin=1,
xmax=100,
xminorticks=true,
xlabel style={font=\color{white!15!black}},
xlabel={$a_1$},
ymode=log,
ymin=1e-07,
ymax=1,
yminorticks=true,
ylabel style={font=\color{white!15!black}},
ylabel={$\Pr(\mathbf{A}_2=[a_1,a_2])$},
axis background/.style={fill=white},
xmajorgrids,
xminorgrids,
ymajorgrids,
yminorgrids,
grid style={dotted},legend style={font=\tiny,draw=none},mark size=1pt
]
\addplot [color=blue, forget plot]
  table[row sep=crcr]{%
1	0.223434566413295\\
2	0.146433881560569\\
3	0.00984193226904524\\
4	0.00697745319060151\\
5	0.00517855799481408\\
6	0.0039864995940912\\
7	0.00315959871915156\\
8	0.0025639335301859\\
9.00000000000001	0.00212124414697195\\
10	0.00178356690113821\\
12	0.00131107998697417\\
14	0.00100384887673886\\
17	0.00071167820390655\\
20	0.000530603841038876\\
24	0.000379689853185142\\
30	0.000250472010824275\\
38	0.000160180753357267\\
49	9.84616193538858e-05\\
66.0000000000001	5.53433256946785e-05\\
93.0000000000001	2.83358977641613e-05\\
100	2.45770712693657e-05\\
};
\addplot [color=blue, dashed, forget plot]
  table[row sep=crcr]{%
1	0.0831612366288608\\
2	0.0332644946515443\\
3	0.0178202649918987\\
4	0.0110881648838481\\
5	0.00756011242080554\\
6	0.00548315845904578\\
7	0.00415806183144304\\
9.00000000000001	0.0026261443145956\\
11	0.00180785297019263\\
14	0.00114705153970843\\
18	0.000709768733674491\\
24	0.00040732034267197\\
34	0.000206611768022015\\
51	9.31604592556307e-05\\
84.0000000000001	3.47349404645433e-05\\
100	2.45784650890676e-05\\
};
\node[right, align=left, font=\color{blue}]
at (axis cs:1.02,0.216) {{\tiny $a_2=1$}};
\addplot [color=red, forget plot]
  table[row sep=crcr]{%
1	0.0833920505273913\\
2	0.056195007953752\\
3	0.00371567376450626\\
4	0.00258999305551124\\
5	0.00189585237868662\\
6	0.00144346864049621\\
7	0.00113399629475263\\
8	0.000913610854977156\\
9.00000000000001	0.000751378698030596\\
10	0.000628617218564745\\
12	0.000458448987708614\\
14	0.000348931604862816\\
17	0.000245744743959702\\
21	0.000166583295106776\\
26	0.000111741618504226\\
33	7.11169308136029e-05\\
43.0000000000001	4.28033598867978e-05\\
58.0000000000001	2.39679467939794e-05\\
83	1.18940203041608e-05\\
100	8.24616676029469e-06\\
};
\addplot [color=red, dashed, forget plot]
  table[row sep=crcr]{%
1	0.0415759024146026\\
2	0.014254595113578\\
3	0.00712729755678901\\
4	0.00426419511944643\\
5	0.00283472061917744\\
6	0.00201988189868515\\
8	0.0011739078328829\\
11	0.000637993387436353\\
15	0.000349867341497357\\
22	0.000165476228515832\\
34	7.019991965319e-05\\
59.0000000000001	2.35535279470887e-05\\
100	8.24632368019157e-06\\
};
\node[right, align=left, font=\color{red}]
at (axis cs:1.02,0.081) {{\tiny $a_2=2$}};
\addplot [color=green, forget plot]
  table[row sep=crcr]{%
1	0.0428729146852638\\
2	0.0296584783610482\\
3	0.0195069490695281\\
4	0.00135067603954049\\
5	0.000983210836386664\\
6	0.000745302583115896\\
7	0.000583452359950111\\
8	0.000468719812277504\\
9.00000000000001	0.000384581182638713\\
11	0.000272094459085355\\
13	0.000202489025651518\\
16	0.000139149700649143\\
19	0.000101442132521922\\
23	7.10326873039371e-05\\
29	4.58316164334238e-05\\
38	2.73192350475038e-05\\
51	1.54597419375942e-05\\
72.0000000000001	7.88465590618802e-06\\
100	4.13324730435533e-06\\
};
\addplot [color=green, dashed, forget plot]
  table[row sep=crcr]{%
1	0.0249445920737287\\
2	0.00791891811864404\\
3	0.00383762954980443\\
4	0.00225742914694377\\
5	0.00148479714724575\\
7	0.000781962133972677\\
10	0.000392519151435546\\
14	0.000203546242951681\\
21	9.17080590945905e-05\\
35	3.33796227401706e-05\\
69.0000000000001	8.65891143908926e-06\\
100	4.13328672897899e-06\\
};
\node[right, align=left, font=\color{green}]
at (axis cs:1.02,0.042) {{\tiny $a_2=3$}};
\addplot [color=black, forget plot]
  table[row sep=crcr]{%
1	0.025959232838127\\
2	0.0183169358333259\\
3	0.0120175491781317\\
4	0.00082910787097755\\
5	0.000601712338383938\\
6	0.000455015797703827\\
7	0.000355518992836229\\
8	0.000285163387228677\\
9.00000000000001	0.000233674744492372\\
11	0.000165000534855845\\
13	0.000122613628596866\\
16	8.41281577732973e-05\\
20	5.56830422212613e-05\\
25	3.66132272271391e-05\\
32	2.28836747733889e-05\\
42	1.35556815757237e-05\\
58.0000000000001	7.23684903351901e-06\\
85	3.42038941632239e-06\\
100	2.48321046392942e-06\\
};
\addplot [color=black, dashed, forget plot]
  table[row sep=crcr]{%
1	0.0166294108106982\\
2	0.00503921539718127\\
3	0.00239847271308146\\
4	0.0013974294798906\\
5	0.0009137038906977\\
7	0.000477856632491331\\
10	0.000238585520956937\\
15	0.000107610509991572\\
25	3.92018170926394e-05\\
48	1.07256535661193e-05\\
100	2.48322469435679e-06\\
};
\node[right, align=left]
at (axis cs:1.02,0.025) {{\tiny $a_2=4$}};
\addplot [color=mycolor1, forget plot]
  table[row sep=crcr]{%
1	0.0173596186537179\\
2	0.0124336854261938\\
3	0.00814621554169855\\
4	0.000560746543067481\\
5	0.000406172246901785\\
6	0.000306678684664647\\
7	0.000239326272448378\\
8	0.00019177551899323\\
9.00000000000001	0.000157021597856821\\
11	0.00011073664966964\\
13	8.22145500952523e-05\\
16	5.63539542521062e-05\\
20	3.72666686656425e-05\\
25	2.44859448845784e-05\\
32	1.52938914199109e-05\\
42	9.05449334950048e-06\\
58.0000000000001	4.83134120099065e-06\\
85	2.2824606235854e-06\\
100	1.65683507864982e-06\\
};
\addplot [color=mycolor1, dashed, forget plot]
  table[row sep=crcr]{%
1	0.011878014440894\\
2	0.00348864759802484\\
3	0.00164104146880774\\
4	0.000950241155271524\\
6	0.000434940371854874\\
9.00000000000001	0.000197184429453575\\
14	8.26638950319032e-05\\
24	2.84341183538064e-05\\
49	6.87441927129169e-06\\
100	1.65684141373664e-06\\
};
\node[right, align=left, font=\color{mycolor1}]
at (axis cs:1.02,0.017) {{\tiny $a_2=5$}};
\addplot [color=blue, forget plot]
  table[row sep=crcr]{%
1	0.012409215215885\\
2	0.00899163205492459\\
3	0.00588588538175524\\
4	0.00040452531062718\\
5	0.000292625821501646\\
6	0.000220712336140405\\
7	0.000172094621806185\\
8	0.000137807945237205\\
9.00000000000001	0.000112771035680728\\
11	7.94613360809858e-05\\
13	5.89576413164054e-05\\
16	4.0385151454116e-05\\
20	2.66902482443598e-05\\
25	1.75278904814682e-05\\
32	1.09429333779094e-05\\
43.0000000000001	6.18719860315369e-06\\
60	3.2326604734186e-06\\
89.0000000000001	1.49011009375716e-06\\
100	1.18412115859838e-06\\
};
\addplot [color=blue, dashed, forget plot]
  table[row sep=crcr]{%
1	0.00890844270560349\\
2	0.00255832200776304\\
3	0.00119347557778418\\
4	0.000688100402087988\\
6	0.000313559265564926\\
9.00000000000001	0.000141725224861869\\
15	5.1718099887392e-05\\
28	1.49843147662818e-05\\
66.0000000000001	2.71405297568632e-06\\
100	1.18412439446876e-06\\
};
\addplot [color=red, forget plot]
  table[row sep=crcr]{%
1	0.00930496503154787\\
2	0.00680466929233564\\
3	0.00445163858284649\\
4	0.000305604376824645\\
5	0.000220853185625181\\
6	0.000166448798035304\\
7	0.000129703857474208\\
8	0.000103810911330403\\
9.00000000000001	8.49157654922563e-05\\
11	5.9796113397836e-05\\
13	4.43462945747249e-05\\
16	3.03615418573948e-05\\
20	2.00567598835275e-05\\
25	1.31667197842275e-05\\
32	8.2174634286671e-06\\
43.0000000000001	4.64476273731978e-06\\
60	2.42614335966938e-06\\
89.0000000000001	1.11809830008412e-06\\
100	8.88456195278552e-07\\
};
\addplot [color=red, dashed, forget plot]
  table[row sep=crcr]{%
1	0.00692875090241757\\
2	0.00195635319597673\\
3	0.000907036481771037\\
4	0.000521285334351178\\
6	0.000236767947306156\\
9.00000000000001	0.000106778695413967\\
15	3.88952179147111e-05\\
29	1.04954570599595e-05\\
74	1.6209342293805e-06\\
100	8.88458016947021e-07\\
};
\node[right, align=left]
at (axis cs:1.02,0.009) {\tiny $\phantom{a_2}~~~~\!\pmb{\pmb{\vdots}}$};
\addplot [color=green, forget plot]
  table[row sep=crcr]{%
1	0.00723266033392212\\
2	0.00532887279813759\\
3	0.00348466512841581\\
4	0.000239014273407441\\
5	0.000172601456130322\\
6	0.00013000605743403\\
7	0.000101258250189732\\
8	8.10129915560573e-05\\
9.00000000000001	6.6246646736093e-05\\
11	4.66272018347524e-05\\
13	3.45677444102361e-05\\
16	2.36577387875908e-05\\
20	1.56229171450197e-05\\
25	1.02531351218909e-05\\
32	6.39745259708987e-06\\
43.0000000000001	3.61518242936408e-06\\
60	1.88797787133857e-06\\
90.0000000000001	8.50982256403725e-07\\
100	6.91238093623881e-07\\
};
\addplot [color=green, dashed, forget plot]
  table[row sep=crcr]{%
1	0.00554297799142305\\
2	0.00154448303166586\\
3	0.000712668598897232\\
4	0.000408573316321113\\
6	0.000185108726986305\\
10	6.7679828955112e-05\\
18	2.11071723811335e-05\\
39	4.52791914053905e-06\\
100	6.91239196326064e-07\\
};
\addplot [color=black, forget plot]
  table[row sep=crcr]{%
1	0.00578153246317426\\
2	0.0042861122713263\\
3	0.00280187639946466\\
4	0.000192050215590845\\
5	0.000138605315667491\\
6	0.000104350676362229\\
7	8.1245595453123e-05\\
8	6.49819997760976e-05\\
9.00000000000001	5.31244971854621e-05\\
11	3.73770615922273e-05\\
13	2.77023329413416e-05\\
16	1.89534781000554e-05\\
20	1.25129885366103e-05\\
25	8.21030319894697e-06\\
32	5.12180507196215e-06\\
43.0000000000001	2.89377871099282e-06\\
60	1.51099854046742e-06\\
90.0000000000001	6.80972288684678e-07\\
100	5.53126995965745e-07\\
};
\addplot [color=black, dashed, forget plot]
  table[row sep=crcr]{%
1	0.00453514934175954\\
2	0.00125029179848006\\
3	0.000574730907365847\\
4	0.000328850644425546\\
6	0.000148693421041304\\
10	5.4277709454197e-05\\
19	1.51852680991532e-05\\
45	2.7244679453084e-06\\
100	5.5312770203597e-07\\
};
\addplot [color=mycolor1, forget plot]
  table[row sep=crcr]{%
1	0.00472635094854534\\
2	0.00352209273774984\\
3	0.00230185321646573\\
4	0.000157689724789656\\
5	0.000113752640485979\\
6	8.56074935790538e-05\\
7	6.66323108688031e-05\\
8	5.32809558893967e-05\\
9.00000000000001	4.35498569852495e-05\\
11	3.06311539514293e-05\\
13	2.26974468817327e-05\\
16	1.55254676285156e-05\\
20	1.02476043290119e-05\\
25	6.7226812612112e-06\\
32	4.19311120997058e-06\\
43.0000000000001	2.36871747389017e-06\\
60	1.23667919346113e-06\\
90.0000000000001	5.57282486257788e-07\\
100	4.52648743487946e-07\\
};
\addplot [color=mycolor1, dashed, forget plot]
  table[row sep=crcr]{%
1	0.00377928146948164\\
2	0.00103284710967201\\
3	0.000473306597695999\\
5	0.000174672672959219\\
8	6.92003265323184e-05\\
14	2.28261337895915e-05\\
30	5.00712784145328e-06\\
96.0000000000001	4.91116804151969e-07\\
100	4.52649216334634e-07\\
};
\addplot [color=blue, forget plot]
  table[row sep=crcr]{%
1	0.00393532612235359\\
2	0.00294560962307866\\
3	0.00192471189679738\\
4	0.000131793321504615\\
5	9.50343008444285e-05\\
6	7.14980483774322e-05\\
7	5.56363457325384e-05\\
8	4.44793074635931e-05\\
9.00000000000001	3.63496982976424e-05\\
11	2.55603595614089e-05\\
13	1.89365138565914e-05\\
16	1.29503371418281e-05\\
20	8.54635088692597e-06\\
25	5.60578594770521e-06\\
32	3.49600940446854e-06\\
43.0000000000001	1.97467443255225e-06\\
60	1.03084638808143e-06\\
90.0000000000001	4.64486894629093e-07\\
100	3.77269379175662e-07\\
};
\addplot [color=blue, dashed, forget plot]
  table[row sep=crcr]{%
1	0.00319784686993153\\
2	0.000867589759494456\\
3	0.00039655334794061\\
5	0.000146037503427801\\
8	5.77857189516224e-05\\
15	1.66033452609136e-05\\
35	3.0698196479476e-06\\
100	3.7726970766049e-07\\
};
\addplot [color=red, forget plot]
  table[row sep=crcr]{%
1	0.00332721229244026\\
2	0.00249993591749798\\
3	0.00163323750942573\\
4	0.000111792007444329\\
5	8.05849443627916e-05\\
6	6.06112120725072e-05\\
7	4.71548120916177e-05\\
8	3.76922197847531e-05\\
9.00000000000001	3.07988185564145e-05\\
11	2.16524665134063e-05\\
13	1.6038832466061e-05\\
16	1.09668277634726e-05\\
20	7.23627643738266e-06\\
25	4.74588137294413e-06\\
32	2.95940531343747e-06\\
43.0000000000001	1.67140636017339e-06\\
60	8.72453807084349e-07\\
90.0000000000001	3.93087616753623e-07\\
100	3.19272035939626e-07\\
};
\addplot [color=red, dashed, forget plot]
  table[row sep=crcr]{%
1	0.00274100682974984\\
2	0.000739056656317733\\
3	0.000337069758793559\\
5	0.00012391039319783\\
9.00000000000001	3.87858220350264e-05\\
18	9.7825912935479e-06\\
48	1.38332961667844e-06\\
100	3.19272271194994e-07\\
};
\addplot [color=green, forget plot]
  table[row sep=crcr]{%
1	0.00284972976859243\\
2	0.00214828315258208\\
3	0.00140331144469059\\
4	9.60227053223097e-05\\
5	6.91981528297931e-05\\
6	5.20350219616478e-05\\
7	4.04753780492118e-05\\
8	3.23484721122485e-05\\
9.00000000000001	2.64292351741504e-05\\
11	1.85771381161898e-05\\
13	1.375898851107e-05\\
16	9.40660235508973e-06\\
20	6.2059904858131e-06\\
25	4.06974119735604e-06\\
32	2.53754091135373e-06\\
43.0000000000001	1.4330194127293e-06\\
61	7.24124146507023e-07\\
92.0000000000001	3.22674988234331e-07\\
100	2.7369392608853e-07\\
};
\addplot [color=green, dashed, forget plot]
  table[row sep=crcr]{%
1	0.00237553575173067\\
2	0.000637116868280254\\
3	0.000290036341781077\\
5	0.000106458068259375\\
9.00000000000001	3.32885698560955e-05\\
18	8.3905896537402e-06\\
50	1.09315528587464e-06\\
100	2.73694098957454e-07\\
};
\addplot [color=black, forget plot]
  table[row sep=crcr]{%
1	0.00246801012197562\\
2	0.00186594488634501\\
3	0.00121874470142336\\
4	8.33701646689787e-05\\
5	6.00655532663644e-05\\
6	4.51587944150439e-05\\
7	3.51212797682885e-05\\
8	2.80659012134677e-05\\
9.00000000000001	2.29279581949695e-05\\
11	1.61135555268807e-05\\
13	1.19329919048268e-05\\
16	8.15722000320595e-06\\
20	5.38111642188607e-06\\
25	3.52848592522835e-06\\
32	2.19988027586247e-06\\
43.0000000000001	1.24223802372982e-06\\
61	6.27677676930063e-07\\
92.0000000000001	2.79682420658486e-07\\
100	2.3722543500885e-07\\
};
\addplot [color=black, dashed, forget plot]
  table[row sep=crcr]{%
1	0.00207859115691619\\
2	0.000554907539109994\\
3	0.000252205195985806\\
5	9.24503109080533e-05\\
9.00000000000001	2.88826932410765e-05\\
19	6.53285505434701e-06\\
57	7.2938994663255e-07\\
100	2.37225564902649e-07\\
};
\addplot [color=mycolor1, forget plot]
  table[row sep=crcr]{%
1	0.00215807899014442\\
2	0.00163582902947906\\
3	0.00106834249147851\\
4	7.30637991802853e-05\\
5	5.26289769644979e-05\\
6	3.95610896146546e-05\\
7	3.0763638121338e-05\\
8	2.45809760743991e-05\\
9.00000000000001	2.00792180648583e-05\\
11	1.41095564674935e-05\\
13	1.04478758605511e-05\\
16	7.14125065264078e-06\\
20	4.71045207292584e-06\\
25	3.0884745610468e-06\\
32	1.9254120269221e-06\\
43.0000000000001	1.08717761517768e-06\\
61	5.49296759931889e-07\\
92.0000000000001	2.44745555342641e-07\\
100	2.07590564214061e-07\\
};
\addplot [color=mycolor1, dashed, forget plot]
  table[row sep=crcr]{%
1	0.00183404901260344\\
2	0.000487645485267024\\
3	0.000221322684750711\\
5	8.10366035458244e-05\\
9.00000000000001	2.52972277600484e-05\\
19	5.71891931019281e-06\\
59.0000000000001	5.95818949000297e-07\\
100	2.07590663658214e-07\\
};
\addplot [color=blue, forget plot]
  table[row sep=crcr]{%
1	0.00190301547486038\\
2	0.00144580553046769\\
3	0.000944162888165552\\
4	6.4557307330559e-05\\
5	4.64929464097663e-05\\
6	3.49434452955663e-05\\
7	2.71696282838621e-05\\
8	2.1707189751463e-05\\
9.00000000000001	1.77303470876518e-05\\
11	1.24575189836049e-05\\
13	9.22376263987039e-06\\
16	6.30396035282028e-06\\
20	4.1578131003162e-06\\
25	2.72593836022615e-06\\
32	1.69929367023558e-06\\
43.0000000000001	9.5944442564646e-07\\
61	4.84734695200974e-07\\
92.0000000000001	2.15970128955395e-07\\
100	1.83182337373878e-07\\
};
\addplot [color=blue, dashed, forget plot]
  table[row sep=crcr]{%
1	0.00163026422690587\\
2	0.000431914158816641\\
3	0.000195785264298759\\
5	7.16136740501436e-05\\
9.00000000000001	2.2340387525003e-05\\
20	4.55743007494113e-06\\
67.0000000000001	4.07825945809587e-07\\
100	1.8318241480926e-07\\
};
\addplot [color=red, forget plot]
  table[row sep=crcr]{%
1	0.00169060232037998\\
2	0.00128707270987533\\
3	0.000840444918623028\\
4	5.74546462665802e-05\\
5	4.13709057379878e-05\\
6	3.10896890366205e-05\\
7	2.41706727039137e-05\\
8	1.93095364883213e-05\\
9.00000000000001	1.57708581329219e-05\\
11	1.10795815151352e-05\\
13	8.2028766690225e-06\\
16	5.60577019071684e-06\\
20	3.69704012596295e-06\\
25	2.42369704044148e-06\\
32	1.5107987376165e-06\\
43.0000000000001	8.52973366049671e-07\\
61	4.30923368162238e-07\\
92.0000000000001	1.919877488631e-07\\
100	1.62839915568113e-07\\
};
\addplot [color=red, dashed, forget plot]
  table[row sep=crcr]{%
1	0.0014586562328663\\
2	0.000385220410533023\\
3	0.000174426724349745\\
5	6.37439856427656e-05\\
9.00000000000001	1.98733340626316e-05\\
20	4.05244824688819e-06\\
70	3.32163731049827e-07\\
100	1.62839976759172e-07\\
};
\addplot [color=green, forget plot]
  table[row sep=crcr]{%
1	0.00151184071322708\\
2	0.00115311841405968\\
3	0.000752927290052427\\
4	5.14630240415129e-05\\
5	3.70511014302913e-05\\
6	2.78401351088176e-05\\
7	2.16422792359059e-05\\
8	1.72883393086996e-05\\
9.00000000000001	1.41191920643253e-05\\
11	9.91828497209194e-06\\
13	7.34259057063195e-06\\
16	5.01748505608936e-06\\
20	3.3088407873173e-06\\
25	2.16908238443525e-06\\
32	1.35201904732422e-06\\
43.0000000000001	7.63293524383901e-07\\
61	3.85601379981209e-07\\
92.0000000000001	1.71789922349549e-07\\
100	1.45707784917097e-07\\
};
\addplot [color=green, dashed, forget plot]
  table[row sep=crcr]{%
1	0.00131278962286195\\
2	0.000345710364994768\\
3	0.000156382462911443\\
5	5.71039442178933e-05\\
10	1.44300152349577e-05\\
24	2.52101038850379e-06\\
100	1.45707833916376e-07\\
};
\addplot [color=black, forget plot]
  table[row sep=crcr]{%
1	0.00135998449340135\\
2	0.0010390407260642\\
3	0.000678403107635062\\
4	4.63621987251024e-05\\
5	3.33743133449635e-05\\
6	2.50747519410143e-05\\
7	1.94908940312664e-05\\
8	1.55687078128988e-05\\
9.00000000000001	1.27140811965312e-05\\
11	8.93047680595784e-06\\
13	6.61089658140257e-06\\
16	4.51718750909564e-06\\
20	2.97873447544641e-06\\
25	1.95258719834614e-06\\
32	1.21702048925815e-06\\
43.0000000000001	6.87050453055181e-07\\
61	3.47072227185046e-07\\
92.0000000000001	1.54620144757384e-07\\
100	1.31144204365157e-07\\
};
\addplot [color=black, dashed, forget plot]
  table[row sep=crcr]{%
1	0.00118776124091171\\
2	0.000311982314686074\\
3	0.000141000486484742\\
5	5.14500537523618e-05\\
10	1.2994184084359e-05\\
25	2.09186551763414e-06\\
100	1.311442440759e-07\\
};
\addplot [color=mycolor1, forget plot]
  table[row sep=crcr]{%
1	0.00122989617183811\\
2	0.00094109271659476\\
3	0.000614421565211521\\
4	4.19839247356364e-05\\
5	3.02189546059383e-05\\
6	2.27019046329213e-05\\
7	1.76451120137224e-05\\
8	1.4093492217149e-05\\
9.00000000000001	1.15087783729374e-05\\
11	8.08323921078056e-06\\
13	5.9833824737695e-06\\
16	4.08816417354488e-06\\
20	2.69568062611832e-06\\
25	1.76696380025172e-06\\
32	1.10127981291609e-06\\
43.0000000000001	6.21687522816474e-07\\
61	3.14043026683075e-07\\
92.0000000000001	1.39901939460975e-07\\
100	1.18660156503711e-07\\
};
\addplot [color=mycolor1, dashed, forget plot]
  table[row sep=crcr]{%
1	0.00107978229266721\\
2	0.00028296053273525\\
3	0.000127781613527713\\
5	4.65962468907403e-05\\
10	1.17625007477413e-05\\
26	1.75046377277896e-06\\
100	1.18660188993894e-07\\
};
\addplot [color=blue, forget plot]
  table[row sep=crcr]{%
1	0.00111760784158044\\
2	0.000856370590985477\\
3	0.000559083664806531\\
4	3.81978862586592e-05\\
5	2.74908794257969e-05\\
6	2.06506575765088e-05\\
7	1.60496689196433e-05\\
8	1.28184684329665e-05\\
9.00000000000001	1.04671141917312e-05\\
11	7.35110812010695e-06\\
13	5.4411666273158e-06\\
16	3.71749017925857e-06\\
20	2.45114253120519e-06\\
25	1.60660885652257e-06\\
32	1.00130037337547e-06\\
43.0000000000001	5.65228559200751e-07\\
61	2.85514517596647e-07\\
92.0000000000001	1.27189783444415e-07\\
100	1.0787771423033e-07\\
};
\addplot [color=blue, dashed, forget plot]
  table[row sep=crcr]{%
1	0.000985887640290853\\
2	0.000257808344179379\\
3	0.000116338420239542\\
5	4.23983635888939e-05\\
10	1.06980151827604e-05\\
26	1.59160754994575e-06\\
100	1.07877741087633e-07\\
};
\addplot [color=red, forget plot]
  table[row sep=crcr]{%
1	0.00102001575572497\\
2	0.00078259583724688\\
3	0.00051089955233792\\
4	3.49018787846027e-05\\
5	2.51162722330785e-05\\
6	1.88654061050665e-05\\
7	1.46612524019291e-05\\
8	1.17089817681069e-05\\
9.00000000000001	9.56074917085457e-06\\
11	6.71413561091831e-06\\
13	4.96945964651974e-06\\
16	3.39504334208013e-06\\
20	2.23843549421271e-06\\
25	1.467135132065e-06\\
32	9.14344671562173e-07\\
43.0000000000001	5.16126568540949e-07\\
61	2.60704518036459e-07\\
92.0000000000001	1.16134961289807e-07\\
100	9.85010733007795e-08\\
};
\addplot [color=red, dashed, forget plot]
  table[row sep=crcr]{%
1	0.000903729886937558\\
2	0.000235867091058852\\
3	0.000106366502684332\\
5	3.87433129535235e-05\\
10	9.7717752363199e-06\\
27	1.34793941362675e-06\\
100	9.85010956830978e-08\\
};
\addplot [color=green, forget plot]
  table[row sep=crcr]{%
1	0.000934664028431342\\
2	0.000717960331042655\\
3	0.000468687107094079\\
4	3.20148331984795e-05\\
5	2.30365975612607e-05\\
6	1.73020651647935e-05\\
7	1.34455288452942e-05\\
8	1.07375652003694e-05\\
9.00000000000001	8.76722414120076e-06\\
11	6.15651552964944e-06\\
13	4.55654428492862e-06\\
16	3.11280518194105e-06\\
20	2.05226473479199e-06\\
25	1.34506794813628e-06\\
33	7.8999148846748e-07\\
45	4.33098332106356e-07\\
64.0000000000001	2.17511363926057e-07\\
95	9.99254552815554e-08\\
};
\addplot [color=green, dashed, forget plot]
  table[row sep=crcr]{%
1	0.000831431117970676\\
2	0.00021661253616264\\
4	5.52997085447797e-05\\
8	1.39716754175189e-05\\
20	2.24973582144267e-06\\
90.0000000000001	1.11466062325871e-07\\
96.0000000000001	9.79740090639256e-08\\
};
\addplot [color=black, forget plot]
  table[row sep=crcr]{%
1	0.000859588939944533\\
2	0.000661014429731044\\
3	0.000431498692539597\\
4	2.94717792634118e-05\\
5	2.12049553004288e-05\\
6	1.5925319434285e-05\\
7	1.23749987149551e-05\\
8	9.88222169187348e-06\\
9.00000000000001	8.06855442159201e-06\\
11	5.66559269349358e-06\\
13	4.19304019115145e-06\\
16	2.86435720254804e-06\\
20	1.88839243636485e-06\\
25	1.23762650261247e-06\\
33	7.26866015798853e-07\\
45	3.98480449758219e-07\\
64.0000000000001	2.00121128176511e-07\\
92.0000000000001	9.79481193431947e-08\\
};
\addplot [color=black, dashed, forget plot]
  table[row sep=crcr]{%
1	0.000767474558258644\\
2	0.000199623234440998\\
4	5.0919512388304e-05\\
8	1.28594968903729e-05\\
20	2.07011533219902e-06\\
92.0000000000001	9.81442836907267e-08\\
};
\addplot [color=mycolor1, forget plot]
  table[row sep=crcr]{%
1	0.000793205194383972\\
2	0.00061058494046795\\
3	0.000398567472061157\\
4	2.72201558355971e-05\\
5	1.95834114417859e-05\\
6	1.47066099780112e-05\\
7	1.14274266719676e-05\\
8	9.12516706050079e-06\\
9.00000000000001	7.45020107940979e-06\\
11	5.23113801294226e-06\\
13	3.87136607613131e-06\\
16	2.6445124861352e-06\\
20	1.74339427274835e-06\\
25	1.14256397940189e-06\\
33	6.71016067653685e-07\\
45	3.67853614538801e-07\\
64.0000000000001	1.84736274640419e-07\\
88.0000000000001	9.87071861094806e-08\\
};
\addplot [color=mycolor1, dashed, forget plot]
  table[row sep=crcr]{%
1	0.00071062431849063\\
2	0.000184557259186281\\
4	4.70399124545467e-05\\
8	1.18750332447808e-05\\
21	1.73382039461321e-06\\
88.0000000000001	9.90173598730147e-08\\
};
\end{axis}
\end{tikzpicture}%

%% file: fig2/b-p.tex
%
%
\definecolor{mycolor1}{rgb}{0.92900,0.69400,0.12500}%
\definecolor{mycolor2}{rgb}{0.46600,0.67400,0.18800}%
\definecolor{mycolor3}{rgb}{0.63500,0.07800,0.18400}%
\definecolor{mycolor4}{rgb}{0.85000,0.32500,0.09800}%
\begin{tikzpicture}

\begin{axis}[%
width=0.951\fwidth,
height=\fheight,
at={(0\fwidth,0\fheight)},
scale only axis,
xmin=1,
xmax=9,
xlabel style={font=\color{white!15!black}},
xlabel={$a$},
ymin=0,
ymax=0.68,
ylabel style={font=\color{white!15!black}},
ylabel={$\Pr(A_{(1)}=a)$},
grid = major,
grid style={dotted},
axis background/.style={fill=white},
xtick={1,2,3,4,5,6,7,8,9},
ytick={0,0.1,0.2,0.3,0.4,0.5,0.6},
legend style={font=\scriptsize, legend cell align=left, align=left, draw={none}}
]
\addplot [color=blue, solid, line width=0.5pt]
  table[row sep=crcr]{%
1	0.301029995663981\\
2	0.176091259055681\\
3	0.1249387366083\\
4	0.0969100130080564\\
5	0.0791812460476248\\
6	0.0669467896306132\\
7	0.0579919469776867\\
8	0.0511525224473813\\
9	0.0457574905606751\\
};
\addlegendentry{Benford}

\addplot [color=blue, line width=0.5pt, only marks, mark=square, mark options={scale=.8, solid, blue}, forget plot]
  table[row sep=crcr]{%
1	0.301455312089801\\
2	0.176096658282995\\
3	0.124831001854061\\
4	0.0967784900576364\\
5	0.0791407103352105\\
6	0.0668697529912082\\
7	0.0581776609240087\\
8	0.0509418252901496\\
9	0.0457085881749302\\
};
\addplot [color=green, solid, line width=0.5pt]
  table[row sep=crcr]{%
1	0.428349096755902\\
2	0.189765693235384\\
3	0.113122857795847\\
4	0.0771988125235018\\
5	0.0569857959998042\\
6	0.0442893348830856\\
7	0.0357006049715592\\
8	0.0295713342177599\\
9	0.025016469617157\\
};
\addlegendentry{Pareto, $s=0.5$, $\rho=0$}

\addplot [color=green, line width=0.5pt, only marks, mark=square, mark options={scale=.8, solid, green}, forget plot]
  table[row sep=crcr]{%
1	0.4282545\\
2	0.1898551\\
3	0.1131265\\
4	0.0772963\\
5	0.0568327\\
6	0.0443452\\
7	0.0357416\\
8	0.0295728\\
9	0.0249753\\
};
\addplot [color=green, dashed, line width=0.5pt]
  table[row sep=crcr]{%
1	0.240878398713803\\
2	0.106713091409651\\
3	0.174473952742775\\
4	0.137281058789998\\
5	0.10133666769113\\
6	0.0787588123069098\\
7	0.063485650746857\\
8	0.052586094766822\\
9	0.0444862728320546\\
};
\addlegendentry{Pareto, $s=0.5$, $\rho=0.5$}

\addplot [color=green, line width=0.5pt, only marks, mark=square, mark options={scale=.8, solid, green}, forget plot]
  table[row sep=crcr]{%
1	0.240994\\
2	0.1066473\\
3	0.1743769\\
4	0.1372341\\
5	0.1013532\\
6	0.0788312\\
7	0.0635164\\
8	0.052499\\
9	0.0445479\\
};
\addplot [color=red, solid, line width=0.5pt]
  table[row sep=crcr]{%
1	0.555555555555556\\
2	0.185185185185185\\
3	0.0925925925925926\\
4	0.0555555555555556\\
5	0.037037037037037\\
6	0.0264550264550265\\
7	0.0198412698412698\\
8	0.0154320987654321\\
9	0.0123456790123457\\
};
\addlegendentry{Pareto, $s=1.0$, $\rho=0$}

\addplot [color=red, line width=0.5pt, only marks, mark=square, mark options={scale=.8, solid, red}, forget plot]
  table[row sep=crcr]{%
1	0.5556575\\
2	0.1850763\\
3	0.0927131\\
4	0.0555032\\
5	0.0370379\\
6	0.0265062\\
7	0.0198132\\
8	0.0153732\\
9	0.0123194\\
};
\addplot [color=red, dashed, line width=0.5pt]
  table[row sep=crcr]{%
1	0.278437352015151\\
2	0.0928124506717171\\
3	0.0464062253358585\\
4	0.0278437352015151\\
5	0.18325043408889\\
6	0.13258921524531\\
7	0.0994419114339825\\
8	0.0773437088930977\\
9	0.0618749671144779\\
};
\addlegendentry{Pareto, $s=1.0$, $\rho=0.7$}

\addplot [color=red, line width=0.5pt, only marks, mark=square, mark options={scale=.8, solid, red}, forget plot]
  table[row sep=crcr]{%
1	0.2784568\\
2	0.0929659\\
3	0.0464011\\
4	0.0277447\\
5	0.1831458\\
6	0.1325341\\
7	0.0995494\\
8	0.0773472\\
9	0.061855\\
};
\end{axis}
\end{tikzpicture}%

%% file: fig2/b-p-2.tex
%
%
\definecolor{mycolor1}{rgb}{0.00000,0.44700,0.74100}%
\definecolor{mycolor2}{rgb}{0.92900,0.69400,0.12500}%
\definecolor{mycolor3}{rgb}{0.46600,0.67400,0.18800}%
\definecolor{mycolor4}{rgb}{0.63500,0.07800,0.18400}%
\definecolor{mycolor5}{rgb}{0.85000,0.32500,0.09800}%
\begin{tikzpicture}

\begin{axis}[%
yticklabel style={
        /pgf/number format/fixed,
        /pgf/number format/precision=5
},
scaled y ticks=false,
width=0.951\fwidth,
height=\fheight,
at={(0\fwidth,0\fheight)},
scale only axis,
xmin=10,
xmax=100,
xlabel style={font=\color{white!15!black}},
xlabel={$a$},
ymin=0,
ymax=0.101128101348877,
ylabel style={font=\color{white!15!black}},
ylabel={$\Pr(A_{(2)}=a)$},
grid = major,
grid style={dotted},
axis background/.style={fill=white},
xtick={10,20,30,40,50,60,70,80,90},
xmajorgrids,
ymajorgrids,
legend style={font=\scriptsize, legend cell align=left, align=left, draw={none}}
]
\addplot [color=blue, solid, line width=0.5pt]
  table[row sep=crcr]{%
10	0.0418511175169273\\
11	0.038172262604448\\
12	0.0350857591109624\\
13	0.0324593112158112\\
14	0.0301972816777539\\
15	0.0282288169886853\\
16	0.0265003080153068\\
17	0.0249704451445166\\
18	0.0236068826973525\\
19	0.0223839329984786\\
20	0.0212809374159352\\
21	0.0202810933157805\\
22	0.0193705946977862\\
23	0.0185379928285448\\
24	0.0177737138567274\\
25	0.01706969021308\\
26	0.0164190756684376\\
27	0.015816022707037\\
28	0.0152555068760561\\
29	0.0147331869407183\\
30	0.0142452926103807\\
31	0.0137885336961488\\
32	0.0133600260750198\\
33	0.0129572309417491\\
34	0.0125779046474565\\
35	0.0122200570332886\\
36	0.0118819166272027\\
37	0.0115619014202012\\
38	0.0112585942058953\\
39	0.0109707216732612\\
40	0.010697136602669\\
41	0.0104368026409634\\
42	0.0101887812299288\\
43	0.00995222034122683\\
44	0.00972634473296978\\
45	0.00951044749355439\\
46	0.00930388267863087\\
47	0.0091060588798616\\
48	0.00891643359071136\\
49	0.00873450825644864\\
50	0.00855982391333438\\
51	0.00839195733671239\\
52	0.00823051763022\\
53	0.00807514319811575\\
54	0.00792549905152869\\
55	0.00778127440649413\\
56	0.00764218053746707\\
57	0.00750794885522969\\
58	0.0073783291822774\\
59	0.00725308820242304\\
60	0.00713200806445909\\
61	0.00701488512221169\\
62	0.00690152879584618\\
63	0.0067917605408875\\
64	0.00668541291333602\\
65	0.00658232872060524\\
66	0.0064823602491221\\
67	0.00638536856077233\\
68	0.00629122285103981\\
69	0.00619979986267319\\
70	0.00611098334931185\\
71	0.00602466358425871\\
72	0.00594073690989119\\
73	0.00585910532399421\\
74	0.00577967609939868\\
75	0.00570236143398462\\
76	0.00562707812807998\\
77	0.00555374728696355\\
78	0.00548229404609752\\
79	0.00541264731714654\\
80	0.0053447395529437\\
81	0.00527850652983344\\
82	0.00521388714582242\\
83	0.00515082323327302\\
84	0.00508925938492838\\
85	0.00502914279216517\\
86	0.00497042309448511\\
87	0.00491305223932478\\
88	0.0048569843514258\\
89	0.00480217561092154\\
90	0.00474858413959208\\
91	0.00469616989450177\\
92	0.00464489456865951\\
93	0.00459472149791642\\
94	0.00454561557396388\\
95	0.00449754316263612\\
96	0.00445047202740352\\
97	0.00440437125759255\\
98	0.0043592112008633\\
99	0.00431496339994805\\
};
\addlegendentry{Benford}

\addplot [color=blue, line width=0.5pt, draw=none, mark=square, mark options={scale=0.3, solid, blue}, forget plot]
  table[row sep=crcr]{%
10	0.0418046843528802\\
11	0.0381639898816777\\
12	0.0350786521914275\\
13	0.0323957455006743\\
14	0.0302425745175407\\
15	0.0283436126213811\\
16	0.0265036992926682\\
17	0.0249131910996795\\
18	0.0236385427013719\\
19	0.0223352707602003\\
20	0.0213064245341794\\
21	0.0202714732867784\\
22	0.0193582221377083\\
23	0.018516229598845\\
24	0.0177902324662035\\
25	0.0170793477635355\\
26	0.0164326158265173\\
27	0.0157619642155673\\
28	0.015208308834016\\
29	0.0147533346177231\\
30	0.0142211969344786\\
31	0.0137583162150869\\
32	0.0133453765722307\\
33	0.0129508520758324\\
34	0.012582949475944\\
35	0.0122532783214194\\
36	0.01187887037088\\
37	0.0115104673594031\\
38	0.0112407455131846\\
39	0.0109769285237108\\
40	0.0106794838754876\\
41	0.0103833402973946\\
42	0.0101747687472947\\
43	0.00992556377620593\\
44	0.00976713346727684\\
45	0.00950952158150078\\
46	0.00933037423280648\\
47	0.00918375363736673\\
48	0.00893284726687701\\
49	0.00871466781427724\\
50	0.0085333186545934\\
51	0.0084333364192048\\
52	0.00825478956441673\\
53	0.00811937818856009\\
54	0.00796805372418814\\
55	0.00784385156791461\\
56	0.00762987557265851\\
57	0.00753079407812926\\
58	0.00741990287011067\\
59	0.00726317396058258\\
60	0.0071302646426686\\
61	0.007006763062619\\
62	0.00688676436368913\\
63	0.00676166146655624\\
64	0.00669800911249503\\
65	0.00659222210267945\\
66	0.00649724398317616\\
67	0.00637624446106923\\
68	0.00628657070440437\\
69	0.00619879851177594\\
70	0.00614295257849581\\
71	0.00601424671792549\\
72	0.00593928506196346\\
73	0.00585091237542879\\
74	0.00580047088730481\\
75	0.00570629342634317\\
76	0.00563483465150086\\
77	0.00552884747704987\\
78	0.00549421899512349\\
79	0.00543697190939548\\
80	0.00534129321366824\\
81	0.00527684020106538\\
82	0.00518576529195263\\
83	0.00510389795606887\\
84	0.00509719244078254\\
85	0.00503183868732032\\
86	0.00499560888831064\\
87	0.00490123126271358\\
88	0.00483157396958999\\
89	0.0047667206277163\\
90	0.00473179189883679\\
91	0.00470246777974884\\
92	0.004643118965349\\
93	0.00461869887982866\\
94	0.00454784059889259\\
95	0.00446176980566516\\
96	0.00445726610136838\\
97	0.00441142839985888\\
98	0.00437750049415644\\
99	0.00431755118585036\\
};
\addplot [color=green, solid, line width=0.5pt]
  table[row sep=crcr]{%
10	0.0680598135483133\\
11	0.0593643000591847\\
12	0.0523754959449399\\
13	0.046658467531892\\
14	0.0419111399622822\\
15	0.0379178395919914\\
16	0.0345208901759904\\
17	0.0316026233202864\\
18	0.0290736321276359\\
19	0.0268648944933855\\
20	0.0249223659535332\\
21	0.0232031897327712\\
22	0.0216729904066564\\
23	0.0203039087255862\\
24	0.019073152851287\\
25	0.0179619154993885\\
26	0.0169545543379628\\
27	0.0160379644536106\\
28	0.0152010927556167\\
29	0.0144345585189714\\
30	0.0137303541676015\\
31	0.0130816073358393\\
32	0.0124823901696734\\
33	0.0119275653663655\\
34	0.0114126610206636\\
35	0.010933768232395\\
36	0.0104874568288922\\
37	0.0100707056022822\\
38	0.00968084425157689\\
39	0.00931550482055721\\
40	0.0089725808833534\\
41	0.00865019308560266\\
42	0.0083466599259352\\
43	0.00806047287926537\\
44	0.00779027513406694\\
45	0.0075348433510708\\
46	0.00729307195857443\\
47	0.00706395958588013\\
48	0.00684659730586268\\
49	0.00664015841388998\\
50	0.00644388951602093\\
51	0.00625710273669644\\
52	0.00607916888673721\\
53	0.00590951145763208\\
54	0.00574760132891139\\
55	0.0055929520926653\\
56	0.00544511591363395\\
57	0.00530367985530833\\
58	0.00516826261254377\\
59	0.00503851159965496\\
60	0.00491410035009954\\
61	0.00479472618990451\\
62	0.00468010815211073\\
63	0.00456998510388096\\
64	0.00446411406163438\\
65	0.00436226867275802\\
66	0.00426423784517904\\
67	0.00416982450842385\\
68	0.00407884449182438\\
69	0.0039911255072701\\
70	0.00390650622542288\\
71	0.00382483543562584\\
72	0.00374597128087728\\
73	0.00366978056023348\\
74	0.00359613809187487\\
75	0.00352492613083233\\
76	0.00345603383602799\\
77	0.0033893567818728\\
78	0.0033247965101818\\
79	0.00326226011861031\\
80	0.00320165988222099\\
81	0.0031429129051429\\
82	0.00308594079959301\\
83	0.00303066938980956\\
84	0.0029770284386974\\
85	0.00292495139518956\\
86	0.00287437516054501\\
87	0.00282523987195371\\
88	0.00277748870199435\\
89	0.00273106767261344\\
90	0.00268592548243285\\
91	0.00264201334628755\\
92	0.00259928484601324\\
93	0.00255769579157628\\
94	0.00251720409173328\\
95	0.00247776963347046\\
96	0.00243935416954155\\
97	0.00240192121348737\\
98	0.0023654359415618\\
99	0.00232986510105271\\
};
\addlegendentry{Pareto, $s=0.5$, $\rho=0$}

\addplot [color=red, line width=0.5pt, draw=none, mark=square, mark options={solid, scale=0.3, red}, forget plot]
  table[row sep=crcr]{%
10	0.1011281\\
11	0.0841849\\
12	0.071365\\
13	0.061048\\
14	0.0529493\\
15	0.0462197\\
16	0.0407671\\
17	0.0362586\\
18	0.0324253\\
19	0.0292143\\
20	0.0264437\\
21	0.0240179\\
22	0.0219697\\
23	0.0200918\\
24	0.0185027\\
25	0.0171814\\
26	0.0157749\\
27	0.0147286\\
28	0.0136332\\
29	0.0127448\\
30	0.0119049\\
31	0.0111643\\
32	0.0105483\\
33	0.0099372\\
34	0.0093507\\
35	0.0088035\\
36	0.0083065\\
37	0.0079379\\
38	0.0075307\\
39	0.007131\\
40	0.0067983\\
41	0.0064694\\
42	0.0061581\\
43	0.0058862\\
44	0.005657\\
45	0.005373\\
46	0.0051105\\
47	0.0049509\\
48	0.0047597\\
49	0.0045583\\
50	0.004374\\
51	0.004159\\
52	0.0040551\\
53	0.0038922\\
54	0.0037379\\
55	0.0036182\\
56	0.0034532\\
57	0.0033484\\
58	0.0032586\\
59	0.0031284\\
60	0.0030268\\
61	0.0029389\\
62	0.0028674\\
63	0.0027274\\
64	0.0026445\\
65	0.0026028\\
66	0.0025276\\
67	0.002422\\
68	0.0023843\\
69	0.0022933\\
70	0.002255\\
71	0.0021687\\
72	0.002096\\
73	0.0020638\\
74	0.0020218\\
75	0.0019454\\
76	0.0018831\\
77	0.0018431\\
78	0.0018135\\
79	0.0017616\\
80	0.0017086\\
81	0.0016858\\
82	0.0016288\\
83	0.0015956\\
84	0.0015671\\
85	0.0015209\\
86	0.0014617\\
87	0.001448\\
88	0.0014144\\
89	0.0013826\\
90	0.0013256\\
91	0.0013388\\
92	0.0012843\\
93	0.001261\\
94	0.0012556\\
95	0.0012033\\
96	0.0012007\\
97	0.0011751\\
98	0.0011305\\
99	0.0011142\\
};
\addplot [color=green, dashed, line width=0.5pt]
  table[row sep=crcr]{%
10	0.0382728457429666\\
11	0.0333829991642795\\
12	0.0294529057971792\\
13	0.026237984463235\\
14	0.023568365986628\\
15	0.0213227681645155\\
16	0.019412523128317\\
17	0.0177714610774212\\
18	0.0163493048187515\\
19	0.0151072403705093\\
20	0.0140148762971887\\
21	0.0130481124629697\\
22	0.0121876181461169\\
23	0.01141772693929\\
24	0.0107256220499509\\
25	0.0101007273648831\\
26	0.00953424655442178\\
27	0.00901881018419906\\
28	0.00854820264453491\\
29	0.00811714876609602\\
30	0.00772114555794192\\
31	0.0132679174460424\\
32	0.0221971774268024\\
33	0.021210543902901\\
34	0.0202949001067999\\
35	0.0194432949218058\\
36	0.0186496285424908\\
37	0.0179085284171022\\
38	0.0172152460043728\\
39	0.0165655704165152\\
40	0.0159557558397762\\
41	0.0153824602569883\\
42	0.0148426934888876\\
43	0.0143337729563747\\
44	0.0138532858694496\\
45	0.0133990567890779\\
46	0.012969119699865\\
47	0.012561693884918\\
48	0.0121751630178435\\
49	0.0118080569868173\\
50	0.0114590360469058\\
51	0.0111268769631652\\
52	0.0108104608614339\\
53	0.0105087625484963\\
54	0.0102208411003154\\
55	0.00994583154772088\\
56	0.00968293751449073\\
57	0.00943142468413311\\
58	0.00919061498956053\\
59	0.00895988143490856\\
60	0.00873864347144709\\
61	0.00852636286028141\\
62	0.00832253996365379\\
63	0.00812671041441621\\
64	0.00793844211986955\\
65	0.00775733256182348\\
66	0.00758300635959064\\
67	0.00741511306680576\\
68	0.00725332517656213\\
69	0.00709733631245962\\
70	0.00694685958585806\\
71	0.00680162610196072\\
72	0.00666138359938134\\
73	0.0065258952096241\\
74	0.0063949383244376\\
75	0.00626830356036745\\
76	0.00614579381100638\\
77	0.00602722337848022\\
78	0.0059124171766255\\
79	0.00580120999911571\\
80	0.00569344584650122\\
81	0.00558897730676139\\
82	0.00548766498451519\\
83	0.00538937697453373\\
84	0.00529398837563577\\
85	0.00520138084143018\\
86	0.00511144216472442\\
87	0.00502406589271642\\
88	0.00493915097037206\\
89	0.00485660140963157\\
90	0.00477632598230915\\
91	0.00469823793475133\\
92	0.00462225472249145\\
93	0.00454829776330318\\
94	0.00447629220719498\\
95	0.00440616672201999\\
96	0.00433785329348851\\
97	0.00427128703848029\\
98	0.00420640603064565\\
99	0.00414315113737007\\
};
\addlegendentry{Pareto, $s=0.5$, $\rho=0.5$}

\addplot [color=green, line width=0.5pt, draw=none, mark=square, mark options={solid, scale=0.3, green}, forget plot]
  table[row sep=crcr]{%
10	0.0680825\\
11	0.059284\\
12	0.0523774\\
13	0.0466136\\
14	0.0419188\\
15	0.0379099\\
16	0.0344532\\
17	0.0315871\\
18	0.0291245\\
19	0.0269088\\
20	0.0249464\\
21	0.0232228\\
22	0.0217313\\
23	0.0203351\\
24	0.0190726\\
25	0.0178964\\
26	0.016919\\
27	0.0159816\\
28	0.015239\\
29	0.0143964\\
30	0.013777\\
31	0.0130824\\
32	0.0125115\\
33	0.0119252\\
34	0.0113951\\
35	0.010953\\
36	0.0104418\\
37	0.0100436\\
38	0.0097409\\
39	0.0093091\\
40	0.0089329\\
41	0.0086745\\
42	0.0083181\\
43	0.0080383\\
44	0.0077772\\
45	0.0075246\\
46	0.0073228\\
47	0.0070805\\
48	0.0068686\\
49	0.0066371\\
50	0.0064482\\
51	0.0062879\\
52	0.0060842\\
53	0.0059126\\
54	0.005739\\
55	0.0056248\\
56	0.0054503\\
57	0.0053249\\
58	0.0052029\\
59	0.0050394\\
60	0.0048978\\
61	0.0048124\\
62	0.0046787\\
63	0.0045753\\
64	0.0044655\\
65	0.0043777\\
66	0.0042489\\
67	0.0041633\\
68	0.0040433\\
69	0.0040105\\
70	0.0039124\\
71	0.0038424\\
72	0.0037012\\
73	0.0036477\\
74	0.0036095\\
75	0.0035415\\
76	0.0034345\\
77	0.0034202\\
78	0.0033268\\
79	0.0032687\\
80	0.0031865\\
81	0.0031483\\
82	0.0030829\\
83	0.0030255\\
84	0.0029827\\
85	0.0029377\\
86	0.0028465\\
87	0.002805\\
88	0.0027913\\
89	0.0027528\\
90	0.0026763\\
91	0.0026469\\
92	0.002609\\
93	0.0025526\\
94	0.0025084\\
95	0.0024827\\
96	0.0024206\\
97	0.0023906\\
98	0.0023703\\
99	0.0023363\\
};
\addplot [color=red, solid, line width=0.5pt]
  table[row sep=crcr]{%
10	0.101010101010101\\
11	0.0841750841750843\\
12	0.0712250712250711\\
13	0.0610500610500612\\
14	0.0529100529100529\\
15	0.0462962962962963\\
16	0.0408496732026143\\
17	0.0363108206245461\\
18	0.0324886289798571\\
19	0.0292397660818713\\
20	0.0264550264550265\\
21	0.024050024050024\\
22	0.0219587176108916\\
23	0.0201288244766505\\
24	0.0185185185185185\\
25	0.0170940170940171\\
26	0.0158277936055714\\
27	0.0146972369194591\\
28	0.0136836343732895\\
29	0.0127713920817371\\
30	0.0119474313022701\\
31	0.0112007168458781\\
32	0.0105218855218856\\
33	0.00990295107942152\\
34	0.00933706816059765\\
35	0.00881834215167543\\
36	0.00834167500834171\\
37	0.00790263948158698\\
38	0.00749737591842847\\
39	0.00712250712250706\\
40	0.00677506775067749\\
41	0.00645244547683566\\
42	0.0061523317337272\\
43	0.00587268029128485\\
44	0.00561167227833903\\
45	0.00536768652710684\\
46	0.00513927433446393\\
47	0.0049251379038614\\
48	0.00472411186696897\\
49	0.00453514739229033\\
50	0.00435729847494548\\
51	0.00418971007206304\\
52	0.00403160780519274\\
53	0.00388228899759302\\
54	0.00374111485222584\\
55	0.00360750360750361\\
56	0.0034809245335562\\
57	0.00336089265308859\\
58	0.00324696408857716\\
59	0.00313873195229131\\
60	0.00303582270795388\\
61	0.00293789294318114\\
62	0.0028446265005404\\
63	0.00275573192239853\\
64	0.00267094017094016\\
65	0.00259000259000264\\
66	0.00251268907985325\\
67	0.00243878645985751\\
68	0.00236809699725304\\
69	0.00230043708304584\\
70	0.00223563603845289\\
71	0.00217353503738493\\
72	0.00211398613225089\\
73	0.00205685137191991\\
74	0.00200200200200196\\
75	0.00194931773879148\\
76	0.0018986861092124\\
77	0.00185000185000195\\
78	0.00180316636012834\\
79	0.00175808720112514\\
80	0.00171467764060362\\
81	0.00167285623473523\\
82	0.00163254644594646\\
83	0.00159367629247142\\
84	0.00155617802676622\\
85	0.00151998784009733\\
86	0.00148504559089968\\
87	0.0014512945547428\\
88	0.00141868119396205\\
89	0.00138715494520736\\
90	0.00135666802333467\\
91	0.00132717524021875\\
92	0.0012986338372033\\
93	0.00127100333002872\\
94	0.00124424536518608\\
95	0.00121832358674467\\
96	0.00119320351279117\\
97	0.00116885242069331\\
98	0.00114523924047738\\
99	0.00112233445566778\\
};
\addlegendentry{Pareto, $s=1.0$, $\rho=0$}

\addplot [color=green, line width=0.5pt, draw=none, mark=square, mark options={solid, scale=0.3, green}, forget plot]
  table[row sep=crcr]{%
10	0.0383526\\
11	0.0334187\\
12	0.0295374\\
13	0.0262659\\
14	0.0235739\\
15	0.0213165\\
16	0.0193853\\
17	0.0176994\\
18	0.016362\\
19	0.0151184\\
20	0.014034\\
21	0.0130282\\
22	0.012212\\
23	0.0114366\\
24	0.0107574\\
25	0.0101228\\
26	0.0095393\\
27	0.0090743\\
28	0.0085744\\
29	0.0081615\\
30	0.0076924\\
31	0.0132911\\
32	0.0221581\\
33	0.0211487\\
34	0.0203024\\
35	0.0193735\\
36	0.0186669\\
37	0.0179563\\
38	0.0172617\\
39	0.016541\\
40	0.0159739\\
41	0.0153742\\
42	0.0148091\\
43	0.0142499\\
44	0.013837\\
45	0.0134268\\
46	0.0129797\\
47	0.0125592\\
48	0.0121865\\
49	0.0118521\\
50	0.0114554\\
51	0.0111128\\
52	0.0107996\\
53	0.0105598\\
54	0.0102246\\
55	0.0099448\\
56	0.0096375\\
57	0.0094686\\
58	0.0091441\\
59	0.0089431\\
60	0.0086905\\
61	0.0085262\\
62	0.0083166\\
63	0.0080891\\
64	0.0078768\\
65	0.00778\\
66	0.0076082\\
67	0.0074021\\
68	0.0072174\\
69	0.0070984\\
70	0.0069666\\
71	0.0067736\\
72	0.0066648\\
73	0.0065688\\
74	0.0064213\\
75	0.0062707\\
76	0.0061158\\
77	0.0060673\\
78	0.0059167\\
79	0.0057986\\
80	0.0056757\\
81	0.0055978\\
82	0.0055111\\
83	0.0054095\\
84	0.0052821\\
85	0.0051821\\
86	0.005117\\
87	0.0050336\\
88	0.0049679\\
89	0.0048567\\
90	0.0047841\\
91	0.0046827\\
92	0.0045891\\
93	0.00452\\
94	0.0044467\\
95	0.0043719\\
96	0.0043478\\
97	0.00424\\
98	0.0041824\\
99	0.0041309\\
};
\addplot [color=red, dashed, line width=0.5pt]
  table[row sep=crcr]{%
10	0.0506249730936639\\
11	0.0421874775780533\\
12	0.0356970964121988\\
13	0.0305975112104562\\
14	0.026517843049062\\
15	0.0232031126679293\\
16	0.0204733347069964\\
17	0.0181985197395524\\
18	0.0162828860827574\\
19	0.0146545974744816\\
20	0.013258921524531\\
21	0.0120535650223009\\
22	0.0110054289334052\\
23	0.0100883098556214\\
24	0.00928124506717169\\
25	0.00856730313892773\\
26	0.00793268809159978\\
27	0.00736606751362834\\
28	0.00685806285751604\\
29	0.00640085866701497\\
30	0.00598790004333659\\
31	0.00561365629062805\\
32	0.00527343469725665\\
33	0.00496323265624156\\
34	0.00467961936159918\\
35	0.00441964050817701\\
36	0.0041807410212485\\
37	0.00396070202013019\\
38	0.00375758909602093\\
39	0.00356970964121988\\
40	0.00339557746359941\\
41	0.00323388329866611\\
42	0.00308347012198395\\
43	0.0029433123891665\\
44	0.00281249850520355\\
45	0.00269021596149905\\
46	0.00257573868654164\\
47	0.00246841624126908\\
48	0.00236766455795197\\
49	0.00227295797563388\\
50	0.0194637564329178\\
51	0.0209982920071758\\
52	0.0202059036295467\\
53	0.0194575368284522\\
54	0.0187499900346903\\
55	0.0180803475334515\\
56	0.0174459493743829\\
57	0.0168443649131972\\
58	0.0162733694924108\\
59	0.0157309238426639\\
60	0.0152151558478225\\
61	0.0147243443688605\\
62	0.0142569048650872\\
63	0.0138113765880531\\
64	0.0133864111545745\\
65	0.0129807623317087\\
66	0.0125932768889712\\
67	0.0122228863922368\\
68	0.0118685998301428\\
69	0.011529496977853\\
70	0.0112047224150966\\
71	0.0108934801257885\\
72	0.0105950286154928\\
73	0.0103086764907497\\
74	0.0100337784509964\\
75	0.00976973164965456\\
76	0.00951597238602697\\
77	0.00927197309407773\\
78	0.00903723959802499\\
79	0.00881130860807433\\
80	0.00859374543256647\\
81	0.00838414188543063\\
82	0.00818211437011895\\
83	0.00798730212321153\\
84	0.00779936560266525\\
85	0.00761798500725452\\
86	0.00744285891513369\\
87	0.00727370303069885\\
88	0.00711024903000901\\
89	0.00695224349600876\\
90	0.00679944693565693\\
91	0.0066516328718384\\
92	0.00650858700362661\\
93	0.00637010642908159\\
94	0.00623599892531135\\
95	0.00610608228103399\\
96	0.00598018367730141\\
97	0.00585813911245853\\
98	0.00573979286776222\\
99	0.00562499701040708\\
};
\addlegendentry{Pareto, $s=1.0$, $\rho=0.7$}

\addplot [color=red, line width=0.5pt, draw=none, mark=square, mark options={solid, scale=0.3, red}, forget plot]
  table[row sep=crcr]{%
10	0.0505919\\
11	0.0421763\\
12	0.035793\\
13	0.0305816\\
14	0.0265118\\
15	0.0232382\\
16	0.0204956\\
17	0.0181749\\
18	0.0162427\\
19	0.014694\\
20	0.0132425\\
21	0.0120256\\
22	0.0109376\\
23	0.0100441\\
24	0.009292\\
25	0.0085619\\
26	0.0079254\\
27	0.0073495\\
28	0.0068527\\
29	0.0063996\\
30	0.0059888\\
31	0.005565\\
32	0.005242\\
33	0.0050142\\
34	0.0046711\\
35	0.0043928\\
36	0.0041586\\
37	0.003953\\
38	0.0037529\\
39	0.0035467\\
40	0.003389\\
41	0.0032511\\
42	0.0031059\\
43	0.0029298\\
44	0.002805\\
45	0.0026812\\
46	0.0025876\\
47	0.0024552\\
48	0.0023412\\
49	0.0022996\\
50	0.0194281\\
51	0.0210452\\
52	0.0201907\\
53	0.0195062\\
54	0.0188631\\
55	0.01811\\
56	0.0174773\\
57	0.0168033\\
58	0.0162375\\
59	0.01581\\
60	0.0152852\\
61	0.0146673\\
62	0.0142581\\
63	0.0137452\\
64	0.0133637\\
65	0.0129957\\
66	0.0125739\\
67	0.0122196\\
68	0.0118744\\
69	0.0114966\\
70	0.011269\\
71	0.0108899\\
72	0.0105647\\
73	0.0102441\\
74	0.0100571\\
75	0.0097815\\
76	0.0095241\\
77	0.0092782\\
78	0.0090234\\
79	0.0087715\\
80	0.0086016\\
81	0.0083949\\
82	0.0081946\\
83	0.0079751\\
84	0.0077529\\
85	0.0076214\\
86	0.0074991\\
87	0.0072945\\
88	0.0071734\\
89	0.0069153\\
90	0.0067981\\
91	0.0066163\\
92	0.0065576\\
93	0.0064177\\
94	0.0062157\\
95	0.0061353\\
96	0.0060182\\
97	0.0058031\\
98	0.0057426\\
99	0.0056564\\
};
\end{axis}

\end{tikzpicture}%

%% file: fig2/msd_real_datasets.tex
%
%
\definecolor{mycolor1}{rgb}{1.00000,0.00000,1.00000}%
\begin{tikzpicture}

\begin{axis}[%
width=0.951\fwidth,
height=\fheight,
at={(0\fwidth,0\fheight)},
scale only axis,
xmin=1,
xmax=9,
xlabel style={font=\color{white!15!black}},
xlabel={$a$},
ymin=0,
ymax=1,
ylabel style={font=\color{white!15!black}},
ylabel={$\Pr(A_{(1)}=a)$},
axis background/.style={fill=white},
xmajorgrids,
ymajorgrids,
legend style={legend cell align=left, align=left, draw=white!15!black},
grid style={dotted},legend style={font=\tiny,draw={none}, fill=none},mark size=1pt,xtick={1,2,3,4,5,6,7,8,9}
]
\addplot [color=blue, dotted, mark=square, mark options={solid, blue}]
  table[row sep=crcr]{%
1	0.30489970486719\\
2	0.178267720474213\\
3	0.12416212295533\\
4	0.0935295883147416\\
5	0.0761217547896553\\
6	0.0664298934520534\\
7	0.0585763593617128\\
8	0.0514544044820169\\
9	0.0465584513030864\\
};
\addlegendentry{US total income per ZIP code, 2016\\$p=159,928$ (NBER)}

\addplot [color=blue]
  table[row sep=crcr]{%
1	0.301029995663981\\
2	0.176091259055681\\
3	0.1249387366083\\
4	0.0969100130080564\\
5	0.0791812460476248\\
6	0.0669467896306132\\
7	0.0579919469776867\\
8	0.0511525224473813\\
9	0.0457574905606751\\
};
\addlegendentry{Benford}

\addplot [color=mycolor1, dotted, mark=square, mark options={solid, mycolor1}]
  table[row sep=crcr]{%
1	0.28296943231441\\
2	0.085589519650655\\
3	0.0436681222707424\\
4	0.022707423580786\\
5	0.197379912663755\\
6	0.136244541484716\\
7	0.0943231441048035\\
8	0.0733624454148472\\
9	0.0637554585152838\\
};
\addlegendentry{World Cities Population, 2020 ($\ge 500$k)\\ $p=1,145$
(World Population Review)}

\addplot [color=mycolor1]
  table[row sep=crcr]{%
1	0.266687984534464\\
2	0.0816447479451935\\
3	0.0387476784454579\\
4	0.0223798625606347\\
5	0.203067273679362\\
6	0.141627807916524\\
7	0.103984198747096\\
8	0.079388177007898\\
9	0.0624722691633703\\
};
\addlegendentry{Pareto, $\hat{s}=1.15$, $\hat{\rho}=0.70$}

\addplot [color=red, dotted, mark=square, mark options={solid, red}]
  table[row sep=crcr]{%
1	0.0647350993377483\\
2	0.371026490066225\\
3	0.211092715231788\\
4	0.135596026490066\\
5	0.0806291390728477\\
6	0.0498344370860927\\
7	0.0390728476821192\\
8	0.027317880794702\\
9	0.0206953642384106\\
};
\addlegendentry{Diameter of Mercurian craters ($\ge 20$ km)\\ $p=6,040$~\cite{fassett11:_mercury}}

\addplot [color=red]
  table[row sep=crcr]{%
1	0.0663830158460497\\
2	0.454611452249524\\
3	0.196026027283916\\
4	0.105549645234553\\
5	0.064599343052013\\
6	0.0429892147161257\\
7	0.0303513093673662\\
8	0.0223938914998626\\
9	0.0170961007505896\\
};
\addlegendentry{Pareto, $\hat{s}=1.42$, $\hat{\rho}=0.30$}

\addplot [color=green, dotted, mark=square, mark options={solid,
  green}] table[row sep=crcr]{%
  1	0.277\\
  2	0.1158\\
  3	0.0616\\
  4	0.123\\
  5	0.1232\\
  6	0.0986\\
  7	0.0794\\
  8	0.0676\\
  9	0.0538\\
}; \addlegendentry{Largest US Cities Population, July 2008\\ $p=5,000$
(US Census Bureau)}

\addplot [color=green]
  table[row sep=crcr]{%
1	0.259131068850512\\
2	0.0979583308046936\\
3	0.0529253627095318\\
4	0.144947677381756\\
5	0.141267933347827\\
6	0.104722206324487\\
7	0.0810718281289231\\
8	0.0648255215044148\\
9	0.0531500709478552\\
};
\addlegendentry{Pareto, $\hat{s}=0.78$, $\hat{\rho}=0.63$}

\addplot [color=black, dotted, mark=square, mark options={solid, black}]
  table[row sep=crcr]{%
1	0.672020888404961\\
2	0.180597410849509\\
3	0.0682289272946554\\
4	0.0329527543268178\\
5	0.0181794206644684\\
6	0.0111775483576445\\
7	0.00740517398264646\\
8	0.00540100370451482\\
9	0.00403687241478228\\
};
\addlegendentry{Lunar Craters (diameter $\ge 1$ km)\\ $p=1,296,796$~\cite{robbins18}}

\addplot [color=black]
  table[row sep=crcr]{%
1	0.685509575812625\\
2	0.161987673515822\\
3	0.0656685134961156\\
4	0.0338157363890522\\
5	0.0199800053434115\\
6	0.0129144097619874\\
7	0.0088939557470511\\
8	0.00642134484889156\\
9	0.00480878508504345\\
};
\addlegendentry{Pareto, $\hat{s}=1.59$, $\hat{\rho}=0.00$}

\end{axis}
\end{tikzpicture}%

%% file: fig2/paretoj_general.tex
%
%
\begin{tikzpicture}

\begin{axis}[%
ticklabel style={
/pgf/number format/fixed,
/pgf/number format/precision=5
}, scaled ticks=false,
width=0.951\fwidth,
height=\fheight,
at={(0\fwidth,0\fheight)},
scale only axis,
xmin=0,
xmax=9,
xlabel style={font=\color{white!15!black}},
xlabel={$a$},
ymin=0.07,
ymax=0.155,
ylabel style={font=\color{white!15!black}},
ylabel={$\Pr(A_{[j]}=a)$},
axis background/.style={fill=white},
legend style={legend cell align=left, align=left, draw=white!15!black},
legend style={font=\tiny,draw={none}, fill=white},xtick={0,1,2,3,4,5,6,7,8,9},
grid=major, grid style={dotted}
]
\addplot [color=blue]
  table[row sep=crcr]{%
0	0.119679268596881\\
1	0.113890103407556\\
2	0.108821499005508\\
3	0.104329560230961\\
4	0.100308202267578\\
5	0.0966772358023224\\
6	0.0933747357830365\\
7	0.0903519892696031\\
8	0.0875700535788619\\
9	0.0849973520576922\\
};
\addlegendentry{Benford, $j=2$}

\addplot [color=blue, line width=0.5pt, draw=none, mark=square, mark options={scale=0.6, solid, blue}, forget plot]
  table[row sep=crcr]{%
0	0.11969732\\
1	0.11389972\\
2	0.10879724\\
3	0.10437456\\
4	0.10034846\\
5	0.09659186\\
6	0.09343752\\
7	0.09026926\\
8	0.08762484\\
9	0.08495922\\
};
\addplot [color=blue, dashed]
  table[row sep=crcr]{%
0	0.101784364644181\\
1	0.101375977447795\\
2	0.10097219813707\\
3	0.100572932110944\\
4	0.100178087627925\\
5	0.0997875756921836\\
6	0.0994013099449835\\
7	0.0990192065618586\\
8	0.0986411841548051\\
9	0.0982671636782246\\
};
\addlegendentry{Benford, $j=3$}

\addplot [color=blue, line width=0.5pt, draw=none, mark=square, mark options={scale=0.6, solid, blue}, forget plot]
  table[row sep=crcr]{%
0	0.1018366\\
1	0.10134766\\
2	0.10095584\\
3	0.1005161\\
4	0.10024772\\
5	0.0997606\\
6	0.09944072\\
7	0.09896686\\
8	0.09868156\\
9	0.09824634\\
};
\addplot [color=blue, dashdotted]
  table[row sep=crcr]{%
0	0.100176146939765\\
1	0.100136888118105\\
2	0.100097672594297\\
3	0.10005850028342\\
4	0.100019371096996\\
5	0.0999802849481267\\
6	0.0999412417495204\\
7	0.0999022414154643\\
8	0.0998632838590611\\
9	0.0998243689957831\\
};
\addlegendentry{Benford, $j=4$}

\addplot [color=blue, line width=0.5pt, draw=none, mark=square, mark options={scale=0.6, solid, blue}, forget plot]
  table[row sep=crcr]{%
0	0.10023196\\
1	0.10011098\\
2	0.10003994\\
3	0.1000473\\
4	0.100085\\
5	0.09998632\\
6	0.09993254\\
7	0.09989098\\
8	0.0997904\\
9	0.09988458\\
};
\addplot [color=green]
  table[row sep=crcr]{%
0	0.105672019043794\\
1	0.108858130840169\\
2	0.119894585279158\\
3	0.112112457056985\\
4	0.105323066403793\\
5	0.0993277591547113\\
6	0.0939805040660286\\
7	0.0891712054049414\\
8	0.0848150127325704\\
9	0.0808452600178522\\
};
\addlegendentry{Pareto, $j=2$, $s=1.0$, $\rho=0.5$}

\addplot [color=green, line width=0.5pt, draw=none, mark=square, mark options={scale=0.6, solid, green}, forget plot]
  table[row sep=crcr]{%
0	0.10570848\\
1	0.10887754\\
2	0.11987438\\
3	0.11212196\\
4	0.10526638\\
5	0.09939294\\
6	0.09396538\\
7	0.08920318\\
8	0.08476378\\
9	0.08082598\\
};
\addplot [color=green, dashed]
  table[row sep=crcr]{%
0	0.101686319036034\\
1	0.101024883938813\\
2	0.10037315133021\\
3	0.0997308665139656\\
4	0.0990977844814012\\
5	0.0984736694392367\\
6	0.100294363713314\\
7	0.100388434085207\\
8	0.0997702232316549\\
9	0.0991603042301641\\
};
\addlegendentry{Pareto, $j=3$, $s=1.0$, $\rho=0.5$}

\addplot [color=green, line width=0.5pt, draw=none, mark=square, mark options={scale=0.6, solid, green}, forget plot]
  table[row sep=crcr]{%
0	0.10171672\\
1	0.10109494\\
2	0.10041776\\
3	0.0996679\\
4	0.09912282\\
5	0.09844198\\
6	0.10020104\\
7	0.1004118\\
8	0.09982696\\
9	0.09909808\\
};
\addplot [color=green, dashdotted]
  table[row sep=crcr]{%
0	0.100040841342661\\
1	0.0999773910804151\\
2	0.100142402795484\\
3	0.100166744461546\\
4	0.10010336451893\\
5	0.100040074300637\\
6	0.0999768735779298\\
7	0.0999137621217869\\
8	0.0998507397034132\\
9	0.0997878060970834\\
};
\addlegendentry{Pareto, $j=4$, $s=1.0$, $\rho=0.5$}

\addplot [color=green, line width=0.5pt, draw=none, mark=square, mark options={scale=0.6, solid, green}, forget plot]
  table[row sep=crcr]{%
0	0.10010338\\
1	0.09992864\\
2	0.10013144\\
3	0.1002705\\
4	0.1001418\\
5	0.10004428\\
6	0.09989362\\
7	0.0998434\\
8	0.0998743\\
9	0.09976864\\
};
\addplot [color=red]
  table[row sep=crcr]{%
0	0.138753025462623\\
1	0.127817969076182\\
2	0.116313685030193\\
3	0.106929781959253\\
4	0.0990983088946269\\
5	0.0924367000017385\\
6	0.0866790624502398\\
7	0.0816354483744322\\
8	0.0771668245345554\\
9	0.0731691942161552\\
};
\addlegendentry{Pareto, $j=2$, $s=1.5$, $\rho=0.7$}

\addplot [color=red, line width=0.5pt, draw=none, mark=square, mark options={scale=0.6, solid, red}, forget plot]
  table[row sep=crcr]{%
0	0.13886278\\
1	0.12778032\\
2	0.11630466\\
3	0.10684262\\
4	0.09906408\\
5	0.0924618\\
6	0.08664702\\
7	0.08172206\\
8	0.07712788\\
9	0.07318678\\
};
\addplot [color=red, dashed]
  table[row sep=crcr]{%
0	0.101022615028278\\
1	0.102607235590384\\
2	0.102326800904258\\
3	0.10150204076281\\
4	0.100693055703438\\
5	0.0998993374432331\\
6	0.0991203997255923\\
7	0.0983557771510419\\
8	0.0976050240811759\\
9	0.0968677136098037\\
};
\addlegendentry{Pareto, $j=3$, $s=1.5$, $\rho=0.7$}

\addplot [color=red, line width=0.5pt, draw=none, mark=square, mark options={scale=0.6, solid, red}, forget plot]
  table[row sep=crcr]{%
0	0.10103406\\
1	0.10260196\\
2	0.10234544\\
3	0.10147546\\
4	0.10069956\\
5	0.09989764\\
6	0.09907006\\
7	0.09830534\\
8	0.0976394\\
9	0.09693108\\
};
\addplot [color=red, dashdotted]
  table[row sep=crcr]{%
0	0.100114867103457\\
1	0.100073196951996\\
2	0.100254460706992\\
3	0.1001747388936\\
4	0.10009516878489\\
5	0.100015749887461\\
6	0.0999364817106425\\
7	0.0998573637659206\\
8	0.0997783955665454\\
9	0.0996995766283817\\
};
\addlegendentry{Pareto, $j=4$, $s=1.5$, $\rho=0.7$}

\addplot [color=red, line width=0.5pt, draw=none, mark=square, mark options={scale=0.6, solid, red}, forget plot]
  table[row sep=crcr]{%
0	0.10008958\\
1	0.1000536\\
2	0.10022886\\
3	0.10015832\\
4	0.10014868\\
5	0.09998746\\
6	0.09989104\\
7	0.0998645\\
8	0.09986642\\
9	0.09971154\\
};
\end{axis}
\end{tikzpicture}%

%% file: fig2/cf_general_benford_pseudo_n1e+08_mincha_110522_092109.tex
%
%
\definecolor{mycolor1}{rgb}{1.00000,0.00000,1.00000}%
\begin{tikzpicture}

\begin{axis}[%
width=0.964\fwidth,
height=\fheight,
at={(0\fwidth,0\fheight)},
scale only axis,
xmode=log,
xmin=1,
xmax=100,
xminorticks=true,
xlabel style={font=\color{white!15!black}},
xlabel={$a_1$},
ymode=log,
ymin=5e-05,
ymax=1,
yminorticks=true,
ylabel style={font=\color{white!15!black}},
ylabel={$\Pr(\mathbf{A}_2=[a_1,a_2])$},
axis background/.style={fill=white},
xmajorgrids,
xminorgrids,
ymajorgrids,
yminorgrids,
grid style={dotted},legend style={font=\tiny,draw=none},mark size=1pt
]
\addplot [color=blue, draw=none, mark=square, mark options={solid, blue}, forget plot]
  table[row sep=crcr]{%
1	0.1666342\\
2	0.0666654\\
3	0.03573436\\
4	0.02223866\\
5	0.01514011\\
6	0.01098554\\
7	0.00833154999999999\\
8	0.00651838999999999\\
9.00000000000001	0.00525805\\
10	0.00432018999999999\\
11	0.00362315\\
12	0.00308237\\
13	0.00264009\\
14	0.00229519\\
15	0.00201437\\
16	0.0017854\\
17	0.00159081\\
18	0.00142709\\
19	0.00128164\\
20	0.00116468\\
21	0.0010565\\
22	0.000966149999999999\\
23	0.000887889999999999\\
24	0.000812639999999999\\
25	0.000752269999999999\\
26	0.00069589\\
27	0.00064831\\
28	0.00060624\\
29	0.000563739999999999\\
30	0.00052964\\
31	0.00049479\\
32	0.00046457\\
33	0.00044083\\
34	0.00041499\\
35	0.00038962\\
36	0.00037069\\
37	0.00034927\\
38	0.00033632\\
39	0.00031612\\
40	0.00029946\\
41	0.00028817\\
42	0.00027327\\
43.0000000000001	0.00026245\\
44	0.00025166\\
45	0.00023739\\
46	0.00022741\\
47.0000000000001	0.00021793\\
48	0.00021127\\
49	0.0002015\\
50	0.00019554\\
51	0.00018732\\
52	0.00018266\\
53	0.00017238\\
54.0000000000001	0.00016624\\
55	0.00016001\\
56	0.00015578\\
57	0.00014927\\
58.0000000000001	0.00014765\\
59.0000000000001	0.0001387\\
60	0.00013715\\
61	0.00013083\\
62.0000000000001	0.00012787\\
63.0000000000001	0.00012432\\
64.0000000000001	0.00011828\\
65	0.00011486\\
66.0000000000001	0.00011289\\
67.0000000000001	0.00010711\\
68	0.00010638\\
69.0000000000001	0.00010149\\
70	9.874e-05\\
71	9.659e-05\\
72.0000000000001	9.521e-05\\
73	9.28100000000001e-05\\
74	8.959e-05\\
75	8.667e-05\\
76	8.573e-05\\
77.0000000000001	8.17399999999999e-05\\
78.0000000000001	8.05500000000001e-05\\
79.0000000000001	7.87e-05\\
80	7.923e-05\\
81	7.449e-05\\
82	7.33999999999999e-05\\
83	7.068e-05\\
84.0000000000001	6.844e-05\\
85	6.777e-05\\
85.9999999999999	6.67499999999999e-05\\
86.9999999999999	6.326e-05\\
88.0000000000001	6.269e-05\\
89.0000000000001	6.13e-05\\
90.0000000000001	6.207e-05\\
91.0000000000001	5.992e-05\\
92.0000000000001	5.784e-05\\
93.0000000000001	5.716e-05\\
93.9999999999999	5.559e-05\\
95	5.423e-05\\
96.0000000000001	5.368e-05\\
97.0000000000001	5.344e-05\\
98.0000000000001	5.145e-05\\
98.9999999999999	5.034e-05\\
100	5.001e-05\\
};
\addplot [color=red, draw=none, mark=square, mark options={solid, red}, forget plot]
  table[row sep=crcr]{%
1	0.0834012099999999\\
2	0.02857529\\
3	0.01429987\\
4	0.00854552\\
5	0.00568414\\
6	0.0040561\\
7	0.00301702\\
8	0.00235493\\
9.00000000000001	0.00187676\\
10	0.0015395\\
11	0.00128752\\
12	0.00108325\\
13	0.00092726\\
14	0.000801089999999999\\
15	0.000702489999999999\\
16	0.00061835\\
17	0.00054629\\
18	0.000487579999999999\\
19	0.00044319\\
20	0.00040052\\
21	0.00036418\\
22	0.00033039\\
23	0.00030742\\
24	0.00028103\\
25	0.00025884\\
26	0.00023938\\
27	0.00022328\\
28	0.00020711\\
29	0.00019439\\
30	0.00018056\\
31	0.0001655\\
32	0.00016039\\
33	0.00014864\\
34	0.00013976\\
35	0.00013308\\
36	0.0001258\\
37	0.00011993\\
38	0.00011314\\
39	0.00010624\\
40	0.00010334\\
41	9.69499999999999e-05\\
42	9.40399999999999e-05\\
43.0000000000001	8.874e-05\\
44	8.445e-05\\
45	8.16600000000001e-05\\
46	7.92799999999999e-05\\
47.0000000000001	7.427e-05\\
48	7.134e-05\\
49	6.781e-05\\
50	6.42e-05\\
51	6.33e-05\\
52	5.934e-05\\
53	5.854e-05\\
54.0000000000001	5.644e-05\\
55	5.31e-05\\
56	5.265e-05\\
57	5.013e-05\\
};
\addplot [color=green, draw=none, mark=square, mark options={solid, green}, forget plot]
  table[row sep=crcr]{%
1	0.05001203\\
2	0.01587363\\
3	0.00769438\\
4	0.00451514\\
5	0.00297354\\
6	0.00210202\\
7	0.00156506\\
8	0.00120611\\
9.00000000000001	0.000964419999999999\\
10	0.00078509\\
11	0.00065347\\
12	0.000550789999999999\\
13	0.00047203\\
14	0.00040855\\
15	0.00035521\\
16	0.00031467\\
17	0.00027682\\
18	0.00024739\\
19	0.00022458\\
20	0.00020372\\
21	0.00018476\\
22	0.00016967\\
23	0.0001539\\
24	0.00014146\\
25	0.00013246\\
26	0.00011889\\
27	0.00011141\\
28	0.00010407\\
29	9.741e-05\\
30	9.172e-05\\
31	8.51200000000001e-05\\
32	7.777e-05\\
33	7.509e-05\\
34	7.042e-05\\
35	6.708e-05\\
36	6.416e-05\\
37	6.111e-05\\
38	5.628e-05\\
39	5.341e-05\\
40	5.064e-05\\
};
\addplot [color=black, draw=none, mark=square, mark options={solid, black}, forget plot]
  table[row sep=crcr]{%
1	0.03331951\\
2	0.01009865\\
3	0.00481614999999999\\
4	0.00280614\\
5	0.00182994\\
6	0.00128839\\
7	0.00095536\\
8	0.00074495\\
9.00000000000001	0.00058781\\
10	0.00047717\\
11	0.00040051\\
12	0.00033523\\
13	0.00028527\\
14	0.00024721\\
15	0.00021572\\
16	0.00018947\\
17	0.00016863\\
18	0.00014862\\
19	0.00013637\\
20	0.00012439\\
21	0.00011197\\
22	0.00010246\\
23	9.214e-05\\
24	8.38e-05\\
25	7.923e-05\\
26	7.194e-05\\
27	6.801e-05\\
28	6.249e-05\\
29	5.87999999999999e-05\\
30	5.515e-05\\
31	5.118e-05\\
};
\addplot [color=mycolor1, draw=none, mark=square, mark options={solid, mycolor1}, forget plot]
  table[row sep=crcr]{%
1	0.02382539\\
2	0.00697884999999999\\
3	0.00328429\\
4	0.00189611\\
5	0.0012404\\
6	0.00086909\\
7	0.000646589999999999\\
8	0.00049604\\
9.00000000000001	0.00039666\\
10	0.00031909\\
11	0.00026756\\
12	0.00022425\\
13	0.00019124\\
14	0.00016564\\
15	0.00014669\\
16	0.00012821\\
17	0.00011287\\
18	0.00010334\\
19	9.105e-05\\
20	8.227e-05\\
21	7.54400000000001e-05\\
22	6.828e-05\\
23	6.197e-05\\
24	5.77299999999999e-05\\
25	5.127e-05\\
};
\addplot [color=blue, draw=none, mark=square, mark options={solid, blue}, forget plot]
  table[row sep=crcr]{%
1	0.01785796\\
2	0.0051257\\
3	0.0023961\\
4	0.00138335\\
5	0.000898\\
6	0.000627879999999999\\
7	0.00046826\\
8	0.00035764\\
9.00000000000001	0.00028405\\
10	0.00022967\\
11	0.00018808\\
12	0.00015974\\
13	0.00013758\\
14	0.00011851\\
15	0.00010409\\
16	9.221e-05\\
17	8.056e-05\\
18	7.233e-05\\
19	6.535e-05\\
20	5.873e-05\\
21	5.276e-05\\
};
\addplot [color=red, draw=none, mark=square, mark options={solid, red}, forget plot]
  table[row sep=crcr]{%
1	0.01389045\\
2	0.00391684\\
3	0.00182242\\
4	0.00104238\\
5	0.00067739\\
6	0.000474769999999999\\
7	0.00035009\\
8	0.00026911\\
9.00000000000001	0.00021488\\
10	0.00017401\\
11	0.00014756\\
12	0.00011996\\
13	0.00010254\\
14	8.926e-05\\
15	7.762e-05\\
16	6.874e-05\\
17	6.175e-05\\
18	5.426e-05\\
};
\addplot [color=green, draw=none, mark=square, mark options={solid, green}, forget plot]
  table[row sep=crcr]{%
1	0.0111142\\
2	0.00310452\\
3	0.0014299\\
4	0.000815069999999999\\
5	0.00053065\\
6	0.00036906\\
7	0.00027561\\
8	0.00021239\\
9.00000000000001	0.00016691\\
10	0.00013372\\
11	0.00011126\\
12	9.57200000000001e-05\\
13	8.198e-05\\
14	7.189e-05\\
15	6.096e-05\\
16	5.405e-05\\
};
\addplot [color=black, draw=none, mark=square, mark options={solid, black}, forget plot]
  table[row sep=crcr]{%
1	0.00908752000000001\\
2	0.00250796\\
3	0.00115323\\
4	0.00066281\\
5	0.00043121\\
6	0.00029768\\
7	0.00022179\\
8	0.00017165\\
9.00000000000001	0.00013215\\
10	0.00010824\\
11	8.914e-05\\
12	7.69299999999999e-05\\
13	6.391e-05\\
14	5.72800000000001e-05\\
};
\addplot [color=mycolor1, draw=none, mark=square, mark options={solid, mycolor1}, forget plot]
  table[row sep=crcr]{%
1	0.00757837999999999\\
2	0.00207023\\
3	0.00095418\\
4	0.00054225\\
5	0.00034689\\
6	0.000247\\
7	0.00018117\\
8	0.00013843\\
9.00000000000001	0.00010939\\
10	8.917e-05\\
11	7.347e-05\\
12	6.243e-05\\
13	5.376e-05\\
};
\addplot [color=blue, draw=none, mark=square, mark options={solid, blue}, forget plot]
  table[row sep=crcr]{%
1	0.00642056999999999\\
2	0.00173966\\
3	0.0007965\\
4	0.00045339\\
5	0.00029265\\
6	0.00020479\\
7	0.00015123\\
8	0.00011673\\
9.00000000000001	8.931e-05\\
10	7.241e-05\\
11	6.199e-05\\
12	5.121e-05\\
};
\addplot [color=red, draw=none, mark=square, mark options={solid, red}, forget plot]
  table[row sep=crcr]{%
1	0.00548590999999999\\
2	0.00148306\\
3	0.00068111\\
4	0.00038442\\
5	0.00024726\\
6	0.00017132\\
7	0.00013023\\
8	9.965e-05\\
9.00000000000001	7.84899999999999e-05\\
10	6.385e-05\\
11	5.314e-05\\
};
\addplot [color=green, draw=none, mark=square, mark options={solid, green}, forget plot]
  table[row sep=crcr]{%
1	0.00476217\\
2	0.00127737\\
3	0.00057774\\
4	0.00032718\\
5	0.00021029\\
6	0.00014946\\
7	0.00010843\\
8	8.463e-05\\
9.00000000000001	6.599e-05\\
10	5.362e-05\\
};
\addplot [color=black, draw=none, mark=square, mark options={solid, black}, forget plot]
  table[row sep=crcr]{%
1	0.00415436999999999\\
2	0.00110993\\
3	0.00050752\\
4	0.00028884\\
5	0.00018662\\
6	0.00012948\\
7	9.515e-05\\
8	7.36500000000001e-05\\
9.00000000000001	5.719e-05\\
};
\addplot [color=mycolor1, draw=none, mark=square, mark options={solid, mycolor1}, forget plot]
  table[row sep=crcr]{%
1	0.00366664\\
2	0.000973749999999999\\
3	0.00044138\\
4	0.00025291\\
5	0.00016272\\
6	0.00011234\\
7	8.198e-05\\
8	6.421e-05\\
9.00000000000001	5.046e-05\\
};
\addplot [color=blue, draw=none, mark=square, mark options={solid, blue}, forget plot]
  table[row sep=crcr]{%
1	0.00326887\\
2	0.000868279999999999\\
3	0.00039148\\
4	0.00022328\\
5	0.00014435\\
6	9.89e-05\\
7	7.352e-05\\
8	5.694e-05\\
};
\addplot [color=red, draw=none, mark=square, mark options={solid, red}, forget plot]
  table[row sep=crcr]{%
1	0.00292811\\
2	0.000775309999999999\\
3	0.00035012\\
4	0.00019808\\
5	0.00012759\\
6	8.94199999999999e-05\\
7	6.514e-05\\
8	5.15e-05\\
};
\addplot [color=green, draw=none, mark=square, mark options={solid, green}, forget plot]
  table[row sep=crcr]{%
1	0.00263104\\
2	0.00068661\\
3	0.00031425\\
4	0.00017548\\
5	0.00011249\\
6	7.962e-05\\
7	5.77e-05\\
};
\addplot [color=black, draw=none, mark=square, mark options={solid, black}, forget plot]
  table[row sep=crcr]{%
1	0.00237771\\
2	0.000624289999999999\\
3	0.00028375\\
4	0.00016086\\
5	0.00010385\\
6	7.201e-05\\
7	5.286e-05\\
};
\addplot [color=mycolor1, draw=none, mark=square, mark options={solid, mycolor1}, forget plot]
  table[row sep=crcr]{%
1	0.002162\\
2	0.00057079\\
3	0.00025758\\
4	0.00014309\\
5	9.371e-05\\
6	6.405e-05\\
};
\addplot [color=blue, draw=none, mark=square, mark options={solid, blue}, forget plot]
  table[row sep=crcr]{%
1	0.00198123\\
2	0.00051789\\
3	0.00023088\\
4	0.00013173\\
5	8.461e-05\\
6	5.898e-05\\
};
\addplot [color=red, draw=none, mark=square, mark options={solid, red}, forget plot]
  table[row sep=crcr]{%
1	0.00180672\\
2	0.00047281\\
3	0.00021149\\
4	0.00012121\\
5	7.71399999999999e-05\\
6	5.393e-05\\
};
\addplot [color=green, draw=none, mark=square, mark options={solid, green}, forget plot]
  table[row sep=crcr]{%
1	0.00166668\\
2	0.00043275\\
3	0.00019377\\
4	0.0001095\\
5	7.238e-05\\
6	5.062e-05\\
};
\addplot [color=black, draw=none, mark=square, mark options={solid, black}, forget plot]
  table[row sep=crcr]{%
1	0.00154958\\
2	0.00039842\\
3	0.00018052\\
4	0.00010158\\
5	6.488e-05\\
};
\addplot [color=mycolor1, draw=none, mark=square, mark options={solid, mycolor1}, forget plot]
  table[row sep=crcr]{%
1	0.00142794\\
2	0.000374059999999999\\
3	0.00016407\\
4	9.26399999999999e-05\\
5	6.066e-05\\
};
\addplot [color=blue, draw=none, mark=square, mark options={solid, blue}, forget plot]
  table[row sep=crcr]{%
1	0.00132191\\
2	0.00034293\\
3	0.000153\\
4	8.783e-05\\
5	5.639e-05\\
};
\addplot [color=red, draw=none, mark=square, mark options={solid, red}, forget plot]
  table[row sep=crcr]{%
1	0.0012349\\
2	0.00031829\\
3	0.00014329\\
4	8.18e-05\\
5	5.082e-05\\
};
\addplot [color=green, draw=none, mark=square, mark options={solid, green}, forget plot]
  table[row sep=crcr]{%
1	0.00115138\\
2	0.00029535\\
3	0.00013152\\
4	7.664e-05\\
};
\addplot [color=black, draw=none, mark=square, mark options={solid, black}, forget plot]
  table[row sep=crcr]{%
1	0.00107638\\
2	0.00027565\\
3	0.00012486\\
4	7.013e-05\\
};
\addplot [color=mycolor1, draw=none, mark=square, mark options={solid, mycolor1}, forget plot]
  table[row sep=crcr]{%
1	0.00101149\\
2	0.00026154\\
3	0.00011741\\
4	6.539e-05\\
};
\addplot [color=blue, draw=none, mark=square, mark options={solid, blue}, forget plot]
  table[row sep=crcr]{%
1	0.00094442\\
2	0.0002445\\
3	0.00010996\\
4	6.145e-05\\
};
\addplot [color=red, draw=none, mark=square, mark options={solid, red}, forget plot]
  table[row sep=crcr]{%
1	0.00088801\\
2	0.0002322\\
3	0.00010261\\
4	5.886e-05\\
};
\addplot [color=green, draw=none, mark=square, mark options={solid, green}, forget plot]
  table[row sep=crcr]{%
1	0.000844439999999999\\
2	0.00021702\\
3	9.681e-05\\
4	5.55e-05\\
};
\addplot [color=black, draw=none, mark=square, mark options={solid, black}, forget plot]
  table[row sep=crcr]{%
1	0.000793429999999999\\
2	0.00020212\\
3	9.374e-05\\
4	5.165e-05\\
};
\addplot [color=mycolor1, draw=none, mark=square, mark options={solid, mycolor1}, forget plot]
  table[row sep=crcr]{%
1	0.00074952\\
2	0.00019318\\
3	8.72299999999999e-05\\
};
\addplot [color=blue, draw=none, mark=square, mark options={solid, blue}, forget plot]
  table[row sep=crcr]{%
1	0.000707079999999999\\
2	0.00018244\\
3	8.15899999999999e-05\\
};
\addplot [color=red, draw=none, mark=square, mark options={solid, red}, forget plot]
  table[row sep=crcr]{%
1	0.000674529999999999\\
2	0.00017363\\
3	7.72400000000001e-05\\
};
\addplot [color=green, draw=none, mark=square, mark options={solid, green}, forget plot]
  table[row sep=crcr]{%
1	0.00064396\\
2	0.00016388\\
3	7.40299999999999e-05\\
};
\addplot [color=black, draw=none, mark=square, mark options={solid, black}, forget plot]
  table[row sep=crcr]{%
1	0.00061256\\
2	0.00015568\\
3	6.97800000000001e-05\\
};
\addplot [color=mycolor1, draw=none, mark=square, mark options={solid, mycolor1}, forget plot]
  table[row sep=crcr]{%
1	0.00057807\\
2	0.00014799\\
3	6.619e-05\\
};
\addplot [color=blue, draw=none, mark=square, mark options={solid, blue}, forget plot]
  table[row sep=crcr]{%
1	0.000550699999999999\\
2	0.00014241\\
3	6.271e-05\\
};
\addplot [color=red, draw=none, mark=square, mark options={solid, red}, forget plot]
  table[row sep=crcr]{%
1	0.00052838\\
2	0.00013403\\
3	6.078e-05\\
};
\addplot [color=green, draw=none, mark=square, mark options={solid, green}, forget plot]
  table[row sep=crcr]{%
1	0.00050214\\
2	0.00013188\\
3	5.657e-05\\
};
\addplot [color=black, draw=none, mark=square, mark options={solid, black}, forget plot]
  table[row sep=crcr]{%
1	0.00047929\\
2	0.00012371\\
3	5.475e-05\\
};
\addplot [color=mycolor1, draw=none, mark=square, mark options={solid, mycolor1}, forget plot]
  table[row sep=crcr]{%
1	0.00046451\\
2	0.00011799\\
3	5.272e-05\\
};
\addplot [color=blue, draw=none, mark=square, mark options={solid, blue}, forget plot]
  table[row sep=crcr]{%
1	0.00043945\\
2	0.00011446\\
};
\addplot [color=red, draw=none, mark=square, mark options={solid, red}, forget plot]
  table[row sep=crcr]{%
1	0.00042193\\
2	0.00010943\\
};
\addplot [color=green, draw=none, mark=square, mark options={solid, green}, forget plot]
  table[row sep=crcr]{%
1	0.00040739\\
2	0.0001037\\
};
\addplot [color=black, draw=none, mark=square, mark options={solid, black}, forget plot]
  table[row sep=crcr]{%
1	0.00039297\\
2	0.00010184\\
};
\addplot [color=mycolor1, draw=none, mark=square, mark options={solid, mycolor1}, forget plot]
  table[row sep=crcr]{%
1	0.0003778\\
2	9.65599999999999e-05\\
};
\addplot [color=blue, dashed, forget plot]
  table[row sep=crcr]{%
1	0.166666666666667\\
2	0.0666666666666667\\
3	0.0357142857142857\\
4	0.0222222222222222\\
5	0.0151515151515152\\
6	0.010989010989011\\
7	0.00833333333333333\\
9.00000000000001	0.00526315789473683\\
11	0.0036231884057971\\
14	0.00229885057471264\\
18	0.00142247510668564\\
24	0.00081632653061224\\
34	0.000414078674948242\\
51	0.000186706497386106\\
84.0000000000001	6.96136442742777e-05\\
100	4.92586572090042e-05\\
};
\node[right, align=left, font=\color{blue}]
at (axis cs:1.02,0.157) {{\tiny $a_2=1$}};
\addplot [color=red, dashed, forget plot]
  table[row sep=crcr]{%
1	0.0833333333333333\\
2	0.0285714285714285\\
3	0.0142857142857143\\
4	0.00854700854700857\\
5	0.00568181818181818\\
6	0.00404858299595142\\
8	0.00235294117647058\\
11	0.00127877237851662\\
15	0.00070126227208976\\
22	0.000331674958540632\\
34	0.000140706345856195\\
58.0000000000001	4.8840048840048e-05\\
};
\node[right, align=left, font=\color{red}]
at (axis cs:1.02,0.078) {{\tiny $a_2=2$}};
\addplot [color=green, dashed, forget plot]
  table[row sep=crcr]{%
1	0.05\\
2	0.0158730158730159\\
3	0.00769230769230772\\
4	0.00452488687782803\\
5	0.00297619047619047\\
7	0.00156739811912224\\
10	0.000786782061369007\\
14	0.000407996736026109\\
21	0.000183823529411764\\
35	6.69075337883068e-05\\
41	4.88758553274689e-05\\
};
\node[right, align=left, font=\color{green}]
at (axis cs:1.02,0.047) {{\tiny $a_2=3$}};
\addplot [color=black, dashed, forget plot]
  table[row sep=crcr]{%
1	0.0333333333333333\\
2	0.0101010101010101\\
3	0.00480769230769229\\
4	0.00280112044817926\\
5	0.00183150183150183\\
7	0.000957854406130276\\
10	0.000478240076518421\\
15	0.000215703192407241\\
25	7.85792865000758e-05\\
32	4.81486831335129e-05\\
};
\node[right, align=left]
at (axis cs:1.02,0.031) {{\tiny $a_2=4$}};
\addplot [color=mycolor1, dashed, forget plot]
  table[row sep=crcr]{%
1	0.0238095238095237\\
2	0.00699300699300703\\
3	0.00328947368421057\\
4	0.00190476190476191\\
6	0.000871839581516986\\
9.00000000000001	0.000395256916996045\\
14	0.000165700082850037\\
24	5.69962952408076e-05\\
26	4.86215782564276e-05\\
};
\node[right, align=left, font=\color{mycolor1}]
at (axis cs:1.02,0.022) {{\tiny $a_2=5$}};
\addplot [color=blue, dashed, forget plot]
  table[row sep=crcr]{%
1	0.0178571428571429\\
2	0.0051282051282051\\
3	0.00239234449760761\\
4	0.00137931034482758\\
6	0.000628535512256444\\
9.00000000000001	0.000284090909090901\\
15	0.000103669914990667\\
22	4.85083676934339e-05\\
};
\addplot [color=red, dashed, forget plot]
  table[row sep=crcr]{%
1	0.0138888888888888\\
2	0.00392156862745096\\
3	0.00181818181818183\\
4	0.00104493207941486\\
6	0.000474608448030378\\
9.00000000000001	0.000214041095890405\\
15	7.79666302822418e-05\\
19	4.8775729197155e-05\\
};
\node[right, align=left]
at (axis cs:1.02,0.013) {\tiny $\phantom{a_2}~~~~\!\pmb{\pmb{\vdots}}$};
\addplot [color=green, dashed, forget plot]
  table[row sep=crcr]{%
1	0.0111111111111111\\
2	0.00309597523219812\\
3	0.00142857142857139\\
4	0.000819000819000825\\
6	0.000371057513914669\\
10	0.000135666802333484\\
17	4.73978576168363e-05\\
};
\addplot [color=black, dashed, forget plot]
  table[row sep=crcr]{%
1	0.00909090909090914\\
2	0.0025062656641604\\
3	0.00115207373271892\\
4	0.000659195781147009\\
6	0.000298062593144577\\
10	0.000108802089000104\\
15	4.86949746786175e-05\\
};
\addplot [color=mycolor1, dashed, forget plot]
  table[row sep=crcr]{%
1	0.00757575757575767\\
2	0.00207039337474124\\
3	0.000948766603415585\\
5	0.00035014005602238\\
8	0.000138715494520713\\
14	4.57561198810297e-05\\
};
\addplot [color=blue, dashed, forget plot]
  table[row sep=crcr]{%
1	0.00641025641025639\\
2	0.00173913043478258\\
3	0.000794912559618388\\
5	0.000292740046838435\\
8	0.000115834588208047\\
13	4.4232130219396e-05\\
};
\addplot [color=red, dashed, forget plot]
  table[row sep=crcr]{%
1	0.00549450549450547\\
2	0.00148148148148147\\
3	0.000675675675675669\\
5	0.000248385494287129\\
9.00000000000001	7.77484061576794e-05\\
12	4.39270810454611e-05\\
};
\addplot [color=green, dashed, forget plot]
  table[row sep=crcr]{%
1	0.00476190476190474\\
2	0.00127713920817368\\
3	0.000581395348837221\\
5	0.000213401621852321\\
9.00000000000001	6.67289470172166e-05\\
11	4.48028673835227e-05\\
};
\addplot [color=black, dashed, forget plot]
  table[row sep=crcr]{%
1	0.00416666666666665\\
2	0.0011123470522803\\
3	0.000505561172901958\\
5	0.000185322461082271\\
9.00000000000001	5.78971746178769e-05\\
10	4.69682025268925e-05\\
};
\addplot [color=mycolor1, dashed, forget plot]
  table[row sep=crcr]{%
1	0.00367647058823528\\
2	0.000977517106549419\\
3	0.000443655723158797\\
5	0.000162443144899271\\
9.00000000000001	5.07099391480748e-05\\
10	4.11336432067644e-05\\
};
\addplot [color=blue, dashed, forget plot]
  table[row sep=crcr]{%
1	0.00326797385620914\\
2	0.000865800865800902\\
3	0.000392464678178994\\
5	0.000143554407120328\\
9.00000000000001	4.47828034035036e-05\\
};
\addplot [color=red, dashed, forget plot]
  table[row sep=crcr]{%
1	0.0029239766081871\\
2	0.000772200772200748\\
3	0.000349650349650343\\
5	0.00012777919754664\\
9.00000000000001	3.98374631503384e-05\\
};
\addplot [color=green, dashed, forget plot]
  table[row sep=crcr]{%
1	0.00263157894736854\\
2	0.0006930006930006\\
3	0.000313479623824442\\
5	0.000114468864468864\\
8	4.50755014649584e-05\\
};
\addplot [color=black, dashed, forget plot]
  table[row sep=crcr]{%
1	0.00238095238095226\\
2	0.000625390869293418\\
3	0.000282645562464712\\
5	0.000103135313531344\\
8	4.0595948524319e-05\\
};
\addplot [color=mycolor1, dashed, forget plot]
  table[row sep=crcr]{%
1	0.00216450216450214\\
2	0.000567214974475316\\
3	0.000256147540983575\\
5	9.34055669717937e-05\\
7	4.79202606862017e-05\\
};
\addplot [color=blue, dashed, forget plot]
  table[row sep=crcr]{%
1	0.00197628458498034\\
2	0.000516795865633024\\
3	0.000233208955223884\\
5	8.49906510283493e-05\\
7	4.35919790758577e-05\\
};
\addplot [color=red, dashed, forget plot]
  table[row sep=crcr]{%
1	0.00181159420289856\\
2	0.000472813238770631\\
3	0.000213219616204685\\
5	7.76638707673249e-05\\
7	3.98247710075572e-05\\
};
\addplot [color=green, dashed, forget plot]
  table[row sep=crcr]{%
1	0.00166666666666659\\
2	0.000434216239687479\\
4	0.000110852455381894\\
6	4.96154800297687e-05\\
};
\addplot [color=black, dashed, forget plot]
  table[row sep=crcr]{%
1	0.00153846153846149\\
2	0.000400160064025545\\
4	0.000102072062876407\\
6	4.56725279744385e-05\\
};
\addplot [color=mycolor1, dashed, forget plot]
  table[row sep=crcr]{%
1	0.00142450142450146\\
2	0.000369959304476586\\
4	9.42951438001049e-05\\
6	4.21816341164838e-05\\
};
\addplot [color=blue, dashed, forget plot]
  table[row sep=crcr]{%
1	0.00132275132275139\\
2	0.000343053173241736\\
4	8.73743993009824e-05\\
6	3.90762377398268e-05\\
};
\addplot [color=red, dashed, forget plot]
  table[row sep=crcr]{%
1	0.00123152709359597\\
2	0.000318979266347774\\
4	8.11886011204132e-05\\
6	3.63015936399869e-05\\
};
\addplot [color=green, dashed, forget plot]
  table[row sep=crcr]{%
1	0.00114942528735629\\
2	0.000297353553374957\\
4	7.56372437788477e-05\\
5	4.85767026134343e-05\\
};
\addplot [color=black, dashed, forget plot]
  table[row sep=crcr]{%
1	0.00107526881720432\\
2	0.000277854959711054\\
4	7.06364342727971e-05\\
5	4.53597024403507e-05\\
};
\addplot [color=mycolor1, dashed, forget plot]
  table[row sep=crcr]{%
1	0.00100806451612911\\
2	0.000260213374967455\\
4	6.61157024793301e-05\\
5	4.24520292069686e-05\\
};
\addplot [color=blue, dashed, forget plot]
  table[row sep=crcr]{%
1	0.000946969696969724\\
2	0.000244200244200243\\
4	6.20155038759674e-05\\
5	3.98152572065602e-05\\
};
\addplot [color=red, dashed, forget plot]
  table[row sep=crcr]{%
1	0.000891265597147916\\
2	0.00022962112514352\\
4	5.82852480037477e-05\\
5	3.74167477363019e-05\\
};
\addplot [color=green, dashed, forget plot]
  table[row sep=crcr]{%
1	0.000840336134453889\\
2	0.000216309755569966\\
4	5.48817298721138e-05\\
5	3.5228633833595e-05\\
};
\addplot [color=black, dashed, forget plot]
  table[row sep=crcr]{%
1	0.000793650793650791\\
2	0.000204123290467428\\
4	5.17678728580994e-05\\
5	3.32270069112184e-05\\
};
\addplot [color=mycolor1, dashed, forget plot]
  table[row sep=crcr]{%
1	0.000750750750750705\\
2	0.000192938452633628\\
4	4.89117143555872e-05\\
};
\addplot [color=blue, dashed, forget plot]
  table[row sep=crcr]{%
1	0.000711237553342792\\
2	0.000182648401826413\\
4	4.62855820412e-05\\
};
\addplot [color=red, dashed, forget plot]
  table[row sep=crcr]{%
1	0.000674763832658409\\
2	0.000173160173160269\\
4	4.38654208887168e-05\\
};
\addplot [color=green, dashed, forget plot]
  table[row sep=crcr]{%
1	0.000641025641025861\\
2	0.000164392569455829\\
4	4.16302402064806e-05\\
};
\addplot [color=black, dashed, forget plot]
  table[row sep=crcr]{%
1	0.000609756097560976\\
2	0.0001562744178778\\
4	3.95616568421575e-05\\
};
\addplot [color=mycolor1, dashed, forget plot]
  table[row sep=crcr]{%
1	0.000580720092915099\\
2	0.000148743120630623\\
4	3.76435159044064e-05\\
};
\addplot [color=blue, dashed, forget plot]
  table[row sep=crcr]{%
1	0.000553709856035534\\
2	0.000141743444365738\\
4	3.5861574323115e-05\\
};
\addplot [color=red, dashed, forget plot]
  table[row sep=crcr]{%
1	0.000528541226215484\\
2	0.000135226504394847\\
4	3.42032356260846e-05\\
};
\addplot [color=green, dashed, forget plot]
  table[row sep=crcr]{%
1	0.000505050505050564\\
2	0.000129148908691645\\
4	3.26573266712648e-05\\
};
\addplot [color=black, dashed, forget plot]
  table[row sep=crcr]{%
1	0.000483091787439771\\
2	0.000123472033584482\\
4	3.12139089178043e-05\\
};
\addplot [color=mycolor1, dashed, forget plot]
  table[row sep=crcr]{%
1	0.00046253469010149\\
2	0.000118161408483919\\
4	2.98641182619186e-05\\
};
\addplot [color=blue, dashed, forget plot]
  table[row sep=crcr]{%
1	0.000443262411347733\\
2	0.000113186191284698\\
4	2.86000286000287e-05\\
};
\addplot [color=red, dashed, forget plot]
  table[row sep=crcr]{%
1	0.000425170068027225\\
2	0.000108518719479089\\
3	4.85672656629355e-05\\
};
\addplot [color=green, dashed, forget plot]
  table[row sep=crcr]{%
1	0.000408163265306016\\
2	0.000104134124752742\\
3	4.65983224604005e-05\\
};
\addplot [color=black, dashed, forget plot]
  table[row sep=crcr]{%
1	0.000392156862745075\\
2	0.000100010001000073\\
3	4.47467334884765e-05\\
};
\addplot [color=mycolor1, dashed, forget plot]
  table[row sep=crcr]{%
1	0.000377073906485781\\
2	9.61261174661222e-05\\
3	4.3003354261606e-05\\
};
\end{axis}
\end{tikzpicture}%

%% file: fig2/cf_general_pareto_pseudo_n1e+08_s1_50_rho0_48_mincha_110522_102323.tex
%
%
\definecolor{mycolor1}{rgb}{1.00000,0.00000,1.00000}%
\begin{tikzpicture}

\begin{axis}[%
width=0.964\fwidth,
height=\fheight,
at={(0\fwidth,0\fheight)},
scale only axis,
xmode=log,
xmin=1,
xmax=100,
xminorticks=true,
xlabel style={font=\color{white!15!black}},
xlabel={$a_1$},
ymode=log,
ymin=5e-05,
ymax=1,
yminorticks=true,
ylabel style={font=\color{white!15!black}},
ylabel={$\Pr(\mathbf{A}_2=[a_1,a_2])$},
axis background/.style={fill=white},
xmajorgrids,
xminorgrids,
ymajorgrids,
yminorgrids,
grid style={dotted},legend style={font=\tiny,draw=none},mark size=1pt
]
\addplot [color=blue, draw=none, mark=square, mark options={solid, blue}, forget plot]
  table[row sep=crcr]{%
1	0.42178694\\
2	0.01114858\\
3	0.00838743999999999\\
4	0.00634795999999999\\
5	0.00490768\\
6	0.00388305\\
7	0.00316702\\
8	0.00260486\\
9.00000000000001	0.00218525\\
10	0.00185326\\
11	0.00159934\\
12	0.00139698\\
13	0.00121023\\
14	0.00107441\\
15	0.000960279999999999\\
16	0.000855329999999999\\
17	0.000775999999999999\\
18	0.000702839999999999\\
19	0.000637489999999999\\
20	0.00058037\\
21	0.00053184\\
22	0.000495149999999999\\
23	0.00045485\\
24	0.0004157\\
25	0.00039288\\
26	0.00036548\\
27	0.00034114\\
28	0.00031528\\
29	0.00029581\\
30	0.00027707\\
31	0.00026615\\
32	0.00025102\\
33	0.00023494\\
34	0.00022381\\
35	0.00020987\\
36	0.00020183\\
37	0.00018931\\
38	0.00018113\\
39	0.00017249\\
40	0.00016425\\
41	0.00015638\\
42	0.00015006\\
43.0000000000001	0.00014259\\
44	0.000138\\
45	0.00012954\\
46	0.00012575\\
47.0000000000001	0.00012183\\
48	0.00011801\\
49	0.00011232\\
50	0.00010796\\
51	0.00010334\\
52	9.994e-05\\
53	9.612e-05\\
54.0000000000001	9.167e-05\\
55	8.99699999999999e-05\\
56	8.709e-05\\
57	8.486e-05\\
58.0000000000001	8.20199999999999e-05\\
59.0000000000001	7.754e-05\\
60	7.50099999999999e-05\\
61	7.45700000000001e-05\\
62.0000000000001	7.045e-05\\
63.0000000000001	6.83e-05\\
64.0000000000001	6.83e-05\\
65	6.46e-05\\
66.0000000000001	6.345e-05\\
67.0000000000001	6.12400000000001e-05\\
68	6.149e-05\\
69.0000000000001	5.975e-05\\
70	5.678e-05\\
71	5.344e-05\\
72.0000000000001	5.314e-05\\
73	5.176e-05\\
};
\addplot [color=red, draw=none, mark=square, mark options={solid, red}, forget plot]
  table[row sep=crcr]{%
1	0.13554337\\
2	0.004052\\
3	0.00307478\\
4	0.00231668\\
5	0.00177443\\
6	0.00139865\\
7	0.00112569\\
8	0.000924029999999999\\
9.00000000000001	0.00076774\\
10	0.00065217\\
11	0.00055985\\
12	0.00048665\\
13	0.000422889999999999\\
14	0.00037605\\
15	0.00033059\\
16	0.00029648\\
17	0.00026546\\
18	0.0002409\\
19	0.00021722\\
20	0.00020032\\
21	0.00018353\\
22	0.00016651\\
23	0.00015353\\
24	0.00014259\\
25	0.00013272\\
26	0.00012326\\
27	0.00011534\\
28	0.00010892\\
29	9.99499999999999e-05\\
30	9.50400000000001e-05\\
31	8.94900000000001e-05\\
32	8.45100000000001e-05\\
33	7.96399999999999e-05\\
34	7.43500000000001e-05\\
35	6.98900000000001e-05\\
36	6.76500000000001e-05\\
37	6.413e-05\\
38	5.99500000000001e-05\\
39	5.687e-05\\
40	5.56400000000001e-05\\
41	5.33e-05\\
42	5.001e-05\\
};
\addplot [color=green, draw=none, mark=square, mark options={solid, green}, forget plot]
  table[row sep=crcr]{%
1	0.0644374199999999\\
2	0.00207981\\
3	0.00159937\\
4	0.00120322\\
5	0.000921819999999999\\
6	0.0007191\\
7	0.00058396\\
8	0.0004723\\
9.00000000000001	0.00039259\\
10	0.00033493\\
11	0.00028601\\
12	0.00024613\\
13	0.00021554\\
14	0.00018942\\
15	0.00016786\\
16	0.00014973\\
17	0.00013509\\
18	0.00012221\\
19	0.00011088\\
20	0.00010225\\
21	9.34400000000001e-05\\
22	8.59999999999999e-05\\
23	7.83299999999999e-05\\
24	7.27400000000001e-05\\
25	6.777e-05\\
26	6.321e-05\\
27	5.868e-05\\
28	5.468e-05\\
29	5.167e-05\\
};
\addplot [color=black, draw=none, mark=square, mark options={solid, black}, forget plot]
  table[row sep=crcr]{%
1	0.03717866\\
2	0.00126872\\
3	0.000974769999999999\\
4	0.00073281\\
5	0.000559449999999999\\
6	0.00043897\\
7	0.00035298\\
8	0.00028759\\
9.00000000000001	0.00023827\\
10	0.00020185\\
11	0.00017512\\
12	0.00014935\\
13	0.00012921\\
14	0.00011306\\
15	0.00010101\\
16	9.28e-05\\
17	8.193e-05\\
18	7.306e-05\\
19	6.76500000000001e-05\\
20	6.123e-05\\
21	5.61100000000001e-05\\
22	5.043e-05\\
};
\addplot [color=mycolor1, draw=none, mark=square, mark options={solid, mycolor1}, forget plot]
  table[row sep=crcr]{%
1	0.02404729\\
2	0.000841679999999999\\
3	0.000659529999999999\\
4	0.000496029999999999\\
5	0.00037951\\
6	0.00029706\\
7	0.00023716\\
8	0.00019332\\
9.00000000000001	0.00016135\\
10	0.00013595\\
11	0.00011504\\
12	0.00010175\\
13	8.637e-05\\
14	7.77399999999999e-05\\
15	6.737e-05\\
16	6.095e-05\\
17	5.426e-05\\
};
\addplot [color=blue, draw=none, mark=square, mark options={solid, blue}, forget plot]
  table[row sep=crcr]{%
1	0.01677925\\
2	0.00061651\\
3	0.000474849999999999\\
4	0.00035758\\
5	0.00027255\\
6	0.0002109\\
7	0.00017026\\
8	0.00013852\\
9.00000000000001	0.00011515\\
10	9.84000000000001e-05\\
11	8.269e-05\\
12	7.15600000000001e-05\\
13	6.21600000000001e-05\\
14	5.565e-05\\
};
\addplot [color=red, draw=none, mark=square, mark options={solid, red}, forget plot]
  table[row sep=crcr]{%
1	0.01236636\\
2	0.00045817\\
3	0.00035753\\
4	0.00026888\\
5	0.00020597\\
6	0.00015917\\
7	0.00012825\\
8	0.00010358\\
9.00000000000001	8.72700000000001e-05\\
10	7.339e-05\\
11	6.377e-05\\
12	5.483e-05\\
};
\addplot [color=green, draw=none, mark=square, mark options={solid, green}, forget plot]
  table[row sep=crcr]{%
1	0.00945971999999999\\
2	0.00035934\\
3	0.00027843\\
4	0.00020948\\
5	0.00015887\\
6	0.0001257\\
7	0.00010051\\
8	8.183e-05\\
9.00000000000001	6.872e-05\\
10	5.819e-05\\
};
\addplot [color=black, draw=none, mark=square, mark options={solid, black}, forget plot]
  table[row sep=crcr]{%
1	0.00750434\\
2	0.00028622\\
3	0.00022234\\
4	0.0001678\\
5	0.00012871\\
6	0.00010162\\
7	7.87999999999999e-05\\
8	6.539e-05\\
9.00000000000001	5.438e-05\\
};
\addplot [color=mycolor1, draw=none, mark=square, mark options={solid, mycolor1}, forget plot]
  table[row sep=crcr]{%
1	0.00605946999999999\\
2	0.00023506\\
3	0.00018322\\
4	0.00013811\\
5	0.00010585\\
6	8.24800000000001e-05\\
7	6.661e-05\\
8	5.456e-05\\
};
\addplot [color=blue, draw=none, mark=square, mark options={solid, blue}, forget plot]
  table[row sep=crcr]{%
1	0.00499982\\
2	0.00019309\\
3	0.00015477\\
4	0.00011484\\
5	8.78500000000001e-05\\
6	6.869e-05\\
7	5.402e-05\\
};
\addplot [color=red, draw=none, mark=square, mark options={solid, red}, forget plot]
  table[row sep=crcr]{%
1	0.00420814\\
2	0.00527723999999999\\
3	0.00013111\\
4	9.824e-05\\
5	7.464e-05\\
6	5.877e-05\\
};
\addplot [color=green, draw=none, mark=square, mark options={solid, green}, forget plot]
  table[row sep=crcr]{%
1	0.00357395\\
2	0.00451795999999999\\
3	0.00011269\\
4	8.499e-05\\
5	6.404e-05\\
};
\addplot [color=black, draw=none, mark=square, mark options={solid, black}, forget plot]
  table[row sep=crcr]{%
1	0.00308835\\
2	0.00392045\\
3	9.59899999999999e-05\\
4	7.16600000000001e-05\\
5	5.593e-05\\
};
\addplot [color=mycolor1, draw=none, mark=square, mark options={solid, mycolor1}, forget plot]
  table[row sep=crcr]{%
1	0.00268259\\
2	0.00343324\\
3	8.681e-05\\
4	6.498e-05\\
};
\addplot [color=blue, draw=none, mark=square, mark options={solid, blue}, forget plot]
  table[row sep=crcr]{%
1	0.00235988\\
2	0.00303217\\
3	7.54e-05\\
4	5.522e-05\\
};
\addplot [color=red, draw=none, mark=square, mark options={solid, red}, forget plot]
  table[row sep=crcr]{%
1	0.00208719\\
2	0.00269006\\
3	6.696e-05\\
};
\addplot [color=green, draw=none, mark=square, mark options={solid, green}, forget plot]
  table[row sep=crcr]{%
1	0.00185546\\
2	0.00241983\\
3	5.919e-05\\
};
\addplot [color=black, draw=none, mark=square, mark options={solid, black}, forget plot]
  table[row sep=crcr]{%
1	0.00167022\\
2	0.00216965\\
3	5.385e-05\\
};
\addplot [color=mycolor1, draw=none, mark=square, mark options={solid, mycolor1}, forget plot]
  table[row sep=crcr]{%
1	0.00149778\\
2	0.00196851\\
};
\addplot [color=blue, draw=none, mark=square, mark options={solid, blue}, forget plot]
  table[row sep=crcr]{%
1	0.0013655\\
2	0.00179585\\
};
\addplot [color=red, draw=none, mark=square, mark options={solid, red}, forget plot]
  table[row sep=crcr]{%
1	0.00124072\\
2	0.00162811\\
};
\addplot [color=green, draw=none, mark=square, mark options={solid, green}, forget plot]
  table[row sep=crcr]{%
1	0.00113295\\
2	0.00149266\\
};
\addplot [color=black, draw=none, mark=square, mark options={solid, black}, forget plot]
  table[row sep=crcr]{%
1	0.00104334\\
2	0.0013756\\
};
\addplot [color=mycolor1, draw=none, mark=square, mark options={solid, mycolor1}, forget plot]
  table[row sep=crcr]{%
1	0.00096539\\
2	0.00127275\\
};
\addplot [color=blue, draw=none, mark=square, mark options={solid, blue}, forget plot]
  table[row sep=crcr]{%
1	0.00088451\\
2	0.00117971\\
};
\addplot [color=red, draw=none, mark=square, mark options={solid, red}, forget plot]
  table[row sep=crcr]{%
1	0.00082032\\
2	0.00109672\\
};
\addplot [color=green, draw=none, mark=square, mark options={solid, green}, forget plot]
  table[row sep=crcr]{%
1	0.00076413\\
2	0.00102095\\
};
\addplot [color=black, draw=none, mark=square, mark options={solid, black}, forget plot]
  table[row sep=crcr]{%
1	0.00070941\\
2	0.00095355\\
};
\addplot [color=mycolor1, draw=none, mark=square, mark options={solid, mycolor1}, forget plot]
  table[row sep=crcr]{%
1	0.00066695\\
2	0.00088851\\
};
\addplot [color=blue, draw=none, mark=square, mark options={solid, blue}, forget plot]
  table[row sep=crcr]{%
1	0.00062301\\
2	0.00083873\\
};
\addplot [color=red, draw=none, mark=square, mark options={solid, red}, forget plot]
  table[row sep=crcr]{%
1	0.00058171\\
2	0.00078589\\
};
\addplot [color=green, draw=none, mark=square, mark options={solid, green}, forget plot]
  table[row sep=crcr]{%
1	0.00054844\\
2	0.0007402\\
};
\addplot [color=black, draw=none, mark=square, mark options={solid, black}, forget plot]
  table[row sep=crcr]{%
1	0.00051691\\
2	0.00069561\\
};
\addplot [color=mycolor1, draw=none, mark=square, mark options={solid, mycolor1}, forget plot]
  table[row sep=crcr]{%
1	0.0004912\\
2	0.00065981\\
};
\addplot [color=blue, draw=none, mark=square, mark options={solid, blue}, forget plot]
  table[row sep=crcr]{%
1	0.000456\\
2	0.00062134\\
};
\addplot [color=red, draw=none, mark=square, mark options={solid, red}, forget plot]
  table[row sep=crcr]{%
1	0.00043701\\
2	0.00059032\\
};
\addplot [color=green, draw=none, mark=square, mark options={solid, green}, forget plot]
  table[row sep=crcr]{%
1	0.00041415\\
2	0.00055925\\
};
\addplot [color=black, draw=none, mark=square, mark options={solid, black}, forget plot]
  table[row sep=crcr]{%
1	0.00038927\\
2	0.00053471\\
};
\addplot [color=mycolor1, draw=none, mark=square, mark options={solid, mycolor1}, forget plot]
  table[row sep=crcr]{%
1	0.0003757\\
2	0.00050859\\
};
\addplot [color=blue, draw=none, mark=square, mark options={solid, blue}, forget plot]
  table[row sep=crcr]{%
1	0.0003552\\
2	0.00048125\\
};
\addplot [color=red, draw=none, mark=square, mark options={solid, red}, forget plot]
  table[row sep=crcr]{%
1	0.00033867\\
2	0.00046089\\
};
\addplot [color=green, draw=none, mark=square, mark options={solid, green}, forget plot]
  table[row sep=crcr]{%
1	0.00032092\\
2	0.00043776\\
};
\addplot [color=black, draw=none, mark=square, mark options={solid, black}, forget plot]
  table[row sep=crcr]{%
1	0.00030833\\
2	0.0004179\\
};
\addplot [color=mycolor1, draw=none, mark=square, mark options={solid, mycolor1}, forget plot]
  table[row sep=crcr]{%
1	0.00029531\\
2	0.0004023\\
};
\addplot [color=blue, draw=none, mark=square, mark options={solid, blue}, forget plot]
  table[row sep=crcr]{%
1	0.00028186\\
2	0.00038352\\
};
\addplot [color=red, draw=none, mark=square, mark options={solid, red}, forget plot]
  table[row sep=crcr]{%
1	0.00027112\\
2	0.00036444\\
};
\addplot [color=green, draw=none, mark=square, mark options={solid, green}, forget plot]
  table[row sep=crcr]{%
1	0.00025853\\
2	0.00035586\\
};
\addplot [color=black, draw=none, mark=square, mark options={solid, black}, forget plot]
  table[row sep=crcr]{%
1	0.0002479\\
2	0.00033786\\
};
\addplot [color=mycolor1, draw=none, mark=square, mark options={solid, mycolor1}, forget plot]
  table[row sep=crcr]{%
1	0.00023846\\
2	0.00032478\\
};
\addplot [color=blue, dashed, forget plot]
  table[row sep=crcr]{%
1	0.421784848588095\\
2	0.0111462841309558\\
3	0.0083867547265466\\
4	0.00634584992507405\\
5	0.00491365165795627\\
6	0.00389690289774613\\
7	0.00315752011763761\\
8	0.00260622283049777\\
9.00000000000001	0.00218558344278982\\
10	0.00185799293192482\\
11	0.00159822802607047\\
13	0.00121799388655511\\
15	0.000958339713907695\\
17	0.000773390025350054\\
20	0.000582049432066898\\
24	0.000420256260601752\\
29	0.000297759382655691\\
35	0.000210252646608387\\
43.0000000000001	0.00014289059130182\\
55	8.95045288471589e-05\\
72.0000000000001	5.33276837744118e-05\\
75	4.92798807030228e-05\\
};
\node[right, align=left, font=\color{blue}]
at (axis cs:1.02,0.381) {{\tiny $a_2=1$}};
\addplot [color=red, dashed, forget plot]
  table[row sep=crcr]{%
1	0.135543651402059\\
2	0.00404521836097951\\
3	0.00307554921500614\\
4	0.00231392191585163\\
5	0.00177732989678694\\
6	0.00139882472030578\\
7	0.00112584264958191\\
8	0.000923930603562592\\
9.00000000000001	0.000770986903295028\\
10	0.00065263849403598\\
11	0.000559322280461609\\
13	0.000423675407328557\\
15	0.000331774361484391\\
18	0.000241128413516318\\
21	0.000183051802473189\\
25	0.000133309653549783\\
30	9.51865798917301e-05\\
37	6.42530596303139e-05\\
43.0000000000001	4.83337739684319e-05\\
};
\node[right, align=left, font=\color{red}]
at (axis cs:1.02,0.122) {{\tiny $a_2=2$}};
\addplot [color=green, dashed, forget plot]
  table[row sep=crcr]{%
1	0.0644572718822734\\
2	0.0020807232231828\\
3	0.00159427256164544\\
4	0.00119764199611859\\
5	0.000917155286062956\\
6	0.000719695315932783\\
7	0.000577725091598428\\
8	0.000473041953344841\\
9.00000000000001	0.000393971316236663\\
10	0.00033293970304448\\
11	0.000284923260670484\\
13	0.00021531426078616\\
15	0.000168298936261922\\
18	0.000122059577991252\\
21	9.25152062976176e-05\\
25	6.72696527419939e-05\\
30	4.79641726549437e-05\\
};
\node[right, align=left, font=\color{green}]
at (axis cs:1.02,0.058) {{\tiny $a_2=3$}};
\addplot [color=black, dashed, forget plot]
  table[row sep=crcr]{%
1	0.0371861274159075\\
2	0.00126592153428805\\
3	0.000975123909952538\\
4	0.000732072077580813\\
5	0.000559735112048824\\
6	0.000438523929040491\\
7	0.000351517651860692\\
8	0.00028746978053343\\
9.00000000000001	0.000239166352497051\\
10	0.000201933445131403\\
11	0.000172675579616328\\
13	0.000130322870958873\\
15	0.000101764553457705\\
18	7.37213793600753e-05\\
21	5.58298057985793e-05\\
23	4.72863502791388e-05\\
};
\node[right, align=left]
at (axis cs:1.02,0.034) {{\tiny $a_2=4$}};
\addplot [color=mycolor1, dashed, forget plot]
  table[row sep=crcr]{%
1	0.0240590411671855\\
2	0.000850891065342439\\
3	0.000657921423935324\\
4	0.000493778791977195\\
5	0.000377168671939986\\
6	0.000295195521412768\\
7	0.000236414269367129\\
8	0.000193188979796213\\
9.00000000000001	0.000160620869745642\\
10	0.000135538442031535\\
11	0.000115843348486425\\
13	8.73597043546597e-05\\
15	6.81733674753753e-05\\
18	4.93516069822606e-05\\
};
\node[right, align=left, font=\color{mycolor1}]
at (axis cs:1.02,0.022) {{\tiny $a_2=5$}};
\addplot [color=blue, dashed, forget plot]
  table[row sep=crcr]{%
1	0.0167894335846358\\
2	0.000611054013684584\\
3	0.000473815066519794\\
4	0.000355541502175191\\
5	0.000271396183868082\\
6	0.000212264634590984\\
7	0.000169891951750156\\
8	0.000138755266456806\\
9.00000000000001	0.000115310884586627\\
10	9.72657759186361e-05\\
11	8.31038440061445e-05\\
13	6.26354532379794e-05\\
15	4.88580663588423e-05\\
};
\addplot [color=red, dashed, forget plot]
  table[row sep=crcr]{%
1	0.0123610356460102\\
2	0.000460028082932781\\
3	0.000357490118737633\\
4	0.00026822423479495\\
5	0.00020464508872988\\
6	0.000159976389800573\\
7	0.000127983620100445\\
8	0.00010448675208664\\
9.00000000000001	8.68033580712782e-05\\
10	7.31983488197031e-05\\
11	6.25250941175755e-05\\
13	4.71061036127075e-05\\
};
\node[right, align=left]
at (axis cs:1.02,0.011) {\tiny $\phantom{a_2}~~~~\!\pmb{\pmb{\vdots}}$};
\addplot [color=green, dashed, forget plot]
  table[row sep=crcr]{%
1	0.00947065149366709\\
2	0.000358804262776196\\
3	0.00027931441759043\\
4	0.000209554141428572\\
5	0.000159823324444705\\
6	0.000124889920604954\\
7	9.98793391777198e-05\\
8	8.151785223729e-05\\
9.00000000000001	6.77043886918277e-05\\
10	5.70802775470591e-05\\
11	4.87479825929376e-05\\
};
\addplot [color=black, dashed, forget plot]
  table[row sep=crcr]{%
1	0.00748288989211032\\
2	0.000287663544144907\\
3	0.000224251847946751\\
4	0.000168235496018867\\
5	0.000128273356542103\\
6	0.000100205551802566\\
7	8.01164166709205e-05\\
8	6.53726609207648e-05\\
9.00000000000001	5.42840804205782e-05\\
10	4.57579298505075e-05\\
};
\addplot [color=mycolor1, dashed, forget plot]
  table[row sep=crcr]{%
1	0.00605874476696377\\
2	0.000235764095378472\\
3	0.000184009151355426\\
4	0.000138040305559982\\
5	0.000105226215931397\\
6	8.21812246623894e-05\\
7	6.56910802566598e-05\\
8	5.35917746320982e-05\\
9.00000000000001	4.44942007056359e-05\\
};
\addplot [color=blue, dashed, forget plot]
  table[row sep=crcr]{%
1	0.00500424957806344\\
2	0.000196743112491914\\
3	0.000153706266893979\\
4	0.000115304705132824\\
5	8.78783969264567e-05\\
6	6.86187382800588e-05\\
7	5.48399783207474e-05\\
8	4.47322258807647e-05\\
};
\addplot [color=red, dashed, forget plot]
  table[row sep=crcr]{%
1	0.0042020503008382\\
2	0.00527045332370953\\
3	0.000130318934820648\\
4	9.77585453148966e-05\\
5	7.44938410728859e-05\\
6	5.81576760426802e-05\\
7	4.64723910892542e-05\\
};
\addplot [color=green, dashed, forget plot]
  table[row sep=crcr]{%
1	0.00357782760368741\\
2	0.00452189938757333\\
3	0.000111891722913731\\
4	8.39341872726897e-05\\
5	6.39507306192097e-05\\
6	4.99193861353332e-05\\
};
\addplot [color=black, dashed, forget plot]
  table[row sep=crcr]{%
1	0.00308269024139606\\
2	0.0039222093050977\\
3	9.71146804862807e-05\\
4	7.28485568982887e-05\\
5	5.54979330040254e-05\\
6	4.3315794217165e-05\\
};
\addplot [color=mycolor1, dashed, forget plot]
  table[row sep=crcr]{%
1	0.00268342241029607\\
2	0.00343437223734996\\
3	8.50834978604707e-05\\
4	6.38230551448412e-05\\
5	4.86171293892446e-05\\
};
\addplot [color=blue, dashed, forget plot]
  table[row sep=crcr]{%
1	0.00235682629140988\\
2	0.00303220445372174\\
3	7.5157569019264e-05\\
4	5.63769986657645e-05\\
5	4.29412768414065e-05\\
};
\addplot [color=red, dashed, forget plot]
  table[row sep=crcr]{%
1	0.00208630823380695\\
2	0.00269675916141348\\
3	6.68728021601526e-05\\
4	5.01621709319317e-05\\
5	3.82045491790526e-05\\
};
\addplot [color=green, dashed, forget plot]
  table[row sep=crcr]{%
1	0.00185974741685513\\
2	0.00241405229550287\\
3	5.98863051838946e-05\\
4	4.4921304677903e-05\\
};
\addplot [color=black, dashed, forget plot]
  table[row sep=crcr]{%
1	0.00166812503671665\\
2	0.00217358049593974\\
3	5.39402898466099e-05\\
4	4.04609916767436e-05\\
};
\addplot [color=mycolor1, dashed, forget plot]
  table[row sep=crcr]{%
1	0.0015046192444067\\
2	0.00196733088734705\\
3	4.88379243108467e-05\\
};
\addplot [color=blue, dashed, forget plot]
  table[row sep=crcr]{%
1	0.00136399353744687\\
2	0.00178910460479977\\
3	4.44268181300587e-05\\
};
\addplot [color=red, dashed, forget plot]
  table[row sep=crcr]{%
1	0.00124217478153021\\
2	0.00163404538571663\\
3	4.05875017716489e-05\\
};
\addplot [color=green, dashed, forget plot]
  table[row sep=crcr]{%
1	0.00113595657831922\\
2	0.00149830513867349\\
3	3.72252459521631e-05\\
};
\addplot [color=black, dashed, forget plot]
  table[row sep=crcr]{%
1	0.00104278719547646\\
2	0.00137880281841052\\
3	3.4264158271518e-05\\
};
\addplot [color=mycolor1, dashed, forget plot]
  table[row sep=crcr]{%
1	0.000960615608952219\\
2	0.00127304799829054\\
3	3.16428603600636e-05\\
};
\addplot [color=blue, dashed, forget plot]
  table[row sep=crcr]{%
1	0.000887778158040197\\
2	0.00117901003474884\\
3	2.93112797579577e-05\\
};
\addplot [color=red, dashed, forget plot]
  table[row sep=crcr]{%
1	0.000822914022983845\\
2	0.00109501983913159\\
3	2.7228239683766e-05\\
};
\addplot [color=green, dashed, forget plot]
  table[row sep=crcr]{%
1	0.000764901448206511\\
2	0.00101969528927619\\
3	2.53596276844097e-05\\
};
\addplot [color=black, dashed, forget plot]
  table[row sep=crcr]{%
1	0.000712809092809908\\
2	0.000951883995423598\\
3	2.36769895445088e-05\\
};
\addplot [color=mycolor1, dashed, forget plot]
  table[row sep=crcr]{%
1	0.000665858544679598\\
2	0.000890618954514894\\
3	2.2156439222521e-05\\
};
\addplot [color=blue, dashed, forget plot]
  table[row sep=crcr]{%
1	0.00062339516520582\\
2	0.000835083879522546\\
3	2.07778061640769e-05\\
};
\addplot [color=red, dashed, forget plot]
  table[row sep=crcr]{%
1	0.000584865215116803\\
2	0.000784585864459009\\
3	1.95239626988871e-05\\
};
\addplot [color=green, dashed, forget plot]
  table[row sep=crcr]{%
1	0.000549797761953629\\
2	0.000738533663371031\\
2.99657829230977	1.85723562146892e-05\\
};
\addplot [color=black, dashed, forget plot]
  table[row sep=crcr]{%
1	0.000517790260551099\\
2	0.000696420303234739\\
2.9773636117174	1.85723562146892e-05\\
};
\addplot [color=mycolor1, dashed, forget plot]
  table[row sep=crcr]{%
1	0.000488496978784947\\
2	0.000657809069925103\\
2.9588128470171	1.85723562146892e-05\\
};
\addplot [color=blue, dashed, forget plot]
  table[row sep=crcr]{%
1	0.000461619644845368\\
2	0.000622322139600479\\
2.94088541044711	1.85723562146892e-05\\
};
\addplot [color=red, dashed, forget plot]
  table[row sep=crcr]{%
1	0.0004368998419636\\
2	0.000589631299804452\\
2.92354421669142	1.85723562146892e-05\\
};
\addplot [color=green, dashed, forget plot]
  table[row sep=crcr]{%
1	0.000414112787319953\\
2	0.000559450332526932\\
2.90675529710256	1.85723562146892e-05\\
};
\addplot [color=black, dashed, forget plot]
  table[row sep=crcr]{%
1	0.000393062214638084\\
2	0.00053152872753947\\
2.89048746515275	1.85723562146892e-05\\
};
\addplot [color=mycolor1, dashed, forget plot]
  table[row sep=crcr]{%
1	0.000373576142312481\\
2	0.000505646466937952\\
2.87471202520438	1.85723562146892e-05\\
};
\addplot [color=blue, dashed, forget plot]
  table[row sep=crcr]{%
1	0.00035550335623058\\
2	0.000481609677274445\\
2.85940251808327	1.85723562146892e-05\\
};
\addplot [color=red, dashed, forget plot]
  table[row sep=crcr]{%
1	0.000338710472653325\\
2	0.000459246988153652\\
2.84453449803665	1.85723562146892e-05\\
};
\addplot [color=green, dashed, forget plot]
  table[row sep=crcr]{%
1	0.000323079474377819\\
2	0.000438406469100274\\
2.83008533657446	1.85723562146892e-05\\
};
\addplot [color=black, dashed, forget plot]
  table[row sep=crcr]{%
1	0.000308505635012613\\
2	0.000418953042065378\\
2.81603404941842	1.85723562146892e-05\\
};
\addplot [color=mycolor1, dashed, forget plot]
  table[row sep=crcr]{%
1	0.000294895763061964\\
2	0.000400766287002273\\
2.8023611433826	1.85723562146892e-05\\
};
\addplot [color=blue, dashed, forget plot]
  table[row sep=crcr]{%
1	0.000282166710732223\\
2	0.000383738573729878\\
2.78904848051497	1.85723562146892e-05\\
};
\addplot [color=red, dashed, forget plot]
  table[row sep=crcr]{%
1	0.0002702441028265\\
2	0.000367773465772054\\
2.77607915722903	1.85723562146892e-05\\
};
\addplot [color=green, dashed, forget plot]
  table[row sep=crcr]{%
1	0.000259061249385128\\
2	0.000352784351867162\\
2.76343739648771	1.85723562146892e-05\\
};
\addplot [color=black, dashed, forget plot]
  table[row sep=crcr]{%
1	0.000248558212335937\\
2	0.000338693268754714\\
2.75110845140867	1.85723562146892e-05\\
};
\addplot [color=mycolor1, dashed, forget plot]
  table[row sep=crcr]{%
1	0.000238681001724195\\
2	0.000325429885302215\\
2.73907851886123	1.85723562146892e-05\\
};
\end{axis}
\end{tikzpicture}%

%% file: fig2/cf_gk_a1_a2_n1e+08_mincha_110522_101623.tex
%
%
\definecolor{mycolor1}{rgb}{1.00000,0.00000,1.00000}%
\begin{tikzpicture}

\begin{axis}[%
width=0.951\fwidth,
height=\fheight,
at={(0\fwidth,0\fheight)},
scale only axis,
xmode=log,
xmin=1,
xmax=100,
xminorticks=true,
xlabel style={font=\color{white!15!black}},
xlabel={$a$},
ymode=log,
ymin=0.0001,
ymax=1,
yminorticks=true,
ylabel style={font=\color{white!15!black}},
ylabel={$\Pr(A_{j}=a)$},
axis background/.style={fill=white},
xmajorgrids,
xminorgrids,
ymajorgrids,
yminorgrids,
legend style={legend cell align=left, align=left,
  draw=white!15!black},
legend pos={south west},
grid style={dotted},legend style={font=\scriptsize,draw={none}, fill=none},mark size=0.5pt
]
\addplot [color=blue]
  table[row sep=crcr]{%
1	0.415037499278844\\
2	0.169925001442312\\
3	0.0931094043914815\\
4	0.0588936890535686\\
5	0.0406419844973458\\
6	0.0297473433940519\\
7	0.0227200765000834\\
8	0.0179219079972623\\
9	0.0144995696951152\\
10	0.0119726416660759\\
11	0.0100536646639231\\
12	0.008562013503424\\
13	0.0073795303655975\\
14	0.00642626915943309\\
15	0.00564656314114204\\
16	0.0050006810583666\\
17	0.00445964819069988\\
18	0.00400193055749626\\
19	0.00361125355237883\\
20	0.00327513203286088\\
21	0.00298385843882141\\
22	0.00272979275557238\\
23	0.00250685561057491\\
24	0.00231016068720104\\
25	0.00213574434399099\\
26	0.00198036412824081\\
27	0.00184134682416746\\
28	0.0017164725890216\\
29	0.00160388570258964\\
30	0.0015020251652318\\
31	0.00140957025467135\\
32	0.00132539746656754\\
33	0.00124854619725898\\
34	0.00117819119728114\\
35	0.00111362031070861\\
36	0.00105421637200174\\
37	0.000999442395972731\\
38	0.000948829393548805\\
39	0.000901966294392579\\
40	0.000858491570044747\\
41	0.000818086237338932\\
42	0.000780467988138446\\
43	0.00074538624282195\\
44	0.000712617965039369\\
45	0.000681964106713503\\
46	0.000653246577105653\\
47	0.000626305649466702\\
48	0.000600997734535647\\
49	0.000577193462745493\\
50	0.000554776027174385\\
51	0.000533639747489665\\
52	0.000513688821837545\\
53	0.000494836239078218\\
54	0.000477002828246661\\
55	0.000460116425806821\\
56	0.000444111144307132\\
57	0.000428926728561379\\
58	0.000414507987591855\\
59	0.000400804292309638\\
60	0.000387769130378782\\
61	0.000375359710947564\\
62	0.000363536612957722\\
63	0.00035226347162902\\
64	0.000341506698455589\\
65	0.000331235230680024\\
66	0.000321420306752048\\
67	0.000312035264737475\\
68	0.000303055361032243\\
69	0.000294457607081855\\
70	0.000286220622085304\\
71	0.000278324499925382\\
72	0.000270750688772219\\
73	0.000263481882001594\\
74	0.00025650191922659\\
75	0.000249795696388886\\
76	0.000243349083968874\\
77	0.00023714885249243\\
78	0.000231182604595318\\
79	0.000225438712997055\\
80	0.000219906263803531\\
81	0.000214575004617836\\
82	0.000209435297005569\\
83	0.000204478072894203\\
84	0.000199694794545239\\
85	0.000195077417765719\\
86	0.000190618358061512\\
87	0.000186310459468345\\
88	0.000182146965823484\\
89	0.000178121494255287\\
90	0.000174228010705049\\
91	0.000170460807298065\\
92	0.00016681448141252\\
93	0.000163283916295482\\
94	0.000159864263102127\\
95	0.000156550924236891\\
96	0.000153339537890928\\
97	0.000150225963679218\\
98	0.000147206269286409\\
99	0.000144276718044916\\
100	0.00014143375736906\\
};
\addlegendentry{Gauss-Kuz'min}

\addplot [color=green]
  table[row sep=crcr]{%
1	0.5\\
2	0.166666666666667\\
3	0.0833333333333333\\
4	0.05\\
5	0.0333333333333333\\
6	0.0238095238095238\\
7	0.0178571428571429\\
8	0.0138888888888889\\
9	0.0111111111111111\\
10	0.00909090909090909\\
11	0.00757575757575758\\
12	0.00641025641025641\\
13	0.00549450549450549\\
14	0.00476190476190476\\
15	0.00416666666666667\\
16	0.00367647058823529\\
17	0.00326797385620915\\
18	0.00292397660818713\\
19	0.00263157894736842\\
20	0.00238095238095238\\
21	0.00216450216450216\\
22	0.00197628458498024\\
23	0.00181159420289855\\
24	0.00166666666666667\\
25	0.00153846153846154\\
26	0.00142450142450142\\
27	0.00132275132275132\\
28	0.00123152709359606\\
29	0.00114942528735632\\
30	0.0010752688172043\\
31	0.00100806451612903\\
32	0.000946969696969697\\
33	0.00089126559714795\\
34	0.000840336134453782\\
35	0.000793650793650794\\
36	0.000750750750750751\\
37	0.000711237553342817\\
38	0.000674763832658569\\
39	0.000641025641025641\\
40	0.000609756097560976\\
41	0.000580720092915215\\
42	0.000553709856035437\\
43	0.000528541226215645\\
44	0.000505050505050505\\
45	0.000483091787439614\\
46	0.000462534690101758\\
47	0.000443262411347518\\
48	0.000425170068027211\\
49	0.000408163265306122\\
50	0.000392156862745098\\
51	0.000377073906485671\\
52	0.000362844702467344\\
53	0.000349406009783368\\
54	0.000336700336700337\\
55	0.000324675324675325\\
56	0.00031328320802005\\
57	0.000302480338777979\\
58	0.000292226767971946\\
59	0.000282485875706215\\
60	0.000273224043715847\\
61	0.000264410364886304\\
62	0.000256016385048643\\
63	0.000248015873015873\\
64	0.000240384615384615\\
65	0.000233100233100233\\
66	0.000226142017186793\\
67	0.000219490781387182\\
68	0.000213128729752771\\
69	0.00020703933747412\\
70	0.000201207243460765\\
71	0.000195618153364632\\
72	0.000190258751902588\\
73	0.000185116623472788\\
74	0.00018018018018018\\
75	0.000175438596491228\\
76	0.000170881749829118\\
77	0.000166500166500167\\
78	0.000162284972411555\\
79	0.000158227848101266\\
80	0.000154320987654321\\
81	0.000150557061126167\\
82	0.000146929180135175\\
83	0.000143430866322433\\
84	0.000140056022408964\\
85	0.000136798905608755\\
86	0.000133654103180968\\
87	0.000130616509926855\\
88	0.000127681307456588\\
89	0.000124843945068664\\
90	0.000122100122100122\\
91	0.000119445771619685\\
92	0.000116877045348294\\
93	0.000114390299702585\\
94	0.000111982082866741\\
95	0.000109649122807018\\
96	0.000107388316151203\\
97	0.000105196717862403\\
98	0.00010307153164296\\
99	0.000101010101010101\\
100	9.9009900990099e-05\\
};
\addlegendentry{Benford, $j=1$}

\addplot [color=green, line width=1.0pt, draw=none, mark=square, mark options={solid, green}, forget plot]
  table[row sep=crcr]{%
1	0.50002927\\
2	0.16667335\\
3	0.08329263\\
4	0.04997441\\
5	0.03333291\\
6	0.02379779\\
7	0.0178598\\
8	0.01389124\\
9	0.01110989\\
10	0.00909978\\
11	0.00759082\\
12	0.00641799\\
13	0.00549178\\
14	0.0047714\\
15	0.00415738\\
16	0.00367186\\
17	0.00326882\\
18	0.002927\\
19	0.00262948\\
20	0.00238214\\
21	0.00216793\\
22	0.0019708\\
23	0.00182137\\
24	0.00166462\\
25	0.00153168\\
26	0.00142965\\
27	0.00132534\\
28	0.00122783\\
29	0.00115342\\
30	0.00107473\\
31	0.00101201\\
32	0.00094492\\
33	0.00089409\\
34	0.00084018\\
35	0.00079021\\
36	0.00075641\\
37	0.00071119\\
38	0.00067725\\
39	0.00063684\\
40	0.00061167\\
41	0.00057783\\
42	0.00055388\\
43	0.00053429\\
44	0.00050467\\
45	0.00048675\\
46	0.00046109\\
47	0.00044314\\
48	0.00042432\\
49	0.00040899\\
50	0.00039342\\
51	0.00037797\\
52	0.00036343\\
53	0.00035073\\
54	0.00033919\\
55	0.0003249\\
56	0.00031237\\
57	0.00030233\\
58	0.00029276\\
59	0.00028089\\
60	0.00027253\\
61	0.00026239\\
62	0.00025718\\
63	0.00025119\\
64	0.00023863\\
65	0.00023599\\
66	0.00022332\\
67	0.00021988\\
68	0.00021246\\
69	0.00021016\\
70	0.00020166\\
71	0.00019367\\
72	0.0001915\\
73	0.00018403\\
74	0.00017956\\
75	0.00017701\\
76	0.00016944\\
77	0.00016797\\
78	0.00016036\\
79	0.00015817\\
80	0.00015256\\
81	0.00014857\\
82	0.00014589\\
83	0.00014315\\
84	0.00013971\\
85	0.000136\\
86	0.00013337\\
87	0.00012994\\
88	0.00012457\\
89	0.00012417\\
90	0.00012227\\
91	0.00011942\\
92	0.00011635\\
93	0.00011507\\
94	0.00011195\\
95	0.00010893\\
96	0.00010721\\
97	0.00010501\\
98	0.00010351\\
99	0.00010146\\
100	9.897e-05\\
};
\addplot [color=red]
  table[row sep=crcr]{%
1	0.386294361119891\\
2	0.168523753999383\\
3	0.095419753355459\\
4	0.0615863632159228\\
5	0.0430876087827263\\
6	0.0318527368494288\\
7	0.0245124210715115\\
8	0.0194506138359767\\
9	0.0158116632792451\\
10	0.0131074815826282\\
11	0.0110429064337878\\
12	0.00943083234071462\\
13	0.00814794987929224\\
14	0.00711029571737543\\
15	0.00625908369522943\\
16	0.00555215412469723\\
17	0.0049586308592362\\
18	0.00445548032333942\\
19	0.00402523593249304\\
20	0.00365445751451787\\
21	0.00333266631355089\\
22	0.00305159483686568\\
23	0.00280464945904013\\
24	0.00258651949298505\\
25	0.0023928888070297\\
26	0.00222022035127667\\
27	0.00206559325794364\\
28	0.00192657834714316\\
29	0.00180114202581239\\
30	0.00168757141173792\\
31	0.00158441548834443\\
32	0.00149043848379149\\
33	0.00140458265548071\\
34	0.00132593837197392\\
35	0.00125371990139511\\
36	0.00118724569526574\\
37	0.00112592223842611\\
38	0.00106923074640664\\
39	0.00101671615053078\\
40	0.00096797793188097\\
41	0.000922662457709222\\
42	0.000880456545289121\\
43	0.000841082033532281\\
44	0.000804291185946493\\
45	0.000769862782582731\\
46	0.000737598785371052\\
47	0.000707321482705137\\
48	0.000678871036123851\\
49	0.000652103365650802\\
50	0.000626888321407559\\
51	0.000603108098026439\\
52	0.000580655855719314\\
53	0.000559434517770185\\
54	0.000539355719129553\\
55	0.00052033888476144\\
56	0.000502310419809415\\
57	0.000485202996297396\\
58	0.000468954923444898\\
59	0.00045350959059931\\
60	0.000438814973343415\\
61	0.000424823194729473\\
62	0.000411490134728076\\
63	0.000398775081909886\\
64	0.000386640422243456\\
65	0.00037505136053273\\
66	0.000363975670663486\\
67	0.000353383471284552\\
68	0.00034324702400701\\
69	0.000333540551598954\\
70	0.000324240073900306\\
71	0.000315323259553413\\
72	0.000306769291802844\\
73	0.000298558746870903\\
74	0.000290673483578807\\
75	0.000283096543033601\\
76	0.000275812057358493\\
77	0.000268805166531827\\
78	0.000262061942537528\\
79	0.000255569320078308\\
80	0.000249315033241881\\
81	0.000243287557514904\\
82	0.000237476056646813\\
83	0.000231870333911033\\
84	0.000226460787336347\\
85	0.000221238368558474\\
86	0.000216194544958359\\
87	0.000211321264772302\\
88	0.000206610924950112\\
89	0.000202056341463308\\
90	0.000197650721907472\\
91	0.000193387640149645\\
92	0.000189261012889741\\
93	0.000185265077932595\\
94	0.000181394374072052\\
95	0.000177643722412135\\
96	0.000174008209034682\\
97	0.000170483168897562\\
98	0.000167064170855991\\
99	0.000163747003736781\\
100	0.000160527663351395\\
};
\addlegendentry{Benford, $j=2$}

\addplot [color=red, line width=1.0pt, draw=none, mark=square, mark options={solid, red}, forget plot]
  table[row sep=crcr]{%
1	0.38634814\\
2	0.16847878\\
3	0.09542806\\
4	0.06157394\\
5	0.04309\\
6	0.03183798\\
7	0.02451577\\
8	0.0194677\\
9	0.01579727\\
10	0.01310286\\
11	0.01105674\\
12	0.00942835\\
13	0.00814547\\
14	0.0070998\\
15	0.00626508\\
16	0.00553878\\
17	0.00496026\\
18	0.00446282\\
19	0.00401788\\
20	0.00366921\\
21	0.00334078\\
22	0.00305883\\
23	0.00280052\\
24	0.0025825\\
25	0.00238823\\
26	0.00222838\\
27	0.00206251\\
28	0.00191771\\
29	0.00179593\\
30	0.00168715\\
31	0.00159268\\
32	0.00149054\\
33	0.00140516\\
34	0.00133036\\
35	0.00125188\\
36	0.00118971\\
37	0.00112663\\
38	0.00106926\\
39	0.00101757\\
40	0.00096959\\
41	0.00092271\\
42	0.00088209\\
43	0.00084149\\
44	0.00079913\\
45	0.00076857\\
46	0.00073576\\
47	0.00070984\\
48	0.00067653\\
49	0.00065544\\
50	0.00062995\\
51	0.00060617\\
52	0.00057496\\
53	0.00056004\\
54	0.0005406\\
55	0.00052201\\
56	0.00050319\\
57	0.00048563\\
58	0.00046618\\
59	0.0004501\\
60	0.00044015\\
61	0.00042274\\
62	0.000412\\
63	0.00039699\\
64	0.00039029\\
65	0.00037548\\
66	0.00036528\\
67	0.00035396\\
68	0.00034328\\
69	0.00033167\\
70	0.00031923\\
71	0.00031958\\
72	0.00030656\\
73	0.00030069\\
74	0.0002883\\
75	0.00028318\\
76	0.00027735\\
77	0.00027026\\
78	0.00026202\\
79	0.00025799\\
80	0.00024631\\
81	0.00024246\\
82	0.00023739\\
83	0.00023032\\
84	0.00022964\\
85	0.0002225\\
86	0.00021484\\
87	0.00021177\\
88	0.00020669\\
89	0.00020425\\
90	0.00019736\\
91	0.00019571\\
92	0.00018831\\
93	0.0001864\\
94	0.0001818\\
95	0.00017508\\
96	0.00017351\\
97	0.00017188\\
98	0.0001673\\
99	0.0001629\\
100	0.00016152\\
};
\addplot [color=mycolor1]
  table[row sep=crcr]{%
1	0.792351979697098\\
2	0.090464404273033\\
3	0.018075036413832\\
4	0.0136229737660956\\
5	0.0104804746788819\\
6	0.00826023876599684\\
7	0.0066565673511907\\
8	0.00546875153356785\\
9	0.00456786424970426\\
10	0.00386995324788507\\
11	0.0033190899420765\\
12	0.00287709595108604\\
13	0.00251729303578713\\
14	0.00222063179255403\\
15	0.00197323989574349\\
16	0.00176483353804578\\
17	0.00158766527757544\\
18	0.00143581263344839\\
19	0.00130468778309762\\
20	0.00119069362733416\\
21	0.00109097857695021\\
22	0.00100325908960951\\
23	0.000925689455462747\\
24	0.000856765028035381\\
25	0.000795249458529891\\
26	0.000740119379281313\\
27	0.000690521923730738\\
28	0.000645741794918052\\
29	0.000605175510558261\\
30	0.000568311094383707\\
31	0.000534711938214937\\
32	0.000504003885199902\\
33	0.000475864820783296\\
34	0.000450016230701923\\
35	0.000426216312851354\\
36	0.00040425432488242\\
37	0.000383945920756385\\
38	0.000365129283519568\\
39	0.000347661902771289\\
40	0.000331417876957426\\
41	0.000316285645103516\\
42	0.000302166071656003\\
43	0.000288970823020786\\
44	0.000276620986139054\\
45	0.000265045888746986\\
46	0.000254182088375583\\
47	0.000243972503075715\\
48	0.000234365661621144\\
49	0.00022531505479242\\
50	0.000216778572471173\\
51	0.000208718013818196\\
52	0.000201098659894578\\
53	0.00019388889979401\\
54	0.000187059902767895\\
55	0.000180585329990546\\
56	0.000174441080582721\\
57	0.00016860506731943\\
58	0.000163057018123526\\
59	0.000157778300012924\\
60	0.000152751762645797\\
61	0.000147961599010817\\
62	0.000143393221148374\\
63	0.000139033149078243\\
64	0.000134868911354402\\
65	0.000130888955876635\\
66	0.000127082569767957\\
67	0.000123439807279874\\
68	0.000119951424819715\\
69	0.000116608822306857\\
70	0.000113403990163797\\
71	0.00011032946133113\\
72	0.000107378267770064\\
73	0.000104543900979267\\
74	0.000101820276107977\\
75	9.9201699296288e-05\\
76	9.66828379148322e-05\\
77	9.42586934137092e-05\\
78	9.19245765223331e-05\\
79	8.96760845706525e-05\\
80	8.75090807264634e-05\\
81	8.54196749666913e-05\\
82	8.34042066185938e-05\\
83	8.14592283246233e-05\\
84	7.95814913003136e-05\\
85	7.77679317668229e-05\\
86	7.60156584527268e-05\\
87	7.43219410696776e-05\\
88	7.26841996763685e-05\\
89	7.10999948531912e-05\\
90	6.95670186178919e-05\\
91	6.80830860187912e-05\\
92	6.66461273491139e-05\\
93	6.52541809293832e-05\\
94	6.39053864122909e-05\\
95	6.25979785662817e-05\\
96	6.13302814995121e-05\\
97	6.01007032892335e-05\\
98	5.89077309834523e-05\\
99	5.77499259469566e-05\\
100	5.66259195235687e-05\\
};
\addlegendentry{Pareto, $j=1$, $s=1.5$, $\rho=0.48$}

\addplot [color=mycolor1, line width=1.0pt, draw=none, mark=square, mark options={solid, mycolor1}, forget plot]
  table[row sep=crcr]{%
1	0.79229645\\
2	0.09050813\\
3	0.01809873\\
4	0.01361983\\
5	0.01046277\\
6	0.00826769\\
7	0.00665963\\
8	0.00545508\\
9	0.00457494\\
10	0.00386929\\
11	0.00332415\\
12	0.00287371\\
13	0.00252342\\
14	0.00221828\\
15	0.00196957\\
16	0.00176693\\
17	0.00158604\\
18	0.00143656\\
19	0.00130322\\
20	0.00119222\\
21	0.00109309\\
22	0.00100682\\
23	0.00092574\\
24	0.00085558\\
25	0.00079634\\
26	0.0007353\\
27	0.00068837\\
28	0.00064087\\
29	0.0006035\\
30	0.00056968\\
31	0.00053709\\
32	0.000504\\
33	0.00047475\\
34	0.00044671\\
35	0.00042484\\
36	0.00040271\\
37	0.00038324\\
38	0.00036612\\
39	0.00034721\\
40	0.00033566\\
41	0.00031596\\
42	0.00030402\\
43	0.00029053\\
44	0.0002785\\
45	0.00026552\\
46	0.00025515\\
47	0.00024555\\
48	0.00023522\\
49	0.00022603\\
50	0.00021709\\
51	0.00020867\\
52	0.00020041\\
53	0.00019477\\
54	0.00018692\\
55	0.00018285\\
56	0.00017667\\
57	0.00016942\\
58	0.00016285\\
59	0.00015907\\
60	0.00015049\\
61	0.0001488\\
62	0.00014561\\
63	0.0001363\\
64	0.0001329\\
65	0.00013204\\
66	0.00012764\\
67	0.00012561\\
68	0.00011914\\
69	0.00011638\\
70	0.00011294\\
71	0.00010902\\
72	0.0001072\\
73	0.00010569\\
74	0.00010168\\
75	9.872e-05\\
76	9.523e-05\\
77	9.428e-05\\
78	9.267e-05\\
79	9.064e-05\\
80	8.771e-05\\
81	8.409e-05\\
82	8.258e-05\\
83	7.97e-05\\
84	8.018e-05\\
85	7.654e-05\\
86	7.621e-05\\
87	7.324e-05\\
88	7.141e-05\\
89	7.21e-05\\
90	6.987e-05\\
91	6.687e-05\\
92	6.644e-05\\
93	6.541e-05\\
94	6.382e-05\\
95	6.242e-05\\
96	6.213e-05\\
97	5.965e-05\\
98	5.941e-05\\
99	5.784e-05\\
100	5.578e-05\\
};
\end{axis}
\end{tikzpicture}%

%% file: fig2/cf_general_benford_US_total_income_per_ZIP_code,_2016_n159924_mincha_110522_092226.tex
%
%
\definecolor{mycolor1}{rgb}{1.00000,0.00000,1.00000}%
\begin{tikzpicture}

\begin{axis}[%
width=0.964\fwidth,
height=\fheight,
at={(0\fwidth,0\fheight)},
scale only axis,
unbounded coords=jump,
xmode=log,
xmin=1,
xmax=100,
xminorticks=true,
xlabel style={font=\color{white!15!black}},
xlabel={$a_1$},
ymode=log,
ymin=5e-05,
ymax=1,
yminorticks=true,
ylabel style={font=\color{white!15!black}},
ylabel={$\Pr(\mathbf{A}_2=[a_1,a_2])$},
axis background/.style={fill=white},
xmajorgrids,
xminorgrids,
ymajorgrids,
yminorgrids,
grid style={dotted},legend style={font=\tiny,draw=none},mark size=1pt
]
\addplot [color=blue, draw=none, mark=square, mark options={solid, blue}, forget plot]
  table[row sep=crcr]{%
1	0.163802806333008\\
2	0.0657437282709287\\
3	0.0355356294239764\\
4	0.022041719816913\\
5	0.0152822590730597\\
6	0.0110489982741802\\
7	0.00889172356869513\\
8	0.00663440134063681\\
9.00000000000001	0.00562142017458293\\
10	0.00443335584402591\\
11	0.00361421675295765\\
12	0.00326405042395138\\
13	0.00259498261674295\\
14	0.00226982816837998\\
15	0.00212600985468097\\
16	0.00163827818213652\\
17	0.00155073659988494\\
18	0.00160701333133238\\
19	0.00131937670393437\\
20	0.0013256296740952\\
21	0.00104424601685801\\
22	0.00100672819589305\\
23	0.000906680673319826\\
24	0.000812886120907431\\
25	0.00076286235962082\\
26	0.000781621270103299\\
27	0.000794127210424951\\
28	0.000619044045921812\\
29	0.00058777919511768\\
30	0.000481478702383632\\
31	0.00043145494109702\\
32	0.000443960881418674\\
33	0.000456466821740327\\
34	0.000493984642705284\\
35	0.000400190090292889\\
36	0.000368925239488757\\
37	0.000331407418523799\\
38	0.00038143117981041\\
39	0.000250118806433056\\
40	0.000262624746754709\\
41	0.000350166329006278\\
42	0.000281383657237187\\
43.0000000000001	0.000281383657237187\\
44	0.000275130687076361\\
45	0.00022510692578975\\
46	0.000187589104824792\\
47.0000000000001	0.000218853955628924\\
48	0.000200095045146445\\
49	0.000212600985468097\\
50	0.000187589104824792\\
51	0.000231359895950576\\
52	0.000150071283859833\\
53	0.000125059403216528\\
54.0000000000001	0.000262624746754709\\
55	0.000181336134663965\\
56	0.00013756534353818\\
57	0.000125059403216528\\
58.0000000000001	0.000200095045146445\\
59.0000000000001	0.000125059403216528\\
60	0.000118806433055701\\
61	0.000181336134663965\\
62.0000000000001	0.000112553462894875\\
63.0000000000001	0.000125059403216528\\
64.0000000000001	0.000112553462894875\\
65	0.000106300492734049\\
66.0000000000001	0.000162577224181486\\
67.0000000000001	0.000100047522573222\\
68	0.000118806433055701\\
69.0000000000001	8.12886120907431e-05\\
70	9.37945524123959e-05\\
71	0.000150071283859833\\
nan	nan\\
73	8.75415822515695e-05\\
74	8.12886120907431e-05\\
75	6.87826717690904e-05\\
76	0.000125059403216528\\
77.0000000000001	9.37945524123959e-05\\
78.0000000000001	9.37945524123959e-05\\
nan	nan\\
80	0.000162577224181486\\
81	8.12886120907431e-05\\
82	5.62767314474376e-05\\
83	8.12886120907431e-05\\
nan	nan\\
85	5.00237612866112e-05\\
85.9999999999999	6.25297016082639e-05\\
86.9999999999999	8.75415822515695e-05\\
88.0000000000001	6.87826717690904e-05\\
89.0000000000001	0.000118806433055701\\
nan	nan\\
91.0000000000001	7.50356419299167e-05\\
92.0000000000001	0.000100047522573222\\
93.0000000000001	5.62767314474376e-05\\
93.9999999999999	5.62767314474376e-05\\
95	6.25297016082639e-05\\
nan	nan\\
98.9999999999999	7.50356419299167e-05\\
};
\addplot [color=red, draw=none, mark=square, mark options={solid, red}, forget plot]
  table[row sep=crcr]{%
1	0.0795752982666767\\
2	0.029388959755884\\
3	0.0150071283859834\\
4	0.0085790750606538\\
5	0.00584027413021185\\
6	0.00400815387308971\\
7	0.00290763112478427\\
8	0.00241989945223981\\
9.00000000000001	0.00190715589905205\\
10	0.00148195392811585\\
11	0.00125059403216528\\
12	0.00108176383782296\\
13	0.00102548710637553\\
14	0.000825392061229083\\
15	0.000681573747530077\\
16	0.000700332658012555\\
17	0.000619044045921812\\
18	0.000387684149971236\\
19	0.000462719791901153\\
20	0.000318901478202146\\
21	0.000281383657237187\\
22	0.000281383657237187\\
23	0.000318901478202146\\
24	0.000200095045146445\\
25	0.000287636627398014\\
26	0.000256371776593882\\
27	0.000218853955628924\\
28	0.000168830194342312\\
29	0.000175083164503139\\
30	0.00013756534353818\\
31	0.000150071283859833\\
32	0.000187589104824792\\
33	0.000150071283859833\\
34	0.000162577224181486\\
35	8.12886120907431e-05\\
36	0.000150071283859833\\
37	0.000125059403216528\\
38	0.000112553462894875\\
39	6.87826717690904e-05\\
40	7.50356419299167e-05\\
41	6.25297016082639e-05\\
42	0.000100047522573222\\
43.0000000000001	7.50356419299167e-05\\
44	8.75415822515695e-05\\
45	0.000131312373377354\\
nan	nan\\
47.0000000000001	9.37945524123959e-05\\
nan	nan\\
50	7.50356419299167e-05\\
nan	nan\\
53	0.000125059403216528\\
nan	nan\\
56	5.62767314474376e-05\\
nan	nan\\
60	0.000125059403216528\\
nan	nan\\
62.0000000000001	5.00237612866112e-05\\
63.0000000000001	6.87826717690904e-05\\
nan	nan\\
69.0000000000001	5.00237612866112e-05\\
70	5.62767314474376e-05\\
nan	nan\\
75	6.25297016082639e-05\\
nan	nan\\
78.0000000000001	5.00237612866112e-05\\
nan	nan\\
81	5.62767314474376e-05\\
nan	nan\\
83	9.37945524123959e-05\\
nan	nan\\
85.9999999999999	5.00237612866112e-05\\
86.9999999999999	5.00237612866112e-05\\
nan	nan\\
90.0000000000001	6.25297016082639e-05\\
nan	nan\\
93.9999999999999	5.00237612866112e-05\\
};
\addplot [color=green, draw=none, mark=square, mark options={solid, green}, forget plot]
  table[row sep=crcr]{%
1	0.0491545984342562\\
2	0.0167329481503714\\
3	0.00768490032765563\\
4	0.00490232860608789\\
5	0.00281383657237188\\
6	0.00208849203371602\\
7	0.00165078412245817\\
8	0.00115054650959205\\
9.00000000000001	0.000944198494284784\\
10	0.000719091568495035\\
11	0.000569020284635202\\
12	0.000631549986243465\\
13	0.000412696030614542\\
14	0.00038143117981041\\
15	0.00036267226932793\\
16	0.000368925239488757\\
17	0.000237612866111403\\
18	0.000237612866111403\\
19	0.000243865836272229\\
20	0.00022510692578975\\
21	0.00022510692578975\\
22	0.000200095045146445\\
23	0.000206348015307271\\
24	7.50356419299167e-05\\
25	0.000112553462894875\\
26	0.000118806433055701\\
27	0.000168830194342312\\
28	8.12886120907431e-05\\
29	8.12886120907431e-05\\
30	8.75415822515695e-05\\
31	5.00237612866112e-05\\
32	0.000100047522573222\\
33	7.50356419299167e-05\\
34	8.75415822515695e-05\\
35	8.12886120907431e-05\\
nan	nan\\
37	7.50356419299167e-05\\
38	6.87826717690904e-05\\
39	5.00237612866112e-05\\
nan	nan\\
41	6.25297016082639e-05\\
42	0.000143818313699007\\
43.0000000000001	5.62767314474376e-05\\
nan	nan\\
46	6.87826717690904e-05\\
nan	nan\\
52	6.87826717690904e-05\\
nan	nan\\
55	5.00237612866112e-05\\
nan	nan\\
57	9.37945524123959e-05\\
nan	nan\\
61	5.62767314474376e-05\\
};
\addplot [color=black, draw=none, mark=square, mark options={solid, black}, forget plot]
  table[row sep=crcr]{%
1	0.0335346789725119\\
2	0.0103111477952027\\
3	0.00503364097946525\\
4	0.00289512518446262\\
5	0.00171956679422726\\
6	0.0012318351216828\\
7	0.000987969285410569\\
8	0.00067532077736925\\
9.00000000000001	0.000619044045921812\\
10	0.000462719791901153\\
11	0.000331407418523799\\
12	0.000375178209649584\\
13	0.000212600985468097\\
14	0.000256371776593882\\
15	0.000306395537880493\\
16	0.000125059403216528\\
17	0.000168830194342312\\
18	0.000143818313699007\\
19	0.000150071283859833\\
20	8.12886120907431e-05\\
21	0.000125059403216528\\
22	7.50356419299167e-05\\
23	6.87826717690904e-05\\
nan	nan\\
25	9.37945524123959e-05\\
26	6.25297016082639e-05\\
27	5.00237612866112e-05\\
28	5.00237612866112e-05\\
29	0.000112553462894875\\
30	6.87826717690904e-05\\
nan	nan\\
33	5.00237612866112e-05\\
nan	nan\\
37	5.00237612866112e-05\\
nan	nan\\
43.0000000000001	5.00237612866112e-05\\
nan	nan\\
93.0000000000001	0.000100047522573222\\
};
\addplot [color=mycolor1, draw=none, mark=square, mark options={solid, mycolor1}, forget plot]
  table[row sep=crcr]{%
1	0.0234236262224556\\
2	0.00667817213176258\\
3	0.00322027963282559\\
4	0.00193842074985618\\
5	0.00134438858457767\\
6	0.000994222255571397\\
7	0.00069407968785173\\
8	0.000525249493509417\\
9.00000000000001	0.000387684149971236\\
10	0.000331407418523799\\
11	0.000331407418523799\\
12	0.000275130687076361\\
13	0.000200095045146445\\
14	0.000100047522573222\\
15	0.000125059403216528\\
16	0.000168830194342312\\
17	8.12886120907431e-05\\
18	0.000200095045146445\\
19	5.00237612866112e-05\\
nan	nan\\
21	6.25297016082639e-05\\
22	7.50356419299167e-05\\
nan	nan\\
40	8.12886120907431e-05\\
41	6.25297016082639e-05\\
nan	nan\\
47.0000000000001	9.37945524123959e-05\\
nan	nan\\
51	0.000106300492734049\\
nan	nan\\
60	5.00237612866112e-05\\
};
\addplot [color=blue, draw=none, mark=square, mark options={solid, blue}, forget plot]
  table[row sep=crcr]{%
1	0.0179460243615717\\
2	0.0057089617568345\\
3	0.00238863460143568\\
4	0.00138190640554263\\
5	0.000700332658012555\\
6	0.000619044045921812\\
7	0.000493984642705284\\
8	0.000331407418523799\\
9.00000000000001	0.00031264850804132\\
10	0.000212600985468097\\
11	0.000187589104824792\\
12	0.00013756534353818\\
13	0.000112553462894875\\
14	0.00013756534353818\\
15	0.000162577224181486\\
16	0.000106300492734049\\
17	6.87826717690904e-05\\
18	5.00237612866112e-05\\
19	6.87826717690904e-05\\
20	0.000150071283859833\\
21	5.62767314474376e-05\\
nan	nan\\
23	6.87826717690904e-05\\
24	9.37945524123959e-05\\
25	5.00237612866112e-05\\
26	6.25297016082639e-05\\
nan	nan\\
31	5.62767314474376e-05\\
nan	nan\\
92.0000000000001	5.00237612866112e-05\\
};
\addplot [color=red, draw=none, mark=square, mark options={solid, red}, forget plot]
  table[row sep=crcr]{%
1	0.013975388309447\\
2	0.00400815387308971\\
3	0.00190090292889122\\
4	0.00110052274830544\\
5	0.000700332658012555\\
6	0.000500237612866112\\
7	0.000275130687076361\\
8	0.000306395537880493\\
9.00000000000001	0.000243865836272229\\
10	0.000256371776593882\\
11	0.00015632425402066\\
12	0.00015632425402066\\
13	5.62767314474376e-05\\
14	0.00013756534353818\\
15	8.75415822515695e-05\\
16	6.25297016082639e-05\\
17	8.75415822515695e-05\\
nan	nan\\
28	5.62767314474376e-05\\
nan	nan\\
50	6.25297016082639e-05\\
};
\addplot [color=green, draw=none, mark=square, mark options={solid, green}, forget plot]
  table[row sep=crcr]{%
1	0.0113303819314174\\
2	0.00319526775218229\\
3	0.00145694204747255\\
4	0.000787874240264125\\
5	0.000506490583026937\\
6	0.000300142567719666\\
7	0.000356419299167104\\
8	0.00022510692578975\\
9.00000000000001	0.000243865836272229\\
10	0.000150071283859833\\
11	0.000125059403216528\\
12	5.00237612866112e-05\\
13	6.25297016082639e-05\\
14	0.000100047522573222\\
15	0.000106300492734049\\
16	7.50356419299167e-05\\
17	8.12886120907431e-05\\
18	6.87826717690904e-05\\
nan	nan\\
25	0.000106300492734049\\
nan	nan\\
33	8.75415822515695e-05\\
nan	nan\\
49	7.50356419299167e-05\\
};
\addplot [color=black, draw=none, mark=square, mark options={solid, black}, forget plot]
  table[row sep=crcr]{%
1	0.00940446712188289\\
2	0.00241364648207899\\
3	0.00126935294264776\\
4	0.000644055926565118\\
5	0.000425201970936194\\
6	0.000318901478202146\\
7	0.000231359895950576\\
8	0.000168830194342312\\
9.00000000000001	0.00015632425402066\\
10	0.000118806433055701\\
nan	nan\\
12	7.50356419299167e-05\\
13	9.37945524123959e-05\\
nan	nan\\
15	5.62767314474376e-05\\
nan	nan\\
17	5.00237612866112e-05\\
nan	nan\\
65	7.50356419299167e-05\\
nan	nan\\
73	9.37945524123959e-05\\
};
\addplot [color=mycolor1, draw=none, mark=square, mark options={solid, mycolor1}, forget plot]
  table[row sep=crcr]{%
1	0.00775368299942472\\
2	0.00210099797403767\\
3	0.000956704434606437\\
4	0.000531502463670243\\
5	0.000350166329006278\\
6	0.000300142567719666\\
7	0.000181336134663965\\
8	0.000193842074985618\\
9.00000000000001	0.000150071283859833\\
10	8.75415822515695e-05\\
11	6.25297016082639e-05\\
nan	nan\\
13	5.00237612866112e-05\\
14	6.25297016082639e-05\\
};
\addplot [color=blue, draw=none, mark=square, mark options={solid, blue}, forget plot]
  table[row sep=crcr]{%
1	0.00643430629549035\\
2	0.00173832570470974\\
3	0.000875415822515695\\
4	0.000462719791901153\\
5	0.000268877716915535\\
6	0.000275130687076361\\
7	0.00022510692578975\\
8	8.75415822515695e-05\\
9.00000000000001	8.12886120907431e-05\\
10	5.00237612866112e-05\\
11	6.25297016082639e-05\\
12	7.50356419299167e-05\\
nan	nan\\
19	6.87826717690904e-05\\
nan	nan\\
97.0000000000001	6.25297016082639e-05\\
};
\addplot [color=red, draw=none, mark=square, mark options={solid, red}, forget plot]
  table[row sep=crcr]{%
1	0.00550886671168805\\
2	0.00151947174908081\\
3	0.000594032165278507\\
4	0.000318901478202146\\
5	0.000275130687076361\\
6	8.75415822515695e-05\\
7	0.00015632425402066\\
8	0.000118806433055701\\
nan	nan\\
10	5.00237612866112e-05\\
11	6.25297016082639e-05\\
12	5.00237612866112e-05\\
nan	nan\\
14	5.62767314474376e-05\\
15	6.25297016082639e-05\\
};
\addplot [color=green, draw=none, mark=square, mark options={solid, green}, forget plot]
  table[row sep=crcr]{%
1	0.00465846276981566\\
2	0.00153197768940247\\
3	0.00051899652334859\\
4	0.000350166329006278\\
5	0.000268877716915535\\
6	0.000181336134663965\\
nan	nan\\
8	0.000112553462894875\\
9.00000000000001	8.12886120907431e-05\\
nan	nan\\
12	5.00237612866112e-05\\
nan	nan\\
16	5.62767314474376e-05\\
};
\addplot [color=black, draw=none, mark=square, mark options={solid, black}, forget plot]
  table[row sep=crcr]{%
1	0.00425827267952277\\
2	0.00108801680798379\\
3	0.000487731672544458\\
4	0.000281383657237187\\
5	0.00022510692578975\\
6	8.12886120907431e-05\\
7	0.000112553462894875\\
8	7.50356419299167e-05\\
9.00000000000001	0.000150071283859833\\
10	0.000150071283859833\\
};
\addplot [color=mycolor1, draw=none, mark=square, mark options={solid, mycolor1}, forget plot]
  table[row sep=crcr]{%
1	0.00376428803681749\\
2	0.00121932918136115\\
3	0.00038143117981041\\
4	0.000237612866111403\\
5	0.000118806433055701\\
6	0.000175083164503139\\
7	7.50356419299167e-05\\
8	7.50356419299167e-05\\
nan	nan\\
10	7.50356419299167e-05\\
nan	nan\\
22	5.00237612866112e-05\\
nan	nan\\
30	5.00237612866112e-05\\
};
\addplot [color=blue, draw=none, mark=square, mark options={solid, blue}, forget plot]
  table[row sep=crcr]{%
1	0.00311397914009154\\
2	0.000919186613641478\\
3	0.000418949000775368\\
4	0.00022510692578975\\
5	0.00022510692578975\\
6	5.62767314474376e-05\\
7	8.12886120907431e-05\\
8	6.25297016082639e-05\\
9.00000000000001	8.12886120907431e-05\\
};
\addplot [color=red, draw=none, mark=square, mark options={solid, red}, forget plot]
  table[row sep=crcr]{%
1	0.00282634251269353\\
2	0.000769115329781646\\
3	0.000337660388684625\\
4	0.000200095045146445\\
5	9.37945524123959e-05\\
6	5.62767314474376e-05\\
7	6.25297016082639e-05\\
nan	nan\\
14	6.25297016082639e-05\\
};
\addplot [color=green, draw=none, mark=square, mark options={solid, green}, forget plot]
  table[row sep=crcr]{%
1	0.00261374152722543\\
2	0.000631549986243465\\
3	0.00038143117981041\\
4	0.000200095045146445\\
5	0.00015632425402066\\
6	7.50356419299167e-05\\
7	0.000106300492734049\\
8	0.000112553462894875\\
nan	nan\\
13	6.25297016082639e-05\\
};
\addplot [color=black, draw=none, mark=square, mark options={solid, black}, forget plot]
  table[row sep=crcr]{%
1	0.00245116430304394\\
2	0.000556514344313548\\
3	0.000325154448362973\\
4	6.25297016082639e-05\\
5	0.000118806433055701\\
6	0.000100047522573222\\
7	6.25297016082639e-05\\
};
\addplot [color=mycolor1, draw=none, mark=square, mark options={solid, mycolor1}, forget plot]
  table[row sep=crcr]{%
1	0.00216978064580676\\
2	0.000456466821740327\\
3	0.000287636627398014\\
4	0.00013756534353818\\
5	0.000118806433055701\\
6	8.75415822515695e-05\\
nan	nan\\
13	6.87826717690904e-05\\
};
\addplot [color=blue, draw=none, mark=square, mark options={solid, blue}, forget plot]
  table[row sep=crcr]{%
1	0.00184462619744378\\
2	0.000562767314474375\\
3	0.000218853955628924\\
4	8.75415822515695e-05\\
5	7.50356419299167e-05\\
nan	nan\\
7	6.25297016082639e-05\\
};
\addplot [color=red, draw=none, mark=square, mark options={solid, red}, forget plot]
  table[row sep=crcr]{%
1	0.00195092669017783\\
2	0.00036267226932793\\
3	0.00022510692578975\\
4	0.000106300492734049\\
5	9.37945524123959e-05\\
6	5.62767314474376e-05\\
};
\addplot [color=green, draw=none, mark=square, mark options={solid, green}, forget plot]
  table[row sep=crcr]{%
1	0.00163827818213652\\
2	0.000456466821740327\\
3	0.000168830194342312\\
4	6.25297016082639e-05\\
5	5.00237612866112e-05\\
6	5.62767314474376e-05\\
nan	nan\\
8	5.00237612866112e-05\\
9.00000000000001	8.12886120907431e-05\\
};
\addplot [color=black, draw=none, mark=square, mark options={solid, black}, forget plot]
  table[row sep=crcr]{%
1	0.00139441234586429\\
2	0.000368925239488757\\
3	0.000218853955628924\\
4	9.37945524123959e-05\\
5	8.75415822515695e-05\\
6	6.87826717690904e-05\\
nan	nan\\
11	5.00237612866112e-05\\
nan	nan\\
31	9.37945524123959e-05\\
};
\addplot [color=mycolor1, draw=none, mark=square, mark options={solid, mycolor1}, forget plot]
  table[row sep=crcr]{%
1	0.00140066531602511\\
2	0.000375178209649584\\
3	0.000231359895950576\\
4	9.37945524123959e-05\\
5	6.25297016082639e-05\\
nan	nan\\
7	6.25297016082639e-05\\
};
\addplot [color=blue, draw=none, mark=square, mark options={solid, blue}, forget plot]
  table[row sep=crcr]{%
1	0.00140066531602511\\
2	0.000300142567719666\\
3	0.000181336134663965\\
4	5.62767314474376e-05\\
5	7.50356419299167e-05\\
nan	nan\\
17	5.00237612866112e-05\\
};
\addplot [color=red, draw=none, mark=square, mark options={solid, red}, forget plot]
  table[row sep=crcr]{%
1	0.0012318351216828\\
2	0.000268877716915535\\
3	0.000193842074985618\\
4	0.000118806433055701\\
5	6.25297016082639e-05\\
nan	nan\\
11	0.000100047522573222\\
};
\addplot [color=green, draw=none, mark=square, mark options={solid, green}, forget plot]
  table[row sep=crcr]{%
1	0.00118181136039619\\
2	0.000356419299167104\\
3	7.50356419299167e-05\\
4	5.00237612866112e-05\\
};
\addplot [color=black, draw=none, mark=square, mark options={solid, black}, forget plot]
  table[row sep=crcr]{%
1	0.00105049898701883\\
2	0.000331407418523799\\
3	0.000162577224181486\\
4	5.00237612866112e-05\\
5	9.37945524123959e-05\\
};
\addplot [color=mycolor1, draw=none, mark=square, mark options={solid, mycolor1}, forget plot]
  table[row sep=crcr]{%
1	0.000894174732998173\\
2	0.000287636627398014\\
3	0.00013756534353818\\
4	0.000131312373377354\\
nan	nan\\
34	9.37945524123959e-05\\
};
\addplot [color=blue, draw=none, mark=square, mark options={solid, blue}, forget plot]
  table[row sep=crcr]{%
1	0.00106925789750131\\
2	0.000206348015307271\\
3	0.000168830194342312\\
nan	nan\\
8	0.000100047522573222\\
nan	nan\\
10	7.50356419299167e-05\\
};
\addplot [color=red, draw=none, mark=square, mark options={solid, red}, forget plot]
  table[row sep=crcr]{%
1	0.000994222255571397\\
2	0.000287636627398014\\
nan	nan\\
4	5.62767314474376e-05\\
};
\addplot [color=green, draw=none, mark=square, mark options={solid, green}, forget plot]
  table[row sep=crcr]{%
1	0.000950451464445612\\
2	0.000181336134663965\\
3	0.00015632425402066\\
nan	nan\\
6	6.87826717690904e-05\\
};
\addplot [color=black, draw=none, mark=square, mark options={solid, black}, forget plot]
  table[row sep=crcr]{%
1	0.000800380180585777\\
2	0.000175083164503139\\
nan	nan\\
4	7.50356419299167e-05\\
nan	nan\\
9.00000000000001	5.62767314474376e-05\\
};
\addplot [color=mycolor1, draw=none, mark=square, mark options={solid, mycolor1}, forget plot]
  table[row sep=crcr]{%
1	0.00083164503138991\\
2	0.000162577224181486\\
3	0.000187589104824792\\
};
\addplot [color=blue, draw=none, mark=square, mark options={solid, blue}, forget plot]
  table[row sep=crcr]{%
1	0.000725344538655861\\
2	0.00022510692578975\\
nan	nan\\
12	8.12886120907431e-05\\
};
\addplot [color=red, draw=none, mark=square, mark options={solid, red}, forget plot]
  table[row sep=crcr]{%
1	0.000706585628173383\\
2	0.000143818313699007\\
3	0.000118806433055701\\
nan	nan\\
5	8.75415822515695e-05\\
nan	nan\\
36	5.00237612866112e-05\\
};
\addplot [color=green, draw=none, mark=square, mark options={solid, green}, forget plot]
  table[row sep=crcr]{%
1	0.000650308896725945\\
2	0.000131312373377354\\
};
\addplot [color=black, draw=none, mark=square, mark options={solid, black}, forget plot]
  table[row sep=crcr]{%
1	0.000525249493509417\\
2	0.000193842074985618\\
};
\addplot [color=mycolor1, draw=none, mark=square, mark options={solid, mycolor1}, forget plot]
  table[row sep=crcr]{%
1	0.000669067807208424\\
2	0.000112553462894875\\
3	0.000175083164503139\\
4	6.87826717690904e-05\\
};
\addplot [color=blue, draw=none, mark=square, mark options={solid, blue}, forget plot]
  table[row sep=crcr]{%
1	0.00053775543383107\\
2	0.000150071283859833\\
};
\addplot [color=red, draw=none, mark=square, mark options={solid, red}, forget plot]
  table[row sep=crcr]{%
1	0.000600285135439334\\
2	0.000193842074985618\\
};
\addplot [color=green, draw=none, mark=square, mark options={solid, green}, forget plot]
  table[row sep=crcr]{%
1	0.000443960881418674\\
2	0.000100047522573222\\
3	8.75415822515695e-05\\
nan	nan\\
7	8.12886120907431e-05\\
nan	nan\\
9.00000000000001	5.00237612866112e-05\\
};
\addplot [color=black, draw=none, mark=square, mark options={solid, black}, forget plot]
  table[row sep=crcr]{%
1	0.000468972762061979\\
2	0.000143818313699007\\
};
\addplot [color=mycolor1, draw=none, mark=square, mark options={solid, mycolor1}, forget plot]
  table[row sep=crcr]{%
1	0.000400190090292889\\
2	0.000106300492734049\\
};
\addplot [color=blue, draw=none, mark=square, mark options={solid, blue}, forget plot]
  table[row sep=crcr]{%
1	0.000350166329006278\\
2	8.75415822515695e-05\\
};
\addplot [color=red, draw=none, mark=square, mark options={solid, red}, forget plot]
  table[row sep=crcr]{%
1	0.000462719791901153\\
2	8.12886120907431e-05\\
};
\addplot [color=green, draw=none, mark=square, mark options={solid, green}, forget plot]
  table[row sep=crcr]{%
1	0.000343913358845452\\
2	0.000131312373377354\\
nan	nan\\
4	6.87826717690904e-05\\
};
\addplot [color=black, draw=none, mark=square, mark options={solid, black}, forget plot]
  table[row sep=crcr]{%
1	0.000412696030614542\\
2	0.000106300492734049\\
3	5.62767314474376e-05\\
};
\addplot [color=mycolor1, draw=none, mark=square, mark options={solid, mycolor1}, forget plot]
  table[row sep=crcr]{%
1	0.000400190090292889\\
2	8.75415822515695e-05\\
nan	nan\\
5	6.25297016082639e-05\\
};
\addplot [color=blue, dashed, forget plot]
  table[row sep=crcr]{%
1	0.166666666666667\\
2	0.0666666666666667\\
3	0.0357142857142857\\
4	0.0222222222222222\\
5	0.0151515151515152\\
6	0.010989010989011\\
7	0.00833333333333333\\
9.00000000000001	0.00526315789473683\\
11	0.0036231884057971\\
14	0.00229885057471264\\
18	0.00142247510668564\\
24	0.00081632653061224\\
34	0.000414078674948242\\
51	0.000186706497386106\\
84.0000000000001	6.96136442742777e-05\\
100	4.92586572090042e-05\\
};
\node[right, align=left, font=\color{blue}]
at (axis cs:1.02,0.157) {{\tiny $a_2=1$}};
\addplot [color=red, dashed, forget plot]
  table[row sep=crcr]{%
1	0.0833333333333333\\
2	0.0285714285714285\\
3	0.0142857142857143\\
4	0.00854700854700857\\
5	0.00568181818181818\\
6	0.00404858299595142\\
8	0.00235294117647058\\
11	0.00127877237851662\\
15	0.00070126227208976\\
22	0.000331674958540632\\
34	0.000140706345856195\\
58.0000000000001	4.8840048840048e-05\\
};
\node[right, align=left, font=\color{red}]
at (axis cs:1.02,0.078) {{\tiny $a_2=2$}};
\addplot [color=green, dashed, forget plot]
  table[row sep=crcr]{%
1	0.05\\
2	0.0158730158730159\\
3	0.00769230769230772\\
4	0.00452488687782803\\
5	0.00297619047619047\\
7	0.00156739811912224\\
10	0.000786782061369007\\
14	0.000407996736026109\\
21	0.000183823529411764\\
35	6.69075337883068e-05\\
41	4.88758553274689e-05\\
};
\node[right, align=left, font=\color{green}]
at (axis cs:1.02,0.047) {{\tiny $a_2=3$}};
\addplot [color=black, dashed, forget plot]
  table[row sep=crcr]{%
1	0.0333333333333333\\
2	0.0101010101010101\\
3	0.00480769230769229\\
4	0.00280112044817926\\
5	0.00183150183150183\\
7	0.000957854406130276\\
10	0.000478240076518421\\
15	0.000215703192407241\\
25	7.85792865000758e-05\\
32	4.81486831335129e-05\\
};
\node[right, align=left]
at (axis cs:1.02,0.031) {{\tiny $a_2=4$}};
\addplot [color=mycolor1, dashed, forget plot]
  table[row sep=crcr]{%
1	0.0238095238095237\\
2	0.00699300699300703\\
3	0.00328947368421057\\
4	0.00190476190476191\\
6	0.000871839581516986\\
9.00000000000001	0.000395256916996045\\
14	0.000165700082850037\\
24	5.69962952408076e-05\\
26	4.86215782564276e-05\\
};
\node[right, align=left, font=\color{mycolor1}]
at (axis cs:1.02,0.022) {{\tiny $a_2=5$}};
\addplot [color=blue, dashed, forget plot]
  table[row sep=crcr]{%
1	0.0178571428571429\\
2	0.0051282051282051\\
3	0.00239234449760761\\
4	0.00137931034482758\\
6	0.000628535512256444\\
9.00000000000001	0.000284090909090901\\
15	0.000103669914990667\\
22	4.85083676934339e-05\\
};
\addplot [color=red, dashed, forget plot]
  table[row sep=crcr]{%
1	0.0138888888888888\\
2	0.00392156862745096\\
3	0.00181818181818183\\
4	0.00104493207941486\\
6	0.000474608448030378\\
9.00000000000001	0.000214041095890405\\
15	7.79666302822418e-05\\
19	4.8775729197155e-05\\
};
\node[right, align=left]
at (axis cs:1.02,0.013) {\tiny $\phantom{a_2}~~~~\!\pmb{\pmb{\vdots}}$};
\addplot [color=green, dashed, forget plot]
  table[row sep=crcr]{%
1	0.0111111111111111\\
2	0.00309597523219812\\
3	0.00142857142857139\\
4	0.000819000819000825\\
6	0.000371057513914669\\
10	0.000135666802333484\\
17	4.73978576168363e-05\\
};
\addplot [color=black, dashed, forget plot]
  table[row sep=crcr]{%
1	0.00909090909090914\\
2	0.0025062656641604\\
3	0.00115207373271892\\
4	0.000659195781147009\\
6	0.000298062593144577\\
10	0.000108802089000104\\
15	4.86949746786175e-05\\
};
\addplot [color=mycolor1, dashed, forget plot]
  table[row sep=crcr]{%
1	0.00757575757575767\\
2	0.00207039337474124\\
3	0.000948766603415585\\
5	0.00035014005602238\\
8	0.000138715494520713\\
14	4.57561198810297e-05\\
};
\addplot [color=blue, dashed, forget plot]
  table[row sep=crcr]{%
1	0.00641025641025639\\
2	0.00173913043478258\\
3	0.000794912559618388\\
5	0.000292740046838435\\
8	0.000115834588208047\\
13	4.4232130219396e-05\\
};
\addplot [color=red, dashed, forget plot]
  table[row sep=crcr]{%
1	0.00549450549450547\\
2	0.00148148148148147\\
3	0.000675675675675669\\
5	0.000248385494287129\\
9.00000000000001	7.77484061576794e-05\\
12	4.39270810454611e-05\\
};
\addplot [color=green, dashed, forget plot]
  table[row sep=crcr]{%
1	0.00476190476190474\\
2	0.00127713920817368\\
3	0.000581395348837221\\
5	0.000213401621852321\\
9.00000000000001	6.67289470172166e-05\\
11	4.48028673835227e-05\\
};
\addplot [color=black, dashed, forget plot]
  table[row sep=crcr]{%
1	0.00416666666666665\\
2	0.0011123470522803\\
3	0.000505561172901958\\
5	0.000185322461082271\\
9.00000000000001	5.78971746178769e-05\\
10	4.69682025268925e-05\\
};
\addplot [color=mycolor1, dashed, forget plot]
  table[row sep=crcr]{%
1	0.00367647058823528\\
2	0.000977517106549419\\
3	0.000443655723158797\\
5	0.000162443144899271\\
9.00000000000001	5.07099391480748e-05\\
10	4.11336432067644e-05\\
};
\addplot [color=blue, dashed, forget plot]
  table[row sep=crcr]{%
1	0.00326797385620914\\
2	0.000865800865800902\\
3	0.000392464678178994\\
5	0.000143554407120328\\
9.00000000000001	4.47828034035036e-05\\
};
\addplot [color=red, dashed, forget plot]
  table[row sep=crcr]{%
1	0.0029239766081871\\
2	0.000772200772200748\\
3	0.000349650349650343\\
5	0.00012777919754664\\
9.00000000000001	3.98374631503384e-05\\
};
\addplot [color=green, dashed, forget plot]
  table[row sep=crcr]{%
1	0.00263157894736854\\
2	0.0006930006930006\\
3	0.000313479623824442\\
5	0.000114468864468864\\
8	4.50755014649584e-05\\
};
\addplot [color=black, dashed, forget plot]
  table[row sep=crcr]{%
1	0.00238095238095226\\
2	0.000625390869293418\\
3	0.000282645562464712\\
5	0.000103135313531344\\
8	4.0595948524319e-05\\
};
\addplot [color=mycolor1, dashed, forget plot]
  table[row sep=crcr]{%
1	0.00216450216450214\\
2	0.000567214974475316\\
3	0.000256147540983575\\
5	9.34055669717937e-05\\
7	4.79202606862017e-05\\
};
\addplot [color=blue, dashed, forget plot]
  table[row sep=crcr]{%
1	0.00197628458498034\\
2	0.000516795865633024\\
3	0.000233208955223884\\
5	8.49906510283493e-05\\
7	4.35919790758577e-05\\
};
\addplot [color=red, dashed, forget plot]
  table[row sep=crcr]{%
1	0.00181159420289856\\
2	0.000472813238770631\\
3	0.000213219616204685\\
5	7.76638707673249e-05\\
7	3.98247710075572e-05\\
};
\addplot [color=green, dashed, forget plot]
  table[row sep=crcr]{%
1	0.00166666666666659\\
2	0.000434216239687479\\
4	0.000110852455381894\\
6	4.96154800297687e-05\\
};
\addplot [color=black, dashed, forget plot]
  table[row sep=crcr]{%
1	0.00153846153846149\\
2	0.000400160064025545\\
4	0.000102072062876407\\
6	4.56725279744385e-05\\
};
\addplot [color=mycolor1, dashed, forget plot]
  table[row sep=crcr]{%
1	0.00142450142450146\\
2	0.000369959304476586\\
4	9.42951438001049e-05\\
6	4.21816341164838e-05\\
};
\addplot [color=blue, dashed, forget plot]
  table[row sep=crcr]{%
1	0.00132275132275139\\
2	0.000343053173241736\\
4	8.73743993009824e-05\\
6	3.90762377398268e-05\\
};
\addplot [color=red, dashed, forget plot]
  table[row sep=crcr]{%
1	0.00123152709359597\\
2	0.000318979266347774\\
4	8.11886011204132e-05\\
6	3.63015936399869e-05\\
};
\addplot [color=green, dashed, forget plot]
  table[row sep=crcr]{%
1	0.00114942528735629\\
2	0.000297353553374957\\
4	7.56372437788477e-05\\
5	4.85767026134343e-05\\
};
\addplot [color=black, dashed, forget plot]
  table[row sep=crcr]{%
1	0.00107526881720432\\
2	0.000277854959711054\\
4	7.06364342727971e-05\\
5	4.53597024403507e-05\\
};
\addplot [color=mycolor1, dashed, forget plot]
  table[row sep=crcr]{%
1	0.00100806451612911\\
2	0.000260213374967455\\
4	6.61157024793301e-05\\
5	4.24520292069686e-05\\
};
\addplot [color=blue, dashed, forget plot]
  table[row sep=crcr]{%
1	0.000946969696969724\\
2	0.000244200244200243\\
4	6.20155038759674e-05\\
5	3.98152572065602e-05\\
};
\addplot [color=red, dashed, forget plot]
  table[row sep=crcr]{%
1	0.000891265597147916\\
2	0.00022962112514352\\
4	5.82852480037477e-05\\
5	3.74167477363019e-05\\
};
\addplot [color=green, dashed, forget plot]
  table[row sep=crcr]{%
1	0.000840336134453889\\
2	0.000216309755569966\\
4	5.48817298721138e-05\\
5	3.5228633833595e-05\\
};
\addplot [color=black, dashed, forget plot]
  table[row sep=crcr]{%
1	0.000793650793650791\\
2	0.000204123290467428\\
4	5.17678728580994e-05\\
5	3.32270069112184e-05\\
};
\addplot [color=mycolor1, dashed, forget plot]
  table[row sep=crcr]{%
1	0.000750750750750705\\
2	0.000192938452633628\\
4	4.89117143555872e-05\\
};
\addplot [color=blue, dashed, forget plot]
  table[row sep=crcr]{%
1	0.000711237553342792\\
2	0.000182648401826413\\
4	4.62855820412e-05\\
};
\addplot [color=red, dashed, forget plot]
  table[row sep=crcr]{%
1	0.000674763832658409\\
2	0.000173160173160269\\
4	4.38654208887168e-05\\
};
\addplot [color=green, dashed, forget plot]
  table[row sep=crcr]{%
1	0.000641025641025861\\
2	0.000164392569455829\\
4	4.16302402064806e-05\\
};
\addplot [color=black, dashed, forget plot]
  table[row sep=crcr]{%
1	0.000609756097560976\\
2	0.0001562744178778\\
4	3.95616568421575e-05\\
};
\addplot [color=mycolor1, dashed, forget plot]
  table[row sep=crcr]{%
1	0.000580720092915099\\
2	0.000148743120630623\\
4	3.76435159044064e-05\\
};
\addplot [color=blue, dashed, forget plot]
  table[row sep=crcr]{%
1	0.000553709856035534\\
2	0.000141743444365738\\
4	3.5861574323115e-05\\
};
\addplot [color=red, dashed, forget plot]
  table[row sep=crcr]{%
1	0.000528541226215484\\
2	0.000135226504394847\\
4	3.42032356260846e-05\\
};
\addplot [color=green, dashed, forget plot]
  table[row sep=crcr]{%
1	0.000505050505050564\\
2	0.000129148908691645\\
4	3.26573266712648e-05\\
};
\addplot [color=black, dashed, forget plot]
  table[row sep=crcr]{%
1	0.000483091787439771\\
2	0.000123472033584482\\
4	3.12139089178043e-05\\
};
\addplot [color=mycolor1, dashed, forget plot]
  table[row sep=crcr]{%
1	0.00046253469010149\\
2	0.000118161408483919\\
4	2.98641182619186e-05\\
};
\addplot [color=blue, dashed, forget plot]
  table[row sep=crcr]{%
1	0.000443262411347733\\
2	0.000113186191284698\\
4	2.86000286000287e-05\\
};
\addplot [color=red, dashed, forget plot]
  table[row sep=crcr]{%
1	0.000425170068027225\\
2	0.000108518719479089\\
3	4.85672656629355e-05\\
};
\addplot [color=green, dashed, forget plot]
  table[row sep=crcr]{%
1	0.000408163265306016\\
2	0.000104134124752742\\
3	4.65983224604005e-05\\
};
\addplot [color=black, dashed, forget plot]
  table[row sep=crcr]{%
1	0.000392156862745075\\
2	0.000100010001000073\\
3	4.47467334884765e-05\\
};
\addplot [color=mycolor1, dashed, forget plot]
  table[row sep=crcr]{%
1	0.000377073906485781\\
2	9.61261174661222e-05\\
3	4.3003354261606e-05\\
};
\end{axis}
\end{tikzpicture}%

%% file: fig2/cf_general_pareto_Lunar_Craters_diameter_gt_1km_n1296785_s1_59_rho0_00_mincha_110522_094958.tex
%
%
\definecolor{mycolor1}{rgb}{1.00000,0.00000,1.00000}%
\begin{tikzpicture}

\begin{axis}[%
width=0.964\fwidth,
height=\fheight,
at={(0\fwidth,0\fheight)},
scale only axis,
unbounded coords=jump,
xmode=log,
xmin=1,
xmax=100,
xminorticks=true,
xlabel style={font=\color{white!15!black}},
xlabel={$a_1$},
ymode=log,
ymin=5e-05,
ymax=1,
yminorticks=true,
ylabel style={font=\color{white!15!black}},
ylabel={$\Pr(\mathbf{A}_2=[a_1,a_2])$},
axis background/.style={fill=white},
xmajorgrids,
xminorgrids,
ymajorgrids,
yminorgrids,
grid style={dotted},legend style={font=\tiny,draw=none},mark size=1pt
]
\addplot [color=blue, draw=none, mark=square, mark options={solid, blue}, forget plot]
  table[row sep=crcr]{%
1	0.076254737678181\\
2	0.0736922465944624\\
3	0.0573572334658405\\
4	0.0431150884687901\\
5	0.0324687592777523\\
6	0.025214665499678\\
7	0.0199632167244378\\
8	0.0161098408757041\\
9	0.013418569770625\\
10	0.0110234156008899\\
11	0.00940402611072768\\
12	0.00804373893899142\\
13	0.00713688082450059\\
14	0.00618761012812455\\
15	0.0053825422101582\\
16	0.00487282008968333\\
17	0.00430603376812656\\
18	0.0038256148860451\\
19	0.00350790609083233\\
20	0.00315626723011139\\
21	0.00288636898175102\\
22	0.00256634677298087\\
23	0.00241443261604661\\
24	0.00224324001280089\\
25	0.00206664944458796\\
26	0.00191242187409632\\
27	0.00172503537594898\\
28	0.00168570734547361\\
29	0.00155230049699835\\
30	0.00146361964396565\\
31	0.00136028717173626\\
32	0.00128240224863798\\
33	0.00126620835373636\\
34	0.00109270233693326\\
35	0.00110272712901522\\
36	0.00106494137424477\\
37	0.000980116210474365\\
38	0.000902231287376088\\
39	0.000863674394753179\\
40	0.00084131139703189\\
41	0.000803525642261439\\
42	0.000747232579031991\\
43	0.000711760237818914\\
44	0.000681685861573044\\
45	0.000680143585868128\\
46	0.00058529363001577\\
47	0.000603029800622308\\
48	0.000578353389343646\\
49	0.000535940807458446\\
50	0.000573726562228897\\
51	0.000563701770146941\\
52	0.000481190019933914\\
53	0.000465767262884749\\
54	0.000458055884360168\\
55	0.00045265791939296\\
56	0.00042721037026184\\
57	0.000424125818852007\\
58	0.000384797788376639\\
59	0.000387111201934014\\
60	0.000387111201934014\\
61	0.00035009658501602\\
62	0.000351638860720937\\
63	0.000339300655081606\\
64	0.000330818138704566\\
65	0.000312310830245569\\
66	0.000302286038163612\\
67	0.00030768400313082\\
68	0.000284549867557074\\
69	0.00026912711050791\\
70	0.000257560042721037\\
71	0.000270669386212826\\
72	0.000247535250639081\\
73	0.000262958007688245\\
74	0.000224401115065335\\
75	0.000232112493589917\\
76	0.000230570217885\\
77	0.000216689736540753\\
78	0.000249848664196455\\
79	0.000200495841639131\\
80	0.000222858839360418\\
81	0.000216689736540753\\
82	0.000218232012245669\\
83.0000000000001	0.000193555600967007\\
84	0.000203580393048963\\
85	0.000193555600967007\\
86	0.000189699911704716\\
87	0.000165023500426054\\
88	0.000168879189688345\\
89	0.000186615360294883\\
90	0.000172734878950636\\
91.0000000000001	0.000171192603245719\\
92	0.000163481224721137\\
93	0.000178903981770301\\
94	0.000158854397606388\\
95.0000000000001	0.000155769846196555\\
96	0.000168879189688345\\
97	0.000138033675590017\\
98.0000000000001	0.000146516191967057\\
99	0.000137262537737559\\
100	0.000151143019081806\\
};
\addplot [color=red, draw=none, mark=square, mark options={solid, red}, forget plot]
  table[row sep=crcr]{%
1	0.0219303893860586\\
2	0.0261485134390049\\
3	0.0211600226714529\\
4	0.015853823108688\\
5	0.0119179355097414\\
6	0.00925288309164587\\
7	0.00724869581310702\\
8	0.00577582251491188\\
9.00000000000001	0.00474018437906052\\
10	0.00400529000566787\\
11	0.00328196270006208\\
12	0.00290256287665265\\
13	0.00244373585444002\\
14	0.00216612622755507\\
15	0.00192707349329303\\
16	0.00169418986185065\\
17	0.00151220132867052\\
18	0.00129782500568714\\
19	0.00115516450298238\\
20	0.00107419502847426\\
21	0.00099939465678582\\
22	0.00093076338791704\\
23	0.000776535817425402\\
24	0.000754172819704114\\
25	0.000724098443458245\\
26	0.000643128968950134\\
27	0.000575268837933813\\
28	0.000565244045851857\\
29	0.000534398531753529\\
30	0.000489672536310954\\
31	0.000463453849327375\\
32	0.000463453849327375\\
33	0.000425668094556923\\
34	0.000400220545425803\\
35	0.000386340064081556\\
36	0.000350867722868479\\
37	0.000331589276557023\\
38	0.000320022208770151\\
39	0.000302286038163612\\
40	0.000289947832524281\\
41	0.000257560042721037\\
42	0.000244450699229247\\
43.0000000000001	0.000244450699229247\\
44	0.000245992974934164\\
45	0.00021360518513092\\
46	0.000235968182852207\\
47.0000000000001	0.000204351530901422\\
48	0.000189699911704716\\
49	0.000205893806606338\\
50	0.000168879189688345\\
51	0.000160396673311304\\
52	0.000187386498147341\\
53	0.000159625535458846\\
54.0000000000001	0.000151143019081806\\
55	0.000157312121901472\\
56	0.000151143019081806\\
57	0.00011104385075398\\
58.0000000000001	0.000142660502704766\\
59.0000000000001	0.000139575951294933\\
60	0.000134949124180184\\
61	0.000120297504983478\\
62.0000000000001	0.000130322297065435\\
63.0000000000001	0.000124153194245769\\
64.0000000000001	8.94519908851505e-05\\
65	0.00011104385075398\\
66.0000000000001	0.000104103610081856\\
67.0000000000001	8.7138577327776e-05\\
68	0.000105645885786773\\
69.0000000000001	0.000109501575049064\\
70	8.09694745081105e-05\\
71	9.33076801474415e-05\\
72.0000000000001	8.55963016228596e-05\\
73	7.01735445736957e-05\\
74	0.000106417023639231\\
75	8.48251637704015e-05\\
76	9.17654044425252e-05\\
77.0000000000001	8.79097151802342e-05\\
78.0000000000001	9.25365422949834e-05\\
79.0000000000001	6.63178553114047e-05\\
80	6.94024067212376e-05\\
81	6.78601310163211e-05\\
82	7.86560609507359e-05\\
nan	nan\\
84.0000000000001	5.47507875245318e-05\\
85	7.09446824261539e-05\\
85.9999999999999	6.70889931638629e-05\\
86.9999999999999	5.70642010819064e-05\\
88.0000000000001	5.62930632294482e-05\\
nan	nan\\
90.0000000000001	6.24621660491137e-05\\
91.0000000000001	6.24621660491137e-05\\
92.0000000000001	6.01487524917392e-05\\
nan	nan\\
95	5.08950982622409e-05\\
96.0000000000001	5.08950982622409e-05\\
};
\addplot [color=green, draw=none, mark=square, mark options={solid, green}, forget plot]
  table[row sep=crcr]{%
1	0.00977571455561253\\
2	0.0131448158330024\\
3	0.0109162274393982\\
4	0.00816249416826998\\
5	0.00608581993160007\\
6	0.00456976291366725\\
7	0.00371148648388129\\
8	0.00306604410137378\\
9.00000000000001	0.00246841226571868\\
10	0.00201112751921097\\
11	0.00176976137139156\\
12	0.0014073265807362\\
13	0.00121299984191674\\
14	0.00110041371545784\\
15	0.000970862556244867\\
16	0.00081971953716306\\
17	0.000818177261458145\\
18	0.00067474562090092\\
19	0.000603029800622307\\
20	0.000511264396179783\\
21	0.000478876606376539\\
22	0.000437235162343796\\
23	0.000411016475360218\\
24	0.000433379473081505\\
25	0.000339300655081605\\
26	0.000337758379376689\\
27	0.000323106760179983\\
28	0.00027760962688495\\
29	0.000268355972655452\\
30	0.000275296213327575\\
31	0.000235968182852207\\
32	0.00022208770150796\\
33	0.000218232012245669\\
34	0.000202809255196505\\
35	0.000193555600967007\\
36	0.000168108051835886\\
37	0.000157312121901472\\
38	0.000150371881229348\\
39	0.000154998708344097\\
40	0.000147287329819515\\
41	0.000128780021360518\\
42	0.0001364913998851\\
43.0000000000001	0.000124924332098227\\
44	0.000118755229278562\\
45	0.000101790196524482\\
46	0.000114128402163813\\
47.0000000000001	0.000101790196524482\\
48	8.3282888065485e-05\\
49	0.000109501575049064\\
50	0.00010256133437694\\
51	8.7138577327776e-05\\
52	9.79345072621908e-05\\
53	8.63674394753179e-05\\
54.0000000000001	8.01983366556523e-05\\
55	6.40044417540301e-05\\
56	8.01983366556523e-05\\
57	7.78849230982777e-05\\
58.0000000000001	6.47755796064884e-05\\
59.0000000000001	5.47507875245318e-05\\
60	5.55219253769901e-05\\
61	5.08950982622409e-05\\
62.0000000000001	6.3233303901572e-05\\
63.0000000000001	7.40292338359868e-05\\
64.0000000000001	5.78353389343646e-05\\
65	6.09198903441974e-05\\
nan	nan\\
67.0000000000001	6.01487524917392e-05\\
68	5.86064767868228e-05\\
nan	nan\\
71	5.62930632294482e-05\\
nan	nan\\
73	5.55219253769901e-05\\
nan	nan\\
75	5.24373739671573e-05\\
nan	nan\\
79.0000000000001	5.24373739671573e-05\\
};
\addplot [color=black, draw=none, mark=square, mark options={solid, black}, forget plot]
  table[row sep=crcr]{%
1	0.00543266617056798\\
2	0.00790647640125387\\
3	0.00656006971086186\\
4	0.00496150094271602\\
5	0.00378011775275007\\
6	0.00296271162914438\\
7	0.00230801559240738\\
8	0.00184147719167017\\
9.00000000000001	0.00139730178865425\\
10	0.00126003925091669\\
11	0.00102407106806448\\
12	0.000864445532605635\\
13	0.000716387064933662\\
14	0.000661636277409131\\
15	0.000596089559950184\\
16	0.000508950982622408\\
17	0.000439548575901171\\
18	0.00046962295214704\\
19	0.000370917307032392\\
20	0.000336216103671773\\
21	0.00030305717601607\\
22	0.000263729145540703\\
23	0.000271440524065284\\
24	0.000247535250639081\\
25	0.000206664944458796\\
26	0.000197411290229298\\
27	0.00020512266875388\\
28	0.000178132843917843\\
29	0.000163481224721137\\
30	0.000120297504983478\\
31	0.000140347089147391\\
32	0.000122610918540853\\
33	0.000118755229278562\\
34	0.000107959299344147\\
35	0.000117212953573646\\
36	0.000109501575049064\\
37	9.63922315572743e-05\\
38	9.02231287376089e-05\\
39	0.000105645885786773\\
40	9.63922315572743e-05\\
41	8.94519908851505e-05\\
42	7.32580959835285e-05\\
43.0000000000001	6.86312688687794e-05\\
44	7.01735445736957e-05\\
45	6.55467174589466e-05\\
46	7.40292338359868e-05\\
47.0000000000001	7.32580959835285e-05\\
48	5.78353389343646e-05\\
49	6.16910281966556e-05\\
50	5.86064767868228e-05\\
nan	nan\\
52	5.62930632294482e-05\\
53	5.55219253769901e-05\\
nan	nan\\
56	7.01735445736957e-05\\
57	5.08950982622409e-05\\
};
\addplot [color=mycolor1, draw=none, mark=square, mark options={solid, mycolor1}, forget plot]
  table[row sep=crcr]{%
1	0.00339917565363572\\
2	0.00531622435484679\\
3	0.00440551055109367\\
4	0.00337218582879969\\
5	0.0025270187425055\\
6	0.00195869014524381\\
7	0.00157080780545734\\
8	0.00124076080460523\\
9.00000000000001	0.000997081243228446\\
10	0.00083668456991714\\
11	0.000700964307884499\\
12	0.000571413148671522\\
13	0.000491985949868328\\
14	0.000417956716032341\\
15	0.000374001858442225\\
16	0.00035472341213077\\
17	0.000305370589573445\\
18	0.000288405556819365\\
19	0.000232883631442374\\
20	0.000224401115065334\\
21	0.000221316563655502\\
22	0.000187386498147341\\
23	0.000164252362573596\\
24	0.000154998708344097\\
25	0.000135720262032642\\
26	0.000130322297065435\\
27	0.000134949124180184\\
28	0.000124153194245769\\
29	0.00010256133437694\\
30	0.00010256133437694\\
31	8.09694745081105e-05\\
32	7.1715820278612e-05\\
33	9.63922315572743e-05\\
34	7.09446824261539e-05\\
35	6.78601310163211e-05\\
36	6.63178553114047e-05\\
37	6.70889931638629e-05\\
38	6.01487524917392e-05\\
39	6.78601310163211e-05\\
40	5.16662361146991e-05\\
41	5.9377614639281e-05\\
42	5.70642010819064e-05\\
43.0000000000001	5.9377614639281e-05\\
};
\addplot [color=blue, draw=none, mark=square, mark options={solid, blue}, forget plot]
  table[row sep=crcr]{%
1	0.00232035379804671\\
2	0.00387265429504505\\
3	0.00319559526058675\\
4	0.00247226795498097\\
5	0.00179212436911284\\
6	0.00140655544288375\\
7	0.0010919311990808\\
8	0.000912256079458043\\
9.00000000000001	0.000700193170032041\\
10	0.000587607043573144\\
11	0.000491985949868328\\
12	0.000431066059524131\\
13	0.000363205928507809\\
14	0.000339300655081605\\
15	0.000276067351180034\\
16	0.000255246629163662\\
17	0.000219003150098127\\
18	0.000191242187409632\\
19	0.000185844222442425\\
20	0.000163481224721137\\
21	0.000130322297065435\\
22	0.000114128402163813\\
23	0.000124153194245769\\
24	0.000114899540016271\\
25	0.000123382056393311\\
26	0.000105645885786773\\
27	9.8705645114649e-05\\
28	7.32580959835285e-05\\
29	8.09694745081105e-05\\
30	7.01735445736957e-05\\
31	6.09198903441974e-05\\
32	6.09198903441974e-05\\
33	5.62930632294482e-05\\
34	5.9377614639281e-05\\
35	5.01239604097827e-05\\
36	5.08950982622409e-05\\
37	5.47507875245318e-05\\
38	5.01239604097827e-05\\
};
\addplot [color=red, draw=none, mark=square, mark options={solid, red}, forget plot]
  table[row sep=crcr]{%
1	0.00176436340642435\\
2	0.00288482670604611\\
3	0.00239823872114498\\
4	0.00189854139275207\\
5	0.00138573472086738\\
6	0.00105260316860543\\
7	0.000795043125884399\\
8	0.000645442382507509\\
9.00000000000001	0.000531313980343696\\
10	0.000451886781540502\\
11	0.000359350239245519\\
12	0.000314624243802943\\
13	0.000290718970376739\\
14	0.000263729145540703\\
15	0.000203580393048963\\
16	0.000163481224721137\\
17	0.000160396673311304\\
18	0.000165794638278512\\
19	0.000125695469950686\\
20	0.000113357264311354\\
21	0.000143431640557224\\
22	9.17654044425252e-05\\
23	8.3282888065485e-05\\
24	9.25365422949834e-05\\
25	7.63426473933612e-05\\
26	6.24621660491137e-05\\
27	6.78601310163211e-05\\
28	7.55715095409031e-05\\
29	5.47507875245318e-05\\
30	6.01487524917392e-05\\
31	5.78353389343646e-05\\
};
\addplot [color=green, draw=none, mark=square, mark options={solid, green}, forget plot]
  table[row sep=crcr]{%
1	0.00127623314581831\\
2	0.00222396156648943\\
3	0.00191627756335861\\
4	0.00144742574906403\\
5	0.00106108568498248\\
6	0.000800441090851606\\
7	0.000659322863851757\\
8	0.000519746912556823\\
9.00000000000001	0.000438777438048712\\
10	0.000329275862999649\\
11	0.000274525075475117\\
12	0.000239823872114498\\
13	0.000220545425803043\\
14	0.000189699911704716\\
15	0.000157312121901472\\
16	0.00014111822699985\\
17	0.000120297504983478\\
18	0.000103332472229398\\
19	8.63674394753179e-05\\
20	9.63922315572743e-05\\
21	9.17654044425252e-05\\
22	8.63674394753179e-05\\
23	7.32580959835285e-05\\
24	6.01487524917392e-05\\
25	6.47755796064884e-05\\
26	5.39796496720736e-05\\
27	5.86064767868228e-05\\
};
\addplot [color=black, draw=none, mark=square, mark options={solid, black}, forget plot]
  table[row sep=crcr]{%
1	0.00106031454713002\\
2	0.00181448736683413\\
3	0.00152839522357214\\
4	0.00118446774137579\\
5	0.000870614635425302\\
6	0.000680143585868127\\
7	0.000515120085442074\\
8	0.000400220545425803\\
9.00000000000001	0.000342385206491439\\
10	0.000283007591852157\\
11	0.000248306388491539\\
12	0.000197411290229298\\
13	0.00016656577613097\\
14	0.000148829605524432\\
15	0.000132635710622809\\
16	0.000117984091426104\\
17	0.000112586126458896\\
18	0.000116441815721187\\
19	8.17406123605687e-05\\
20	9.25365422949834e-05\\
21	6.3233303901572e-05\\
22	7.78849230982777e-05\\
23	5.16662361146991e-05\\
24	6.16910281966556e-05\\
nan	nan\\
26	5.55219253769901e-05\\
};
\addplot [color=mycolor1, draw=none, mark=square, mark options={solid, mycolor1}, forget plot]
  table[row sep=crcr]{%
1	0.000861360981195804\\
2	0.00139730178865425\\
3	0.00126852176729373\\
4	0.00095621093704816\\
5	0.000724869581310702\\
6	0.000553676978064983\\
7	0.000421812405294633\\
8	0.000350867722868479\\
9.00000000000001	0.000269898248360368\\
10	0.000238281596409582\\
11	0.000198182428081756\\
12	0.000171963741098177\\
13	0.000126466607803144\\
14	0.000129551159212977\\
15	0.000104874747934315\\
16	0.000108730437196605\\
17	9.17654044425252e-05\\
18	8.48251637704015e-05\\
19	6.3233303901572e-05\\
20	7.24869581310704e-05\\
nan	nan\\
22	5.24373739671573e-05\\
};
\addplot [color=blue, draw=none, mark=square, mark options={solid, blue}, forget plot]
  table[row sep=crcr]{%
1	0.000672432207343546\\
2	0.00122765146111345\\
3	0.00100093693249074\\
4	0.00081971953716306\\
5	0.00064467124465505\\
6	0.000430294921671672\\
7	0.000361663652802893\\
8	0.000279923040442325\\
9.00000000000001	0.000221316563655502\\
10	0.000191242187409632\\
11	0.000159625535458846\\
12	0.0001364913998851\\
13	0.000131864572770351\\
14	0.00012800888350806\\
15	9.71633694097326e-05\\
16	0.000101019058672024\\
17	7.09446824261539e-05\\
18	5.24373739671573e-05\\
19	5.32085118196155e-05\\
nan	nan\\
21	5.24373739671573e-05\\
nan	nan\\
23	5.01239604097827e-05\\
};
\addplot [color=red, draw=none, mark=square, mark options={solid, red}, forget plot]
  table[row sep=crcr]{%
1	0.000561388356589565\\
2	0.000997852381080904\\
3	0.000876783738244967\\
4	0.000660094001704215\\
5	0.000497383914835535\\
6	0.000361663652802893\\
7	0.000286863281114448\\
8	0.000201266979491589\\
9.00000000000001	0.000210520633721087\\
10	0.000137262537737559\\
11	0.000125695469950686\\
12	0.000121068642835937\\
13	9.8705645114649e-05\\
14	8.01983366556523e-05\\
15	7.63426473933612e-05\\
16	5.39796496720736e-05\\
17	5.55219253769901e-05\\
18	5.62930632294482e-05\\
19	5.01239604097827e-05\\
};
\addplot [color=green, draw=none, mark=square, mark options={solid, green}, forget plot]
  table[row sep=crcr]{%
1	0.000488130260606037\\
2	0.00092228087154\\
3	0.000753401681851655\\
4	0.000567557459409231\\
5	0.000427981508114298\\
6	0.000323877898032442\\
7	0.000246764112786622\\
8	0.000193555600967007\\
9.00000000000001	0.000155769846196555\\
10	0.000151143019081806\\
11	0.00011952636713102\\
12	0.000125695469950686\\
13	8.63674394753179e-05\\
14	7.40292338359868e-05\\
15	5.62930632294482e-05\\
16	5.08950982622409e-05\\
17	5.70642010819064e-05\\
};
\addplot [color=black, draw=none, mark=square, mark options={solid, black}, forget plot]
  table[row sep=crcr]{%
1	0.00046962295214704\\
2	0.000742605751917241\\
3	0.000636959866130468\\
4	0.000447259954425752\\
5	0.000372459582737308\\
6	0.00030768400313082\\
7	0.000215918598688295\\
8	0.000174277154655552\\
9.00000000000001	0.000155769846196555\\
10	0.000113357264311354\\
11	0.000101790196524482\\
12	7.24869581310704e-05\\
13	6.01487524917392e-05\\
14	6.78601310163211e-05\\
15	6.86312688687794e-05\\
nan	nan\\
17	5.47507875245318e-05\\
};
\addplot [color=mycolor1, draw=none, mark=square, mark options={solid, mycolor1}, forget plot]
  table[row sep=crcr]{%
1	0.000350867722868479\\
2	0.000671661069491087\\
3	0.000593776146392809\\
4	0.000404847372540552\\
5	0.000315395381655402\\
6	0.000245221837081706\\
7	0.000187386498147341\\
8	0.000161938949016221\\
9.00000000000001	0.00012800888350806\\
10	9.71633694097326e-05\\
11	9.25365422949834e-05\\
12	8.48251637704015e-05\\
13	6.78601310163211e-05\\
14	5.47507875245318e-05\\
15	5.16662361146991e-05\\
};
\addplot [color=blue, draw=none, mark=square, mark options={solid, blue}, forget plot]
  table[row sep=crcr]{%
1	0.000349325447163562\\
2	0.000608427765589516\\
3	0.000460369297917542\\
4	0.000355494549983228\\
5	0.000271440524065284\\
6	0.000215147460835836\\
7	0.000165794638278512\\
8	0.000138804813442475\\
9.00000000000001	0.000114899540016271\\
10	0.00010256133437694\\
11	8.09694745081105e-05\\
12	6.63178553114047e-05\\
13	5.08950982622409e-05\\
};
\addplot [color=red, draw=none, mark=square, mark options={solid, red}, forget plot]
  table[row sep=crcr]{%
1	0.000289176694671823\\
2	0.000506637569065034\\
3	0.00044417540301592\\
4	0.000366290479917642\\
5	0.000254475491311204\\
6	0.000189699911704716\\
7	0.000152685294786722\\
8	0.000116441815721187\\
9.00000000000001	8.01983366556523e-05\\
10	6.94024067212376e-05\\
11	6.78601310163211e-05\\
12	5.32085118196155e-05\\
};
\addplot [color=green, draw=none, mark=square, mark options={solid, green}, forget plot]
  table[row sep=crcr]{%
1	0.000241366147819415\\
2	0.000460369297917542\\
3	0.000370917307032392\\
4	0.000289176694671823\\
5	0.000219003150098127\\
6	0.00018353080888505\\
7	0.000109501575049064\\
8	9.79345072621908e-05\\
9.00000000000001	7.86560609507359e-05\\
10	6.78601310163211e-05\\
11	5.32085118196155e-05\\
12	5.55219253769901e-05\\
};
\addplot [color=black, draw=none, mark=square, mark options={solid, black}, forget plot]
  table[row sep=crcr]{%
1	0.000242908423524331\\
2	0.000427210370261839\\
3	0.000347012033606188\\
4	0.000256788904868579\\
5	0.000201266979491589\\
6	0.00017504829250801\\
7	0.000127237745655602\\
8	8.79097151802342e-05\\
9.00000000000001	8.94519908851505e-05\\
10	8.17406123605687e-05\\
11	5.01239604097827e-05\\
};
\addplot [color=mycolor1, draw=none, mark=square, mark options={solid, mycolor1}, forget plot]
  table[row sep=crcr]{%
1	0.000226714528622709\\
2	0.000407160786097927\\
3	0.000329275862999649\\
4	0.000232883631442374\\
5	0.000195097876671923\\
6	0.000121839780688395\\
7	0.00010256133437694\\
8	0.000109501575049064\\
9.00000000000001	8.09694745081105e-05\\
10	6.16910281966556e-05\\
};
\addplot [color=blue, draw=none, mark=square, mark options={solid, blue}, forget plot]
  table[row sep=crcr]{%
1	0.000180446257475217\\
2	0.000355494549983228\\
3	0.000284549867557074\\
4	0.000237510458557124\\
5	0.000169650327540803\\
6	0.000128780021360518\\
7	0.000104103610081856\\
8	7.78849230982777e-05\\
9.00000000000001	6.40044417540301e-05\\
nan	nan\\
11	5.86064767868228e-05\\
};
\addplot [color=red, draw=none, mark=square, mark options={solid, red}, forget plot]
  table[row sep=crcr]{%
1	0.000170421465393261\\
2	0.000320022208770151\\
3	0.000293803521786572\\
4	0.00021360518513092\\
5	0.000161938949016221\\
6	0.000107959299344147\\
7	8.63674394753179e-05\\
8	7.09446824261539e-05\\
nan	nan\\
11	5.08950982622409e-05\\
};
\addplot [color=green, draw=none, mark=square, mark options={solid, green}, forget plot]
  table[row sep=crcr]{%
1	0.000139575951294933\\
2	0.000289947832524281\\
3	0.000237510458557124\\
4	0.000186615360294883\\
5	0.000138804813442475\\
6	0.000105645885786773\\
7	8.17406123605687e-05\\
8	7.63426473933612e-05\\
9.00000000000001	5.55219253769901e-05\\
};
\addplot [color=black, draw=none, mark=square, mark options={solid, black}, forget plot]
  table[row sep=crcr]{%
1	0.000137262537737559\\
2	0.000269898248360368\\
3	0.000234425907147291\\
4	0.000171192603245719\\
5	0.000142660502704766\\
6	0.00011104385075398\\
7	8.17406123605687e-05\\
8	5.78353389343646e-05\\
9.00000000000001	5.70642010819064e-05\\
};
\addplot [color=mycolor1, draw=none, mark=square, mark options={solid, mycolor1}, forget plot]
  table[row sep=crcr]{%
1	0.000142660502704766\\
2	0.000215918598688295\\
3	0.000197411290229298\\
4	0.000171192603245719\\
5	0.000105645885786773\\
6	8.86808530326924e-05\\
7	8.09694745081105e-05\\
8	5.62930632294482e-05\\
};
\addplot [color=blue, draw=none, mark=square, mark options={solid, blue}, forget plot]
  table[row sep=crcr]{%
1	0.000109501575049064\\
2	0.000219003150098127\\
3	0.000212062909426004\\
4	0.000141889364852308\\
5	9.63922315572743e-05\\
6	8.25117502130268e-05\\
7	7.24869581310704e-05\\
8	5.47507875245318e-05\\
};
\addplot [color=red, draw=none, mark=square, mark options={solid, red}, forget plot]
  table[row sep=crcr]{%
1	0.000109501575049064\\
2	0.000212834047278461\\
3	0.000181988533180134\\
4	0.000161938949016221\\
5	9.4849955852358e-05\\
6	8.01983366556523e-05\\
7	5.47507875245318e-05\\
8	6.01487524917392e-05\\
};
\addplot [color=green, draw=none, mark=square, mark options={solid, green}, forget plot]
  table[row sep=crcr]{%
1	0.000106417023639231\\
2	0.000198953565934214\\
3	0.000154227570491639\\
4	0.000130322297065435\\
5	0.000116441815721187\\
6	6.40044417540301e-05\\
7	5.86064767868228e-05\\
};
\addplot [color=black, draw=none, mark=square, mark options={solid, black}, forget plot]
  table[row sep=crcr]{%
1	8.40540259179432e-05\\
2	0.000187386498147341\\
3	0.000159625535458846\\
4	0.000116441815721187\\
5	0.000104874747934315\\
6	8.63674394753179e-05\\
};
\addplot [color=mycolor1, draw=none, mark=square, mark options={solid, mycolor1}, forget plot]
  table[row sep=crcr]{%
1	0.000106417023639231\\
2	0.000162710086868679\\
3	0.00014111822699985\\
4	0.000115670677868729\\
5	7.86560609507359e-05\\
6	6.47755796064884e-05\\
7	5.24373739671573e-05\\
};
\addplot [color=blue, draw=none, mark=square, mark options={solid, blue}, forget plot]
  table[row sep=crcr]{%
1	8.17406123605687e-05\\
2	0.000165023500426054\\
3	0.000145745054114599\\
4	0.000123382056393311\\
5	8.3282888065485e-05\\
6	5.01239604097827e-05\\
};
\addplot [color=red, draw=none, mark=square, mark options={solid, red}, forget plot]
  table[row sep=crcr]{%
1	8.55963016228596e-05\\
2	0.000158854397606388\\
3	0.000124153194245769\\
4	9.63922315572743e-05\\
5	7.86560609507359e-05\\
6	6.16910281966556e-05\\
};
\addplot [color=green, draw=none, mark=square, mark options={solid, green}, forget plot]
  table[row sep=crcr]{%
1	7.1715820278612e-05\\
2	0.000146516191967057\\
3	0.000121068642835937\\
4	0.000104103610081856\\
5	7.78849230982777e-05\\
6	5.70642010819064e-05\\
};
\addplot [color=black, draw=none, mark=square, mark options={solid, black}, forget plot]
  table[row sep=crcr]{%
1	7.78849230982777e-05\\
2	0.000140347089147391\\
3	0.000112586126458896\\
4	9.79345072621908e-05\\
5	6.70889931638629e-05\\
};
\addplot [color=mycolor1, draw=none, mark=square, mark options={solid, mycolor1}, forget plot]
  table[row sep=crcr]{%
1	6.16910281966556e-05\\
2	0.0001364913998851\\
3	0.000106417023639231\\
4	8.3282888065485e-05\\
5	6.63178553114047e-05\\
};
\addplot [color=blue, draw=none, mark=square, mark options={solid, blue}, forget plot]
  table[row sep=crcr]{%
1	5.86064767868228e-05\\
2	0.000125695469950686\\
3	0.000118755229278562\\
4	8.63674394753179e-05\\
5	7.01735445736957e-05\\
};
\addplot [color=red, draw=none, mark=square, mark options={solid, red}, forget plot]
  table[row sep=crcr]{%
1	6.63178553114047e-05\\
2	0.000123382056393311\\
3	0.000101019058672024\\
4	6.86312688687794e-05\\
};
\addplot [color=green, draw=none, mark=square, mark options={solid, green}, forget plot]
  table[row sep=crcr]{%
1	5.9377614639281e-05\\
2	0.000114899540016271\\
3	8.40540259179432e-05\\
4	6.94024067212376e-05\\
};
\addplot [color=black, draw=none, mark=square, mark options={solid, black}, forget plot]
  table[row sep=crcr]{%
1	5.39796496720736e-05\\
2	0.000115670677868729\\
3	8.79097151802342e-05\\
4	8.09694745081105e-05\\
5	5.08950982622409e-05\\
};
\addplot [color=mycolor1, draw=none, mark=square, mark options={solid, mycolor1}, forget plot]
  table[row sep=crcr]{%
1	5.62930632294482e-05\\
2	0.000106417023639231\\
3	8.94519908851505e-05\\
4	6.01487524917392e-05\\
};
\addplot [color=blue, draw=none, mark=square, mark options={solid, blue}, forget plot]
  table[row sep=crcr]{%
1	5.47507875245318e-05\\
2	9.94767829671071e-05\\
3	7.78849230982777e-05\\
4	6.86312688687794e-05\\
};
\addplot [color=red, draw=none, mark=square, mark options={solid, red}, forget plot]
  table[row sep=crcr]{%
2	9.33076801474415e-05\\
3	7.1715820278612e-05\\
4	5.62930632294482e-05\\
};
\addplot [color=green, draw=none, mark=square, mark options={solid, green}, forget plot]
  table[row sep=crcr]{%
2	9.4849955852358e-05\\
3	6.09198903441974e-05\\
4	5.24373739671573e-05\\
};
\addplot [color=black, draw=none, mark=square, mark options={solid, black}, forget plot]
  table[row sep=crcr]{%
2	6.01487524917392e-05\\
3	6.16910281966556e-05\\
4	5.01239604097827e-05\\
};
\addplot [color=mycolor1, draw=none, mark=square, mark options={solid, mycolor1}, forget plot]
  table[row sep=crcr]{%
1	5.39796496720736e-05\\
2	9.17654044425252e-05\\
3	6.01487524917392e-05\\
};
\addplot [color=blue, draw=none, mark=square, mark options={solid, blue}, forget plot]
  table[row sep=crcr]{%
2	7.09446824261539e-05\\
3	5.62930632294482e-05\\
};
\addplot [color=red, draw=none, mark=square, mark options={solid, red}, forget plot]
  table[row sep=crcr]{%
2	6.24621660491137e-05\\
3	6.55467174589466e-05\\
};
\addplot [color=green, draw=none, mark=square, mark options={solid, green}, forget plot]
  table[row sep=crcr]{%
2	6.09198903441974e-05\\
3	5.55219253769901e-05\\
4	5.70642010819064e-05\\
};
\addplot [color=black, draw=none, mark=square, mark options={solid, black}, forget plot]
  table[row sep=crcr]{%
2	6.70889931638629e-05\\
3	5.01239604097827e-05\\
};
\addplot [color=mycolor1, draw=none, mark=square, mark options={solid, mycolor1}, forget plot]
  table[row sep=crcr]{%
2	6.55467174589466e-05\\
3	5.01239604097827e-05\\
};
\addplot [color=blue, dashed, forget plot]
  table[row sep=crcr]{%
1	0.0751419115715888\\
2	0.0655949060171695\\
3	0.0503702005116061\\
4	0.0385632759454372\\
5	0.0300892524271979\\
6	0.0239917629799622\\
7	0.0195173172100009\\
8	0.0161592217304967\\
9	0.0135842892881073\\
10	0.0115711527310951\\
11	0.00996981910638883\\
13	0.00761750719362459\\
15	0.00600512461820213\\
17	0.0048534175117185\\
20	0.00365885318401957\\
24	0.00264608033925351\\
29	0.00187744734521801\\
35	0.00132725204037999\\
43	0.000902978884618003\\
54	0.000586359734480887\\
70	0.000356488156740449\\
95.0000000000001	0.000197281821319853\\
100	0.000178524788008373\\
};
\node[right, align=left, font=\color{blue}]
at (axis cs:1.02,0.09) {{\tiny $a_2=1$}};
\addplot [color=red, dashed, forget plot]
  table[row sep=crcr]{%
1	0.0234977666851366\\
2	0.0235658964858113\\
3	0.0183745085477111\\
4	0.0140162558049039\\
5	0.0108599783990348\\
6	0.00859852548935973\\
7	0.00695083364676152\\
8	0.00572329101589114\\
9.00000000000001	0.00478843185549948\\
10	0.00406200026284476\\
11	0.00348729731527869\\
13	0.00264874146187077\\
15	0.00207836847925103\\
17	0.00167346724459899\\
20	0.00125600881300442\\
24	0.000904375156147774\\
29	0.000639151607854173\\
36	0.000427150744617798\\
45	0.000280158774324052\\
58.0000000000001	0.000172406312820235\\
78.0000000000001	9.72194024365936e-05\\
100	5.98964512813285e-05\\
};
\node[right, align=left, font=\color{red}]
at (axis cs:1.02,0.03) {{\tiny $a_2=2$}};
\addplot [color=green, dashed, forget plot]
  table[row sep=crcr]{%
1	0.0110174021686988\\
2	0.0120651877901196\\
3	0.00950296682941698\\
4	0.00724466980796346\\
5	0.00559902188435135\\
6	0.00442110808155233\\
7	0.00356510738204518\\
8	0.00292917304288764\\
9.00000000000001	0.0024461459519623\\
10	0.0020717049375415\\
11	0.00177610020069058\\
13	0.00134590986226965\\
15	0.00105417648229267\\
18	0.00076625364216959\\
21	0.000581714047583476\\
25	0.000423629818267823\\
30	0.0003024640455567\\
37	0.000204150792059016\\
47.0000000000001	0.000129605305418765\\
61	7.85362553832671e-05\\
78.0000000000001	4.87613700692752e-05\\
};
\node[right, align=left, font=\color{green}]
at (axis cs:1.02,0.015) {{\tiny $a_2=3$}};
\addplot [color=black, dashed, forget plot]
  table[row sep=crcr]{%
1	0.00630068447121745\\
2	0.00732064851589653\\
3	0.00580484971540684\\
4	0.00442500984364801\\
5	0.00341534499197771\\
6	0.00269291306626451\\
7	0.00216862491011438\\
8	0.00177971075451699\\
9.00000000000001	0.0014847307487672\\
10	0.00125635805646662\\
11	0.00107627446289665\\
13	0.00081457223535089\\
15	0.000637386184910604\\
18	0.000462781610553631\\
21	0.000351034040421028\\
25	0.000255425053165531\\
30	0.000182231086713905\\
37	0.000122907902398328\\
47.0000000000001	7.79743325215188e-05\\
60	4.87540854745214e-05\\
};
\node[right, align=left]
at (axis cs:1.02,0.008) {{\tiny $a_2=4$}};
\addplot [color=mycolor1, dashed, forget plot]
  table[row sep=crcr]{%
1	0.00405219139007595\\
2	0.00491183016325702\\
3	0.00391326452484675\\
4	0.00298318509560102\\
5	0.00230064098186766\\
6	0.00181234608332476\\
7	0.00145827203033238\\
8	0.00119586928592508\\
9.00000000000001	0.000997023069902346\\
10	0.000843201842027408\\
11	0.000721993448112246\\
13	0.000546007423257056\\
15	0.000426977106691095\\
18	0.000309793704543787\\
21	0.000234864358412468\\
25	0.000170806204249585\\
31	0.000114560686450618\\
39	7.43173611762738e-05\\
49	4.80738837293077e-05\\
};
\node[right, align=left, font=\color{mycolor1}]
at (axis cs:1.02,0.005) {{\tiny $a_2=5$}};
\addplot [color=blue, dashed, forget plot]
  table[row sep=crcr]{%
1	0.002815524876295\\
2	0.0035229053835224\\
3	0.0028165490309121\\
4	0.00214728512472473\\
5	0.00165508560041638\\
6	0.00130299048672122\\
7	0.00104782192749303\\
8	0.000858839850069514\\
9.00000000000001	0.00071571951106747\\
10	0.000605067976382822\\
11	0.000517919888029516\\
13	0.00039146482134117\\
15	0.000305995406037984\\
18	0.00022190784276338\\
21	0.000168174402490534\\
25	0.000122261524017397\\
31	8.19703800688936e-05\\
39	5.31576993251931e-05\\
41	4.8335808967237e-05\\
};
\addplot [color=red, dashed, forget plot]
  table[row sep=crcr]{%
1	0.00206605156481476\\
2	0.00264970021988742\\
3	0.00212413580291004\\
4	0.00161952600842175\\
5	0.00124780573844252\\
6	0.000981905452569931\\
7	0.000789282028771043\\
8	0.00064668916086552\\
9.00000000000001	0.000538749131184457\\
10	0.000455330892696004\\
11	0.000389655384867904\\
13	0.000294400681181084\\
15	0.000230052256327662\\
18	0.000166774837188244\\
21	0.000126358020436507\\
25	9.1837128737532e-05\\
31	6.15552435502403e-05\\
35	4.89885885474051e-05\\
};
\node[right, align=left]
at (axis cs:1.02,0.003) {\tiny $\phantom{a_2}~~~~\!\pmb{\pmb{\vdots}}$};
\addplot [color=green, dashed, forget plot]
  table[row sep=crcr]{%
1	0.00157882971556622\\
2	0.00206515628501458\\
3	0.00165907085128114\\
4	0.00126503334574901\\
5	0.000974386517896053\\
6	0.00076648379374028\\
7	0.00061592134938998\\
8	0.000504504882769263\\
9.00000000000001	0.000420193825574621\\
10	0.000355056876996218\\
11	0.000303788546447137\\
13	0.000229455323097538\\
15	0.000179259894761558\\
18	0.000129918147178859\\
21	9.84134310182405e-05\\
25	7.15125960127072e-05\\
31	4.79222789142448e-05\\
};
\addplot [color=black, dashed, forget plot]
  table[row sep=crcr]{%
1	0.00124483380455176\\
2	0.00165472973114211\\
3	0.00133165292953652\\
4	0.00101544544637666\\
5	0.000781959693792355\\
6	0.00061494624264613\\
7	0.000494024850992036\\
8	0.000404568115637271\\
9.00000000000001	0.000336892803537769\\
10	0.000284621157421616\\
11	0.000243487932928998\\
13	0.000183865616603829\\
15	0.000143616559719034\\
18	0.000104063553134276\\
21	7.88159335937678e-05\\
25	5.72629025313576e-05\\
27	4.96558450260587e-05\\
};
\addplot [color=mycolor1, dashed, forget plot]
  table[row sep=crcr]{%
1	0.00100616992563629\\
2	0.00135554261232817\\
3	0.00109244493420169\\
4	0.000833087513658071\\
5	0.000641411582241624\\
6	0.00050430522313083\\
7	0.000405056484887312\\
8	0.000331649789703347\\
9.00000000000001	0.00027612876052846\\
10	0.000233253490342644\\
11	0.000199520377280326\\
13	0.000150635191918991\\
15	0.000117642676423113\\
18	8.52283636514905e-05\\
21	6.45421698957927e-05\\
25	4.68864357044407e-05\\
};
\addplot [color=blue, dashed, forget plot]
  table[row sep=crcr]{%
1	0.000829842340530739\\
2	0.00113074056609796\\
3	0.000912374115759505\\
4	0.00069580377836309\\
5	0.000535631287732929\\
6	0.000421059336880525\\
7	0.000338136064730931\\
8	0.000276815625596949\\
9.00000000000001	0.000230444425157961\\
10	0.000194640864072559\\
11	0.000166475708320457\\
13	0.000125666860185974\\
15	9.81307416675625e-05\\
18	7.10824685164479e-05\\
21	5.38239672518981e-05\\
22	4.94477542878624e-05\\
};
\addplot [color=red, dashed, forget plot]
  table[row sep=crcr]{%
1	0.000695953083759804\\
2	0.000957559879832147\\
3	0.000773432621973935\\
4	0.000589870411514925\\
5	0.000454024973411071\\
6	0.000356853987127703\\
7	0.000286534322426313\\
8	0.000234542272809168\\
9.00000000000001	0.00019523129071014\\
10	0.000164883173858654\\
11	0.000141012523963202\\
13	0.000106431300967019\\
15	8.31013767897126e-05\\
18	6.01885261552749e-05\\
20	4.9790285015755e-05\\
};
\addplot [color=green, dashed, forget plot]
  table[row sep=crcr]{%
1	0.000591935989370574\\
2	0.000821323540331153\\
3	0.000663981765691158\\
4	0.000506417181800112\\
5	0.000389748120325328\\
6	0.000306293788374823\\
7	0.000245907285042646\\
8	0.000201265497222769\\
9.00000000000001	0.000167516395451317\\
10	0.000141465123248918\\
11	0.00012097635212649\\
13	9.1298276959445e-05\\
15	7.12791804007296e-05\\
18	5.16207077193853e-05\\
19	4.68449049645659e-05\\
};
\addplot [color=black, dashed, forget plot]
  table[row sep=crcr]{%
1	0.000509544203582845\\
2	0.000712223372064001\\
3	0.000576227509721306\\
4	0.000439503722534429\\
5	0.000338218469873031\\
6	0.00026576798966077\\
7	0.000213348837115213\\
8	0.000174601546503561\\
9.00000000000001	0.000145311946531113\\
10	0.000122705303543358\\
11	0.000104927266562804\\
13	7.91785788696597e-05\\
15	6.18122322046111e-05\\
17	4.95670620316302e-05\\
};
\addplot [color=mycolor1, dashed, forget plot]
  table[row sep=crcr]{%
1	0.000443186701996501\\
2	0.000623502837613032\\
3	0.000504790961457502\\
4	0.000385030218719717\\
5	0.000296274394354978\\
6	0.000232786064857604\\
7	0.000186855064213675\\
8	0.000152907124432389\\
9.00000000000001	0.000127247916143935\\
10	0.000107445075539947\\
11	9.18731959393405e-05\\
13	6.93219747652795e-05\\
15	5.41138665110796e-05\\
16	4.83120296546624e-05\\
};
\addplot [color=blue, dashed, forget plot]
  table[row sep=crcr]{%
1	0.000388965932122343\\
2	0.00055038463464574\\
3	0.000445862757999848\\
4	0.000340093020211984\\
5	0.000261677234649871\\
6	0.000205585024215915\\
7	0.000165007797117867\\
8	0.000135019528995028\\
9.00000000000001	0.000112355137170894\\
10	9.48649456367078e-05\\
11	8.11125592058427e-05\\
13	6.11980344523145e-05\\
15	4.77693742415019e-05\\
};
\addplot [color=red, dashed, forget plot]
  table[row sep=crcr]{%
1	0.000344098497880529\\
2	0.000489413569096042\\
3	0.000396683814479215\\
4	0.000302588944131932\\
5	0.000232805809425529\\
6	0.000182888457861685\\
7	0.000146780503356703\\
8	0.000120097308745959\\
9.00000000000001	9.99323373445105e-05\\
10	8.4372028363977e-05\\
11	7.21378416257967e-05\\
13	5.44231303296894e-05\\
14	4.79016298283174e-05\\
};
\addplot [color=green, dashed, forget plot]
  table[row sep=crcr]{%
1	0.000306554246189639\\
2	0.000438040580220567\\
3	0.000355216009007897\\
4	0.000270964311080868\\
5	0.000208462722311342\\
6	0.000163753812880529\\
7	0.000131415296626691\\
8	0.000107519327963712\\
9.00000000000001	8.94619334585455e-05\\
10	7.55287924644189e-05\\
11	6.45745542321441e-05\\
13	4.87142394108298e-05\\
};
\addplot [color=black, dashed, forget plot]
  table[row sep=crcr]{%
1	0.000274824579833188\\
2	0.000394351828557626\\
3	0.000319927347688709\\
4	0.000244051218305286\\
5	0.000187748031710511\\
6	0.000147472814587879\\
7	0.000118342766574575\\
8	9.68189893242322e-05\\
9.00000000000001	8.05551676335558e-05\\
10	6.80066443258242e-05\\
11	5.81414767937827e-05\\
13	4.38588613700663e-05\\
};
\addplot [color=mycolor1, dashed, forget plot]
  table[row sep=crcr]{%
1	0.000247769751238636\\
2	0.000356887782648419\\
3	0.00028964828068001\\
4	0.00022095803339961\\
5	0.000169974776404014\\
6	0.000133504898051806\\
7	0.000107128398141254\\
8	8.76402820374489e-05\\
9.00000000000001	7.29154546752487e-05\\
10	6.1554923612317e-05\\
11	5.26241054554788e-05\\
12	4.54855539689936e-05\\
};
\addplot [color=blue, dashed, forget plot]
  table[row sep=crcr]{%
1	0.000224515829359451\\
2	0.000324519687078806\\
3	0.000263473455761276\\
4	0.00020099451463218\\
5	0.000154611224462253\\
6	0.00012143170922696\\
7	9.74359577885498e-05\\
8	7.97077535884538e-05\\
9.00000000000001	6.63133448660647e-05\\
10	5.59797274050644e-05\\
11	4.78565198446245e-05\\
};
\addplot [color=red, dashed, forget plot]
  table[row sep=crcr]{%
1	0.000204383712636425\\
2	0.000296363515616682\\
3	0.000240693198064523\\
4	0.000183619585478634\\
5	0.000141240601686965\\
6	0.000110925383986715\\
7	8.90019720757862e-05\\
8	7.28055823268625e-05\\
9.00000000000001	6.05690775431488e-05\\
10	5.11291563464189e-05\\
11	4.37087507513812e-05\\
};
\addplot [color=green, dashed, forget plot]
  table[row sep=crcr]{%
1	0.000186839292769747\\
2	0.000271718951192601\\
3	0.000220744863033251\\
4	0.000168404254975218\\
5	0.000129532503872709\\
6	0.000101726034547167\\
7	8.16176173569772e-05\\
8	6.6762729442893e-05\\
9.00000000000001	5.55402045083752e-05\\
10	4.68828493886964e-05\\
};
\addplot [color=black, dashed, forget plot]
  table[row sep=crcr]{%
1	0.000171457867838121\\
2	0.000250025406365872\\
3	0.000203177737472359\\
4	0.000155004863721553\\
5	0.00011922227304928\\
6	9.36255105448147e-05\\
7	7.51156590179658e-05\\
8	6.14422324018982e-05\\
9.00000000000001	5.11126666453143e-05\\
10	4.31444392614e-05\\
};
\addplot [color=mycolor1, dashed, forget plot]
  table[row sep=crcr]{%
1	0.000157898335359286\\
2	0.000230829854816816\\
3	0.000187627341639574\\
4	0.000143143495103248\\
5	0.000110095899194751\\
6	8.64555126950384e-05\\
7	6.93608928432621e-05\\
8	5.67333702614082e-05\\
9.00000000000001	4.71942648334056e-05\\
};
\addplot [color=blue, dashed, forget plot]
  table[row sep=crcr]{%
1	0.00014588421224293\\
2	0.000213762983009458\\
3	0.000173796360745437\\
4	0.000132593439441497\\
5	0.000101978817946269\\
6	8.00787748997944e-05\\
7	6.42430548863926e-05\\
8	5.25458529667678e-05\\
9.00000000000001	4.37098180502415e-05\\
};
\addplot [color=red, dashed, forget plot]
  table[row sep=crcr]{%
1	0.000135189494485988\\
2	0.000198521290460496\\
3	0.000161440320833206\\
4	0.000123168274416385\\
5	9.47274962008748e-05\\
6	7.43824417262507e-05\\
7	5.96714972104784e-05\\
8	4.88054583360209e-05\\
};
\addplot [color=green, dashed, forget plot]
  table[row sep=crcr]{%
1	0.000125627996862443\\
2	0.000184853501315918\\
3	0.000150356706290378\\
4	0.000114713568864055\\
5	8.82230399008589e-05\\
6	6.92730380526529e-05\\
7	5.55711441310169e-05\\
8	4.54507186941067e-05\\
};
\addplot [color=black, dashed, forget plot]
  table[row sep=crcr]{%
1	0.000117045227769563\\
2	0.000172550140028704\\
3	0.000140376605320069\\
4	0.000107100513563996\\
5	8.23662870346187e-05\\
6	6.46726079391869e-05\\
7	5.18793886907787e-05\\
8	4.24303781467739e-05\\
};
\addplot [color=mycolor1, dashed, forget plot]
  table[row sep=crcr]{%
1	0.000109312133327433\\
2	0.000161435456268566\\
3	0.00013135823381627\\
4	0.0001002209839156\\
5	7.70740041378115e-05\\
6	6.0515722116461e-05\\
7	4.85436875511392e-05\\
};
\addplot [color=blue, dashed, forget plot]
  table[row sep=crcr]{%
1	0.000102320235222719\\
2	0.000151361112921911\\
3	0.000123181869994604\\
4	9.39836780582405e-05\\
5	7.22759115149296e-05\\
6	5.67471377175833e-05\\
7	4.5519680270445e-05\\
};
\addplot [color=red, dashed, forget plot]
  table[row sep=crcr]{%
1	9.59778185980475e-05\\
2	0.000142201210563741\\
3	0.000115745859167293\\
4	8.8311070797209e-05\\
5	6.79123370496138e-05\\
6	5.33199526940309e-05\\
7	4.27697062502231e-05\\
};
\addplot [color=green, dashed, forget plot]
  table[row sep=crcr]{%
1	9.02069187338484e-05\\
2	0.000133848334517461\\
3	0.000108963437849385\\
4	8.31369921515668e-05\\
5	6.39323512557911e-05\\
6	5.01941389114968e-05\\
7	4.0261626069793e-05\\
};
\addplot [color=black, dashed, forget plot]
  table[row sep=crcr]{%
1	8.49409208985018e-05\\
2	0.000126210391170889\\
3	0.000102760190633488\\
4	7.84046883087753e-05\\
5	6.02922739635359e-05\\
6	4.73353676557764e-05\\
};
\addplot [color=mycolor1, dashed, forget plot]
  table[row sep=crcr]{%
1	8.01226348828843e-05\\
2	0.000119208058454112\\
3	9.70719997372971e-05\\
4	7.40652581904255e-05\\
5	5.69544703541425e-05\\
6	4.47140627820031e-05\\
};
\addplot [color=blue, dashed, forget plot]
  table[row sep=crcr]{%
1	7.57027399265198e-05\\
2	0.000112772717907217\\
3	9.18433810853658e-05\\
4	7.00763847151231e-05\\
5	5.38863739879547e-05\\
6	4.23046324340981e-05\\
};
\addplot [color=red, dashed, forget plot]
  table[row sep=crcr]{%
1	7.16385208332504e-05\\
2	0.000106844767102876\\
3	8.70261258324457e-05\\
4	6.64012989388673e-05\\
5	5.10596892028259e-05\\
6	4.00848418621314e-05\\
};
\addplot [color=green, dashed, forget plot]
  table[row sep=crcr]{%
1	6.78928346078412e-05\\
2	0.000101372234525665\\
3	8.25781848866613e-05\\
4	6.30079294703378e-05\\
5	4.84497362062131e-05\\
};
\addplot [color=black, dashed, forget plot]
  table[row sep=crcr]{%
1	6.44332608100706e-05\\
2	9.63096365047545e-05\\
3	7.84627479973876e-05\\
4	5.98682002321826e-05\\
5	4.60349104113301e-05\\
};
\addplot [color=mycolor1, dashed, forget plot]
  table[row sep=crcr]{%
1	6.1231399238082e-05\\
2	9.16170290333752e-05\\
3	7.46474795612073e-05\\
4	5.69574477193857e-05\\
5	4.37962338000197e-05\\
};
\addplot [color=blue, dashed, forget plot]
  table[row sep=crcr]{%
1	5.8262286460442e-05\\
2	8.72592174107067e-05\\
3	7.1103881396195e-05\\
4	5.42539350664173e-05\\
5	4.17169808276663e-05\\
};
\addplot [color=red, dashed, forget plot]
  table[row sep=crcr]{%
1	5.5503908767407e-05\\
2	8.32050943782026e-05\\
3	6.7806758927369e-05\\
4	5.1738444969085e-05\\
5	3.97823650531981e-05\\
};
\addplot [color=green, dashed, forget plot]
  table[row sep=crcr]{%
1	5.29367937563974e-05\\
2	7.94270834200139e-05\\
3	6.47337720470765e-05\\
4	4.93939371691843e-05\\
};
\addplot [color=black, dashed, forget plot]
  table[row sep=crcr]{%
1	5.05436663756909e-05\\
2	7.59006685572272e-05\\
3	6.18650556312162e-05\\
4	4.7205259058325e-05\\
};
\addplot [color=mycolor1, dashed, forget plot]
  table[row sep=crcr]{%
1	4.83091580623097e-05\\
2	7.26039956061574e-05\\
3	5.91828976412856e-05\\
4	4.51589001899677e-05\\
};
\addplot [color=blue, dashed, forget plot]
  table[row sep=crcr]{%
1	4.62195598107673e-05\\
2	6.95175327595221e-05\\
3	5.66714650318545e-05\\
4	4.32427832472992e-05\\
};
\addplot [color=red, dashed, forget plot]
  table[row sep=crcr]{%
1	4.42626117523881e-05\\
2	6.66237806056236e-05\\
3	5.43165695207071e-05\\
4	4.14460854100129e-05\\
};
\addplot [color=green, dashed, forget plot]
  table[row sep=crcr]{%
1	4.24273232081385e-05\\
2	6.39070235317619e-05\\
3	5.2105466730401e-05\\
4	3.97590851667351e-05\\
};
\addplot [color=black, dashed, forget plot]
  table[row sep=crcr]{%
1	4.07038182721475e-05\\
2	6.13531158909497e-05\\
3	5.00266833750595e-05\\
4	3.81730305186734e-05\\
};
\addplot [color=mycolor1, dashed, forget plot]
  table[row sep=crcr]{%
1	3.90832028703825e-05\\
2	5.89492974886197e-05\\
3	4.80698681067926e-05\\
};
\end{axis}
\end{tikzpicture}%